# A priori bounds for the dynamic fractional $\Phi^4$ model on $\mathbb{T}^3$ in the full subcritical regime

S. Esquivel & H. Weber

Institut für Analysis und Numerik
Universität Münster

*Email:* salvador.esquivel@uni-muenster.de

*Email:* hendrik.weber@uni-muenster.de

**Abstract:** We show a priori bounds for the dynamic fractional $\Phi^4$ model on $\mathbb{T}^3$ in the full subcritical regime using the framework of Hairer's regularity structures theory [Hai14]. Assuming the model bounds our estimates imply global existence of solutions and existence of an invariant measure. We extend the method developed for the usual heat operator in [CMW23] to the fractional heat operator, thereby treating a more physically relevant model. A key ingredient in this work is the development of localised multilevel Schauder estimates for the fractional heat operator which is not covered by Hairer's original work. Furthermore, the algebraic arguments from [CMW23] are streamlined significantly.

## 1. Introduction

In this work we show a priori bounds for the dynamic fractional $\Phi^4$ model on $\mathbb{T}^3$, which is formally given as the solution to the non-linear stochastic partial differential equation (SPDE):

$$(\partial_t + (-\Delta)^s)\,\varphi = -\lambda\,\varphi^3 + m\,\varphi + \xi \qquad \text{on} \quad \mathbb{R}^+ \times \mathbb{T}^3, \tag{1.1}$$

where $\xi$ is a space-time white noise on $\mathbb{R}^+ \times \mathbb{T}^3$, $\lambda > 0$, $m \in \mathbb{R}$ and $(-\Delta)^s$ is the fractional Laplacian on the torus for $s \in (3/4, 1)$ defined as the Fourier multiplier with symbol $|\cdot|^{2s}$.

In the regime $s \in (3/4, 1)$ the SPDE (1.1) is highly singular but subcritical, and therefore solutions can be described using Hairer's theory of regularity structures [Hai14]. For simplicity, we consider $\lambda = 1$ and $m = 0$. In Section 2 we will describe the regularity structure $(T, A, G)$ where a precise meaning of (renormalised) solution of (1.1) is given. We state our main result:

**Theorem 1.1.** *Fix $s \in (3/4, 1)$, $0 < \kappa \ll 1$ and $\gamma \in (3 - 2s, 2s)$. Let $(\Pi, \Gamma)$ be a smooth, 1-periodic in space and weakly admissible model (Definition 2.11) on the regularity structure described in Section 2, $\Phi \in \mathcal{D}^\gamma(\Gamma)$ a 1-periodic modelled distribution which solves (2.9) and (2.10). Consider $v \colon \mathbb{R} \times \mathbb{T}^3 \to T$ the coefficient at $\mathbf{1}$ of $\Phi$, then for all $t \in (0,1)$ we have the bound*

$$\|v\|_{(t^{2s},1)\times\mathbb{R}^3} \lesssim \max\left\{ t^{-s}, \max_{\substack{\tau \in \mathcal{V}_{0,2s} \\ |\mathfrak{e}(\tau)|=0}} [\Pi;\tau]_K^{\frac{2s}{\mathfrak{l}(\tau)(4-3s-2\kappa)}}, \max_{\substack{\tau \in \mathcal{V}_{0,\gamma} \\ k \in \{0,e_1,\ldots,e_d\}}} [\Gamma\mathcal{I}(\tau), \boldsymbol{X}^k]_K^{\frac{2s}{\mathfrak{l}(\tau)(4-3s-2\kappa)}}, \right.$$

$$\left. \max_{\substack{\tau \in \mathcal{V}_{0,\gamma} \\ |\mathfrak{e}(\tau)|=0}} \sup_{x \in (0,1]\times(0,1]^3} \|\Pi_x \mathcal{I}(\tau)\|_{(0,1]\times\mathbb{R}^3}^{\frac{2s}{\mathfrak{l}(\tau)(4-3s-2\kappa)}} \right\}, \tag{1.2}$$

*where $\|\cdot\|_P$ denotes the supremum norm on the set $P$, $K = [0,1] \times \overline{B_4} \subset \mathbb{R} \times \mathbb{R}^3$, the sets of trees $\mathcal{V}_{0,2s}, \mathcal{V}_{0,\gamma}$ are defined in (2.4), $\mathfrak{l}(\tau)$ is the number of noises in the tree $\tau$, $|\mathfrak{e}(\tau)|$ and $|\mathfrak{n}(\tau)|$ are the total edge and total node decorations of the tree $\tau$ as defined in (A.13) and $\mathfrak{m}(\tau)$ is defined in (2.2). Moreover, the implicit proportionality constant depends only on $s$ and $\gamma$.*

Theorem 1.1 takes as an input a model $(\Pi; \Gamma)$. In the applications of our result to equation (1.1) this model would have to be constructed as a polynomial random variable in the noise $\xi$. However, we stress, that this construction is not (fully) contained in the existing literature [CH16, LOTT24, HS24] because the fractional heat operator is not covered there. Still, for any given $\tau$ one would expect to construct $[\Pi; \tau]_K^{\overline{\frac{2s}{\mathfrak{l}(\tau)(4-3s-2\kappa)}}}$, $[\Gamma \mathcal{I}(\tau), \boldsymbol{X}^k]_K^{\overline{\frac{2s}{\mathfrak{l}(\tau)(4-3s-2\kappa)}}}$ and $\sup_{x \in (0,1] \times (0,1]^3} \|\Pi_x \mathcal{I}(\tau)\|_{(0,1] \times \mathbb{R}^3}^{\overline{\frac{2s}{\mathfrak{l}(\tau)(4-3s-2\kappa)}}}$ as random variables in the $\mathfrak{l}(\tau)$-th inhomogeneous Wiener chaos over $\xi$. Therefore, just as in [CMW23], all terms on the right hand side of (1.2) would have the same stretched exponential stochastic integrability.

Even though we consider the equation (1.1) on the torus the linear polynomials that appear in the framework of regularity structures force us to leave the periodic setting. In general, the only part in our argument that relies on us working in the torus is in the large scale control where we use that $\|v\|_{(t^{2s},1] \times \mathbb{R}^3}$ is finite and allows us to avoid the use of weights.

The SPDE (1.1) arises as the stochastic quantisation of the *fractional* $\Phi_3^4$ measure which was introduced in [BMS03] in the context of Euclidean quantum field theory (EQFT). This measure on fields $\phi \colon \mathbb{T}^3 \longrightarrow \mathbb{R}$ can formally be written as

$$\mu(\mathrm{d}\phi) \propto \exp\left\{-\int_{\mathbb{T}^3} (\lambda \, |\phi(x)|^4 + m \, |\phi(x)|^2) \, \mathrm{d}x\right\} \nu_s(\mathrm{d}\phi), \tag{1.3}$$

where $\nu$ is a centred Gaussian measure on fields with covariance given by $(-\Delta)^{-s}$. In [BMS03] it is shown, under an ultraviolet cutoff of the measure, that for $0 < \varepsilon := 4s - 3 \ll 1$ in the infrared limit this measure converges to a non-Gaussian fixed point of the renormalisation group flow. The relevance of this measure in EQFT comes from it being a candidate to satisfy the Osterwalder-Schrader axioms, in particular reflection positive, precisely in the regime $s \in (3/4, 1]$.

In the last years a lot of work has been developed around the dynamical $\Phi^4$ model. Since the solutions for these equations are distributional valued as soon as $s < 3/2$, there is no canonical meaning for the cubic term in the equation and therefore giving a meaning to this equation was an open problem for a long time. The first works to obtain local in time existence and uniqueness results for this equation when $s = 1$ were Hairer's theory of regularity structures [Hai14], the theory of paracontrolled distributions by Gubinelli, Imkeller, and Perkowski [GIP15], and Kupiainen's renormalisation group approach [Kup15]. A priori estimates which lead to global in time existence of solutions were obtained in $d = 2$ in [MW17a], and for $d = 3$ in [MW17b, MW20, GH21] to mention some.

For the equation (1.1) for $s = 1$, the dimension of the space only plays a role via the *effective* dimension of the driving noise. In [BCCH20, Section 2.8.2] a way to cover the full subcritical regime of fractional dimension by fixing $s = 1$ and instead considering driving noises of regularity $-3 + \delta$ for $\delta \in (0, 1/2)$ is proposed, and in [CMW23] a priori estimates for this model are obtained.

In this work we follow the strategy originally developed in [MW20] by treating small and large scales using different arguments. The small scale behaviour is treated using the local approximations to the solution in the framework of regularity structures and the regularisation properties of the heat operator for them. The large scales are controlled by considering a (random) regularisation of the solution and using the strong dampening effect of the $-\varphi^3$ non-linearity. The large scale control relies on the maximum principle, which still holds true for the fractional heat operator considered in (1.1) as long as $s \in (0, 1)$, which allows us to replicate the strategy used in [MW20]. This maximum principle is the reason we consider the three-dimensional torus and $s \in (3/4, 1)$ instead of the four-dimensional torus and $s > 1$, where the *effective* dimension of the equation is the same (it is singular and subcritical) but the operator no longer satisfies a maximum principle. For the small scales, it is well-known that the fractional heat operator has some smoothing properties similar to those of the classical heat operator in the form of Schauder estimates. However, the *multilevel* Schauder estimates one would need to use in the framework of regularity structures are not contained in the general framework developed in [Hai14] since the kernel associated to the fractional heat operator does not satisfy the smoothness conditions there. This lack of smoothness was pointed out in [CL22] and can be seen from the non-smoothness of the symbol $|\cdot|^{2s}$ at the origin, or equivalently by the lack of decay of its Green's kernel.

Section 3 is dedicated to developing a novel multilevel Schauder estimates for the fractional heat operator. The formulation of our multilevel Schauder estimates follows that one introduced in [OSSW18] where we consider base-point dependent functions or *germs*. Our proof is strongly inspired by the work [FR17], and relies on the Liouville principle for the fractional heat operator developed in there and the scaling properties of the operator; see also [SS24] for a proof of multilevel Schauder estimates for the local heat operator using a similar strategy. This part of the work is done in a framework independent of equation (1.1), and we cover the full regime of $s \in (0,1)$ without restricting ourselves to the periodic setting.

As in [CMW23], as the effective dimension of our equation approaches criticality the number of terms needed to describe the solution diverges and the algebraic relationships are untreatable "by hand". In the previously mentioned work, some algebraic components of the theory of regularity structures are re-developed based on non-commutative product of trees. We clean up the algebraic arguments by putting ourselves into the standard framework of decorated non-planar trees as developed in [BHZ19, BCCH20]. In order to deal with the combinatorial factors arising in this framework our proofs rely on a duality formula obtained in [BM23].

Independently to our work, the recent work of Duch, Gubinelli and Rinaldi [DGR23] constructs the *fractional* $\Phi^4$ measure (1.3) in the full space via stochastic quantisation covering the full subcritical regime $s \in (3/4, 1)$. This last work follows the approach introduced in [AK20] of using finite dimensional approximations of the equation in stationarity combined with a priori estimates, similar to the ones we obtain in Theorem 1.1, to show tightness of the corresponding finite dimensional measures. Moreover, Duch et al show that their constructed measure is reflection positive, translation invariant and, following the *tilting* argument from [HS22], they show the required exponential integrability to apply the Osterwalder-Schrader reconstruction. The main differences between our approach and [DGR23] are that in the later, the *flow equation* approach to SPDEs developed by one the authors in [Duc21] is used for the small scales analysis, which allows them to avoid the algebra required to describe the recentering procedure using regularity structures. The large scales analysis in [DGR23] and our work are similar in the sense that a maximum principle that exploits the damping of the negative sign in the cube is used on a regularised equation, and the (random) scales of the regularisation are chosen as a function of the solution.

## 1.1. Notation

Given $x \in \mathbb{R}^{1+d}$ we will denote $x_0 \in \mathbb{R}$ to the time component and $x_{1:d} := (x_1, \ldots, x_d) \in \mathbb{R}^d$ the space component. We fix $s \in (0,1)$ and work with the fractional parabolic metric

$$d(x,y) := \max\left\{|x_0 - y_0|^{\frac{1}{2s}}, |x_{1:d} - y_{1:d}|\right\},$$

where $|\cdot|$ denotes the Euclidean norm in $\mathbb{R}^d$. This metric reflects the scaling of the operator $\mathscr{L} := (\partial_t + (-\Delta)^s)$. For $r > 0$ we define $r\,x := (r^{2s}\,x_0, r\,x_{1:d}) \in \mathbb{R}^{1+d}$, and with this definition the metric $d$ is 1-homogeneous, i.e., $d(r\,x, r\,y) = r\,d(x,y)$. Clearly $d$ is also translation invariant.

Given $r > 0$ and $x \in \mathbb{R}^{1+d}$ we define $B_r(x_{1:d}) \subset \mathbb{R}^d$ as the Euclidean ball of radius $r > 0$ centred at $x_{1:d} \in \mathbb{R}^d$ and $B_r(x) \subset \mathbb{R}^{1+d}$ as the half-parabolic ball of radius $r$ centred at $x \in \mathbb{R}^d$, i.e.,

$$B_r(x) := \{y \in \mathbb{R}^{1+d} : d(x,y) < r, y_0 \leqslant x_0\} = (x_0 - r^{2s}, x_0] \times B_r(x_{1:d}) \subset \mathbb{R}^{1+d}.$$

It will be clear from the context if we are considering an Euclidean ball or a half-parabolic ball.

Consider a space-time function $v : \mathbb{R}^{1+d} \to \mathbb{R}$. Since our analysis will not require time derivatives to expand our solutions we will denote by $\nabla$ the gradient in the spatial components, i.e., $\nabla v(x) = (\partial_1 v(x), \ldots, \partial_d v(x)) \in \mathbb{R}^d$. Moreover, given a space-time vector $x \in \mathbb{R}^{1+d}$ and a space vector $w \in \mathbb{R}^d$ we will write $w \cdot x := w \cdot x_{1:d}$ to denote the spatial inner product.

**Fractional Laplacian**

We consider $\mathscr{S}(\mathbb{R}^d)$ the space of Schwartz functions and its dual $\mathscr{S}'(\mathbb{R}^d)$ the space of tempered distributions. The fractional Laplacian $(-\Delta)^s$ can be defined on Schwartz functions as the Fourier multiplier with symbol $\zeta \mapsto |\zeta|^{2s}$. Even though $(-\Delta)^s(\mathscr{S}(\mathbb{R}^d)) \subset C^\infty(\mathbb{R}^d)$, due to the symbol not being smooth at the origin one has in general that $(-\Delta)^s \varphi \notin \mathscr{S}(\mathbb{R}^d)$ for $\varphi \in \mathscr{S}(\mathbb{R}^d)$ since it does not have the required decay at infinity. However, it can be shown that

$$(-\Delta)^s \varphi \in \mathscr{S}_s(\mathbb{R}^d) := \{\varphi \in C^\infty(\mathbb{R}^d) : (1+|\cdot|^{d+2s}) D^k \psi \in L^\infty(\mathbb{R}^d) \, \forall k \in \mathbb{N}^d\} \supset \mathscr{S}(\mathbb{R}^d),$$

which allows to extend the definition by duality to the dual space $\mathscr{S}_s'(\mathbb{R}^d) \subset \mathscr{S}'(\mathbb{R}^d)$. We refer the reader to [Sti19, Lemma 1] for the proof and details on the topologies involved in the previous statement. It can be seen that $(-\Delta)^{2s}$ is translation invariant and $2s$-homogeneous, i.e., $(-\Delta)^s(f(z + \sigma \cdot)) = \sigma^{2s}((-\Delta)^s f)(z + \sigma \cdot)$.

Distributions on the torus can be interpreted as 1-periodic tempered distributions, and the fractional Laplacian defined as the corresponding Fourier multiplier agree under this identification. Let $f: \mathbb{T}^d \to \mathbb{R}$ be a continuous function and consider its extension $f: \mathbb{R}^d \to \mathbb{R}$ as a continuous 1-periodic function. Since $f$ is periodic, in particular it has no growth at infinity and therefore if $f$ is smooth enough we can write the action of $-(-\Delta)^s$ to $f$ as a singular integral

$$-((-\Delta)^s f)(x) = c_{d,s} \int_{\mathbb{R}^d} \frac{f(x+y) + f(x-y) - 2f(x)}{|y|^{d+2s}} \, \mathrm{d}y, \tag{1.4}$$

where the constant is explicitly given by

$$c_{d,s} := 2^{2s-1} s \frac{\Gamma\left(\frac{d+2s}{2}\right)}{\Gamma(1-s)} \pi^{-\frac{d}{2}},$$

and via the Bochner's (or semigroup) representation, i.e.,

$$-(-\Delta)^s f = \frac{1}{|\Gamma(-s)|} \int_0^\infty (\mathrm{e}^{t\Delta} f - f) \frac{\mathrm{d}t}{t^{1+s}}, \tag{1.5}$$

where $\{\mathrm{e}^{t\Delta}\}_{t \geqslant 0}$ is the heat semigroup, i.e., $\mathrm{e}^{t\Delta} f$ is the convolution of $f$ with the heat kernel $k_t(z) = (4\pi t)^{-d/2} \exp\{-|z|^2/(4t)\}$. For $f: \mathbb{R}^d \to \mathbb{R}$, not necessarily periodic, to make sense of the pointwise definition of $((-\Delta)^s f)(x)$ as in (1.4) it is enough to assume (see [FR24, Section 1.10]) that $f$ is $2s + \varepsilon$-Hölder continuous around $x$ for some $0 < \varepsilon \ll 1$, and that $f$ satisfies the global integrability condition

$$\|f\|_{L^1_{w_s}(\mathbb{R}^d)} := \int_{\mathbb{R}^d} \frac{|f(x)|}{1+|x|^{d+2s}} < +\infty. \tag{1.6}$$

In particular if $f \in C^2(\mathbb{R}^d) \cap L^\infty(\mathbb{R}^d)$ then one can define $(-\Delta)^s f$ via (1.4). This is the representation we use for the Schauder estimates in Section 3. We refer to the survey [Kwa17] for the equivalence between different representations of the fractional Laplacian.

## 2. Modelled Solution

In this section we define the regularity structure where a renormalised solution to

$$\mathscr{L}\varphi := (\partial_t + (-\Delta)^s)\varphi = -\varphi^3 + \xi$$

can be defined in terms of modelled distributions. We perform the Da Prato-Debussche trick to work with a function-valued remainder and show some algebraic identities that will be needed in the proof of the main result. General definitions of the theory of regularity structures used in this section are recalled in the Appendix A. We work on arbitrary spatial dimensions, i.e., $\mathbb{R}^{1+d}$ since Theorem 1.1 still holds true provided we fix the regularity of the driving noise to be that of a 3-dimensional white noise (see (A.9)).

## 2.1. Definition of our regularity structure

We work with decorated trees $(\tau, \rho_\tau, \mathfrak{l}, \mathfrak{n}, \mathfrak{e})$ where $\tau$ is a non-planar tree with root $\rho_\tau$, node decorations $\mathfrak{n}, \mathfrak{l}$ which represent polynomials and noises respectively, and edge decorations $\mathfrak{e}$ which encode derivatives on the kernel. We refer the reader to Section A.2 for details on decorated trees. Since the non-linearity in our PDE is cubic, in order to describe the solution we will only need to consider a subset $\mathcal{T}$ of decorated trees that contains the noise $\Xi$, the monomials $\{\boldsymbol{X}^k\}_{k \in \mathbb{N}^{1+d}}$ and which is closed under, the following rule:

$$\tau_1, \tau_2, \tau_3 \in \mathcal{T}, j \in \{1, \ldots d\}, k \in \mathbb{N}^{1+d} \Rightarrow \boldsymbol{X}^k \mathcal{I}(\tau_1), \boldsymbol{X}^k \mathcal{I}_j(\tau_1), \boldsymbol{X}^k \mathcal{I}(\tau_1) \mathcal{I}(\tau_2), \mathcal{I}(\tau_1) \mathcal{I}(\tau_2) \mathcal{I}(\tau_3) \in \mathcal{T}. \quad (2.1)$$

This means that $\mathcal{T}$ consists of decorated trees with at most three branches leaving every internal node (with respect to the natural direction associated to the tree which goes from the root to the leaves). We call this trees sub-ternaries. The set $\mathcal{T}$ is stable when removing node decorations $\mathfrak{n}$ and edge decorations $\mathfrak{e}$. Moreover, we follow the convention of setting $\mathcal{I}(\boldsymbol{X}^k) := 0$ (see [FH20, Remark 14.26]).

**Remark 2.1.** The trees with decorated edges are not needed to describe the equation, however it will be convenient to include them since it makes the associated subcritical rule complete in the sense of [BHZ19]. See the discussion in [BCCH20, Remark 2.27] where the authors discuss the associated subcritical rule and how this term relates to the renormalised equation.

We say that a tree $\tau \in \mathcal{T}$ is *full* if every inner node $v \in N_\tau \setminus L_\tau$ contains exactly three outgoing branches. From (2.1) it can be seen inductively that full trees $\tau$ in $\mathcal{T}$ have zero node polynomial decorations and trivial edge decorations. It will be important to see how far a tree is from being full, for this we define the counting function $\mathfrak{m} : \mathcal{T} \to \mathbb{N}$ recursively defined as $\mathfrak{m}(\Xi) := 0$ and

$$\mathfrak{m}(\tau) := \begin{cases} 2 + \mathfrak{m}(\tau_1) & \text{if } \tau = \boldsymbol{X}^k \mathcal{I}(\tau_1) \text{ or } \tau = \boldsymbol{X}^k \mathcal{I}_j(\tau), k \in \mathbb{N}^{1+d} \\ 1 + \mathfrak{m}(\tau_1) + \mathfrak{m}(\tau_2) & \text{if } \tau = \boldsymbol{X}^k \mathcal{I}(\tau_1) \mathcal{I}(\tau_2), k \in \mathbb{N}^{1+d} \\ \mathfrak{m}(\tau_1) + \mathfrak{m}(\tau_2) + \mathfrak{m}(\tau_3) & \text{if } \tau = \mathcal{I}(\tau_1) \mathcal{I}(\tau_2) \mathcal{I}(\tau_3) \end{cases}. \quad (2.2)$$

The mnemonic is that $\mathfrak{m}(\tau)$ counts the number of "missing" branches from $\tau$, and its definition is independent of node or edge decorations of the tree.

We consider homogeneity function $|\cdot|$ as defined in Section A.2, and we omit the sub-index $\mathfrak{s}$ from our notation. In particular $|\Xi| = -s - 3/2 - \kappa$ which is the regularity of a 3-dimensional white noise (see (A.9)), $\mathcal{I}$ and $\mathcal{I}_j$ improve homogeneity by $2s$ and $2s-1$ respectively (see (A.11)), monomials satisfy (A.10) and the homogeneity of a product is the sum of their homogeneities (see (A.12)).

**Lemma 2.2**. *We have that $|\Xi| < |\tau|$ for all $\tau \in \mathcal{T} \setminus \{\Xi\}$ if and only if $s > \frac{3}{4}$ and $\kappa > 0$ is small enough. Moreover, in this case we have that for any $\beta \in \mathbb{R}$ the set $\mathcal{T}_{<\beta} := \{\tau \in \mathcal{T} : |\tau| < \beta\}$ is finite.*

**Proof.** First assume that $s \leqslant \frac{3}{4}$, then $4s \leqslant 3$ and $3s - 9/2 \leqslant s - 3/2$. Consider the tree $\tau := \mathcal{I}(\Xi) \mathcal{I}(\Xi) \mathcal{I}(\Xi) \in \mathcal{T}$, then we have

$$|\tau| = 3(|\Xi| + 2s) = 3s - \frac{9}{2} - 3\kappa \leqslant 3s - \frac{9}{2} - \kappa \leqslant s - \frac{3}{2} - \kappa = |\Xi|,$$

which shows the necessity of condition $s > 3/4$. On the other hand, assume $s > 3/4$ and $\kappa > 0$ is small enough. Since $2s - 1 > 0$ abstract integration by $\mathcal{I}, \mathcal{I}_j$ and multiplication by polynomials increase homogeneity, then from the recursive definition (2.1) it is enough to prove that if the result is true for $\tau_1, \tau_2, \tau_3 \in \mathcal{T}$ then it is also true for $\mathcal{I}(\tau_1)\mathcal{I}(\tau_2)$ and for $\mathcal{I}(\tau_1)\mathcal{I}(\tau_2)\mathcal{I}(\tau_3)$. Since $\kappa \ll 1$ is assumed to be small enough (depending on $s$) we have

$$|\mathcal{I}(\tau_1)\mathcal{I}(\tau_2)| = 4s + |\tau_1| + |\tau_2| > 4s + 2|\Xi| > |\Xi| \Leftrightarrow 4s + |\Xi| > 0 \Leftrightarrow s > 1/2,$$

and

$$|\mathcal{I}(\tau_1)\mathcal{I}(\tau_2)\mathcal{I}(\tau_3)| = 6s + \sum_{i=1}^{3} |\tau_i| > 6s + 3|\Xi| > |\Xi| \Leftrightarrow 6s + 2|\Xi| > 0 \Leftrightarrow s > 3/4,$$

where we used that $|\Xi| = -\left(\frac{2s+d}{2}\right) - \kappa$. □

Subcriticality only restricts us to work on the regime $s > 3/4$, however we focus on the regime $s \in (3/4, 1)$ since for $s > 1$ there is no maximum principle for the fractional Laplacian $-(-\Delta)^s$, and the case $s = 1$ corresponds to the usual Laplacian $\Delta$ and the well-known $\Phi_3^4$ equation which is treated in [MW20].

We define $T := \langle \mathcal{T} \rangle$ as the vector space spanned by the trees in $\mathcal{T}$, and the grading induced by the set of homogeneities $A := \{|\tau| : \tau \in \mathcal{T}\}$, which decomposes $T$ as

$$T = \bigoplus_{\beta \in A} T_\beta, \qquad T_\beta := \langle \{\tau \in \mathcal{T} : |\tau| = \beta\} \rangle. \tag{2.3}$$

Similarly, we define the space $T_{<\gamma} := \bigoplus_{\beta < \gamma} T_\beta$. We consider the inner product on $T$ defined on (A.22) using the symmetry factor of a tree (see (A.6)). This makes $\mathcal{T}$ an orthogonal (but not orthonormal) basis of $T$.

The graded vector space (2.3) will be the basis for our regularity structure. For convenience to the reader we recall in Appendix A some notions of the theory of regularity structures. In particular the description of the associated structure group $G$ (and the associated *coaction* $\Delta$) is presented in Section A.3.

**Remark 2.3.** Since the set $\mathcal{T}_{<2s} := \{\tau \in \mathcal{T} : |\tau| < 2s\}$ is finite in the subcritical regime $s \in (3/4, 1)$ (Lemma 2.2) we can choose $0 < \kappa \ll 1$ small enough so that all the non-polynomial symbols here have non-integer homogeneity. This restriction is only to avoid technical limitations of the Schauder estimates. We consider this restriction on objects with homogeneity up to $2s$ since we only need to expand the right-hand side of (1.1) up to trees of negative homogeneity, which implies that the solution itself will only contain trees with homogeneity less than $2s$.

We define the following sets:

$$\mathcal{P} := \{\boldsymbol{X}^k\}_{k \in \mathbb{N}^{1+d}}, \quad \mathcal{W} := \{\tau \in \mathcal{T} \setminus \mathcal{P} : |\mathcal{I}(\tau)| < 0\}, \quad \mathcal{V} := \{\tau \in \mathcal{T} \setminus \mathcal{P} : |\mathcal{I}(\tau)| > 0\},$$

and denote by $\mathcal{I}(\mathcal{T})$ the set of trees of the form $\mathcal{I}(\tau)$ with $\tau \in \mathcal{T}$, and similarly for other sets. For $\alpha, \beta \in \mathbb{R}$ with $\alpha < \beta$ we define

$$\mathcal{V}_{\alpha,\beta} := \{\tau \in \mathcal{V} : |\mathcal{I}(\tau)| \in [\alpha, \beta)\}. \tag{2.4}$$

Then $\mathcal{W} \cup \mathcal{V}_{0,2s}$ consists of all the trees of negative homogeneity, which are precisely the ones needed to give a good enough local description to the right-hand side of (1.1). Moreover, all trees in $\mathcal{I}(\mathcal{W}) \cup \mathcal{I}(\mathcal{V}_{0,2s}) \cup \{\mathbf{1}\} \cup \{\mathbf{X}_j\}_{j=1}^d$ have homogeneity strictly less than $2s$, and these are precisely the trees needed to give a description of the solution $\varphi$ to (1.1).

**Remark 2.4.** The elements in $\mathcal{W}$ are the ones that also appear in the "Wild expansion", and consists of the trees with the worst homogeneity. Since the structure group acts trivially on $\mathcal{I}(\mathcal{W})$ (see Lemma 2.6) we can perform the De Prato-Debussche trick and work with a remainder $V$ (see (2.13)) which will satisfy an equation that no longer involves trees in $\mathcal{W}$ (see (2.14)). The mnemonics of $\mathcal{V}$ is that these are the non-polynomial trees needed for the expansion of $V$.

The next results give us some characterisation of these trees.

**LEMMA 2.5.** *We have that $\tau \in \mathcal{W} \Longrightarrow \mathfrak{m}(\tau) = 0$, and $\tau$ has zero edge and polynomial decorations, i.e., $\tau$ is full. Moreover, $\tau \in \mathcal{W} \setminus \{\Xi\}$ if and only if there exists $\tau_1, \tau_2, \tau_3 \in \mathcal{W}$ such that $\tau = \mathcal{I}(\tau_1) \mathcal{I}(\tau_2) \mathcal{I}(\tau_3)$ and $|\tau| < -2s$.*

**Proof.** Assume that $\tau \in \mathcal{T}$ is such that $\mathfrak{m}(\tau) > 0$ and consider the node $v \in N_\tau$ where a branch is missing. We consider $\tau'$ the decorated tree obtained by grafting a noise on $\tau$, with an edge $\mathcal{I}$, on the vertex $v$ where it is missing a branch, i.e., $\Xi \curvearrowright_0^v \tau$ and by setting trivial decorations to $\tau'$ (see (A.17) for the definition of $\curvearrowright_0^v$). In this way $\tau' \in \mathcal{T}$ and it satisfies:

$$|\tau'| = |\Xi \curvearrowright_0^v \tau| - |\mathfrak{n}(\tau)| = |\tau| + |\Xi| + 2s - |\mathfrak{n}(\tau)| \geqslant |\tau| + |\Xi| + 2s,$$

and therefore $|\tau| - |\tau'| \geqslant -(|\Xi| + 2s)$. By definition and Lemma 2.2 the set of homogeneities of $\mathcal{W}$ is contained in the interval $[|\Xi|, -2s)$ which has length $-(|\Xi| + 2s) = 3/2 - s + k > 3/4$ in the subcritical regime $s > 3/4$. Since $\tau' \in \mathcal{T}$ by construction, Lemma 2.2 implies:

$$|\tau| = |\tau'| + |\tau| - |\tau'| \geqslant |\Xi| - (|\Xi| + 2s) = 2s,$$

and therefore $\tau \notin \mathcal{W}$. Now assume that $\tau \in \mathcal{T}$ has non-zero polynomial decoration, i.e., $|\mathfrak{n}(\tau)| \neq 0$ and consider $\tau' \in \mathcal{T}$ the decorated tree obtained by setting all decorations in $\tau$ to $0 \in \mathbb{N}^{d+1}$, then we have that $|\tau'| - |\tau| = |\mathfrak{n}(\tau)| \geqslant 1$. Since in the subcritical regime $s > 3/4$ the length of the interval $[|\Xi|, -2s)$ is bounded from below by $3/4$ we conclude that $\tau' \notin \mathcal{W}$. On the other hand, by definitions (2.1) and (2.2) we have that if $\tau \in \mathcal{T}$ has a decorated edge, then $\mathfrak{m}(\tau) = 2$, and by the previous part $\tau \notin \mathcal{W}$. The last part follows immediately from this by the recursion (2.1), the definition (2.2) of $\mathfrak{m}$, and the previous part of this result. $\square$

**LEMMA 2.6.** *If $\tau \in \mathcal{T}$ is full, then the coaction $\Delta$ and the structure group $G$ act trivially on it, i.e., for all $\Gamma \in G$*

$$\Delta \tau = \tau \otimes \mathbf{1}, \qquad \Delta \mathcal{I}(\tau) = \mathcal{I}(\tau) \otimes \mathbf{1}, \qquad \Gamma \tau = \tau, \qquad \Gamma \mathcal{I}(\tau) = \mathcal{I}(\tau).$$

*In particular by Lemma 2.5 this holds for all $\tau \in \mathcal{W}$.*

**Proof.** We consider $\mathscr{T}$ the set of decorated trees (not necessarily sub-ternary) as defined in Section A.2. By (A.25) we can write for $\tau \in \mathcal{T}$

$$\Delta \tau = \sum_{\sigma \in \mathcal{T}, \mu \in \mathcal{T}^+} \frac{\langle \mu \star \sigma, \tau \rangle}{\sigma! \, \mu!} \sigma \otimes \mu, \tag{2.5}$$

where the $\star$-product is defined in (A.19) and the set $\mathcal{T}^+$ is defined in Section A.3. By definition (see (2.1)) if $\tau \in \mathcal{T}$ is full then it has no node or edge decoration. Moreover, given any $\sigma \in \mathcal{T}$ and $v \in N_\tau \setminus L_\tau$ the tree $\sigma \curvearrowright_\rho^v \tau$ is not sub-ternary since it will have 4 outgoing branches at the node $v$, and therefore it will not an element of the regularity structure $T$. This implies that for any $\mathcal{I}_\rho(\sigma) \in \mathcal{T}^+$ we have $\mathcal{I}_\rho(\sigma) \star \tau = \sigma \curvearrowright_\rho \tau \notin \mathcal{T}$. This generalises to any $\mu \in \mathcal{T}^+$ of the form $\mu = \prod_{i \in I} \mathcal{I}_{\rho_i}(\sigma_i)$, i.e., $\mu \star \tau \notin T$. Similarly, we have that for any $k \in \mathbb{N}^{1+d} \setminus \{0\}$ the tree $\uparrow_v^k \tau$ (see (A.4) for the definition) is still sub-ternary, but has non-zero total polynomial decoration. Since full trees in $\mathcal{T}$ do not have polynomial decorations by the rule (2.1) defining $\mathcal{T}$, we have that $\uparrow_v^k \tau \notin \mathcal{T}$. This implies that $\boldsymbol{X}^k \star \tau = \uparrow_{N_\tau \setminus L_\tau}^k \tau \notin T$ since each tree in the linear combination has non-zero total node decorations. We conclude that $\mu \star \tau \notin T$ for any $\mu \in \mathcal{T}^+ \setminus \{\mathbf{1}\}$.

Now, since $|\mu \star \sigma| = |\mu| + |\sigma|$ by (A.20) and $|\mu| > 0$ for $\mu \in \mathcal{T}^+$ we have by (A.25) that $\Delta \tau$ can have a non-trivial component $\sigma \otimes \mu$ only for $\sigma \in \mathcal{T}$ such that $|\sigma| \leqslant |\tau|$. Since $|\mu| > 0$ for every $\mu \neq \mathbf{1} \in \mathcal{T}^+$ we have that $\mu \star \tau = \tau$ if and only if $\mu = \mathbf{1}$, and it is enough to consider $\sigma \in \mathcal{T}$ with $|\sigma| < |\tau|$. Since $\tau \in \mathcal{W}$ by definition $|\tau| < -2s$ and $|\sigma| < |\tau|$ implies that $|\sigma| < -2s$, i.e., $\sigma \in \mathcal{W}$, and the previous argument applies, i.e., $\mu \star \sigma \notin T$ for all $\mu \in \mathcal{T}^+ \setminus \{\mathbf{1}\}$ and in particular $\langle \mu \star \sigma, \tau \rangle = 0$. We conclude on (2.5)

$$\Delta \tau = \tau \otimes \mathbf{1} + \sum_{\substack{\sigma \in \mathcal{T} \\ |\sigma| < |\tau|}} \sum_{\mu \in \mathcal{T}^+} \frac{\langle \mu \star \sigma, \tau \rangle}{\sigma! \, \mu!} \sigma \otimes \mu = \tau \otimes \mathbf{1}.$$

Now, by definition of $\mathcal{W}$ we have that $|\mathcal{I}(\tau)| < 0$ and by definition of $\Delta$ on planted trees (A.14):

$$\Delta \mathcal{I}(\tau) = (\mathcal{I} \otimes \mathrm{Id}_{T^+}) \Delta \tau = \mathcal{I}(\tau) \otimes \mathbf{1}.$$

The triviality of $\Gamma$ follows immediately by this and (A.16). □

**LEMMA 2.7.** *For every $\tau \in \mathcal{T}_{<0}$ we have that $(\mathfrak{m}(\tau), |\mathfrak{n}(\tau)|) \in \{(0,0), (1,0), (2,0), (1,1)\}$. Moreover, if $\tau \in \mathcal{T}_{<0}$ is such that $|\mathcal{I}(\tau)| < 1$, then $\mathfrak{n}(\tau) = 0$.*

**Proof.** Consider $\tau \in \mathcal{T}$ such that $\mathfrak{m}(\tau) \geqslant 3$ and let $(v_1, v_2, v_3) \in (N_\tau)^3$ be three vertex where branches are missing (there can be repetition between two of them) and let $\tau' \in \mathcal{T}$ be the tree obtained by grafting $\Xi$ with an edge $\mathcal{I}$ into this missing spots, i.e., $\tau' = \Xi \curvearrowright_0^{v_1} (\Xi \curvearrowright_0^{v_2} (\Xi \curvearrowright_0^{v_3} \tau)) \in \mathcal{T}$. By Lemma 2.2 and Remark 2.3 we have

$$|\tau| = |\tau'| - 3(|\Xi| + 2s) > |\Xi| - 3(|\Xi| + 2s) = 4s - 3 - 2\kappa > 0 \, \left( \Longleftrightarrow s > \tfrac{3}{4} \right).$$

On the other hand, if $\tau \in \mathcal{T}$ is such that $|\mathfrak{n}(\tau)| > 1$ then we have, by the choice of the scaling $\mathfrak{s} = (2s, 1, \ldots, 1)$, that $|\mathfrak{n}(\tau)| \geqslant 2s$. By the recursive definition of $\mathcal{T}$ in (2.1) we have that if a tree has non-trivial decorations then it can not be full, and therefore we can consider $\tau' \in \mathcal{T}$ the decorated tree obtained by grafting a noise into $\tau$ and setting trivial decorations to $\tau'$. Using Lemma 2.2 we conclude that

$$|\tau| = |\tau'| - (|\Xi| + 2s) + |\mathfrak{n}(\tau)| \geqslant |\Xi| - |\Xi| - 2s + 2s = 0.$$

Assume that $\tau \in \mathcal{T}$ is such that $\mathfrak{m}(\tau) = 2$ and $\mathfrak{n}(\tau) > 0$, then $|\mathfrak{n}(\tau)| > 1$ and for some $(v_1, v_2) \in (N_\tau)^2$ we can set $\tau' \in \mathcal{T}$ as the tree $\Xi \curvearrowright_0^{v_1} (\Xi \curvearrowright_0^{v_2} \tau)$ without decorations, and by Lemma 2.2 we have that $|\tau'| > |\Xi|$ since $\tau' \neq \Xi$, and since the next tree with the worst homogeneity is $\mathcal{I}(\Xi) \mathcal{I}(\Xi) \mathcal{I}(\Xi)$ we conclude that $|\tau'| \geqslant 3 |\mathcal{I}(\Xi)|$, and therefore

$$|\tau| = |\tau'| - 2|\mathcal{I}(\Xi)| + |\mathfrak{n}(\tau)| \geqslant 3|\Xi| - 2|\mathcal{I}(\Xi)| + 1 = |\Xi| + 1 = s - \tfrac{1}{2} > 0.$$

At last, the scenario where $(\mathfrak{m}(\tau), |\mathfrak{n}(\tau)|) = (0, 1)$ is ruled out since by the recursion (2.1) which defines $\mathcal{T}$, full trees in $\mathcal{T}$ do not have polynomial decorations. For the last part observe that $\mathfrak{n}(\tau) \neq 0$ implies by the shown part of this result that $(\mathfrak{m}(\tau), \mathfrak{n}(\tau)) = (1, 1)$ and therefore $\tau$ can be written for some $j \in \{1, \ldots d\}$ as $\uparrow_v^{e_j} \sigma$ where $\sigma$ is the same tree as $\tau$ but without node decorations. Since $(\mathfrak{m}(\tau), \mathfrak{n}(\tau)) = (1, 0)$ then by Lemma 2.5 we have that $\sigma \notin \mathcal{W}$ and therefore $0 < |\mathcal{I}(\sigma)| = |\mathcal{I}(\tau)| - 1$ and the result follows. □

**LEMMA 2.8**. *The trees $\mathcal{I}(\mathcal{V}) \cup \mathcal{P}$ span a function-like sector.*

**Proof.** By (A.24) we have that for $\tau \in \mathcal{V}$

$$\Gamma \mathcal{I}(\tau) = \sum_{\sigma \in \mathcal{T}} \left\{ \sum_{\mu \in \mathcal{T}^+} \langle \mu \star \sigma, \tau \rangle \frac{\gamma(\mu)}{\mu!} \right\} \frac{\mathcal{I}(\sigma)}{\sigma!} + \sum_{k \in \mathbb{N}^{1+d}} \gamma(\mathcal{I}_{\rho+k}(\tau)) \frac{\boldsymbol{X}^k}{k!}.$$

The second term belongs to the span of $\mathcal{P}$, and the first in the span of planted trees, and therefore it is enough to see that if $\sigma \in \mathcal{W}$ then $\langle \mu \star \sigma, \tau \rangle = 0$ for all $\mu \in \mathcal{T}^+$. However, the argument in the proof of Lemma 2.5 tells us that $\mu \star \sigma \notin T$ for all $\mu \in \mathcal{T}^+ \setminus \{\boldsymbol{1}\}$, and in particular $\langle \mu \star \sigma, \tau \rangle = 0$ for all $\mu \in \mathcal{T}^+ \setminus \{\boldsymbol{1}\}$. At last, for $\mu = \boldsymbol{1}$ we have that $\boldsymbol{1} \star \sigma = \sigma \in \mathcal{W}$ and therefore $\langle \mu \star \sigma, \tau \rangle = 0$. □

For the next result we assume we are given a model $(\Pi; \Gamma)$ (see Definition A.2) in our regularity structure.

**LEMMA 2.9**. *For all $\tau_1, \tau_2 \in \mathcal{V}$ we have that*

$$(\Pi_x \mathcal{I}(\tau_1))(y) = \gamma_{y,x}(\mathcal{I}(\tau_1)), \tag{2.6}$$

$$(\Pi_x(\mathcal{I}(\tau_1) \mathcal{I}(\tau_2)))(y) = (\Pi_x \mathcal{I}(\tau_1))(y) \, (\Pi_x \mathcal{I}(\tau_2))(y). \tag{2.7}$$

**Proof.** By [Hai14, Proposition 3.31] we know that the action of $\Pi$ on an element $\tau \in \mathcal{T}$ of positive homogeneity is characterised by the action of $\Pi$ on $\bigoplus_{\beta < |\tau|} T_\beta$ and the action of $\Gamma$ on $\bigoplus_{\beta \leqslant |\tau|} T_\beta$. More specifically, for every $\tau \in \mathcal{T}$ with $|\tau| > 0$, and $x \in \mathbb{R}^{1+d}$ the map $y \mapsto F_{\tau,x}(y) := \Gamma_{yx} \tau - \tau \in \mathcal{D}^{|\tau|}$ is a modelled distribution of positive order which satisfies $\Pi_x \tau = \mathcal{R} F_{\tau,x}$. In particular, for planted trees $\mathcal{I}(\tau)$ with $\tau \in \mathcal{V}$ we have by (A.24) a more explicit representation of $F_{\mathcal{I}(\tau),x}$ given by

$$F_{\mathcal{I}(\tau),x}(y) = \sum_{\sigma \in \mathcal{T}} \left\{ \sum_{\mu \in \mathcal{T}^+ \setminus \{0\}} \langle \mu \star \sigma, \tau \rangle \frac{\gamma_{yx}(\mu)}{\mu!} \right\} \frac{\mathcal{I}(\sigma)}{\sigma!} + \sum_{k \in \mathbb{N}^{1+d}} \gamma_{yx}(\mathcal{I}_k(\tau)) \frac{\boldsymbol{X}^k}{k!}. \tag{2.8}$$

Moreover, by Lemma 2.8 $F_{\mathcal{I}(\tau),x}$ is a function-like modelled distribution and by [Hai14, Proposition 3.28] its reconstruction is given by the coefficient at $\boldsymbol{1}$. At last $\langle \boldsymbol{1}, F_{\mathcal{I}(\tau),x}(y) \rangle = \gamma_{yx}(\mathcal{I}(\tau))$ by (2.8), and we obtain

$$(\Pi_x \mathcal{I}(\tau))(y) = (\mathcal{R} F_{\mathcal{I}(\tau),x})(y) = \langle \boldsymbol{1}, F_{\mathcal{I}(\tau),x}(y) \rangle = \gamma_{yx}(\mathcal{I}(\tau)).$$

For the second part it is easy to see, by multiplying two instances of (A.24), that

$$F_{\mathcal{I}(\tau_1) \mathcal{I}(\tau_2),x}(y) := \Gamma_{yx}(\mathcal{I}(\tau_1) \mathcal{I}(\tau_2)) - \mathcal{I}(\tau_1) \mathcal{I}(\tau_2) = \Gamma_{yx}(\mathcal{I}(\tau_1)) \, \Gamma_{yx}(\mathcal{I}(\tau_2)) - \mathcal{I}(\tau_1) \mathcal{I}(\tau_2)$$

contains only symbols of non-negative homogeneity, and therefore its reconstruction is given by the coefficient at $\boldsymbol{1}$, which is equal to the product $\gamma_{yx}(\mathcal{I}(\tau_1)) \, \gamma_{yx}(\mathcal{I}(\tau_2))$, i.e.,

$$\Pi_x(\mathcal{I}(\tau_1) \mathcal{I}(\tau_2)) = \mathcal{R}(F_{\mathcal{I}(\tau_1) \mathcal{I}(\tau_2),x}) = \gamma_{\cdot x}(\mathcal{I}(\tau_1)) \, \gamma_{\cdot x}(\mathcal{I}(\tau_2)) = \Pi_x(\mathcal{I}(\tau_1)) \, \Pi_x(\mathcal{I}(\tau_2)). \quad \square$$

**Remark 2.10.** For trees $\tau \in \mathcal{T}$ where the structure group acts trivially, e.g., full trees by Lemma 2.6, then one has for all $x, y \in \mathbb{R}^{1+d}$ that $\Pi_x \tau = \Pi_y \Gamma_{yx} \tau = \Pi_y \tau$, i.e., $\Pi.\tau$ is independent of the base point. To emphasise this we will denote by $\mathbf{\Pi}\tau := \Pi_x \tau$ for such trees.

## 2.2. Weakly admissible model

One of the main results of the theory of regularity structures are the multilevel Schauder estimates found in [Hai14, Section 5] which encode the regularisation properties at the level of modelled distributions of some convolution operators which include the Green's kernel of differential operators. A compatibility condition between the regularity structure and the operator is encoded via convolutions in the notion of *admissible* model (see [FH20, Definition 14.23]). Since the kernel of the fractional heat operator does not fit into Hairer's framework we follow instead [CMW23, Section 3.3] (see also [MW20, Remark 2.4]) and consider a *weak* admissibility condition that suits better the multilevel Schauder estimates we prove in Section 3.

**Definition 2.11.** *A model $(\Pi; \Gamma)$ is called weakly admissible for $\mathscr{L}$ if:*

1. *On polynomial symbols it acts as the polynomial model, i.e., for all $k \in \mathbb{N}^{1+d}$:*

$$(\Pi_x \boldsymbol{X}^k)(y) = (y-x)^k, \qquad \Gamma_{xy} \boldsymbol{X}^k = (\boldsymbol{X} + (y-x)\mathbf{1})^k.$$

2. *For all $\tau \in \mathcal{T} \setminus \mathcal{P}$, $x \in \mathbb{R}^{1+d}$ and $t \in \mathbb{R}$ we have $(\Pi_x \mathcal{I}(\tau))(t, \cdot) \in L^\infty(\mathbb{R}^d)$, and weakly*

$$\mathscr{L}(\Pi_x \mathcal{I}(\tau)) = \begin{cases} \Pi_x \tau & \text{if } \Gamma \mathcal{I}(\tau) = \mathcal{I}(\tau) \\ \rho\, \Pi_x \tau & \text{if } \Gamma \mathcal{I}(\tau) \neq \mathcal{I}(\tau) \end{cases},$$

*i.e.,*

$$\langle \mathscr{L}(\Pi_x \mathcal{I}(\tau)), \psi \rangle = \begin{cases} \langle \Pi_x \tau, \psi \rangle & \text{if } \Gamma \mathcal{I}(\tau) = \mathcal{I}(\tau) \\ \langle \Pi_x \tau, \rho\, \psi \rangle & \text{if } \Gamma \mathcal{I}(\tau) \neq \mathcal{I}(\tau) \end{cases} \qquad \forall\, \psi \in \mathscr{D}(\mathbb{R}^{1+d}),$$

*where $\rho \in C_c^\infty(\mathbb{R}^d)$ is some fixed test function with $\mathrm{supp}(\rho) \subset B_1$ and extended to $\rho \in C^\infty(\mathbb{R}^{1+d})$ by setting $\rho(x) := \rho(x_{1:d})$.*

The localiser $\rho$ in the previous definition is the simplest way to guarantee that the $L^\infty$-terms appearing in Theorem 1.1 will be finite.

**Remark 2.12.** The positive renormalisation encoded by structure group is the same for our weak admissibility as the one for Hairer's admissibility condition, which we include in Section A.3.

## 2.3. Modelled solution

From now on we will assume we have a weakly admissible model $(\Pi, \Gamma)$ on the regularity structure $(T, G)$ as described in Section 2.1, and assume that we are given a modelled distribution $\Phi \in \mathcal{D}^\gamma(\Gamma)$ that satisfies the (algebraic) equation

$$\Phi = \mathcal{I}(-\Phi^3 + \Xi) + P(\Phi) \quad \text{in} \quad \mathcal{D}^\gamma(\Gamma), \tag{2.9}$$

where $P(\Phi)$ is the polynomial part of $\Phi$. Moreover, we will assume that the reconstruction $\mathcal{R}\,\Phi$ satisfies the equation:

$$\mathscr{L}(\mathcal{R}\,\Phi) = \mathcal{R}(-\Phi^3 + \Xi) = -\mathcal{R}(\Phi^3) + \mathcal{R}\,\Xi. \tag{2.10}$$

In order to have a well-defined reconstruction for the cube we need to impose some conditions on the exponent $\gamma$. By Lemma 2.2 we have that the homogeneity of a term of the form $\mathcal{I}(\tau)$ in our regularity structure satisfies $|\mathcal{I}(\tau)| \geqslant |\mathcal{I}(\Xi)|$. Also, by [FH20, Theorem 14.5] we know that the cube $\Phi^3$ is a modelled distribution in $\mathcal{D}^\eta(\Gamma)$ with $\eta = \gamma + 2|\mathcal{I}(\Xi)| = \gamma + 2s - 3 - 2\kappa$. If we want to have a well-defined reconstruction of the cube then we need to impose the condition $\gamma + 2s - 3 > 0$ and $0 < \kappa \ll 1$. Since $s \in (3/4, 1)$ then $3 - 2s \in (1, 3/2)$, and therefore it is enough to consider $\Phi \in \mathcal{D}^\gamma$ with $\gamma \in (3 - 2s, 2s) \subset (1, 2s)$. Observe that the interval where $\gamma$ has to be is non-empty precisely on the subcritical regime. In this case $P(\Phi) = \langle \mathbf{1}, \Phi \rangle \mathbf{1} + \sum_{i=1}^d \langle \mathbf{X}_j, \Phi \rangle \mathbf{X}_j = \langle \mathbf{1}, \Phi \rangle \mathbf{1} + \langle \mathbf{X}, \Phi \rangle \cdot \mathbf{X}$, where $\mathbf{X} = (\mathbf{X}_1, \ldots, \mathbf{X}_d)$ denotes the vector of abstract linear monomials in space, and similar notation is used for $\langle \mathbf{X}, \Phi \rangle$.

**Remark 2.13.** The main difference between this notion of solution to the one of Hairer (see e.g. [Hai14, Theorem 7.8]) is that in Hairer's one the notion of admissible model allows to encode the mild formulation of (2.10) directly into the polynomial part of equation (2.9). On the other hand, our formulation in terms of the PDE allows us to ignore the initial conditions which appear in Hairer's mild formulation.

We recall that a modelled distribution $\Phi$ solves the algebraic equation (2.9) if and only if it is coherent, i.e., if

$$\Phi(x) = v(x)\mathbf{1} + v_{\mathbf{X}}(x) \cdot \mathbf{X} + \sum_{\tau \in \mathcal{T}} \frac{\Upsilon_x[\tau]}{\tau!} \mathcal{I}(\tau), \qquad (2.11)$$

where $v(x) := \langle \mathbf{1}, \Phi(x) \rangle$, $v_{\mathbf{X}}(x) := \langle \mathbf{X}, \Phi(x) \rangle$, $\tau!$ is the symmetry factor as defined in (A.6), and $\Upsilon_x[\tau] := \Upsilon[\tau](P(\Phi(x))) = \Upsilon[\tau](v(x), v_{\mathbf{X}}(x))$ where $\Upsilon := \Upsilon^F$ is the coherence map for the particular choice of non-linearity $F(\varphi, D\varphi, \ldots) = -\varphi^3$. We refer to Section A.5 for some details on the coherence map. The next lemma gives us an explicit form of these coefficients in terms of $v$ and $v_{\mathbf{X}}$.

**LEMMA 2.14.** *For every $\tau \in \mathcal{T}_{<0}$ such that $\Upsilon_\cdot[\tau] \neq 0$ there exists a constant $c_\tau \in \mathbb{Z} \setminus \{0\}$ such that*

$$\Upsilon_x[\tau] = \begin{cases} c_\tau v(x)^{\mathfrak{m}(\tau)} & \text{if } \mathfrak{n}(\tau) = 0 \\ c_\tau v_{\mathbf{X}_j}(x) & \text{if } \mathfrak{n}(\tau) = e_j \end{cases}.$$

*Moreover, if $\tau \in \mathcal{T}$ has a non-zero edge decoration then $\Upsilon[\tau] = 0$.*

**Proof.** We show the result by induction. The result is true for $\Xi$ by definition (see Definition A.13). By Lemma 2.7 and the recursive definition of $\mathcal{T}$ in (2.1) we have that every tree $\tau \in \mathcal{T}_{<0}$ is of the form $\mathcal{I}(\tau_1), \mathcal{I}_j(\tau), \mathcal{I}(\tau_1)\mathcal{I}(\tau_2)\mathbf{X}^k, \mathcal{I}(\tau_1)\mathcal{I}(\tau_2)\mathcal{I}(\tau_3)$ for $\tau_1, \tau_2, \tau_3 \in \mathcal{T}_{<0}$ and $k \in \{0, e_1, \ldots, e_j\}$ for which we obtain using the definition of the coherence map that

$$\Upsilon[\tau](\varphi, D\varphi, \ldots) = \begin{cases} \Upsilon[\tau_1] D_0 F(\varphi, D\varphi, \ldots) & \text{if } \tau = \mathcal{I}(\tau_1) \\ \Upsilon[\tau_1] D_{e_j} F(\varphi, D\varphi, \ldots) & \text{if } \tau = \mathcal{I}_j(\tau_1) \\ \Upsilon[\tau_1] \Upsilon[\tau_2] (\partial^k D_0^2 F)(\varphi, D\varphi, \ldots) & \text{if } \tau = \mathcal{I}(\tau_1)\mathcal{I}(\tau_2)\mathbf{X}^k \\ \Upsilon[\tau_1] \Upsilon[\tau_2] \Upsilon[\tau_3] D_0^3 F(\varphi, D\varphi, \ldots) & \text{if } \tau = \mathcal{I}(\tau_1)\mathcal{I}(\tau_2)\mathcal{I}(\tau_3) \end{cases}.$$

We refer to Section A.5 for the definitions of $D$ and $\partial$ in the previous expression. The non-linearity in our case is given by $F(\varphi, D\varphi, \ldots) = -\varphi^3$, and since it does not depend on $\partial^k \varphi$ then $D_k F = 0$ for all $k \in \mathbb{N}^{1+d} \setminus \{0\}$, and in particular $\Upsilon[\mathcal{I}_j(\tau)] = 0$. We have then

$$D_0 F(\varphi, D\varphi, \ldots) = -3\varphi^2, \qquad D_0^2 F(\varphi, D\varphi, \ldots) = -6\varphi, \qquad D_0^3 F(\varphi, D\varphi, \ldots) = -6.$$

Moreover, by definition (A.26) of $\partial^{e_j}$ one has that

$$\partial^k D_0^2 F(\varphi, D\varphi, \ldots) = \sum_{k \in \mathbb{N}^{1+d}} (\partial^{k+e_j}\varphi) D_k D_0^2 F(\varphi, D\varphi, \ldots) = -6\partial_j \varphi,$$

and in general for $k \in \mathbb{N}^{1+d}$ one has $\partial^k D_0^2 F(\varphi, D\varphi, \dots) = -6\, \partial^k \varphi$. We conclude that

$$\Upsilon[\tau](\varphi, D\varphi, \dots) = \begin{cases} -3\, \Upsilon[\tau_1]\, \varphi^2 & \text{if} \quad \tau = \mathcal{I}(\tau_1) \\ -6\, \Upsilon[\tau_1]\, \Upsilon[\tau_2]\, \partial^k \varphi & \text{if} \quad \tau = \mathcal{I}(\tau_1)\, \mathcal{I}(\tau_2)\, \boldsymbol{X}^k \\ -6\, \Upsilon[\tau_1]\, \Upsilon[\tau_2]\, \Upsilon[\tau_3] & \text{if} \quad \tau = \mathcal{I}(\tau_1)\, \mathcal{I}(\tau_2)\, \mathcal{I}(\tau_3) \end{cases}, \tag{2.12}$$

and the result follows by a simple induction taking Lemma 2.7 in consideration. $\square$

**Remark 2.15.** In particular the previous lemma tells us that the trees with edge decorations do not appear in the modelled distribution $\Phi$ given by (2.11). Moreover, the coefficients $\Upsilon[\tau]$ follow the same recursive structure as in [CMW23, Lemma 6.8] with the only difference being the combinatorial factors $c_\tau$ which appear in our setting.

**Corollary 2.16.** *For every $\tau \in \mathcal{W}$ the map $x \mapsto \Upsilon_x[\tau]$ is constant.*

We use the notation $\Upsilon[\tau] \sim v, v^2, v_{\boldsymbol{X}}, 1$ to say that $\tau \in \mathcal{T}$ is such that $\Upsilon[\tau]$ is given by the corresponding expression from Lemma 2.14 up to a proportionality constant that depends on $\tau$.

**Lemma 2.17.** *For all $\tau \in \mathcal{T}_{<0}$ and $\mu \in \mathcal{T}^+ \setminus \{\mathbf{1}\}$ such that $|\mu \star \tau| < 0$ we have*

i. *If $\Upsilon[\tau] \sim v^2$, then:*

$$\Upsilon.[\mu \star \tau] = \begin{cases} 2\, c_\tau\, \Upsilon.[\sigma_1]\, \Upsilon.[\sigma_2] & \mu = \mathcal{I}(\sigma_1)\, \mathcal{I}(\sigma_2) \\ 2\, c_\tau\, v(\cdot)\, \Upsilon.[\sigma] & \mu = \mathcal{I}(\sigma) \\ 0 & \text{else} \end{cases}.$$

ii. *If $\Upsilon[\tau] \sim v$, then:*

$$\Upsilon.[\mu \star \tau] = \begin{cases} c_\tau\, \Upsilon.[\sigma] & \mu = \mathcal{I}_j(\sigma) \\ 0 & \text{else} \end{cases}.$$

iii. *If $\Upsilon[\tau] \sim v_{\boldsymbol{X}}$, then:*

$$\Upsilon.[\mu \star \tau] = \begin{cases} 2\, c_\tau\, \Upsilon.[\sigma] & \mu = \mathcal{I}(\sigma) \\ 2\, c_\tau\, v_{\boldsymbol{X}_j}(\cdot) & \mu = \boldsymbol{X}_j \\ 0 & \text{else} \end{cases}.$$

iv. *If $\Upsilon[\tau] \sim 1$, then $\Upsilon.[\mu \star \tau] = 0$.*

**Proof.** The proof is an application of the morphism property between the $\star$-product and the coherence map showed in [BB21, Proposition 2] which for convenience to the reader we include as Lemma A.14 in Appendix A. Let $\mu = \boldsymbol{X}^k \prod_{j \in J} \mathcal{I}_{\rho_j}(\sigma_j) \in \mathcal{T}^+$ such that $\mu \star \tau$ belongs to the span of $\mathcal{T}_{<0}$ and $\Upsilon[\mu \star \tau] \neq 0$. We consider the cases as in the statement:

i. We have that $\boldsymbol{X}^k \star \tau = \uparrow^k \tau \notin \mathcal{T}_{<0}$ by Lemma 2.7 since $(\mathfrak{m}(\uparrow^{e_j} \tau), \mathfrak{n}(\uparrow^{e_j} \tau)) = (2, k)$, and since $|\boldsymbol{X}^k \star \tau| \leqslant |\boldsymbol{X}^k \prod_{j \in J} \mathcal{I}_{\rho_j}(\sigma_j)|$ (all symbols $\mathcal{I}_{\rho_j}(\sigma_j)$ in $\mathcal{T}^+$ have positive homogeneity) we conclude that $\mu$ has no polynomial factor, i.e., $\mu = \prod_{j \in J} \mathcal{I}_{\rho_j}(\sigma_j)$. Since we have the explicit formula $\Upsilon[\tau] = c_\tau \mathcal{X}_0^2$ then $D_k \Upsilon[\tau] = 0$ for $k \neq 0$ which shows that $\rho_j = 0$ for all $j \in J$. Moreover, since $D_0 \Upsilon[\tau] = 2\, c_\tau \mathcal{X}_0, D_0^2 \Upsilon[\tau] = 2\, c_\tau$ and $D_0^n \Upsilon[\tau] = 0$ for all $n \geqslant 3$, we conclude that $|J| \leqslant 2$. In particular

$$\Upsilon[\mathcal{I}(\sigma) \star \tau] = 2\, c_\tau \mathcal{X}_0\, \Upsilon[\sigma] \quad \text{and} \quad \Upsilon[(\mathcal{I}(\sigma_1)\, \mathcal{I}(\sigma_2)) \star \tau] = 2\, c_\tau\, \Upsilon[\sigma_1]\, \Upsilon[\sigma_2].$$

ii. We have that $\boldsymbol{X}^k \star \tau = \uparrow^k \tau \notin \mathcal{T}_{<0}$ by Lemma 2.7 since $|\mathfrak{n}(\uparrow^k \tau)| > 1$, and we can conclude that $\mu = \prod_{j \in J} \mathcal{I}_{\rho_j}(\sigma_j)$. Since $\Upsilon[\tau] = c_\tau \mathcal{X}_j$ then $D_k \Upsilon[\tau] = 0$ for $k \neq e_j$, $D_{e_j} \Upsilon[\tau] = c_\tau$ and $D_{e_j}^n \Upsilon[\tau] = 0$ for all $n \geq 2$. We conclude that $\mu = \mathcal{I}_j(\sigma)$ and in this case $\Upsilon[\mathcal{I}_j(\sigma) \star \tau] = c_\tau \Upsilon[\sigma]$.

iii. In this case, by Lemma 2.7 we can only raise decoration by some $e_j$ which shows that $k \in \{0\} \cup \{e_j\}_{j=1}^d$. Since $\Upsilon[\tau] = 2 c_\tau \mathcal{X}_0$ then $D_k \Upsilon[\tau] = 0$ for all $k \neq 0$, $D_0 \Upsilon[\tau] = 2 c_\tau$ and $D_0^n \Upsilon[\tau] = 0$ for all $n \geq 2$. We conclude that $\mu = \boldsymbol{X}^k \prod_{j \in J} \mathcal{I}(\sigma_j)$ with $|J| \leq 1$ and $|k| \leq 1$, and then

$$\begin{aligned}
\Upsilon[\boldsymbol{X}_j \star \tau] &= \partial^{e_j} \Upsilon[\tau] = \partial^{e_j}(2 c_\tau \mathcal{X}_0) = 2 c_\tau \mathcal{X}_j, \\
\Upsilon[\mathcal{I}(\sigma) \star \tau] &= \Upsilon[\sigma] \, D_0 \Upsilon[\tau] = 2 c_\tau \Upsilon[\sigma], \\
\Upsilon[(\boldsymbol{X}_j \mathcal{I}^+(\sigma)) \star \tau] &= \Upsilon[\sigma] \, \partial^{e_j}(D_0 \Upsilon[\tau]) = \Upsilon[\sigma] \, \partial^{e_j}(2 c_\tau) = 0.
\end{aligned}$$

iv. Since $\Upsilon[\tau] = c_\tau$ is constant, all derivatives $D_k$ or $\partial^k$ are zero which shows the result. □

### 2.4. Da Prato-Debussche Remainder

A consequence of the action of structure group acting trivially on $\mathcal{W}$ and $\mathcal{I}(\mathcal{W})$ (Lemma 2.6) is that it allows us to perform the Da Prato-Debussche trick, which consists in subtracting the worst terms of $\Phi$ and work instead with the remainder $V$ defined as:

$$V(x) := \Phi(x) - \sum_{\tau \in \mathcal{W}} \frac{\Upsilon[\tau]}{\tau!} \mathcal{I}(\tau) = v(x) \mathbf{1} + v_{\boldsymbol{X}}(x) \cdot \boldsymbol{X} + \sum_{\tau \in \mathcal{V}_{0,\gamma}} \frac{\Upsilon_x[\tau]}{\tau!} \mathcal{I}(\tau), \quad (2.13)$$

where for $\tau \in \mathcal{W}$ we wrote $\Upsilon[\tau]$ without an evaluation on $x$ since by Lemma 2.14 it is constant. The advantage of working with the remainder $V$ is that, by definition of the set $\mathcal{V}$, all trees in $\mathcal{I}(\mathcal{V})$ have positive homogeneity, and therefore $V$ takes values in a function-like sector of the regularity structure ($\mathcal{I}(\mathcal{V}) \cup \mathcal{P}$ spans a sector by Lemma 2.8). Since the structure group acts trivially on $\mathcal{I}(\mathcal{W})$ then $\mathcal{I}(\tau)$ is a modelled distribution of any order, i.e., $\mathcal{I}(\tau) \in \mathcal{D}^\gamma$ for all $\gamma > 0$. Since $\mathcal{D}^\gamma$ is a linear space this implies that $V$ is also a modelled distribution of the same order as $\Phi$, i.e., $V \in \mathcal{D}^\gamma$.

As a consequence of $V$ being function-like we have that the reconstruction $\mathcal{R}V$ is given by the coefficient at $\mathbf{1}$, i.e., $\mathcal{R}V = v$ (see [Hai14, Proposition 3.28]). The following result characterises the relationship between $\Phi$ and $V$.

**Lemma 2.18.** *Given $\Phi \in \mathcal{D}^\gamma$ which satisfies the algebraic equation (2.9) and whose reconstruction satisfies (2.10), then the remainder $V \in \mathcal{D}^\gamma$ defined as in (2.13) satisfies the algebraic equation in $\mathcal{D}^\gamma$*

$$\begin{aligned}
V = &-3 \sum_{\tau \in \mathcal{W}} \frac{\Upsilon[\tau]}{\tau!} \mathcal{I}(\mathcal{I}(\tau) V_{\beta_\tau}^2) - 3 \sum_{\tau_1, \tau_2 \in \mathcal{W}} \frac{\Upsilon[\tau_1] \Upsilon[\tau_2]}{\tau_1! \, \tau_2!} \mathcal{I}(\mathcal{I}(\tau_1) \mathcal{I}(\tau_2) V_{\beta_{\tau_1}, \beta_{\tau_2}}) \\
&+ \sum_{\tau \in \partial \mathcal{W} \cap \mathcal{V}_{0,\gamma}} \frac{\Upsilon[\tau]}{\tau!} \mathcal{I}(\tau) + v(x) \mathbf{1} + v_{\boldsymbol{X}}(x) \cdot \boldsymbol{X},
\end{aligned} \quad (2.14)$$

*and its reconstruction $v = \mathcal{R}V$ satisfies the PDE*

$$\begin{aligned}
\mathscr{L} v = &-v^3 - 3 \sum_{\tau \in \mathcal{W}} \frac{\Upsilon[\tau]}{\tau!} \mathcal{R}(\mathcal{I}(\tau) V^2) \\
&- 3 \sum_{\tau_1, \tau_2 \in \mathcal{W}} \frac{\Upsilon[\tau_1] \Upsilon[\tau_2]}{\tau_1! \, \tau_2!} \mathcal{R}(\mathcal{I}(\tau_1) \mathcal{I}(\tau_2) V) + \sum_{\tau \in \partial \mathcal{W}} \frac{\Upsilon[\tau]}{\tau!} \boldsymbol{\Pi} \tau,
\end{aligned} \quad (2.15)$$

*where $\partial \mathcal{W} := \{\mathcal{I}(\tau_1) \mathcal{I}(\tau_2) \mathcal{I}(\tau_3) \notin \mathcal{W} : \tau_1, \tau_2, \tau_3 \in \mathcal{W}\}$, $\beta_{\tau_1, \tau_2} := \varepsilon - |\mathcal{I}(\tau_1) \mathcal{I}(\tau_2)|$, $\beta_\tau := \varepsilon - |\mathcal{I}(\tau)| \leq 1$ for some small $\varepsilon > 0$ which satisfies (2.18) and (2.19), and $V_\beta := \mathcal{Q}_{<\beta}(V)$, $V_\beta^2 := \mathcal{Q}_{<\beta}(V^2)$ are truncations of $V$ and $V^2$, i.e., $\mathcal{Q}_{<\beta}: T \to T$ is the projection to trees of homogeneity less than $\beta$.*

*Reciprocally, given $V \in \mathcal{D}^\gamma$ which is function-like, satisfies the algebraic equation* (2.14) *and its reconstruction satisfies* (2.15) *then $\Phi$ defined by* (2.13) *satisfies the algebraic equation* (2.9) *and its reconstruction satisfies* (2.10).

**Proof.** We only prove the first part of the statement, the second one being analogous. Let $P(\Phi) = v\,\mathbf{1} + v_{\boldsymbol{X}} \cdot \boldsymbol{X}$ be the polynomial part of $\Phi$ (and of $V$). Then

$$\begin{aligned}
V - P(\Phi) &= \Phi - \sum_{\tau \in \mathcal{W}} \frac{\Upsilon[\tau]}{\tau!} \mathcal{I}(\tau) \\
&= -\mathcal{I}(\Phi^3) - \sum_{\tau \in \mathcal{W}\setminus\{\Xi\}} \frac{\Upsilon[\tau]}{\tau!} \mathcal{I}(\tau) \\
&= -\mathcal{I}\!\left(\!\left(V + \sum_{\tau \in \mathcal{W}} \frac{\Upsilon[\tau]}{\tau!} \mathcal{I}(\tau)\right)^{\!3}\right) - \sum_{\tau \in \mathcal{W}\setminus\{\Xi\}} \frac{\Upsilon[\tau]}{\tau!} \mathcal{I}(\tau) \\
&= -\mathcal{I}(V^3) - 3 \sum_{\tau \in \mathcal{W}} \frac{\Upsilon[\tau]}{\tau!} \mathcal{I}(\mathcal{I}(\tau)\,V^2) - 3 \sum_{\tau_1, \tau_2 \in \mathcal{W}} \frac{\Upsilon[\tau_1]\,\Upsilon[\tau_2]}{\tau_1!\,\tau_2!} \mathcal{I}(\mathcal{I}(\tau_1)\,\mathcal{I}(\tau_2)\,V) \\
&\quad - \sum_{\tau_1, \tau_2, \tau_3 \in \mathcal{W}} \frac{\Upsilon[\tau_1]\,\Upsilon[\tau_2]\,\Upsilon[\tau_3]}{\tau_1!\,\tau_2!\,\tau_3!} \mathcal{I}(\mathcal{I}(\tau_1)\,\mathcal{I}(\tau_2)\,\mathcal{I}(\tau_3)) - \sum_{\tau \in \mathcal{W}\setminus\{\Xi\}} \frac{\Upsilon[\tau]}{\tau!} \mathcal{I}(\tau). \quad (2.16)
\end{aligned}$$

The first thing to observe is that since $V$ is function-like, also $V^3$ is and therefore $\mathcal{I}(V^3)$ takes values in trees of homogeneity $\geqslant 2s$, and since $\gamma < 2s$ we have that $\mathcal{I}(V^3) = 0$ in $\mathcal{D}^\gamma$. On the other hand, by Lemma 2.7 for each $\tau \in \mathcal{W} \setminus \{\Xi\}$ there exists $\tau_1, \tau_2, \tau_3 \in \mathcal{W}$ such that $\tau = \mathcal{I}(\tau_1)\,\mathcal{I}(\tau_2)\,\mathcal{I}(\tau_3)$ and using coherence of $\Upsilon$ as in (2.12) we have that $\Upsilon[\tau] = -6\,\Upsilon[\tau_1]\,\Upsilon[\tau_2]\,\Upsilon[\tau_3]$. Recalling the definition (A.8) of the symmetry factor $\tau!$ we have that $\tau!\,\delta(\tau_1, \tau_2, \tau_3) = 3!\,\tau_1!\,\tau_2!\,\tau_3!$, and putting things together we conclude that

$$\begin{aligned}
-\sum_{\tau \in \mathcal{W}\setminus\{\Xi\}} \frac{\Upsilon[\tau]}{\tau!} \mathcal{I}(\tau) &= \sum_{\tau \in \mathcal{W}\setminus\{\Xi\}} \frac{6\,\Upsilon[\tau_1]\,\Upsilon[\tau_2]\,\Upsilon[\tau_3]}{3!\,\tau_1!\,\tau_2!\,\tau_3!}\,\delta(\tau_1,\tau_2,\tau_3)\,\mathcal{I}(\mathcal{I}(\tau_1)\,\mathcal{I}(\tau_2)\,\mathcal{I}(\tau_3)) \\
&= \sum_{\substack{\tau_1,\tau_2,\tau_3 \in \mathcal{W} \\ \mathcal{I}(\tau_1)\mathcal{I}(\tau_2)\mathcal{I}(\tau_3) \in \mathcal{W}}} \frac{\Upsilon[\tau_1]\,\Upsilon[\tau_2]\,\Upsilon[\tau_3]}{\tau_1!\,\tau_2!\,\tau_3!} \mathcal{I}(\mathcal{I}(\tau_1)\,\mathcal{I}(\tau_2)\,\mathcal{I}(\tau_3)),
\end{aligned}$$

where the factor $\delta(\tau_1, \tau_2, \tau_3)$ is precisely what allows us to re-index the sum. This allows two reduce the last two terms of (2.16) into

$$\begin{aligned}
&-\sum_{\tau_1,\tau_2,\tau_3 \in \mathcal{W}} \frac{\Upsilon[\tau_1]\,\Upsilon[\tau_2]\,\Upsilon[\tau_3]}{\tau_1!\,\tau_2!\,\tau_3!} \mathcal{I}(\mathcal{I}(\tau_1)\,\mathcal{I}(\tau_2)\,\mathcal{I}(\tau_3)) - \sum_{\tau \in \mathcal{W}\setminus\{\Xi\}} \frac{\Upsilon[\tau]}{\tau!} \mathcal{I}(\tau) \\
&= -\sum_{\substack{\tau_1,\tau_2,\tau_3 \in \mathcal{W} \\ \mathcal{I}(\tau_1)\mathcal{I}(\tau_2)\mathcal{I}(\tau_3) \notin \mathcal{W}}} \frac{\Upsilon[\tau_1]\,\Upsilon[\tau_2]\,\Upsilon[\tau_3]}{\tau_1!\,\tau_2!\,\tau_3!} \mathcal{I}(\mathcal{I}(\tau_1)\,\mathcal{I}(\tau_2)\,\mathcal{I}(\tau_3)) \\
&= \sum_{\tau \in \delta\mathcal{W}} \frac{\Upsilon[\tau]}{\tau!} \mathcal{I}(\tau), \quad (2.17)
\end{aligned}$$

where we used the same arguments in the last step with the additional observation that trees of the form $\tau = \mathcal{I}(\tau_1)\,\mathcal{I}(\tau_2)\,\mathcal{I}(\tau_3) \notin \mathcal{V}_{0,\gamma}$ satisfy $|\mathcal{I}(\tau)| > \gamma$ and therefore $\mathcal{I}(\tau) = 0$ in $\mathcal{D}^\gamma$.

To get the truncations in the other terms, observe that by assumption we have that $\gamma > 3 - 2s$, which guarantees by Lemma 2.2 that

$$\gamma + |\mathcal{I}(\tau_1)\,\mathcal{I}(\tau_2)| \geqslant \gamma + |\mathcal{I}(\Xi)\,\mathcal{I}(\Xi)| > 3 - 2s + 2\left(s - \frac{3}{2} - \kappa\right) > 0, \quad \forall\, \tau_1, \tau_2 \in \mathcal{W},$$

and since $\mathcal{W}$ is finite by Lemma 2.2 then we can choose $0 < \varepsilon \ll 1$ which satisfies

$$0 < \varepsilon < \min\{\gamma + |\mathcal{I}(\tau_1)\,\mathcal{I}(\tau_2)| : \tau_1, \tau_2 \in \mathcal{W}\}. \quad (2.18)$$

Now, since $\gamma < 2s$ by assumption and $\varepsilon > 0$ then we have that $\gamma - 2s < \varepsilon$ and therefore $\gamma - 2s - |\mathcal{I}(\tau_1)\mathcal{I}(\tau_2)| < \beta_{\tau_1,\tau_2}$ which implies that $\mathcal{Q}_{<\gamma - 2s - |\mathcal{I}(\tau_1)\mathcal{I}(\tau_2)|} \circ \mathcal{Q}_{<\beta_{\tau_1,\tau_2}} = \mathcal{Q}_{<\gamma - 2s - |\mathcal{I}(\tau_1)\mathcal{I}(\tau_2)|}$ and therefore $\mathcal{Q}_{<\gamma}(\mathcal{I}(\mathcal{I}(\tau_1)\mathcal{I}(\tau_2)V)) = \mathcal{Q}_{<\gamma}(\mathcal{I}(\mathcal{I}(\tau_1)\mathcal{I}(\tau_2)V_{\beta_{\tau_1,\tau_2}}))$, or equivalently we have that $\mathcal{I}(\mathcal{I}(\tau_1)\mathcal{I}(\tau_2)V) = \mathcal{I}(\mathcal{I}(\tau_1)\mathcal{I}(\tau_2)V_{\beta_{\tau_1,\tau_2}})$ in $\mathcal{D}^\gamma$. For the other terms, since $s > 3/4$ we have

$$1 + |\mathcal{I}(\tau)| \geq 1 + |\mathcal{I}(\Xi)| = 1 + s - \frac{3}{2} - \kappa = s - \frac{1}{2} - \kappa > 0,$$

and therefore we can also choose $0 < \varepsilon \ll 1$ such that

$$0 < \varepsilon < \min\{1 + |\mathcal{I}(\tau)| : \tau \in \mathcal{W}\}. \tag{2.19}$$

This implies that $\beta_\tau := \varepsilon - |\mathcal{I}(\tau)| < 1$, and arguing in the same way as before we can conclude that $\mathcal{I}(\mathcal{I}(\tau)V) = \mathcal{I}(\mathcal{I}(\tau)V_{\beta_\tau})$ in $\mathcal{D}^\gamma$. At last, since we are working in $\mathcal{D}^\gamma$ then the last sum in (2.17) gets truncated to $\delta\mathcal{W} \cap \mathcal{V}_{0,\gamma}$ which concludes the proof of (2.14).

Now we show that $v = \mathcal{R}V$ satisfies the right equation. Since $\mathcal{R}\Phi$ satisfies (2.10) we have

$$\begin{aligned}
\mathscr{L}(\mathcal{R}V) &= \mathscr{L}(\mathcal{R}\Phi) - \sum_{\tau \in \mathcal{W}} \frac{\Upsilon[\tau]}{\tau!} \mathscr{L}(\mathcal{R}(\mathcal{I}(\tau))) \\
&= -\mathcal{R}(\Phi^3) + \mathcal{R}\Xi - \sum_{\tau \in \mathcal{W}} \frac{\Upsilon[\tau]}{\tau!} \mathscr{L}(\mathbf{\Pi}(\mathcal{I}(\tau))) \\
&= -\mathcal{R}(V^3) - 3\sum_{\tau \in \mathcal{W}} \frac{\Upsilon[\tau]}{\tau!} \mathcal{R}(\mathcal{I}(\tau)V^2) - 3\sum_{\tau_1,\tau_2 \in \mathcal{W}} \frac{\Upsilon[\tau_1]\Upsilon[\tau_2]}{\tau_1!\tau_2!} \mathcal{R}(\mathcal{I}(\tau_1)\mathcal{I}(\tau_2)V) \\
&\quad - \sum_{\tau_1,\tau_2,\tau_3 \in \mathcal{W}} \frac{\Upsilon[\tau_1]\Upsilon[\tau_2]\Upsilon[\tau_3]}{\tau_1!\tau_2!\tau_3!} \mathcal{R}(\mathcal{I}(\tau_1)\mathcal{I}(\tau_2)\mathcal{I}(\tau_3)) + \mathbf{\Pi}\Xi - \sum_{\tau \in \mathcal{W}} \frac{\Upsilon[\tau]}{\tau!} \mathbf{\Pi}\tau \\
&= -\mathcal{R}(V^3) - 3\sum_{\tau \in \mathcal{W}} \frac{\Upsilon[\tau]}{\tau!} \mathcal{R}(\mathcal{I}(\tau)V^2) - 3\sum_{\tau_1,\tau_2 \in \mathcal{W}} \frac{\Upsilon[\tau_1]\Upsilon[\tau_2]}{\tau_1!\tau_2!} \mathcal{R}(\mathcal{I}(\tau_1)\mathcal{I}(\tau_2)V) \\
&\quad + \sum_{\tau_1,\tau_2,\tau_3 \in \mathcal{W}} \frac{\Upsilon[\mathcal{I}(\tau_1)\mathcal{I}(\tau_2)\mathcal{I}(\tau_3)]}{(\mathcal{I}(\tau_1)\mathcal{I}(\tau_2)\mathcal{I}(\tau_3))!\,\delta(\tau_1,\tau_2,\tau_3)} \mathcal{R}(\mathcal{I}(\tau_1)\mathcal{I}(\tau_2)\mathcal{I}(\tau_3)) \\
&\quad - \sum_{\tau \in \mathcal{W}\setminus\{\Xi\}} \frac{\Upsilon[\tau]}{\tau!} \mathbf{\Pi}\tau \\
&= -\mathcal{R}(V^3) - 3\sum_{\tau \in \mathcal{W}} \frac{\Upsilon[\tau]}{\tau!} \mathcal{R}(\mathcal{I}(\tau)V^2) - 3\sum_{\tau_1,\tau_2 \in \mathcal{W}} \frac{\Upsilon[\tau_1]\Upsilon[\tau_2]}{\tau_1!\tau_2!} \mathcal{R}(\mathcal{I}(\tau_1)\mathcal{I}(\tau_2)V) \\
&\quad + \sum_{\substack{\tau = \mathcal{I}(\tau_1)\mathcal{I}(\tau_2)\mathcal{I}(\tau_3) \\ \tau_1,\tau_2,\tau_3 \in \mathcal{W}}} \frac{\Upsilon[\tau]}{\tau!} \mathbf{\Pi}\tau - \sum_{\tau_1,\tau_2 \in \mathcal{W}\setminus\{\Xi\}} \frac{\Upsilon[\tau]}{\tau!} \mathbf{\Pi}\tau. \\
&= -\mathcal{R}(V^3) - 3\sum_{\tau \in \mathcal{W}} \frac{\Upsilon[\tau]}{\tau!} \mathcal{R}(\mathcal{I}(\tau)V^2) - 3\sum_{\tau_1,\tau_2 \in \mathcal{W}} \frac{\Upsilon[\tau_1]\Upsilon[\tau_2]}{\tau_1!\tau_2!} \mathcal{R}(\mathcal{I}(\tau_1)\mathcal{I}(\tau_2)V) \\
&\quad + \sum_{\tau \in \partial\mathcal{W}} \frac{\Upsilon[\tau]}{\tau!} \mathbf{\Pi}\tau,
\end{aligned}$$

where we used that for tress $\tau \in \mathcal{T}$ where the structure group acts trivially $\mathcal{R}\tau = \mathbf{\Pi}\tau$, and that $\mathscr{L}(\mathbf{\Pi}_x\mathcal{I}(\tau)) = \mathscr{L}(\mathbf{\Pi}\mathcal{I}(\tau)) = \mathbf{\Pi}\tau$ by our assumption of weak admissibility of the model for $\mathscr{L}$ (Definition 2.11). Moreover, since $V$ is function-like then $V^3$ is also function-like, and therefore its reconstruction satisfies $\mathcal{R}(V^3) = \langle \mathbf{1}, V^3 \rangle = v^3$. At last, observe that by construction for all $\tau, \tau_1, \tau_2 \in \mathcal{W}$ we have $\mathcal{I}(\tau_1)\mathcal{I}(\tau_2)V_{\beta_{\tau_1,\tau_2}} = \mathcal{Q}_{<\varepsilon}(\mathcal{I}(\tau_1)\mathcal{I}(\tau_2)V) \in \mathcal{D}^\varepsilon$ and $\mathcal{I}(\tau)V^2_{\beta_\tau} = \mathcal{Q}_{<\varepsilon}(\mathcal{I}(\tau)V^2) \in \mathcal{D}^\varepsilon$, i.e., modelled distributions of positive order $\varepsilon$. By uniqueness of the reconstruction operator on modelled distributions of positive order we have $\mathcal{R}(\mathcal{I}(\tau)V^2) = \mathcal{R}(\mathcal{I}(\tau)V^2_{\beta_\tau})$ and $\mathcal{R}(\mathcal{I}(\tau_1)\mathcal{I}(\tau_2)V) = \mathcal{R}(\mathcal{I}(\tau_1)\mathcal{I}(\tau_2)V_{\beta_{\tau_1,\tau_2}})$, which concludes the proof of (2.15). □

From Lemma 2.18 we see that the equation that $v$ solves involves reconstructions of the modelled distributions $V\mathcal{I}(\tau_1)\mathcal{I}(\tau_2)$ and $V^2\mathcal{I}(\tau)$ and in order to bound this terms with the Reconstruction theorem we need to understand the change of base point formulas for their truncations. This is done in the next lemmas.

**LEMMA 2.19.** *For any $\beta \in (0,\gamma]$ consider $V_\beta$ the truncation of $V$ given by*

$$V_\beta(x) = v(x)\mathbf{1} + \mathbb{1}_{\beta>1} v_{\boldsymbol{X}}(x)\cdot \boldsymbol{X} + \sum_{\sigma \in \mathcal{V}_{0,\beta}} \frac{\Upsilon_x[\sigma]}{\sigma!}\mathcal{I}(\sigma). \tag{2.20}$$

*Then for any $\Gamma \in G$ we have*

$$\langle \mathbf{1}, \Gamma V_\beta(x)\rangle = v(x) + \mathbb{1}_{\beta>1} v_{\boldsymbol{X}}(x)\cdot \gamma(\boldsymbol{X}) + \sum_{\sigma \in \mathcal{V}_{0,\beta}} \frac{\Upsilon_x[\sigma]}{\sigma!}\gamma(\mathcal{I}(\sigma)), \tag{2.21}$$

$$\langle \boldsymbol{X}_j, \Gamma V_\beta(x)\rangle = \mathbb{1}_{\beta>1} v_{\boldsymbol{X}}(x) + \sum_{\sigma \in \mathcal{V}_{1,\beta}} \frac{\Upsilon_x[\sigma]}{\sigma!}\gamma(\mathcal{I}_j(\sigma)), \tag{2.22}$$

$$\langle \mathcal{I}(\tau), \Gamma V_\beta(x)\rangle = \sum_{\substack{\mu \in \mathcal{T} \\ \mu \star \tau \in \mathcal{V}_{0,\beta}}} \frac{\Upsilon_x[\mu \star \tau]}{\mu!}\gamma(\mu). \tag{2.23}$$

**Proof.** We start by writing $V$ in the following form

$$V(x) = \sum_{k \in \mathbb{N}^{1+d}} \frac{v_k(x)}{k!} \boldsymbol{X}^k + \sum_{\sigma \in \mathcal{V}} \frac{\Upsilon_x[\sigma]}{\sigma!}\mathcal{I}(\sigma),$$

and by (A.24) we have that $\mathcal{I}(\mathcal{V}) \cup \mathcal{P}$ spans a sector we only need to look at these coefficients for the change of base point. Using the linearity of the representation (A.23) we have for any $\Gamma \in G$ that

$$\langle \tau, \Gamma V(x)\rangle = \sum_{\mu \in \mathcal{T}^+} \langle \mu \star \tau, V(x)\rangle \frac{\gamma(\mu)}{\mu!},$$

which allows us to deduce the form of the coefficient at $\mathbf{1}$ as

$$\begin{aligned}
\langle \mathbf{1}, \Gamma V(x)\rangle &= \sum_{\mu \in \mathcal{T}^+} \langle \mu, V(x)\rangle \frac{\gamma(\mu)}{\mu!} \\
&= \sum_{k \in \mathbb{N}^{1+d}} \langle \boldsymbol{X}^k, V(x)\rangle \frac{\gamma(\boldsymbol{X}^k)}{\boldsymbol{X}^k!} + \sum_{\sigma \in \mathcal{V}} \langle \mathcal{I}(\sigma), V(x)\rangle \frac{\gamma(\mathcal{I}(\sigma))}{\mathcal{I}(\sigma)!} \\
&= \sum_{k \in \mathbb{N}^{1+d}} v_k(x) \frac{\gamma(\boldsymbol{X}^k)}{k!} + \sum_{\sigma \in \mathcal{V}} \frac{\Upsilon_x[\sigma]}{\sigma!}\gamma(\mathcal{I}(\sigma)).
\end{aligned}$$

To compute the coefficient at polynomial symbols first we see that to obtain a polynomial from $\mu \star \boldsymbol{X}^k$ necessarily $\mu$ needs to be a polynomial, and to obtain a planted tree from $\mu \star \boldsymbol{X}^k$ we have that $\mu$ need to be a planted tree $\mathcal{I}_\rho(\sigma)$. To see which $\rho$ are allowed, we recall that by (A.21)

$$\mathcal{I}_\rho(\sigma) \star \boldsymbol{X}^n = \sum_{m \in \mathbb{N}^{1+d}} \binom{n}{m}\mathcal{I}_{\rho-m}(\sigma)\,\boldsymbol{X}^{n-m},$$

and the only planted term appearing in that sum occurs when the decoration is $\boldsymbol{X}^{n-m} = \boldsymbol{1}$, i.e., for $m = n$ and in which case is given of $\mathcal{I}_{\rho-n}(\sigma)$. Since we are only considering edges without decorations in $\mathcal{T}$ we conclude that $\rho = n$. With this in mind we can write

$$\begin{aligned}
\langle \boldsymbol{X}^n, \Gamma V(x) \rangle &= \sum_{\mu \in \mathcal{T}^+} \langle \mu \star \boldsymbol{X}^n, V(x) \rangle \frac{\gamma(\mu)}{\mu!} \\
&= \sum_{k \in \mathbb{N}^{1+d}} \langle \boldsymbol{X}^k \star \boldsymbol{X}^n, V(x) \rangle \frac{\gamma(\boldsymbol{X}^k)}{k!} + \sum_{\sigma \in \mathcal{V}} \sum_{\rho \in \mathbb{N}^{1+d}} \langle \mathcal{I}_\rho(\sigma) \star \boldsymbol{X}^n, V(x) \rangle \frac{\gamma(\mathcal{I}_\rho(\sigma))}{\mathcal{I}_\rho(\sigma)!} \\
&= \sum_{k \in \mathbb{N}^{1+d}} \langle \boldsymbol{X}^{k+n}, V(x) \rangle \frac{\gamma(\boldsymbol{X})^k}{k!} + \sum_{\substack{\sigma \in \mathcal{V} \\ |\mathcal{I}_n(\sigma)| > 0}} \langle \mathcal{I}_n(\sigma) \star \boldsymbol{X}^n, V(x) \rangle \frac{\gamma(\mathcal{I}_n(\sigma))}{\mathcal{I}_n(\sigma)!} \\
&= \sum_{k \in \mathbb{N}^{1+d}} v_{k+n}(x) \frac{\gamma(\boldsymbol{X})^k}{k!} + \sum_{\substack{\sigma \in \mathcal{V} \\ |\mathcal{I}(\sigma)| > n}} \frac{\langle \mathcal{I}(\sigma), V(x) \rangle}{\mathcal{I}(\sigma)!} \gamma(\mathcal{I}_n(\sigma)) \\
&= \sum_{k \in \mathbb{N}^{1+d}} v_{k+n}(x) \frac{\gamma(\boldsymbol{X})^k}{k!} + \sum_{\substack{\sigma \in \mathcal{V} \\ |\mathcal{I}(\sigma)| > n}} \frac{\Upsilon_x[\sigma]}{\sigma!} \gamma(\mathcal{I}_n(\sigma)).
\end{aligned}$$

At last, for planted trees we have:

$$\begin{aligned}
\langle \mathcal{I}(\tau), \Gamma V(x) \rangle &= \sum_{\mu \in \mathcal{T}^+} \langle \mu \star \mathcal{I}(\tau), V(x) \rangle \frac{\gamma(\mu)}{\mu!} = \sum_{\substack{\mu \in \mathcal{T}^+ \\ \mu \star \tau \in \mathcal{V}_{0,\beta}}} \langle \mathcal{I}(\mu \star \tau), V(x) \rangle \frac{\gamma(\mu)}{\mu!} \\
&= \sum_{\substack{\mu \in \mathcal{T}^+ \\ \mu \star \tau \in \mathcal{V}_{0,\beta}}} \Upsilon_x[\mu \star \tau] \frac{\gamma(\mu)}{\mu!},
\end{aligned}$$

Where we have used that since $V$ takes values in planted trees, and therefore $\mu \star \mathcal{I}(\tau)$ cannot graft (or decorate) the root of $\mathcal{I}(\tau)$ and this is equivalent to only consider the part of the grafting (or decorating) that occurs in $\tau$, which can be equivalently written as $\mathcal{I}(\mu \star \tau)$, i.e., we have $\langle \mu \star \mathcal{I}(\tau), \mathcal{I}(\sigma) \rangle = \langle \mathcal{I}(\mu \star \tau), \mathcal{I}(\sigma) \rangle$. The result now follows by (2.20). □

**Corollary 2.20.** *For any $\beta \in (0, \gamma]$ and $\Gamma \in G$ we have*

$$\langle \boldsymbol{1}, \Gamma_{yx} V_\beta(x) \rangle = (\Pi_x V_\beta(x))(y). \tag{2.24}$$

**Proof.** By (2.6) $\gamma_{yx}(\mathcal{I}(\tau)) = (\Pi_x \mathcal{I}(\tau))(y)$ for all $\tau \in \mathcal{V}$ and therefore applying $\Pi_x(\cdot)(y)$ to (2.20) and comparing it to (2.21) from Lemma 2.19 we conclude the result. □

**Analysis of the square**

For $\beta \in (0, 1]$ the modelled distribution $V_\beta$ simplifies to

$$V_\beta(x) = v(x) \boldsymbol{1} + \sum_{\sigma \in \mathcal{V}_{0,\beta}} \frac{\Upsilon_x[\sigma]}{\sigma!} \mathcal{I}(\sigma),$$

and therefore we can write $V_\beta^2$ as

$$\begin{aligned}
V_\beta^2(x) &= v^2(x) \boldsymbol{1} + 2 v(x) \sum_{\sigma \in \mathcal{V}_{0,\beta}} \frac{\Upsilon_x[\sigma]}{\sigma!} \mathcal{I}(\sigma) + \sum_{\sigma_1, \sigma_2 \in \mathcal{V}_{0,\beta}} \frac{\Upsilon_x[\sigma_1] \Upsilon_x[\sigma_2]}{\sigma_1! \sigma_2!} \mathcal{Q}_{<\beta}(\mathcal{I}(\sigma_1) \mathcal{I}(\sigma_2)) \\
&= v^2(x) \boldsymbol{1} + \sum_{\substack{\mathcal{I}(\sigma) \in \mathcal{T}^+ \\ |\mathcal{I}(\sigma)| < \beta}} \frac{2 v(x) \Upsilon_x[\sigma]}{\sigma!} \mathcal{I}(\sigma) + \sum_{\substack{\mathcal{I}(\sigma_1) \mathcal{I}(\sigma_2) \in \mathcal{T}^+ \\ |\mathcal{I}(\sigma_1) \mathcal{I}(\sigma_2)| < \beta}} \frac{2 \Upsilon_x[\sigma_1] \Upsilon_x[\sigma_2]}{(\mathcal{I}(\sigma_1) \mathcal{I}(\sigma_2))!} \mathcal{I}(\sigma_1) \mathcal{I}(\sigma_2). \tag{2.25}
\end{aligned}$$

The next result relates the action of the structure group to $V^2$ to that one of $V$.

**LEMMA 2.21**. *For any $\beta \in (0,1]$ and $\Gamma \in G$ we have that all non-zero coefficients of $\Gamma V_\beta^2$ are*

$$3\langle \mathbf{1}, \Gamma V_\beta^2(x)\rangle = \langle \mathcal{I}(\mathcal{I}(\Xi)), \Gamma V_{\beta+|\mathcal{I}(\mathcal{I}(\Xi))|}(x)\rangle, \qquad (2.26)$$

$$3\langle \mathcal{I}(\sigma), \Gamma V_\beta^2(x)\rangle = \langle \mathcal{I}(\sigma)\mathcal{I}(\Xi), \Gamma V_{\beta+|\mathcal{I}(\mathcal{I}(\Xi))|}(x)\rangle, \qquad (2.27)$$

$$3\langle \mathcal{I}(\sigma_1)\mathcal{I}(\sigma_2), \Gamma V_\beta^2(x)\rangle = \langle \mathcal{I}(\sigma_1)\mathcal{I}(\sigma_2)\mathcal{I}(\Xi), \Gamma V_{\beta+|\mathcal{I}(\mathcal{I}(\Xi))|}(x)\rangle, \qquad (2.28)$$

*for $\sigma, \sigma_1, \sigma_2 \in \mathcal{V}_{0,\beta}$ with $|\mathcal{I}(\sigma_1)\mathcal{I}(\sigma_2)| < \beta$.*

**Proof.** It is easy to see from (A.24) that, since $\beta \leqslant 1$, the symbols appearing in the description of $V^2$ span a sector, and therefore to describe $\Gamma V^2(x)$ it is enough to look at these coefficients. Consider first $\Gamma = \mathrm{Id}_T \in G$. By coherence (2.12), and using that $\Upsilon[\Xi] = -1$, we can write

$$3\langle \mathbf{1}, V_\beta^2(y)\rangle = 3v^2(y) = -3\Upsilon_y[\Xi]v^2(y) = \Upsilon_y[\mathcal{I}(\Xi)] = \langle \mathcal{I}(\mathcal{I}(\Xi)), V_{\beta+|\mathcal{I}(\mathcal{I}(\Xi))|}(y)\rangle.$$

The reason of the shift of $\beta \leftrightarrow \beta + |\mathcal{I}(\mathcal{I}(\Xi))|$ is only to guarantee that the coefficient $\mathcal{I}(\mathcal{I}(\Xi))$ is not truncated. Similarly, we have

$$3\langle \mathcal{I}(\sigma), V_\beta^2(y)\rangle = 2v(y)\Upsilon_y[\sigma] = \Upsilon_y[\mathcal{I}(\mathcal{I}(\sigma)\mathcal{I}(\Xi))] = \langle \mathcal{I}(\mathcal{I}(\sigma)\mathcal{I}(\Xi)), V_{\beta+|\mathcal{I}(\mathcal{I}(\Xi))|}(y)\rangle, \qquad (2.29)$$

and

$$\begin{aligned}3\langle \mathcal{I}(\sigma_1)\mathcal{I}(\sigma_2), V_\beta^2(y)\rangle &= 2\Upsilon_y[\sigma_1]\Upsilon_y[\sigma_2] = \Upsilon_y[\mathcal{I}(\sigma_1)\mathcal{I}(\sigma_2)\mathcal{I}(\Xi)] \\ &= \langle \mathcal{I}(\mathcal{I}(\sigma_1)\mathcal{I}(\sigma_2)\mathcal{I}(\Xi)), V_{\beta+|\mathcal{I}(\mathcal{I}(\Xi))|}(y)\rangle.\end{aligned}$$

For general elements $\Gamma \in G$ of the structure group we use (A.23) to write the action of $\Gamma$ on $V_\beta^2$ as

$$\langle \tau, \Gamma V_\beta^2(x)\rangle = \sum_{\substack{\mu \in \mathcal{T}^+ \\ \mu \star \tau \in \mathcal{V}_\beta}} \langle \mu \star \tau, V_\beta^2(x)\rangle \frac{\gamma(\mu)}{\mu!},$$

and from here we see that to have non-trivial components we need $\langle \mu \star \tau, V_\beta^2(x)\rangle$ to not be zero. By the explicit form of $V_\beta^2(x)$ in (2.25) we have that $\mu \star \tau$ should be either $\mathbf{1}, \mathcal{I}(\sigma), \mathcal{I}(\sigma_1)\mathcal{I}(\sigma_2)$ for $\sigma \in \mathcal{V}_{0,\beta}$ or $\sigma_1, \sigma_2 \in \mathcal{V}_{0,\beta}$ with $|\mathcal{I}(\sigma_1)\mathcal{I}(\sigma_2)| < \beta$. With this in mind we have

$$\begin{aligned}\langle \tau, \Gamma V_\beta^2(x)\rangle &= \langle \tau, V_\beta^2(x)\rangle + \sum_{\mathcal{I}(\sigma) \in \mathcal{T}^+} \langle \mathcal{I}(\sigma) \star \tau, V_\beta^2(x)\rangle \frac{\gamma(\mathcal{I}(\sigma))}{\mathcal{I}(\sigma)!} \\ &\quad + \sum_{\mathcal{I}(\sigma_1)\mathcal{I}(\sigma_2) \in \mathcal{T}^+} \langle (\mathcal{I}(\sigma_1)\mathcal{I}(\sigma_2)) \star \tau, V_\beta^2(x)\rangle \frac{\gamma(\mathcal{I}(\sigma_1)\mathcal{I}(\sigma_2))}{(\mathcal{I}(\sigma_1)\mathcal{I}(\sigma_2))!},\end{aligned} \qquad (2.30)$$

and in particular we obtain for $\tau = \mathbf{1}$, using that $\Upsilon_\cdot[\Xi] = -1$ and (2.12), that

$$\begin{aligned}&3\langle \mathbf{1}, \Gamma V_\beta^2(x)\rangle \\ &= 3\langle \mathbf{1}, V_\beta^2(x)\rangle + 3\sum_{\mathcal{I}(\sigma) \in \mathcal{T}^+} \langle \mathcal{I}(\sigma), V_\beta^2(x)\rangle \frac{\gamma(\mathcal{I}(\sigma))}{\mathcal{I}(\sigma)!} + 3\sum_{\mathcal{I}(\sigma_1)\mathcal{I}(\sigma_2) \in \mathcal{T}^+} \langle \mathcal{I}(\sigma_1)\mathcal{I}(\sigma_2), V_\beta^2(x)\rangle \frac{\gamma(\mathcal{I}(\sigma_1)\mathcal{I}(\sigma_2))}{(\mathcal{I}(\sigma_1)\mathcal{I}(\sigma_2))!} \\ &= 3v^2(x) + 3\sum_{\sigma \in \mathcal{V}_{0,\beta}} 2v(x)\Upsilon_x[\sigma]\frac{\gamma(\mathcal{I}(\sigma))}{\sigma!} + 3\sum_{\substack{\mathcal{I}(\sigma_1)\mathcal{I}(\sigma_2) \in \mathcal{T}^+ \\ |\mathcal{I}(\sigma_1)\mathcal{I}(\sigma_2)| < \beta}} 2\Upsilon_x[\sigma_1]\Upsilon_x[\sigma_2]\frac{\gamma(\mathcal{I}(\sigma_1)\mathcal{I}_2)}{(\mathcal{I}(\sigma_1)\mathcal{I}(\sigma_2))!} \\ &= 3v^2(x) + \sum_{\sigma \in \mathcal{V}_{0,\beta}} \frac{6v(x)\Upsilon_x[\sigma]}{\sigma!}\gamma(\mathcal{I}(\sigma)) + \sum_{\substack{\mathcal{I}(\sigma_1)\mathcal{I}(\sigma_2) \in \mathcal{T}^+ \\ |\mathcal{I}(\sigma_1)\mathcal{I}(\sigma_2)| < \beta}} \frac{6\Upsilon_x[\sigma_1]\Upsilon_x[\sigma_2]}{(\mathcal{I}(\sigma_1)\mathcal{I}(\sigma_2))!}\gamma(\mathcal{I}(\sigma_1)\mathcal{I}(\sigma_2))\end{aligned}$$

$$
\begin{aligned}
&= -3v^2(x)\Upsilon_x[\Xi] - \sum_{\sigma \in \mathcal{V}_{0,\beta}} \frac{6v(x)\Upsilon_x[\sigma]\Upsilon_x[\Xi]}{\sigma!}\gamma(\mathcal{I}(\sigma)) - \sum_{\substack{\mathcal{I}(\sigma_1)\mathcal{I}(\sigma_2)\in\mathcal{T}^+\\|\mathcal{I}(\sigma_1)\mathcal{I}(\sigma_2)|<\beta}} \frac{6\Upsilon_x[\sigma_1]\Upsilon_x[\sigma_2]\Upsilon_x[\Xi]}{(\mathcal{I}(\sigma_1)\mathcal{I}(\sigma_2))!}\gamma(\mathcal{I}(\sigma_1)\mathcal{I}(\sigma_2)) \\
&= \Upsilon_x[\mathcal{I}(\Xi)] + \sum_{\sigma \in \mathcal{V}_{0,\beta}} \frac{\Upsilon_x[\mathcal{I}(\sigma)\,\mathcal{I}(\Xi)]}{\sigma!}\gamma(\mathcal{I}(\sigma)) + \sum_{\substack{\mathcal{I}(\sigma_1)\mathcal{I}(\sigma_2)\in\mathcal{T}^+\\|\mathcal{I}(\sigma_1)\mathcal{I}(\sigma_2)|<\beta}} \frac{\Upsilon_x[\mathcal{I}(\sigma_1)\,\mathcal{I}(\sigma_2)\,\mathcal{I}(\Xi)]}{(\mathcal{I}(\sigma_1)\mathcal{I}(\sigma_2))!}\gamma(\mathcal{I}(\sigma_1)\mathcal{I}(\sigma_2)).
\end{aligned}
$$

On the other hand we have that $\mathcal{I}(\Xi)\in\mathcal{V}$ since $0<|\mathcal{I}(\mathcal{I}(\Xi))|=3(s-1/2)-\kappa$, and by (2.23) we can consider the component

$$
\begin{aligned}
\langle\mathcal{I}(\mathcal{I}(\Xi)),\Gamma V_{\beta+|\mathcal{I}(\mathcal{I}(\Xi))|}\rangle &= \Upsilon_x[\mathcal{I}(\Xi)] + \sum_{\substack{\mathcal{I}(\sigma)\in\mathcal{T}^+\\\mathcal{I}(\sigma)\star\mathcal{I}(\Xi)\in\mathcal{V}_{0,\beta+|\mathcal{I}(\mathcal{I}(\Xi))|}}} \frac{\Upsilon_x[\mathcal{I}(\sigma)\star\mathcal{I}(\Xi)]}{\mathcal{I}(\sigma)!}\gamma(\mathcal{I}(\sigma)) \\
&\quad + \sum_{\substack{\mathcal{I}(\sigma_1)\mathcal{I}(\sigma_2)\in\mathcal{T}^+\\(\mathcal{I}(\sigma_1)\mathcal{I}(\sigma_2))\star\mathcal{I}(\Xi)\in\mathcal{V}_{0,\beta+|\mathcal{I}(\mathcal{I}(\Xi))|}}} \frac{\Upsilon_x[(\mathcal{I}(\sigma_1)\mathcal{I}(\sigma_2))\star\mathcal{I}(\Xi)]}{(\mathcal{I}(\sigma_1)\mathcal{I}(\sigma_2))!}\gamma(\mathcal{I}(\sigma_1)\mathcal{I}(\sigma_2)) \\
&= \Upsilon_x[\mathcal{I}(\Xi)] + \sum_{\substack{\mathcal{I}(\sigma)\in\mathcal{T}^+\\\mathcal{I}(\sigma)\star\mathcal{I}(\Xi)\in\mathcal{V}_{0,\beta+|\mathcal{I}(\mathcal{I}(\Xi))|}}} \frac{\Upsilon_x[\mathcal{I}(\sigma)\,\mathcal{I}(\Xi)]}{\mathcal{I}(\sigma)!}\gamma(\mathcal{I}(\sigma)) \\
&\quad + \sum_{\substack{\mathcal{I}(\sigma_1)\mathcal{I}(\sigma_2)\in\mathcal{T}^+\\(\mathcal{I}(\sigma_1)\mathcal{I}(\sigma_2))\star\mathcal{I}(\Xi)\in\mathcal{V}_{0,\beta+|\mathcal{I}(\mathcal{I}(\Xi))|}}} \frac{\Upsilon_x[\mathcal{I}(\sigma_1)\,\mathcal{I}(\sigma_2)\,\mathcal{I}(\Xi)]}{(\mathcal{I}(\sigma_1)\mathcal{I}(\sigma_2))!}\gamma(\mathcal{I}(\sigma_1)\mathcal{I}(\sigma_2)) \\
&= 3\,\langle\mathbf{1},\Gamma V_\beta^2(x)\rangle,
\end{aligned}
$$

where we used that the only way to graft into $\mathcal{I}(\Xi)$ is on the root, in which case it is the same as the tree product, which shows (2.26). We proceed similarly for the component of the type $\mathcal{I}(\tau)$. By (2.30) we have that

$$
\begin{aligned}
\langle\mathcal{I}(\tau),\Gamma V_\beta^2(x)\rangle &= \langle\mathcal{I}(\tau),V_\beta^2(x)\rangle + \sum_{\mathcal{I}(\sigma)\in\mathcal{T}^+}\langle\mathcal{I}(\sigma)\star\mathcal{I}(\tau),V_\beta^2(x)\rangle\frac{\gamma(\mathcal{I}(\sigma))}{\mathcal{I}(\sigma)!} \\
&\quad + \sum_{\mathcal{I}(\sigma_1)\mathcal{I}(\sigma_2)\in\mathcal{T}^+}\langle(\mathcal{I}(\sigma_1)\mathcal{I}(\sigma_2))\star\mathcal{I}(\tau),V_\beta^2(x)\rangle\frac{\gamma(\mathcal{I}(\sigma_1)\mathcal{I}(\sigma_2))}{(\mathcal{I}(\sigma_1)\mathcal{I}(\sigma_2))!},
\end{aligned}
$$

and taking into account (2.29), we see that to prove identity (2.27) it will be enough to prove the following identities:

$$
\begin{align}
3\,\langle\mathcal{I}(\sigma)\star\mathcal{I}(\tau),V_\beta^2(x)\rangle &= \langle\mathcal{I}(\mathcal{I}(\sigma)\star(\mathcal{I}(\tau)\mathcal{I}(\Xi))),V_{\beta+|\mathcal{I}(\mathcal{I}(\Xi))|}(x)\rangle, \tag{2.31}\\
3\,\langle(\mathcal{I}(\sigma_1)\mathcal{I}(\sigma_2))\star\mathcal{I}(\tau),V_\beta^2(x)\rangle &= \langle\mathcal{I}((\mathcal{I}(\sigma_1)\mathcal{I}(\sigma_2))\star(\mathcal{I}(\tau)\mathcal{I}(\Xi))),V_{\beta+|\mathcal{I}(\mathcal{I}(\Xi))|}(x)\rangle. \tag{2.32}
\end{align}
$$

Since $\beta\leqslant 1$ all the involved trees have no decorations by Lemma 2.7, and the $\star$-product is precisely the grafting of $\sigma\curvearrowright_0\mathcal{I}(\tau)$ as defined in (A.18) This grafting happens either at the root or in $\tau$, and therefore we can write

$$
\mathcal{I}(\sigma)\star\mathcal{I}(\tau)=\mathcal{I}(\sigma)\,\mathcal{I}(\tau)+\mathcal{I}(\mathcal{I}(\sigma)\star\tau). \tag{2.33}
$$

By a direct examination of (2.25) and coherence (2.12) we obtain

$$
\begin{aligned}
3\,\langle\mathcal{I}(\sigma)\star\mathcal{I}(\tau),V_\beta^2(x)\rangle &= 3\,\langle\mathcal{I}(\sigma)\,\mathcal{I}(\tau),V_\beta^2(x)\rangle + 3\langle\mathcal{I}(\mathcal{I}(\sigma)\star\tau),V_\beta^2(x)\rangle \\
&= 6\,\Upsilon_x[\sigma]\,\Upsilon_x[\tau] + 6\,v(x)\,\Upsilon_x[\mathcal{I}(\sigma)\star\tau] \\
&= -6\,\Upsilon_x[\sigma]\,\Upsilon_x[\tau]\,\Upsilon_x[\Xi] - 6\,v(x)\,\Upsilon_x[\mathcal{I}(\sigma)\star\tau]\,\Upsilon_x[\Xi] \\
&= \Upsilon_x[\mathcal{I}(\sigma)\mathcal{I}(\tau)\mathcal{I}(\Xi)] + \Upsilon_x[(\mathcal{I}(\mathcal{I}(\sigma)\star\tau))\,\mathcal{I}(\Xi)] \\
&= \Upsilon_x[(\mathcal{I}(\sigma)\star\mathcal{I}(\tau))\mathcal{I}(\Xi)] \\
&= \Upsilon_x[\mathcal{I}(\sigma)\star(\mathcal{I}(\tau)\mathcal{I}(\Xi))] \\
&= \langle\mathcal{I}(\mathcal{I}(\sigma)\star(\mathcal{I}(\tau)\mathcal{I}(\Xi))),V_{\beta+|\mathcal{I}(\mathcal{I}(\Xi))|}(x)\rangle,
\end{aligned}
$$

where the identity $(\mathcal{I}(\sigma) \star \mathcal{I}(\tau)) \mathcal{I}(\Xi) = \mathcal{I}(\sigma) \star (\mathcal{I}(\tau) \mathcal{I}(\Xi))$ follows since the grafting in the definition of $\star$ is not allowed on noises (leaves). This proves (2.31). Analogously, we can write

$$\begin{aligned}(\mathcal{I}(\sigma_1)\mathcal{I}(\sigma_2)) \star \mathcal{I}(\tau) &= \mathcal{I}(\sigma_1)\mathcal{I}(\sigma_2)\mathcal{I}(\tau) + \mathcal{I}(\sigma_1)\mathcal{I}(\mathcal{I}(\sigma_2) \star \tau) \\ &\quad + \mathcal{I}(\sigma_2)\mathcal{I}(\mathcal{I}(\sigma_1) \star \tau) + \mathcal{I}((\mathcal{I}(\sigma_1)\mathcal{I}(\sigma_2)) \star \tau),\end{aligned} \qquad (2.34)$$

using that the $\star$-product of each $\mathcal{I}(\sigma_1), \mathcal{I}(\sigma_2)$ is done independently of each other and only on top of $\mathcal{I}(\tau)$, and by splitting this action on the root of $\mathcal{I}(\tau)$ and on $\tau$. It is clear by a direct examination of (2.25) that $V_\beta^2(x)$ has no component in $\mathcal{I}(\sigma_1)\mathcal{I}(\sigma_2)\mathcal{I}(\tau)$ and therefore

$$\begin{aligned}3\langle(\mathcal{I}(\sigma_1)\mathcal{I}(\sigma_2))\star\mathcal{I}(\tau),V_\beta^2(x)\rangle &= 3\langle\mathcal{I}(\sigma_1)\mathcal{I}(\mathcal{I}(\sigma_2)\star\tau),V_\beta^2(x)\rangle + 3\langle\mathcal{I}(\sigma_2)\mathcal{I}(\mathcal{I}(\sigma_1)\star\tau),V_\beta^2(x)\rangle \\ &\quad + 3\langle\mathcal{I}((\mathcal{I}(\sigma_1)\mathcal{I}(\sigma_2))\star\tau),V_\beta^2(x)\rangle \\ &= 6\,\Upsilon_x[\sigma_1]\,\Upsilon_x[\mathcal{I}(\sigma_2)\star\tau] + 6\,\Upsilon_x[\sigma_2]\,\Upsilon_x[\mathcal{I}(\sigma_1)\star\tau] \\ &\quad + 6\,v(x)\,\Upsilon_x[(\mathcal{I}(\sigma_1)\mathcal{I}(\sigma_2))\star\tau] \\ &= \Upsilon_x[\mathcal{I}(\sigma_1)\mathcal{I}(\mathcal{I}(\sigma_2)\star\tau)\mathcal{I}(\Xi)] + \Upsilon_x[\mathcal{I}(\sigma_2)\mathcal{I}(\mathcal{I}(\sigma_1)\star\tau)\mathcal{I}(\Xi)] \\ &\quad + \Upsilon_x[\mathcal{I}((\mathcal{I}(\sigma_1)\mathcal{I}(\sigma_2))\star\tau)\mathcal{I}(\Xi)] \\ &= \Upsilon_x[((\mathcal{I}(\sigma_1)\mathcal{I}(\sigma_2))\star\mathcal{I}(\tau))\mathcal{I}(\Xi)] \\ &= \Upsilon_x[((\mathcal{I}(\sigma_1)\mathcal{I}(\sigma_2))\star(\mathcal{I}(\tau)\mathcal{I}(\Xi)))] \\ &= \langle\mathcal{I}((\mathcal{I}(\sigma_1)\mathcal{I}(\sigma_2))\star(\mathcal{I}(\tau)\mathcal{I}(\Xi))),V_{\beta+|\mathcal{I}(\mathcal{I}(\Xi))|}(x)\rangle,\end{aligned}$$

where we used that since $\mathcal{I}(\sigma_1)\mathcal{I}(\sigma_2)\mathcal{I}(\tau)\mathcal{I}(\Xi)$ is not sub-ternary, then in $T$ we have the identity

$$\begin{aligned}((\mathcal{I}(\sigma_1)\mathcal{I}(\sigma_2))\star\mathcal{I}(\tau))\mathcal{I}(\Xi) &= \mathcal{I}(\sigma_1)\mathcal{I}(\mathcal{I}(\sigma_2)\star\tau)\mathcal{I}(\Xi) + \mathcal{I}(\sigma_2)\mathcal{I}(\mathcal{I}(\sigma_1)\star\tau)\mathcal{I}(\Xi) \\ &\quad + \mathcal{I}((\mathcal{I}(\sigma_1)\mathcal{I}(\sigma_2))\star\tau)\mathcal{I}(\Xi).\end{aligned}$$

This shows (2.32) and therefore (2.27). At last, for the component of the type $\mathcal{I}(\tau_1)\mathcal{I}(\tau_2)$ we use (2.30) to write

$$\begin{aligned}&\langle\mathcal{I}(\tau_1)\,\mathcal{I}(\tau_2),\Gamma V_\beta^2(x)\rangle \\ &= \langle\mathcal{I}(\tau_1)\mathcal{I}(\tau_2),V_\beta^2(x)\rangle + \sum_{\mathcal{I}(\sigma)\in\mathcal{T}^+}\langle\mathcal{I}(\sigma)\star(\mathcal{I}(\tau_1)\mathcal{I}(\tau_2)),V_\beta^2(x)\rangle\frac{\gamma(\mathcal{I}(\sigma))}{\mathcal{I}(\sigma)!} \\ &\quad + \sum_{\mathcal{I}(\sigma_1)\mathcal{I}(\sigma_2)\in\mathcal{T}^+}\langle(\mathcal{I}(\sigma_1)\mathcal{I}(\sigma_2))\star(\mathcal{I}(\tau_1)\mathcal{I}(\tau_2)),V_\beta^2(x)\rangle\frac{\gamma(\mathcal{I}(\sigma_1)\mathcal{I}(\sigma_2))}{(\mathcal{I}(\sigma_1)\mathcal{I}(\sigma_2))!}.\end{aligned} \qquad (2.35)$$

The same argument as (2.33) allows us to write

$$\mathcal{I}(\sigma)\star(\mathcal{I}(\tau_1)\mathcal{I}(\tau_2)) = \mathcal{I}(\sigma)\mathcal{I}(\tau_1)\mathcal{I}(\tau_2) + \mathcal{I}(\sigma\star\mathcal{I}(\tau_1))\mathcal{I}(\tau_2) + \mathcal{I}(\tau_1)\mathcal{I}(\sigma\star\mathcal{I}(\tau_2)), \qquad (2.36)$$

and by a direct examination of (2.25) we see that $V_\beta^2$ has no component in $\mathcal{I}(\sigma)\mathcal{I}(\tau_1)\mathcal{I}(\tau_2)$ and then

$$\begin{aligned}3\langle\mathcal{I}(\sigma)\star(\mathcal{I}(\tau_1)\mathcal{I}(\tau_2)),V_\beta^2(x)\rangle &= 3\langle\mathcal{I}(\sigma\star\mathcal{I}(\tau_1))\mathcal{I}(\tau_2),V_\beta^2(x)\rangle + 3\langle\mathcal{I}(\tau_1)\mathcal{I}(\sigma\star\mathcal{I}(\tau_2)),V_\beta^2(x)\rangle \\ &= 6\,\Upsilon_x[\sigma\star\mathcal{I}(\tau_1)]\,\Upsilon_x[\tau_2] + 6\,\Upsilon_x[\tau_1]\,\Upsilon_x[\sigma\star\mathcal{I}(\tau_2)] \\ &= \Upsilon_x[\mathcal{I}(\sigma\star\mathcal{I}(\tau_1))\mathcal{I}(\tau_2)\mathcal{I}(\Xi)] + \Upsilon_x[\mathcal{I}(\tau_1)\mathcal{I}(\sigma\star\mathcal{I}(\tau_2))\mathcal{I}(\Xi)] \\ &= \Upsilon_x[(\mathcal{I}(\sigma)\star(\mathcal{I}(\tau_1)\mathcal{I}(\tau_2)))\mathcal{I}(\Xi)] \\ &= \Upsilon_x[(\mathcal{I}(\sigma)\star(\mathcal{I}(\tau_1)\mathcal{I}(\tau_2)))\star\mathcal{I}(\Xi)] \\ &= \Upsilon_x[\mathcal{I}(\sigma)\star((\mathcal{I}(\tau_1)\mathcal{I}(\tau_2))\star\mathcal{I}(\Xi))] \\ &= \Upsilon_x[\mathcal{I}(\sigma)\star(\mathcal{I}(\tau_1)\mathcal{I}(\tau_2)\mathcal{I}(\Xi))] \\ &= \langle\mathcal{I}(\mathcal{I}(\sigma)\star(\mathcal{I}(\tau_1)\mathcal{I}(\tau_2)\mathcal{I}(\Xi))),V_{\beta+|\mathcal{I}(\mathcal{I}(\Xi))|}(x)\rangle,\end{aligned}$$

where we used that when multiplying (2.36) by $\mathcal{I}(\Xi)$ the term $\mathcal{I}(\sigma)\mathcal{I}(\tau_1)\mathcal{I}(\tau_2)\mathcal{I}(\Xi)$ appearing is not in $T$. Analogously, by (2.25) we see that $V_\beta^2$ has no component on ternary trees, and therefore we can write modulo components in $T$ or ternary trees

$$(\mathcal{I}(\sigma_1)\mathcal{I}(\sigma_2)) \star (\mathcal{I}(\tau_1)\mathcal{I}(\tau_2)) = \mathcal{I}((\mathcal{I}(\sigma_1)\mathcal{I}(\sigma_2)) \star \tau_1)\mathcal{I}(\tau_2) + \mathcal{I}(\tau_1)\mathcal{I}((\mathcal{I}(\sigma_1)\mathcal{I}(\sigma_2)) \star \tau_2)$$
$$+ \mathcal{I}(\mathcal{I}(\sigma_1) \star \tau_1)\mathcal{I}(\mathcal{I}(\sigma_2) \star \tau_2).$$

Observing that none of these components are planted, we can write by (2.25)

$$\begin{aligned}
& 3 \langle (\mathcal{I}(\sigma_1)\mathcal{I}(\sigma_2)) \star (\mathcal{I}(\tau_1)\mathcal{I}(\tau_2)), V_\beta^2(x) \rangle \\
&= 3 \langle \mathcal{I}((\mathcal{I}(\sigma_1)\mathcal{I}(\sigma_2)) \star \tau_1)\mathcal{I}(\tau_2), V_\beta^2(x) \rangle + 3 \langle \mathcal{I}(\tau_1)\mathcal{I}((\mathcal{I}(\sigma_1)\mathcal{I}(\sigma_2)) \star \tau_2), V_\beta^2(x) \rangle \\
&\quad + 3 \langle \mathcal{I}(\mathcal{I}(\sigma_1) \star \tau_1)\mathcal{I}(\mathcal{I}(\sigma_2) \star \tau_2), V_\beta^2(x) \rangle \\
&= 6\, \Upsilon_x[(\mathcal{I}(\sigma_1)\mathcal{I}(\sigma_2)) \star \tau_1]\, \Upsilon_x[\tau_2] + 6\, \Upsilon_x[\tau_2]\, \Upsilon_x[(\mathcal{I}(\sigma_1)\mathcal{I}(\sigma_2)) \star \tau_2] \\
&\quad + 6\, \Upsilon_x[\mathcal{I}(\sigma_1) \star \tau_1]\, \Upsilon_x[\mathcal{I}(\sigma_2) \star \tau_2] \\
&= \Upsilon_x[\mathcal{I}((\mathcal{I}(\sigma_1)\mathcal{I}(\sigma_2)) \star \tau_1)\mathcal{I}(\tau_2)\mathcal{I}(\Xi)] + \Upsilon_x[\mathcal{I}(\tau_1)\mathcal{I}((\mathcal{I}(\sigma_1)\mathcal{I}(\sigma_2)) \star \tau_2)\mathcal{I}(\Xi)] \\
&\quad + \Upsilon_x[\mathcal{I}(\mathcal{I}(\sigma_1) \star \tau_1)\mathcal{I}(\mathcal{I}(\sigma_2) \star \tau_2)\mathcal{I}(\Xi)] \\
&= \Upsilon_x[((\mathcal{I}(\sigma_1)\mathcal{I}(\sigma_2)) \star (\mathcal{I}(\tau_1)\mathcal{I}(\tau_2)))\,\mathcal{I}(\Xi)] \\
&= \Upsilon_x[(\mathcal{I}(\sigma_1)\mathcal{I}(\sigma_2)) \star (\mathcal{I}(\tau_1)\mathcal{I}(\tau_2)\mathcal{I}(\Xi))] \\
&= \langle \mathcal{I}((\mathcal{I}(\sigma_1)\mathcal{I}(\sigma_2)) \star (\mathcal{I}(\tau_1)\mathcal{I}(\tau_2)\mathcal{I}(\Xi))), V_{\beta+|\mathcal{I}(\mathcal{I}(\Xi))|}(x) \rangle,
\end{aligned}$$

which concludes the proof of (2.28). $\square$

**Analysis of the generalised gradient**

It will be convenient to define a modelled distribution which models the generalised gradient $v_{\boldsymbol{X}}$ appearing as the coefficient at $\boldsymbol{X}$ of $V$. For $j \in \{1, \ldots, d\}$ we define

$$V^{(j)}(x) := v_{\boldsymbol{X}_j}(x)\,\mathbf{1} + \sum_{\sigma \in \mathcal{V}_{1,\gamma}} \frac{\Upsilon_x[\sigma]}{\sigma!}\mathcal{I}_j(\sigma). \qquad (2.37)$$

Due to the restriction in the sum to trees $\sigma \in \mathcal{V}_{1,\gamma}$ we have that $|\mathcal{I}_j(\sigma)| > 0$ and therefore $V^{(j)}$ takes values in a function-like sector.

**Remark 2.22.** The modelled distribution $V^{(j)}$ defined by (2.37) is not the same as $\partial_j V$ where $\partial_j$ is a *realisation* of a partial derivative as defined in [FH20, Section 14.1]. The difference is that $\partial_j V$ will also include trees $\tau \in \mathcal{V}_{0,1}$ and therefore will not be a function-like modelled distribution. Moreover, $\partial_j V$ satisfies that $\mathcal{R}(\partial_j V) = \partial_j \mathcal{R} V = \partial_j v$ which is not the same as $v_{\boldsymbol{X}_j} = \mathcal{R} V^{(j)}$.

Next result describes the action of the structure group on $V^{(j)}$ in terms of the action on $V$, and in particular tells us that $V^{(j)} \in \mathcal{D}^{\gamma-1}$.

**Lemma 2.23**. *For $\Gamma \in G$ all non-zero coefficients of $\Gamma V^{(j)}$ are*

$$\langle \mathbf{1}, \Gamma V^{(j)}(x) \rangle = \langle \boldsymbol{X}_j, \Gamma V(x) \rangle, \qquad (2.38)$$
$$\langle \mathcal{I}_j(\sigma), \Gamma V^{(j)}(x) \rangle = \langle \mathcal{I}(\sigma), \Gamma V(x) \rangle, \qquad (2.39)$$

*for $\sigma \in \mathcal{V}_{1,\gamma}$.*

**Proof.** Using the linearity of the representation (A.23) we have for any $\Gamma \in G$ that

$$\langle \tau, \Gamma V^{(j)}(x) \rangle = \sum_{\mu \in \mathcal{T}^+} \langle \mu \star \tau, V^{(j)}(x) \rangle \frac{\gamma(\mu)}{\mu!}.$$

Analogously to Lemma 2.8 one can see that by (A.24) $\{\mathbf{1}\}\cup\{\mathcal{I}_j(\sigma)\}_{\sigma\in\mathcal{V}_{1,\beta}}$ spans a sector, and therefore it is enough to look at these components. By (2.37) the form of the coefficient at $\mathbf{1}$ is

$$\begin{aligned}\langle \mathbf{1}, \Gamma V^{(j)}(x)\rangle &= \sum_{\mu\in\mathcal{T}^+}\langle \mu, V^{(j)}(x)\rangle\frac{\gamma(\mu)}{\mu!}\\ &= \langle \mathbf{1}, V^{(j)}(x)\rangle + \sum_{\sigma\in\mathcal{V}_{1,\gamma}}\langle \mathcal{I}_j(\sigma), V^{(j)}(x)\rangle\frac{\gamma(\mathcal{I}_j(\sigma))}{\mathcal{I}(\sigma)!}\\ &= v_{\boldsymbol{X}_j}(x) + \sum_{\sigma\in\mathcal{V}_{1,\gamma}}\frac{\Upsilon_x[\sigma]}{\sigma!}\gamma(\mathcal{I}_j(\sigma)),\end{aligned}$$

and (2.38) follows from (2.22) in Lemma 2.19. On the other hand, we have that for $\sigma\in\mathcal{V}_{1,\gamma}$

$$\begin{aligned}\langle \mathcal{I}_j(\sigma), \Gamma V^{(j)}(x)\rangle &= \sum_{\mu\in\mathcal{T}^+}\langle \mu\star\mathcal{I}_j(\sigma), V^{(j)}(x)\rangle\frac{\gamma(\mu)}{\mu!} = \sum_{\substack{\mu\in\mathcal{T}^+ \\ \mu\star\tau\in\mathcal{V}_{1,\gamma}}}\langle \mathcal{I}_j(\mu\star\sigma), V^{(j)}(x)\rangle\frac{\gamma(\mu)}{\mu!}\\ &= \sum_{\substack{\mu\in\mathcal{T}^+ \\ \mu\star\tau\in\mathcal{V}_{1,\gamma}}}\Upsilon_x[\mu\star\sigma]\frac{\gamma(\mu)}{\mu!},\end{aligned}$$

Where we have used that $V^{(j)}$ takes values in trees planted by $\mathcal{I}_j$, and therefore $\mu\star\mathcal{I}_j(\sigma)$ cannot graft (or decorate) the root of $\mathcal{I}_j(\sigma)$ and this is equivalent to only consider the part of the grafting (and decorating) in $\tau$, which can be equivalently written as $\mathcal{I}_j(\mu\star\sigma)$. Equivalently, we have that $\langle \mu\star\mathcal{I}_j(\sigma),\mathcal{I}_j(\tau)\rangle = \langle \mathcal{I}_j(\mu\star\tau),\mathcal{I}_j(\tau)\rangle$. At last, observe that if $\mu\in\mathcal{T}^+$ is such that $\mu\star\sigma\in\mathcal{V}_{0,1}$ then $0<|\mathcal{I}(\mu\star\sigma)|<1$ and therefore $|\mathcal{I}(\sigma)|<1-|\mu|<1$, i.e., $\tau\notin\mathcal{V}_{1,\gamma}$. With this last observation and (2.23) from Lemma 2.19 we conclude (2.39). $\square$

**Back to the coefficients of $V$**

**LEMMA 2.24.** *For any $\beta\in(1,\gamma]$ and $x,y\in\mathbb{R}^{1+d}$ we have that*

$$\langle \mathcal{I}(\tau), \Gamma_{yx}V_\beta(x)\rangle = \tau!\, c_\tau\begin{cases}\Pi_x(V_{\beta-|\mathcal{I}(\tau)|}(x))(y) & \text{if } \Upsilon.[\tau]\sim v\\ \Pi_x(V^2_{\beta-|\mathcal{I}(\tau)|}(x))(y) & \text{if } \Upsilon.[\tau]\sim v^2\\ \gamma_{yx}(V^{(j)}_{\beta-|\mathcal{I}_j(\tau)|})(y) & \text{if } \Upsilon.[\tau]\sim v_{\boldsymbol{X}_j}\\ 1 & \text{if } \Upsilon.[\tau]\sim 1\end{cases}.$$

**Proof.** If $\Upsilon.[\tau]\sim v$ then by Lemma 2.17 and (2.23) we have

$$\begin{aligned}&\frac{\langle \mathcal{I}(\tau), \Gamma_{\cdot x}V_\beta(x)\rangle}{\tau!}\\ &= \Upsilon_x[\tau] + \mathbb{1}_{\mathcal{V}_{0,\beta}}(\boldsymbol{X}_j\star\tau)\sum_{j=1}^d\frac{\Upsilon_x[\boldsymbol{X}_j\star\tau]}{\boldsymbol{X}_j!}\gamma_{\cdot,x}(\boldsymbol{X}_j) + \sum_{\sigma\in\mathcal{V}_{0,\beta-|\mathcal{I}(\tau)|}}\frac{\Upsilon_x[\mathcal{I}(\sigma)\star\tau]}{\mathcal{I}(\sigma)!}\gamma_{\cdot,x}(\mathcal{I}(\sigma))\\ &= 2c_\tau v(x) + \mathbb{1}_{(1,2s)}(|\beta|-|\mathcal{I}(\tau)|)\sum_{j=1}^d 2c_\tau v_{\boldsymbol{X}_j}(x)((\cdot)_j - x_j) + \sum_{\sigma\in\mathcal{V}_{0,\beta-|\mathcal{I}(\tau)|}}\frac{2c_\tau \Upsilon_x[\sigma]}{\sigma!}\gamma_{\cdot,x}(\mathcal{I}(\sigma))\\ &= 2c_\tau\left(v(x)\Pi_x\mathbf{1} + \mathbb{1}_{(1,2s)}(|\beta|-|\mathcal{I}(\tau)|)\, v_{\boldsymbol{X}}(x)\cdot\Pi_x\boldsymbol{X} + \sum_{\sigma\in\mathcal{V}_{0,\beta-|\mathcal{I}(\tau)|}}\frac{\Upsilon_x[\sigma]}{\sigma!}\Pi_x\mathcal{I}(\sigma)\right)\end{aligned}$$

$$
\begin{aligned}
&= 2\, c_\tau\, \Pi_x\!\left( v(x)\, \mathbf{1} + \mathbb{1}_{(1,2s)}(|\beta| - |\mathcal{I}(\tau)|)\, v_{\boldsymbol{X}}(x)\cdot \boldsymbol{X} + \sum_{\sigma\in\mathcal{V}_{0,\beta-|\mathcal{I}(\tau)|}} \frac{2\, c_\tau\, \Upsilon_x[\sigma]}{\sigma!} \mathcal{I}(\sigma) \right) \\
&= 2\, c_\tau\, \Pi_x(V_{\beta-|\mathcal{I}(\tau)|}(x))(y),
\end{aligned}
$$

where we used that $\boldsymbol{X}_j \star \tau \in \mathcal{V}_{0,\beta} \Longleftrightarrow 1 < |\beta| - |\mathcal{I}(\tau)|$. If $\Upsilon.[\tau] \sim v^2$ then by Lemma 2.17, (2.7) and (2.23) we have

$$
\begin{aligned}
\frac{\langle \mathcal{I}(\tau), \Gamma_{\cdot,x} V_\beta(x) \rangle}{\tau!} &= \frac{\Upsilon_x[\mathbf{1}\star\tau]}{\mathbf{1}!}\gamma_{\cdot,x}(\mathbf{1}) + \sum_{\substack{\mathcal{I}(\sigma)\in\mathcal{T}^+ \\ \mathcal{I}(\sigma)\star\tau\in\mathcal{V}_\beta}} \frac{\Upsilon_x[\mathcal{I}(\sigma)\star\tau]}{\mathcal{I}(\sigma)!}\gamma_{\cdot,x}(\mathcal{I}(\sigma)) \\
&\quad + \sum_{\substack{\mathcal{I}(\sigma_1)\mathcal{I}(\sigma_2)\in\mathcal{T}^+ \\ (\mathcal{I}(\sigma_1)\mathcal{I}(\sigma_2))\star\tau\in\mathcal{V}_\beta}} \frac{\Upsilon_x[(\mathcal{I}(\sigma_1)\mathcal{I}(\sigma_2))\star\tau]}{(\mathcal{I}(\sigma_1)\mathcal{I}(\sigma_2))!}\gamma_{\cdot,x}(\mathcal{I}(\sigma_1)\mathcal{I}(\sigma_2)) \\
&= \Upsilon_x[\tau] + \sum_{\substack{\mathcal{I}(\sigma)\in\mathcal{T}^+ \\ |\mathcal{I}(\sigma)|<\beta-|\tau|}} \frac{2\, c_\tau v(x)\, \Upsilon_x[\sigma]}{\sigma!}\gamma_{\cdot,x}(\mathcal{I}(\sigma)) \\
&\quad + \sum_{\substack{\mathcal{I}(\sigma_1)\mathcal{I}(\sigma_2)\in\mathcal{T}^+ \\ |\mathcal{I}(\sigma_1)\mathcal{I}(\sigma_2)|<\beta}} \frac{2\, c_\tau\, \Upsilon_x[\sigma_1]\, \Upsilon_x[\sigma_2]}{(\mathcal{I}(\sigma_1)\mathcal{I}(\sigma_2))!}\gamma_{\cdot,x}(\mathcal{I}(\sigma_1))\, \gamma_{\cdot,x}(\mathcal{I}(\sigma_2)) \\
&= c_\tau v^2(x) + \sum_{\substack{\mathcal{I}(\sigma)\in\mathcal{T}^+ \\ |\mathcal{I}(\sigma)|<\beta-|\tau|}} \frac{2\, c_\tau v(x)\, \Upsilon_x[\sigma]}{\sigma!} \Pi_x \mathcal{I}(\sigma) \\
&\quad + \sum_{\substack{\mathcal{I}(\sigma_1)\mathcal{I}(\sigma_2)\in\mathcal{T}^+ \\ |\mathcal{I}(\sigma_1)\mathcal{I}(\sigma_2)|<\beta}} \frac{2\, c_\tau\, \Upsilon_x[\sigma_1]\, \Upsilon_x[\sigma_2]}{(\mathcal{I}(\sigma_1)\mathcal{I}(\sigma_2))!}\Pi_x\mathcal{I}(\sigma_1)\,\Pi_x\mathcal{I}(\sigma_2) \\
&= c_\tau\, \Pi_x(V^2_{\beta-|\mathcal{I}(\tau)|}(x)),
\end{aligned}
$$

where we used that since $(\mathfrak{m}(\boldsymbol{X}_j\tau), |\mathfrak{n}(\boldsymbol{X}_j\tau)|) = (2,1)$ then Lemma 2.7 implies that $1 + |\tau| = |\boldsymbol{X}_j\tau| > 0$, and therefore $\beta - |\mathcal{I}(\tau)| \leqslant 2\,s - |\mathcal{I}(\tau)| = -|\tau| < 1$, and we can use expression (2.25) for $V^2_{\beta-|\mathcal{I}(\tau)|}$. If $\Upsilon[\tau] \sim v_{\boldsymbol{X}_j}$ by Lemma 2.17 and (2.23) we have

$$
\begin{aligned}
\frac{\langle \mathcal{I}(\tau), \Gamma_{\cdot x} V_\beta(x) \rangle}{\tau!} &= \Upsilon_x[\tau] + \sum_{\substack{\mathcal{I}_j(\sigma)\in\mathcal{T}^+ \\ \mathcal{I}_j(\sigma)\star\tau\in\mathcal{V}_{0,\beta}}} \frac{\Upsilon_x[\mathcal{I}_j(\sigma)\star\tau]}{\mathcal{I}_j(\sigma)!}\gamma_{\cdot,x}(\mathcal{I}_j(\sigma)) \\
&= c_\tau \gamma_{yx}\bigl(V^{(j)}_{\beta-|\mathcal{I}_j(\tau)|}(x)\bigr),
\end{aligned}
$$

where we used that in definition (2.37) all trees in the expansion of $V^{(j)}$ can be interpreted as elements of the structure group and the action of $\gamma$ on $V^{(j)}$ is well-defined. At last if $\Upsilon[\tau] \sim 1$ then by Lemma 2.17 and (2.23) we have that

$$
\frac{\langle \mathcal{I}(\tau), \Gamma_{\cdot x} V_\beta(x) \rangle}{\tau!} = \Upsilon_x[\tau] = c_\tau. \qquad \square
$$

From the results in this section we see that the change of base points for $V^2_\beta$ and $V^{(j)}_\beta$ can be written in terms of the change of base points for a suitable truncation of $V$, and that this change of base points for $V_\beta$ can only be three possibilities which correspond to local expansions of $v$, its square $v^2$, and its generalised derivative $v_{\boldsymbol{X}_j}$ given by $\Pi_x V(x), \Pi_x V^2(x)$ and $\Pi_x V^{(j)}_\beta(x)$ respectively.

**Remark 2.25.** In contrast to (2.6), in general we have that $\Pi_x(\mathcal{I}_j(\tau)) \neq \gamma_{\cdot,x}(\mathcal{I}_j(\tau))$ since in (A.24) we cannot guarantee that trees $\sigma \in \mathcal{T}$ such that $|\mathcal{I}_j(\sigma)| < 0$ will not contribute to the sum.

## 3. Multilevel Schauder Estimate

For this section we drop any assumption of periodicity, consider $s \in (0,1)$, and fix $\gamma \in (0, 2s) \setminus \{1\}$. We consider *germs* of functions i.e., functions of two space-time variables $U: \mathbb{R}^{1+d} \times \mathbb{R}^{1+d} \to \mathbb{R}$. The first argument of $U$ is called the base point and the second the running variable. We think of a germ as a family of space-time functions indexed by the base point. We consider a diagonal $\gamma$-Hölder type seminorm on germs given by

$$[U]_{\gamma;B} := \sup_{x \in B} \inf_{\nu(x) \in \mathbb{R}^d} \sup_{\substack{y \in B \\ y_0 \leq x_0}} \frac{|U(x,y) - U(x,x) - \mathbb{1}_{\gamma>1} \nu(x) \cdot (y-x)|}{d(x,y)^\gamma}, \qquad (3.1)$$

where $B \subset \mathbb{R}^{1+d}$ is a half-parabolic ball, and recall our convention $\nu(x) \cdot (y-x) := \nu(x) \cdot (y-x)_{1:d}$. The seminorm is effectively measuring the error in the Taylor expansion around $x$ of the function $U(x, \cdot)$. If $[U]_{\gamma;B} < +\infty$ for $\gamma \in (1, 2s)$ then the infimum in (3.1) is achieved for all $x \in B$ at the point

$$\nu(x) := \nabla U(x, \cdot)|_x \in \mathbb{R}^d. \qquad (3.2)$$

We call $\nu$ defined by (3.2) the generalised gradient of $U$. The seminorm (3.1) is $\gamma$-homogeneous in the following sense: fix $z \in \mathbb{R}^{1+d}$, $\sigma > 0$ and consider the germ $U_{z,\sigma}(x,y) := U(z + \sigma x, z + \sigma y)$, then we have that $[U_{z,\sigma}]_{\gamma;B_1} = \sigma^\gamma [U]_{\gamma;B_\sigma(z)}$. Moreover, this seminorm is invariant under the recentering given by

$$U_c(x,y) := U(x,y) - U(x,x) - \nu(x) \cdot (y-x), \qquad (3.3)$$

and for this new germ one would have that $U_c(x,x) = 0$ and $\nu_c(x) := \nabla U_c(x, \cdot)|_x = 0$, and the definition of the seminorm simplifies to

$$[U_c]_{\gamma;B} = \sup_{\substack{x,y \in B \\ y_0 \leq x_0}} \frac{|U_c(x,y)|}{d(x,y)^\gamma}.$$

Our Schauder estimate aims to measure the $2s$ improvement of regularity of the germ $U$ compared to the germ $\mathscr{L} U(x, \cdot)$, where the operator $\mathscr{L} = (\partial_t + (-\Delta)^s)$ is understood to act on the running variable of $U(x, \cdot)$. Since $\gamma - 2s < 0$ we have to measure regularity in a distributional sense. Given a half-parabolic ball $B \subset \mathbb{R}^{1+d}$ and a distributional germ $F: \mathbb{R}^{1+d} \to \mathscr{D}'(B)$, we consider the seminorm

$$\|F\|_{\gamma-2s;B} := \sup_{x \in B} \sup_{\psi \in \mathcal{B}_r} \sup_{\substack{\lambda \in (0,1) \\ B_\lambda(x) \subset B}} |\langle F(x), \psi_x^\lambda \rangle| \lambda^{-(\gamma-2s)}, \qquad (3.4)$$

where $r \in \mathbb{Z}^+$ is the smallest integer such that $r > -(\gamma - 2s)$ and

$$\mathcal{B}_r := \left\{ \psi \in \mathscr{D}(B_1) : \operatorname{supp}(\psi) \subset B_1, \|\psi\|_{C^r} := \sum_{k \in \mathbb{N}^{1+d}, k \leq r} \|\partial^k \psi\|_{B_1} \leq 1 \right\}, \qquad (3.5)$$

where $\|\psi\|_B := \sup_{x \in B} |\psi(x)|$. Since $2 > \lfloor -(\gamma - 2s) \rfloor$ it is enough to consider $r = 2$ fixed.

The fractional heat operator $\mathscr{L}$ annihilates linear in space polynomials, which implies then that the seminorm $U \mapsto \|\mathscr{L} U\|_{\gamma-2s;B}$ is invariant under the recentering (3.3). Moreover, due to the $2s$-homogeneity of $\mathscr{L}$ (see Section 1.1) and the $\gamma - 2s$-homogeneity of (3.4) (see Lemma 3.15) one can see that for $z \in \mathbb{R}^{1+d}$ and $\sigma > 0$:

$$\|\mathscr{L}(U_{z,\sigma})\|_{\gamma-2s;B_1} = \sigma^{2s} \|(\mathscr{L} U)_{z,\sigma}\|_{\gamma-2s;B_1} = \sigma^\gamma \|\mathscr{L} U\|_{\gamma-2s;B_\sigma(z)}.$$

Fix a finite set $A \subset (0, \gamma)$. We say that a germ $U : \mathbb{R}^{1+d} \times \mathbb{R}^{1+d} \to \mathbb{R}$ satisfies a *3-point continuity* condition of order $\gamma$ in $B \subset \mathbb{R}^{1+d}$ if there exists some vector-valued germ $\Lambda : \mathbb{R}^{1+d} \times \mathbb{R}^{1+d} \to \mathbb{R}^d$ such that for all $x, y, z \subset B$ with $z_0 \leqslant y_0 \leqslant x_0$ we have that

$$|U(x,z) - U(x,y) - U(y,z) + U(y,y) + \mathbb{1}_{\gamma > 1} \Lambda(x,y) \cdot (z-y)| \leqslant C \sum_{\beta \in A} d(x,y)^\beta d(y,z)^{\gamma - \beta}.$$

We define $[U]_{\gamma\text{-3pt};B}$ as the optimal constant $C$ for which this estimate holds. In particular we have:

$$|U(x,z) - U(x,y) - U(y,z) + U(y,y) + \mathbb{1}_{\gamma > 1} \Lambda(x,y) \cdot (z-y)| \leqslant [U]_{\gamma\text{-3pt};B} \sum_{\beta \in A} d(x,y)^\beta d(y,z)^{\gamma - \beta}. \quad (3.6)$$

In contrast to $\nu$, this property does not characterises the germ $\Lambda$ uniquely, and therefore our definition of $[U]_{\gamma\text{-3pt};B}$ depends on the choice of $\Lambda$. However, it will be clear from the context which $\Lambda$ we are considering. The inclusion of the linear term $\Lambda(x,y) \cdot (z-y)$ is an extra degree of freedom in the definition, and it is related to the recentering (3.3), since for this germ we have that

$$U_c(x,z) - U_c(x,y) - U_c(y,z) = U(x,z) - U(x,y) - U(y,z) + U(y,y) + (\nu(y) - \nu(x)) \cdot (z-y),$$

which produces a linear term. The inclusion of $\Lambda$ allows us to have invariance under the recentering since $U_c$ will then satisfy a *3-point continuity* condition with respect to the germ

$$\Lambda_c(x,y) := \nu(y) - \nu(x) - \Lambda(x,y). \quad (3.7)$$

Then we get for the centred germ $U_c$ that the *3-point continuity* condition reads as

$$|U_c(x,z) - U_c(x,y) - U_c(y,z) + \Lambda_c(x,y) \cdot (z-y)| \leqslant [U]_{\gamma\text{-3pt};B} \sum_{\beta \in A} d(x,y)^\beta d(y,z)^{\gamma - \beta}.$$

The seminorm $[\,\cdot\,]_{\gamma\text{-3pt};B}$ is also $\gamma$-homogeneous (see Lemma 3.15). We will need to consider a *3-point continuity* condition of order $\gamma - 1$ on the germ $\Lambda$. Since $\gamma - 1 \in (0,1)$ there is no need of an extra linear term and this condition reads as

$$|\Lambda(x,z) - \Lambda(x,y) - \Lambda(y,z) + \Lambda(y,y)| \leqslant [\Lambda]_{(\gamma-1)\text{-3pt};B} \sum_{\beta \in A \cap (1,\gamma)} d(x,y)^{\beta - 1} d(y,z)^{\gamma - \beta}. \quad (3.8)$$

Moreover, since

$$\Lambda_c(x,z) - \Lambda_c(x,y) - \Lambda_c(y,z) + \Lambda_c(y,y) = \Lambda(x,z) - \Lambda(x,y) - \Lambda(y,z) + \Lambda(y,y), \quad (3.9)$$

we see that the *3-point continuity* on the germ $\Lambda$ is also invariant under the recentering (3.3).

At last, we consider the following uniform norm for $D_1, D_2 \subset \mathbb{R}^{1+d}$:

$$\|U\|_{D_1 \times D_2} := \sup_{x \in D_1} \sup_{\substack{y \in D_2 \\ y_0 \leqslant x_0}} |U(x,y)|,$$

and when $D_1 = D_2$ we simplify notation and denote $\|U\|_{D_1} := \|U\|_{D_1 \times D_1}$. Observe that this norm is neither invariant under the recentering (3.3) nor $\gamma$-homogeneous.

Our main result is the following Schauder estimate:

**Theorem 3.1**. *Fix $\gamma \in (1, 2s)$ and a finite set $A \subset (0, \gamma)$ and $B_r \subset \mathbb{R}^{1+d}$ the half-parabolic ball of radius $r > 0$. Let $U: B_1 \times ((-1, 0] \times \mathbb{R}^d) \to \mathbb{R}$ be a bounded germ such that $[U]_{\gamma; B_1} < +\infty$. Assume that there exists a germ $\Lambda: B_1 \times ((-1, 0] \times \mathbb{R}^d) \to \mathbb{R}^d$ such that (3.6) holds in $B_1$ and that this germ $\Lambda$ satisfies a 3-point continuity condition of order $\gamma - 1 \in (0, 1)$ in $B_1$ as in (3.8). Then, if $\|\mathscr{L} U\|_{\gamma - 2s; B_1} < +\infty$, we can conclude the estimate*

$$[U]_{\gamma; B_{1/2}} \lesssim \|\mathscr{L} U\|_{\gamma - 2s; B_1} + [U]_{\gamma\text{-3pt}; B_1} + [\Lambda]_{(\gamma - 1)\text{-3pt}; B_1} + \|U\|_{B_1 \times ((-1, 0] \times \mathbb{R}^d)}, \tag{3.10}$$

*for some implicit constant that only depends on $s, d, \gamma, A$.*

**Remark 3.2.** Since the estimate (3.9) is invariant under the shift $U(x, y) \mapsto U(x, y) - U(x, x)$ we will always consider germs that are zero on the diagonal, i.e., $U(x, x) = 0$ for all $x$.

The proof of this result is inspired by the work [FR17] which follows the kernel-free approach to Schauder estimates by scaling in [Sim97]. A crucial element of the proof is the following Liouville-type theorem which we take from [FR17, Theorem 2.1].

**Theorem 3.3**. (Liouville-type theorem) *Let $s \in (0, 1)$, and $u: (-\infty, 0) \times \mathbb{R}^d \to \mathbb{R}$ be any weak solution to:*

$$\mathscr{L} u = 0 \quad in \quad (-\infty, 0) \times \mathbb{R}^d,$$

*such that*

$$\|u\|_{L^\infty(B_R)} \leqslant C (R^\gamma + 1)$$

*for some $\gamma \in [0, 2s)$, with $B_R \subset \mathbb{R}^{1+d}$ the half-parabolic ball of radius $B_R$. Then $u(t, x)$ is a polynomial in the $x$ variable of degree at most $\lfloor \gamma \rfloor \in \{0, 1\}$ in $(-\infty, 0) \times \mathbb{R}^d$.*

For $s = 1$ Liouville-type theorems are available for all $\gamma \geqslant 0$. The restriction of $\gamma \in [0, 2s)$ in Theorem 3.3 is enough to guarantee that the action of the fractional Laplacian is well-defined on regular enough functions.

The first step in the proof of Theorem 3.1 is the following lemma. We postpone its proof to Section 3.1.

**Lemma 3.4**. *Fix $\gamma \in (1, 2s)$. Let $U: B_1 \times ((-1, 0] \times \mathbb{R}^d) \to \mathbb{R}$ be a germ with $[U]_{\gamma; B_1} < +\infty$, $\nu$ defined by (3.2) and $\Lambda$ some germ such that (3.6) holds. Then for every $\delta > 0$ there exists $C(\delta, s, d, \gamma, A) > 0$ such that*

$$[U]_{\gamma; B_{1/2}} \leqslant \delta [U]_{\gamma; (-1, 0] \times \mathbb{R}^d} + C \left( \|\mathscr{L} U\|_{\gamma - 2s; B_1} + [U]_{\gamma\text{-3pt}; B_1} + \|U\|_{B_1} + \|\nu\|_{B_1} \right). \tag{3.11}$$

In order to prove Theorem 3.1 from Lemma 3.4 we need to first do a localisation of the first seminorm appearing in the right hand side of Corollary 3.12, this localisation introduces a reminder term on which $\mathscr{L}$ acts, and which can be controlled using the non-local $L^\infty$ norm appearing in the final estimate. After this localisation, the standard way to conclude would be to use the absorption lemma found in [Sim97, Section 4 p.398 ]. However, in contrast to classical Hölder seminorms $[u]_\alpha$ for functions and $\alpha \in (0, 1)$ where the sub-additivity of the mapping $B \mapsto [u]_{\alpha; B}$ follows from the additivity of the increments, i.e., $u(z) - u(x) = u(z) - u(y) + u(y) - u(x)$, the seminorm (3.1) defined on germs the mapping $B \mapsto [U]_{\gamma; B}$ is not sub-additive. The three point continuity seminorm is precisely a way to control how non-additive our germ is, and we can show the following almost subadditivity property:

**Lemma 3.5**. *Let $B, \{B_i\}_{i \in I}$ be a finite number of half-parabolic balls such that $B \subset \bigcup_{i \in I} B_i$, and let $U$ be a germ such that $[U]_{\gamma;B_i}, [U]_{\gamma\text{-}3\text{pt};B}, [\Lambda]_{(\gamma-1)\text{-}3\text{pt};B} < +\infty$. Then*

$$[U]_{\gamma;B} \lesssim \sum_{i \in I} [U]_{\gamma;B_i} + [U]_{\gamma\text{-}3\text{pt};B} + [\Lambda]_{(\gamma-1)\text{-}3\text{pt};B},$$

*for an implicit constant that depends only on $s, d, \gamma, A$ and $\#I$ the cardinality of $I$.*

With this property in mind we prove the following generalisation of Simon's abstract absorption:

**Lemma 3.6**. *Let $S$ and $D$ be non-negative, monotone functions defined on the convex subsets of a half-parabolic ball $B := B_{\rho_0}(y_0) \subset \mathbb{R}^{1+d}$. Assume that for all half-parabolic balls $B_0, B_1, \ldots, B_n \subset B$ such that $B_0 \subset \bigcup_{i=1}^n B_j$ there exists $C(n) > 0$ such that we have the almost subadditivity property:*

$$S(B_0) \leqslant C(n) \left( \sum_{i=1}^n S(B_j) + D(B) \right). \tag{3.12}$$

*Then for any given constants $\theta_0 \in (0, 1/2], \gamma > 0$, there exists $\varepsilon = \varepsilon(\theta_0, \gamma, s, d, \rho_0) \in (0, 1)$ such that if for some $E \geqslant 0$ the following condition is satisfied:*

$$\sigma^\gamma S(B_{\theta_0 \sigma}(y)) \leqslant \varepsilon \, \sigma^\gamma S(B_\sigma(y)) + E \qquad \forall B_\sigma(y) \subset B,$$

*then for each $\theta \in (0, 1)$ there exists a constant $C = C(d, \theta_0, \theta, \gamma, s, \rho_0) > 0$ such that*

$$\rho^\gamma S(B_{\theta\rho}(y)) \leqslant C(E + D(B)) \qquad \forall B_\sigma(y) \subset B.$$

## 3.1. Proof of Lemma 3.4

The following cone condition (Definition 3.7) is necessary to obtain some control on gradients (Lemma 3.9 and Lemma 3.11).

**Definition 3.7**. *Let $\beta \in (0, 1)$ and $r_0 > 0$. We say that a subset $D \subset \mathbb{R}^{1+d}$ satisfies a $(\beta, r_0)$-spatial cone condition if for every $x \in D, r \in [0, r_0]$ and $\theta \in \mathbb{R}^d$ there exists $y_{1:d} \in \mathbb{R}^d$ such that $x_0 = y_0$, $y \in D$, $d(x, y) = r$ and $\beta r |\theta| \leqslant |\theta \cdot (y_{1:d} - x_{1:d})|$.*

**Remark 3.8.** We only use the cone condition on half-parabolic balls, for which we have that $B_r(z)$ satisfies a $(\beta, r)$-spatial cone condition for $\beta = 1/\sqrt{2}$ (see discussion after [MW20, Corollary 2.12]).

**Lemma 3.9**. *Fix $\gamma \in (1, 2s)$ and a finite set $A \subset (0, \gamma)$. Let $B \subset \mathbb{R}^{1+d}$ be some half-parabolic ball, $U: B \times B \to \mathbb{R}$ a germ such that $[U]_{\gamma;B} < +\infty$, and $\nu: B \to \mathbb{R}^d$ defined by (3.2). Assume that there exists a vector-valued germ $\Lambda: B \times B \to \mathbb{R}^d$ such that (3.6) holds in $B$, and consider the germ $(\nu - \Lambda)(x, y) := \nu(y) - \nu(x) - \Lambda(x, y)$. Then we have that*

$$[\nu - \Lambda]_{\gamma-1;B_\sigma(z)} \lesssim [U]_{\gamma\text{-}3\text{pt};B_\sigma(z)} + [U]_{\gamma;B_\sigma(z)}, \tag{3.13}$$

*where the implicit proportionality constant depends only on $\gamma$ and $A$.*

**Proof.** Given $x, y, z$ we have:

$$\begin{aligned}-(\nu(y) - \nu(x)) \cdot (z - y) &= U(x, z) - U(x, y) - U(y, z) - (U(x, z) - \nu(x) \cdot (z - x)) \\ &\quad + U(x, y) - \nu(x) \cdot (y - x) + U(y, z) - \nu(y) \cdot (z - y),\end{aligned}$$

and therefore

$$\begin{aligned}-(\nu - \Lambda)(x, y) \cdot (z - y) &= U(x, z) - U(x, y) - U(y, z) + \Lambda(x, y) \cdot (z - y) \\ &\quad - (U(x, z) - \nu(x) \cdot (z - x)) + U(x, y) - \nu(x) \cdot (y - x) \\ &\quad + U(y, z) - \nu(y) \cdot (z - y).\end{aligned}$$

Let $r_0$ be the radius of the half-parabolic ball $B$, then $B$ satisfies a $(\beta, r_0)$-spatial cone condition for $\beta = 1/\sqrt{2}$. Then for any $x, y \in B$ with $y_0 \leqslant x_0$ we have $d(x, y) \leqslant 2 r_0$, and therefore we can apply the cone condition to $\theta = (\nu - \lambda)(x, y) \in \mathbb{R}^d$ and $r = \frac{1}{2} d(x, y) \leqslant r_0$ to obtain that there exists $z = (y_0, z_{1:d}) \in B$ with $d(y, z) = r$ and

$$\beta \frac{d(x, y)}{2} |(\nu - \Lambda)(x, y)| \leqslant |(\nu - \Lambda)(x, y) \cdot (z - y)|. \tag{3.14}$$

On the other hand, using that $d(x, z) \leqslant \frac{3}{2} d(x, y)$ and $d(y, z) = \frac{1}{2} d(x, y)$ we have that

$$\begin{aligned}&|(\nu - \Lambda)(x, y) \cdot (z - y)| \\ &\leqslant |U(x, z) - U(x, y) - U(y, z) + \Lambda(x, y) \cdot (z - y)| \\ &\quad + |U(x, z) - \nu(x) \cdot (z - x)| + |U(x, y) - \nu(x) \cdot (y - x)| + |U(y, z) - \nu(y) \cdot (z - y)| \\ &\leqslant [U]_{\gamma\text{-3pt};B} \sum_{\beta \in A} d(x, y)^\beta d(y, z)^{\gamma - \beta} + [U]_{\gamma;B} (d(x, z)^\gamma + d(x, y)^\gamma + d(y, z)^\gamma) \\ &\lesssim d(x, y)^\gamma ([U]_{\gamma\text{-3pt};B} + [U]_{\gamma;B}),\end{aligned}$$

which allows us to conclude on (3.14) for any $x, y \in B$ with $y_0 \leqslant x_0$:

$$\frac{|(\nu - \Lambda)(x, y)|}{d(x, y)^{\gamma - 1}} \lesssim [U]_{\gamma\text{-3pt};B} + [U]_{\gamma;B},$$

which implies (3.13). □

We need the following minor modification of [FR17, Theorem 3.1] that allows us to exchange limits with the fractional heat operator. We omit the proof since it is analogous to the referenced one.

**LEMMA 3.10.** *Let $s \in (0, 1)$, $\{u_k\}_{k \in \mathbb{Z}^+} \subset C((-\infty, 0) \times \mathbb{R}^d)$ and $\{f_k\}_{k \in \mathbb{Z}^+} \subset \mathscr{D}'((-\infty, 0) \times \mathbb{R}^d)$ be such that*

$$\mathscr{L} u_k = f_k \quad \text{in} \quad \mathscr{D}'(I \times D) \qquad \forall k \in \mathbb{Z}^+$$

*for a given open bounded interval $I \subset (-\infty, 0]$ and a bounded open domain $D \subset \mathbb{R}^d$. Suppose that for some function $u: (-\infty, 0] \times \mathbb{R}^d \to \mathbb{R}$ and distribution $f \in \mathscr{D}'(I \times D)$ the following hypotheses hold:*

  i. *$u_k \to u$ uniformly in compact sets of $(-\infty, 0) \times \mathbb{R}^d$;*

  ii. *$f_k \to f$ weakly in $\mathscr{D}'(I \times D)$;*

  iii. *$\sup_{t \in I} |u_k(t, x)| \lesssim (1 + |x|^{2s - \epsilon})$ for some $\epsilon > 0$ and for all $x \in \mathbb{R}^d$.*

*Then u satisfies*

$$\mathscr{L} u = f \quad in \quad \mathscr{D}'(I \times D).$$

**Proof of Lemma 3.4.** We consider $\gamma \in (1, 2s)$ since the case $\gamma \in (0,1)$ is simpler and follows with the same argument. We argue by contradiction using a blow-up argument. Suppose there exists $\delta > 0$ and a sequence of germs $\{U_k\}_{k \in \mathbb{N}}$ (and their respective $\{\nu_k\}_{k \in \mathbb{N}}$ and $\{\Lambda_k\}_{k \in \mathbb{N}}$) such that

$$\delta [U_k]_{\gamma;(-1,0]\times\mathbb{R}^d} + k \left( \|\mathscr{L} U_k\|_{\gamma-2s;B_1} + [U_k]_{\gamma\text{-3pt};B_1} + \|U_k\|_{B_1} + \|\nu_k\|_{B_1} \right) < [U_k]_{\gamma;B_{1/2}}. \tag{3.15}$$

We recall that by Remark 3.2 all the germs we consider satisfy $U_k(x,x) = 0$.

*Step 1.(Defining the parameters)* In particular, $0 < [U_k]_{\gamma;B_{1/2}}$ and therefore there exists $x_k, y_k \in B_{1/2}$ such that $y_{0,k} \leq x_{0,k}$, and since $U_k(x,x) = 0$,

$$0 < \frac{1}{2}[U_k]_{\gamma;B_{1/2}} < \frac{|U_k(x_k, y_k) - \nu_k(x_k) \cdot (y_k - x_k)|}{d(x_k, y_k)^\gamma} \leq r_k^{-\gamma} \left( \|U_k\|_{B_{1/2}} + r_k \|\nu_k\|_{B_{1/2}} \right) \tag{3.16}$$

for $r_k := d(x_k, y_k) \leq 1$. Using (3.15) we conclude that

$$\frac{1}{2}[U_k]_{\gamma;B_{1/2}} < r_k^{-\gamma} (1 + r_k) k^{-1} [U_k]_{\gamma;B_{1/2}} \implies \frac{1}{2} k < r_k^{-\gamma} + r_k^{-(\gamma-1)},$$

and either $r_k^{-\gamma} \to +\infty$ or $r_k^{-(\gamma-1)} \to +\infty$, and since $\gamma > 1$ any of these implies that $r_k \to 0$.

*Step 2.(Constructing the blowing-up sequence)* We define the function $v_k: \mathbb{R}^{1+d} \to \mathbb{R}$ as

$$v_k(y) := \frac{U_k(x_k, x_k + r_k y) - r_k \nu_k(x_k) \cdot y}{r_k^\gamma [U_k]_{\gamma;(-1,0]\times\mathbb{R}^d}}.$$

Define $\xi_k := r_k^{-1}(y_k - x_k)$, and observe that $d(0, \xi_k) = r_k^{-1} d(x_k, y_k) = 1$. To conclude that $\xi_k \in B_1$, the parabolic ball looking into the past, it is enough to observe that the definition of the seminorm (3.1) guarantees that $y_{0,k} \leq x_{0,k}$, which implies that $\xi_{0,k} \leq 0$. Moreover, by (3.15) we have that

$$v_k(\xi_k) = v_k(r_k^{-1}(y_k - x_k)) = \frac{U_k(x_k, y_k) - \nu_k(x_k) \cdot (y_k - x_k)}{d(x_k, y_k)^\gamma [U_k]_{\gamma;(-1,0]\times\mathbb{R}^d}},$$

and therefore, by definition of $x_k, y_k$ in (3.16), and (3.15) we conclude that

$$|v_k(\xi_k)| > \frac{[U_k]_{\gamma;B_{1/2}}}{2 [U_k]_{\gamma;(-1,0]\times\mathbb{R}^d}} > \frac{\delta}{2}. \tag{3.17}$$

*Step 3.(Showing uniform convergence)* We will show that, up to subsequences, $\{v_k\}_{k \in \mathbb{Z}^+}$ converges uniformly in compact subsets of $(-\infty, 0] \times \mathbb{R}^d$, and for this we use the Arzelà-Ascoli theorem. Fix $R > 0$ and let $k_R \in \mathbb{Z}^+$ be such that for all $k \geq k_R$:

$$\overline{B_R} \subset r_k^{-1}((-1 + 2^{-2s}, 0] \times B_{1/2}),$$

which is possible since $r_k \to 0$. Given $y \in \overline{B_R}$ we have that $r_k y \in ((-1 + 2^{-2s}, 0] \times B_{1/2})$ and since $x_k \in (-2^{-2s}, 0] \times B_{1/2}$:

$$x_{0,k} + r_k^{2s} y_0 \in (x_{0,k} - 1 + 2^{-2s}, x_{0,k}] \subset (-2^{-2s} - 1 + 2^{-2s}, 0] = (-1, 0],$$

and
$$|(x_k + r_k y)_{1:d}| \leqslant |(x_k)_{1:d}| + |r_k y_{1:d}| < \frac{1}{2} + \frac{1}{2} = 1,$$

and $x_k + r_k y \in B_1$. Moreover, since $y \in \overline{B_R}$ implies that $y_0 \leqslant 0$ then we have that the time component of $x_k + r_k y$ is in the past of $x_k$ for all $y \in \overline{B_R}$, and hence the pair $(x_k, x_k + r_k y)$ satisfy the condition of definition (3.1) of $[U_k]_{\gamma;B_1}$ and we can conclude that

$$\begin{aligned}|v_k(y)| &\leqslant \frac{|U_k(x_k, x_k + r_k y) - r_k \nu_k(x_k) \cdot y|}{r_k^\gamma [U_k]_{\gamma;(-1,0]\times\mathbb{R}^d}} \leqslant \frac{[U_k]_{\gamma;B_1} d(x_k, x_k + r_k y)^\gamma}{r_k^\gamma [U_k]_{\gamma;(-1,0]\times\mathbb{R}^d}} \\ &= \frac{[U_k]_{\gamma;B_1} r_k^\gamma d(0,y)^\gamma}{r_k^\gamma [U_k]_{\gamma;(-1,0]\times\mathbb{R}^d}} \leqslant d(0,y)^\gamma \leqslant R^\gamma,\end{aligned} \quad (3.18)$$

which implies that the sequence $\{v_k\}_{k \geq k_R}$ is uniformly bounded on $\overline{B_R}$. Given $y, z \in \overline{B_R}$, assume without loss of generality that $z_0 \leqslant y_0$, then by the same arguments as before we have that $x_k$, $x_k + r_k y, x_k + r_k z \in B_1$ and their time components satisfy $(x_k + r_k z)_0 \leqslant (x_k + r_k y)_0 \leqslant x_{0,k}$ we have (omitting the sub-index $k$ for readability)

$$\begin{aligned}&(v(z) - v(y)) r^\gamma [U]_{\gamma;(-1,0]\times\mathbb{R}^d} \\ &= U(x, x+rz) - r\nu(x)\cdot z - (U(x, x+ry) - r\nu(x)\cdot y) \\ &= U(x, x+rz) - U(x, x+ry) - U(x+ry, x+rz) + \Lambda(x, x+ry)\cdot(rz-ry) \\ &\quad + \nu(x)\cdot(ry-rz) + U(x+ry, x+rz) - \Lambda(x, x+ry)\cdot(rz-ry) \\ &= U(x, x+rz) - U(x, x+ry) - U(x+ry, x+rz) + \Lambda(x, x+ry)\cdot(rz-ry) \\ &\quad + U(x+ry, x+rz) - \nu(x+ry)(rz-ry) \\ &\quad + (\nu(x+ry) - \nu(x) - \Lambda(x, x+ry))\cdot(rz-ry),\end{aligned}$$

and we can conclude

$$\begin{aligned}&|v(z)-v(y)|\, r^\gamma [U]_{\gamma;(-1,0]\times\mathbb{R}^d} \\ &\leqslant [U]_{\gamma\text{-3pt};B_1} \sum_{\beta \in A} d(x,x+ry)^\beta d(x+ry, x+rz)^{\gamma-\beta} \\ &\quad + [U]_{\gamma;B_1} d(x+ry, x+rz)^\gamma + r|\nu(x+ry)-\nu(x)-\Lambda(x,x+ry)|\,|y_{1:d}-z_{1:d}| \\ &\leqslant [U]_{\gamma\text{-3pt};B_1} r^\gamma \sum_{\beta \in A} d(0,y)^\beta d(y,z)^{\gamma-\beta} + r^\gamma [U]_{\gamma;B_1} d(y,z)^\gamma + r^\gamma [\delta\nu - \Lambda]_{\gamma-1;B_1} d(0,y)^{\gamma-1} d(y,z) \\ &\lesssim r^\gamma ([U]_{\gamma\text{-3pt};B_1} + [U]_{\gamma;B_1} + [\delta\nu-\Lambda]_{\gamma-1;B_1}) \sum_{\beta \in A \cup \{0,\gamma-1\}} R^\beta d(y,z)^{\gamma-\beta} \\ &\lesssim r^\gamma ([U]_{\gamma\text{-3pt};B_1} + [U]_{\gamma;B_1}) \sum_{\beta \in A \cup \{0,\gamma-1\}} R^\beta d(y,z)^{\gamma-\beta},\end{aligned}$$

where the last inequality follows from Lemma 3.9. Using Assumption (3.15) we have for all $k \in \mathbb{Z}^+$:

$$\begin{aligned}|v_k(z)-v_k(y)|\, r_k^\gamma [U_k]_{\gamma;(-1,0]\times\mathbb{R}^d} &\lesssim r_k^\gamma \left(3k^{-1}[U_k]_{\gamma;B_{1/2}} + [U_k]_{\gamma;B_1}\right) \sum_{\beta \in A \cup \{0,\gamma-1\}} R^\beta d(y,z)^{\gamma-\beta} \\ &\lesssim r_k^\gamma [U_k]_{\gamma;(-1,0]\times\mathbb{R}^d} \sum_{\beta \in A \cup \{0,\gamma-1\}} R^\beta d(y,z)^{\gamma-\beta}.\end{aligned}$$

Since $A$ is finite and $A \cup \{0, \gamma-1\} \subset [0, \gamma)$ all powers $\gamma - \beta$ appearing in the sum are strictly positive and this implies equicontinuity on $\overline{B_R}$ of $\{v_k\}_{k\in\mathbb{Z}^+}$. We conclude by the Arzelà-Ascoli theorem that there exists $v_R \in C(\overline{B_R})$ such that a subsequence of $\{v_k\}_{k \geq k_R}$ converges uniformly on $\overline{B_R}$ to $v_R$. Since this is for all $R > 0$ and

$$\bigcup_{R \in \mathbb{Z}^+} \overline{B_R} = (-\infty, 0] \times \mathbb{R}^d,$$

by a diagonal argument we can conclude that there exists $v \in C((-\infty, 0] \times \mathbb{R}^d)$ such that $\{v_k\}_{k \in \mathbb{Z}^+}$ converges uniformly on compact subsets of $(-\infty, 0] \times \mathbb{R}^d$ to $v$.

*Step 4. (Harmonicity of the limit)* We prove that $\{\mathscr{L} v_k\}_{k \in \mathbb{Z}^+}$ converges weakly to 0 in $(-\infty, 0] \times \mathbb{R}^d$. Since $\mathscr{L}$ is translation invariant, $2s$-homogeneous and annihilates linear in space polynomials, we can conclude that

$$r_k^\gamma [U_k]_{\gamma; (-1,0] \times \mathbb{R}^d} (\mathscr{L} v_k) = \mathscr{L}(U_k(x_k, x_k + r_k \cdot) - r_k \nu(x_k) \cdot (\cdot)_{1:d})$$
$$= r_k^{2s} \mathscr{L}(U_k(x_k, \cdot))(x_k + r_k \cdot),$$

and hence for all $\psi \in \mathcal{D}((0, \infty) \times \mathbb{R}^d)$ we have that

$$r_k^{\gamma - 2s}[U_k]_{\gamma; (-1,0] \times \mathbb{R}^d} |\langle \mathscr{L} v_k, \psi \rangle| = r_k^{-2s} |\langle r_k^{2s} \mathscr{L}(U_k(x_k, \cdot))(x_k + r_k \cdot), \psi \rangle|$$
$$= |\langle \mathscr{L}(U_k(x_k, \cdot)), \psi_{x_k}^{r_k} \rangle|.$$

For $k \in \mathbb{Z}^+$ big enough we have that $\mathrm{supp}(\psi_{x_k}^{r_k}) \subset B_1$ since $x_k \in B_{1/2}$ and $r_k \to 0$, and therefore by definition of the seminorm $\|\mathscr{L} U_k\|_{\gamma - 2s; B_1}$ and (3.15) we have that

$$|\langle \mathscr{L}(U_k(x_k, \cdot)), \psi_{x_k}^{r_k} \rangle| \leqslant r_k^{\gamma - 2s} \|\mathscr{L} U_k\|_{\gamma - 2s; B_1} \leqslant \frac{r_k^{\gamma - 2s}}{k} [U_k]_{\gamma; B_{1/2}}.$$

We conclude that as $k \to \infty$

$$|\langle \mathscr{L} v_k, \psi \rangle| \leqslant \frac{1}{k} \frac{[U_k]_{\gamma; B_{1/2}}}{[U_k]_{\gamma; (-1,0] \times \mathbb{R}^d}} \leqslant \frac{1}{k} \frac{[U_k]_{\gamma; (-1,0] \times \mathbb{R}^d}}{[U_k]_{\gamma; (-1,0] \times \mathbb{R}^d}} \leqslant \frac{1}{k} \to 0,$$

which concludes the proof of the claimed weak convergence.

To apply Lemma 3.10 to $\{v_k\}_{k \in \mathbb{Z}^+}$ we need to show the required growth condition. Let $I \subset (-\infty, 0]$ be a bounded interval and $k_I \in \mathbb{Z}^+$ such that $I \subset r_k^{-2s}(-1 + 2^{-2s}, 0]$. Analogously to (3.18) this condition implies that for all $k_I \leqslant k \in \mathbb{Z}^+$ and $y \in I \times \mathbb{R}^d$ we have $x_k + r_k y \in (-1, 0) \times \mathbb{R}^d$ and $|v_k(y)| \leqslant d(0, y)^\gamma$, and therefore for all $y_{1:d} \in \mathbb{R}^d$:

$$\sup_{y_0 \in I} |v_k(y_0, y_{1:d})| \leqslant \sup_{y_0 \in I} d(0, (y_0, y_{1:d}))^\gamma = \sup_{t \in I} \max\left\{|y_0|^{\frac{\gamma}{2s}}, |y_{1:d}|^\gamma\right\} \leqslant C(1 + |y_{1:d}|^\gamma), \quad (3.19)$$

for some constant $C > 0$ that depends only on the bounded interval $I \subset (-\infty, 0]$. We see that for $\varepsilon = -(\gamma - 2s) > 0$ this is the required growth condition. Therefore, for every bounded open domain $I \times D \subseteq (-\infty, 0] \times \mathbb{R}^d$ we can conclude that $\mathscr{L} v = 0$ in $\mathscr{D}'(I \times D)$ and therefore $\mathscr{L} v = 0$ in $\mathscr{D}'((-\infty, 0] \times \mathbb{R}^d)$.

*Step 5. (Triviality of the limit)* Since $\{v_k\}_{k \in \mathbb{Z}^+}$ converges uniformly to $v$ on $\overline{B_R}$ we can take limits in (3.19) to see that $v$ satisfies the growth condition of the Liouville-type Theorem (Theorem 3.3) for $\beta = \gamma$. This allows us to conclude that $v(x_0, x_{1:d})$ is a polynomial in its spatial variables of degree at most $\lfloor \gamma \rfloor = 1$ in $(-\infty, 0) \times \mathbb{R}^d$ and by continuity in $(-\infty, 0] \times \mathbb{R}^d$.

Since $U(z_k, z_k) = 0$ by assumption, then $v_k(0) = 0$ for all $k \in \mathbb{Z}^+$ and taking the limit we conclude that $v(0) = 0$. We claim that also $\nabla v(0) = 0$ which will imply that $v \equiv 0$ in $(-\infty, 0] \times \mathbb{R}^d$. Given $j \in \{1, \ldots d\}$ and $e_j \in \mathbb{R}^{1+d}$ the canonical vector, observe that for all $h \in [-1, 1]$ we have that the time component of $x_k + r_k h e_j$ is equal to $x_0 \in (-2^{-2s}, 0]$ and therefore $x_k + r_k h e_j \in \overline{B_1}$. Arguing as in (3.18) we have that $|v_k(h e_j)| \leqslant d(0, h e_j)^\gamma$, and using that $v_k$ converges uniformly on $\overline{B_1}$ we obtain

$$h^{-1}|v(h e_j) - v(0)| = h^{-1}|v(h e_j)| = \lim_{k \to \infty} h^{-1}|v_k(h e_j)| \leqslant \lim_{k \to \infty} h^{-1} d(0, h e_j)^\gamma = \lim_{k \to \infty} h^{\gamma - 1} = h^{\gamma - 1},$$

which, since $\gamma - 1 > 0$, allows us to conclude that

$$|\nabla v(0)| = \lim_{h \to 0} h^{-1} |v(h\,e_j) - v(0)| \leqslant \lim_{h \to 0} h^{\gamma - 1} = 0.$$

Therefore, $v$ is a linear in space polynomial in $(-\infty, 0] \times \mathbb{R}^d$ which satisfies $v(0) = 0$ and $\nabla v(0) = 0$, which implies that $v \equiv 0$ in $(-\infty, 0] \times \mathbb{R}^d$.

*Step 6.* (*Reaching the contradiction*) Taking limits in (3.17), justified by the uniform convergence and $\xi_k \to \xi \in \overline{B_1} \subset (-\infty, 0] \times \mathbb{R}^d$, implies that $0 = v(\xi) = \lim_{k \to \infty} v_k(\xi_k) \geqslant \delta/2 > 0$, a contradiction. □

The following result that allows us to control some gradient terms appearing in the estimates, and will also be useful in the following section.

**LEMMA 3.11.** *Let $B \subset \mathbb{R}^{1+d}$ be a domain satisfying a $(\beta, r_0)$-spatial cone condition. Let $U$ be a germ such that for some $\gamma > 0$ we have $[U]_{\gamma;B} < +\infty$, and let $\nu$ be as in (3.2). Then for all $r \in [0, r_0]$*

$$\|\nu\|_B \lesssim r^{\gamma - 1} [U]_{\gamma;B} + r^{-1} \|U\|_B, \qquad (3.20)$$

*with the proportionality constant depending on $\beta$ and $\gamma$. Moreover, assuming that*

$$\left(\frac{\|U\|_B}{[U]_{\gamma;B}}\right)^{\gamma^{-1}} \leqslant r_0 \qquad (3.21)$$

*we can conclude that for every $\delta > 0$ there exists $C = C(d, \delta, \beta, \gamma) > 0$ such that*

$$\|\nu\|_B \leqslant \delta \, [U]_{\gamma;B} + C \, \|U\|_B.$$

**Proof.** Fix $x \in B$. Since $B$ satisfies a spatial cone condition for $\beta \in (0,1)$ and $r_0 > 0$ for all $r \in [0, r_0]$ there exists $y_{1:d} \in \mathbb{R}^d$ such that $y := (x_0, y_{1:d}) \in B$, $d(x,y) = |x_{1:d} - y_{1:d}| = r$ and $|\nu(x) \cdot (y_{1:d} - x_{1:d})| \geqslant \beta r \, |\nu(x)|$. Hence,

$$\beta r \, |\nu(x)| \leqslant |\nu(x) \cdot (y - x)| \leqslant |U(x,y) - \nu(x) \cdot (y - x)| + |U(x,y)| \leqslant r^\gamma [U]_{\gamma;B} + \|U\|_B.$$

Since this last bound is independent of $x \in B$ we conclude that for all $r \in [0, r_0]$:

$$\beta r \, \|\nu\|_B \leqslant r^\gamma [U]_{\gamma;B} + \|U\|_B. \qquad (3.22)$$

Choose $r > 0$ such that both terms in the right hand side are comparable, i.e.,

$$r^\gamma [U]_{\gamma;B} = \|U\|_B \iff r := \left(\frac{\|U\|_B}{[U]_{\gamma;B}}\right)^{\gamma^{-1}}.$$

By assumption $r \leqslant r_0$, and then we can use (3.22) and obtain $\beta r \, \|\nu\|_B \leqslant 2 \, \|U\|_B$, and

$$\|\nu\|_B \leqslant 2 \beta^{-1} r^{-1} \|U\|_B = 2 \beta^{-1} \left(\frac{\|U\|_B}{[U]_{\gamma;B}}\right)^{-\gamma^{-1}} \|U\|_B = 2 \beta^{-1} \|U\|_B^{1 - \gamma^{-1}} [U]_{\gamma;B}^{\gamma^{-1}}.$$

Since $\gamma > 1$ then $\gamma^{-1} \in (0,1)$ and we can use an interpolation between these terms to obtain that for all $\delta > 0$ there exists $C(\delta, \beta) > 0$ such that $\|\nu\|_B \leqslant \delta [U]_{\gamma;B} + C \|U\|_B$, which proves the result. □

Using Lemma 3.11 we are able to bound the gradient and conclude the following:

**Corollary 3.12**. *If condition* (3.21) *of Lemma 3.11 holds, then Lemma 3.4 can be improved to*

$$[U]_{\gamma;B_{1/2}} \leqslant \delta [U]_{\gamma;(-1,0]\times\mathbb{R}^d} + C\left(\|\mathscr{L}U\|_{\gamma-2s;B_1} + [U]_{\gamma\text{-3pt};B_1} + \|U\|_{B_1}\right). \tag{3.23}$$

## 3.2. Proof of the abstract absorption Lemma 3.6

Let $B, \{B_i\}_{i\in I}$ be convex subsets of $\mathbb{R}^{1+d}$ such that $B \subset \bigcup_{i\in I} B_i$. For $x, y \in B$ denote by $\overline{xy}$ the line segment between $x$ and $y$, then by convexity $\overline{xy} \subset B$. We define

$$m_{xy} := \min\left\{n \in \mathbb{N} : \begin{array}{l}\exists \{x_0, \ldots, x_n\} \subset \overline{xy} \quad \text{such that}\quad x_0 = x, x_n = y,\\ \text{and}\ \forall j \in \{1, \ldots n\}\ \exists i_j \in I \quad \text{such that}\quad \overline{x_i x_{i+1}} \subset B_{i_j}\end{array}\right\}.$$

Then $m_{xy}$ is the minimum number of segments in which $\overline{xy}$ has to be split in order for each of these segments to be contained completely in one of the sets $B_i$. Observe that by convexity of the sets and minimality of $m_{xy}$ each segment $\overline{x_i x_{1+i}}$ is contained in different sets $B_{j_i}$. We have the following result.

**Lemma 3.13**. *Fix $\gamma \in (0, 2s)\setminus\{1\}$ and $A \subset (0, \gamma)$ finite. Let $B, \{B_i\}_{i\in I}$ be a finite number of half-parabolic balls such that $B \subset \bigcup_{i\in I} B_i$, and let $U$ be a germ such that $[U]_{\gamma;B} < +\infty$ and $\Lambda$ another germ such that $U$ satisfies the 3-point continuity condition* (3.6) *on $B$ with the set $A$. Given $x, y \in B$ let $\{x_0, \ldots, x_{m_{xy}}\} \subset \overline{xy}$ be such that $x_0 = x, x_{m_{xy}} = y$ and for each $j \in \{1, \ldots n\}$ let $i_j \in I$ be such that $x_j, x_{j+1} \in B_{i_j}$, then we have that*

$$|U(x,y) - \mathbb{1}_{\gamma>1}\nu(x)\cdot(y-x)| \lesssim \left(\sum_{j=0}^{m_{xy}-1}[U]_{\gamma;B_{i_j}} + [U]_{\gamma\text{-3pt};B} + \mathbb{1}_{\gamma>1}[\Lambda]_{(\gamma-1)\text{-3pt};B}\right)d(x,y)^\gamma \tag{3.24}$$

*for an implicit proportionality constant that depends on $m_{x,y}, s, d, \gamma$ and $A$.*

**Proof.** *Step 1.* We show the result for $\gamma \in (0,1)$ first, recalling that for this case there are no $\nu$ or $\Lambda$. The proof follows by induction over the sets

$$D_n = \{(x,y) \in B\times B : x_0 \leqslant y_0, m_{xy} \leqslant n\}.$$

Given $(x,y) \in D_1$ we have that $m_{xy} = 1$ and therefore there exists $i \in I$ such that $\overline{xy} \subset B_i$ and therefore

$$|U(x,y)| \leqslant [U]_{\gamma;B_i} d(x,y)^\gamma.$$

Assume the result is true for $D_n$ and let $(x,y) \in D_{n+1}\setminus D_n$. Let $\{x_0, \ldots, x_{n+1}\} \subset \overline{xy}$ be such that $x_0 = x, x_{n+1} = y$ and $\overline{x_j, x_{j+1}} \in B_{i_j}$. We have then that $(x_0, x_n) \in D_n$, since otherwise one would violate the minimality of the definition of $m_{xy}$. Moreover, their time components satisfy the right order since they are in the line segment $\overline{xy}$. We consider the induction hypothesis

$$|U(x_0, x_n)| \leqslant \left(\sum_{j=0}^{n-1}[U]_{\gamma;B_{i_j}} + \#A\,(n-1)\,[U]_{\gamma\text{-3pt};B}\right)d(x,y)^\gamma. \tag{3.25}$$

Now we use the three point continuity to bound the term $U(x_0, x_{n+1})$:

$$\begin{aligned}
|U(x_0,x_{n+1})| &\leq |U(x_0,x_n)| + |U(x_n,x_{n+1})| + |U(x_0,x_{n+1}) - U(x_0,x_n) - U(x_n,x_{n+1})| \\
&\leq \left(\sum_{j=0}^{n-1} [U]_{\gamma;B_{i_j}} + \#A\,(n-1)\,[U]_{\gamma\text{-3pt};B}\right) d(x_0,x_n)^\gamma + [U]_{\gamma;B_{i_n}} d(x_n,x_{n+1})^\gamma \\
&\quad + [U]_{\gamma\text{-3pt};B} \sum_{\beta \in A} d(x_0,x_n)^\beta d(x_n,x_{n+1})^{\gamma-\beta} \\
&\leq \left(\sum_{j=0}^{n-1} [U]_{\gamma;B_{i_j}} + \#A\,(n-1)\,[U]_{\gamma\text{-3pt};B}\right) d(x_0,x_{n+1})^\gamma + [U]_{\gamma;B_{i_n}} d(x_0,x_{n+1})^\gamma \\
&\quad + [U]_{\gamma\text{-3pt};B} \sum_{\beta \in A} d(x_0,x_{n+1})^\beta d(x_0,x_{n+1})^{\gamma-\beta} \\
&\leq \left(\sum_{j=0}^{n} [U]_{\gamma;B_{i_j}} + \#A\,n\,[U]_{\gamma\text{-3pt};B}\right) d(x_0,x_{n+1})^\gamma,
\end{aligned}$$

which concludes the result for $\gamma \in (0,1)$ by ignoring the constant $\#A\,n$.

*Step 2.* Assume $\gamma \in (1, 2s)$. We proceed analogously by induction, the case $n=1$ following by definition. Assume the result is true for $D_n$ and let $(x,y) \in D_{n+1} \setminus D_n$. Let $\{x_0, \ldots, x_{n+1}\} \subset \overline{xy}$ be such that $x_0 = x$, $x_{n+1} = y$ and $\overline{x_j, x_{j+1}} \in B_{i_j}$. We have then that $(x_0, x_n) \in D_n$ and by induction hypothesis

$$|U(x_0,x_n) - \nu(x_n) \cdot (x_n - x_0)| \lesssim \left(\sum_{j=0}^{n-1} [U]_{\gamma;B_{i_j}} + [U]_{\gamma\text{-3pt};B} + [\nu - \Lambda]_{(\gamma-1)\text{-3pt};B}\right) d(x_0,x_n)^\gamma.$$

Moreover, we have the following identity:

$$\begin{aligned}
U(x,y) - \nu(x) \cdot (y-x) &= U(x_0,x_{n+1}) - \nu(x_0) \cdot (x_{n+1} - x_0) \\
&= U(x_0,x_n) - \nu(x_0) \cdot (x_n - x_0) + U(x_n,x_{n+1}) - \nu(x_n) \cdot (x_{n+1} - x_n) \\
&\quad + U(x_0,x_{n+1}) - U(x_0,x_n) - U(x_n,x_{n+1}) + \Lambda(x_0,x_n) \cdot (x_{n+1} - x_n) \\
&\quad + (\nu(x_n) - \nu(x_0) - \Lambda(x_0,x_n)) \cdot (x_{n+1} - x_n),
\end{aligned}$$

and therefore we have the estimate

$$\begin{aligned}
|U(x,y) - \nu(x) \cdot (y-x)| &\leq |U(x_0,x_n) - \nu(x_0) \cdot (x_n - x_0)| + |U(x_n,x_{n+1}) - \nu(x_n) \cdot (x_{n+1} - x_n)| \\
&\quad + |U(x_0,x_{n+1}) - U(x_0,x_n) - U(x_n,x_{n+1}) + \Lambda(x_0,x_n) \cdot (x_{n+1} - x_n)| \\
&\quad + |\nu(x_n) - \nu(x_0) - \Lambda(x_0,x_n)|\,|x_{n+1} - x_n|.
\end{aligned}$$

The first term can be bounded using the induction hypothesis ($m_{x_0,x_n} = n$), the second the definition of $[U]_{\gamma;B_{i_n}}$ since $x_n, x_{n+1} \in B_{i_n}$, and the third using the three point continuity of $U$. For the last term observe that $\nu - \Lambda$ is a germ and since $\gamma - 1 \in (0, 2s-1) \subset (0,1)$ we can use (3.25) to obtain

$$|\nu(x_n) - \nu(x_0) - \Lambda(x_0,x_n)| \lesssim \left(\sum_{j=0}^{n-1} [\nu-\Lambda]_{\gamma-1;B_{i_j}} + [\nu-\Lambda]_{(\gamma-1)\text{-3pt};B}\right) d(x_0,x_n)^{\gamma-1}.$$

By Lemma 3.9 we have that

$$[\nu - \Lambda]_{\gamma-1;B_i} \lesssim [U]_{\gamma;B_i} + [U]_{\gamma\text{-3pt};B_i} \leq [U]_{\gamma;B_i} + [U]_{\gamma\text{-3pt};B},$$

and therefore

$$|\nu(x_n) - \nu(x_0) - \Lambda(x_0,x_n)| \lesssim \left(\sum_{j=0}^{n-1} [U]_{\gamma;B_{i_j}} + [U]_{\gamma\text{-3pt};B} + [\nu-\Lambda]_{(\gamma-1)\text{-3pt};B}\right) d(x_0,x_n)^{\gamma-1}.$$

We conclude that

$$
\begin{aligned}
&|U(x,y) - \nu(x) \cdot (y-x)| \\
&\lesssim \left( \sum_{j=0}^{n-1} [U]_{\gamma;B_{i_j}} + [U]_{\gamma\text{-3pt};B} + [\nu - \Lambda]_{(\gamma-1)\text{-3pt};B} \right) d(x_0, x_n)^\gamma \\
&\quad + [U]_{\gamma;B_{i_n}} d(x_n, x_{n+1})^\gamma + [U]_{\gamma\text{-3pt};B} \sum_{\beta \in A} d(x_0, x_n)^\beta d(x_n, x_{n+1})^{\gamma-\beta} \\
&\quad + \left( \sum_{j=0}^{n-1} [U]_{\gamma;B_{i_j}} + [U]_{\gamma\text{-3pt};B} + [\nu - \Lambda]_{(\gamma-1)\text{-3pt};B} \right) d(x_0, x_n)^{\gamma-1} |x_{n+1} - x_n| \\
&\lesssim \left( \sum_{j=1}^{n} [U]_{\gamma;B_{i_j}} + [U]_{\gamma\text{-3pt};B} + [\nu - \Lambda]_{(\gamma-1)\text{-3pt};B} \right) d(x,y)^\gamma.
\end{aligned}
$$

We conclude by recalling that $[\nu - \Lambda]_{(\gamma-1)\text{-3pt};B} = [\Lambda]_{(\gamma-1)\text{-3pt};B}$ by (3.9) and (3.7). $\square$

**Proof of Lemma 3.5.** Given $x, y \in B$ let $\{x_0, \ldots, x_{m_{xy}}\} \subset \overline{xy}$ be such that $x_0 = x$, $x_{m_{xy}} = y$ and for each $j \in \{1, \ldots, n\}$ let $i_j \in I$ be such that $x_j, x_{j+1} \in B_{i_j}$. Recall that by convexity of the sets $B_{i_j}$ and minimality of $m_{xy}$ we have that the segments $\overline{x_i x_{1+i}}$ are contained in different sets $B_{j_i}$. This implies that $m_{xy} \leq |I|$ and that the sum in (3.24) is over different subsets $B_{j_i}$ and therefore:

$$
\begin{aligned}
|U(x,y)| &\lesssim \left( \sum_{j=1}^{m_{xy}} [U]_{\gamma;B_{i_j}} + [U]_{\gamma\text{-3pt};B} + [\Lambda]_{(\gamma-1)\text{-3pt};B} \right) d(x,y)^\gamma \\
&\lesssim \left( \sum_{i \in I} [U]_{\gamma;B_i} + [U]_{\gamma\text{-3pt};B} + [\Lambda]_{(\gamma-1)\text{-3pt};B} \right) d(x,y)^\gamma,
\end{aligned}
$$

and the result follows. $\square$

Now we prove our generalisation of the abstract absorption Lemma.

**Proof of Lemma 3.6.** The proof follows the one in [Sim97] with the only difference being the change of the sub-additivity condition (3.12). Let $B = B_{\rho_0}(y_0)$. It is enough to show the result for $\theta = 1/2$, because for any $\theta \in (1/2, 1)$ and any ball $B_\sigma(y) \subset B$ we can find $\{y_i\}_{i=1}^N \subset B_{\theta\sigma}(y)$ such that $B_{(1-\theta)\sigma}(y_i) \subset B_\sigma(y)$ and $B_{\theta\sigma}(y) \subset \bigcup_{i=1}^N B_{(1-\theta)\sigma/2}(y)$, with $N \in \mathbb{N}$ depending only on the ratio of the radii of the balls $(\theta\sigma)/((1-\theta)\sigma/2) = 2\theta(1-\theta)^{-1}$ and the geometry of the space which is determined by the metric, which depends on $s$ and the dimension $d+1$, i.e., $N(\theta, s, d)$.

We proceed to prove the result for $\theta = 1/2$. Define

$$
Q := \sup_{B_\sigma(y) \subset B} \sigma^\gamma S(B_{\sigma/2}(y)) \leq \rho_0^\gamma S(B) < +\infty.
$$

Then we have for every $B_\sigma(y) \subset B$

$$
\left( \frac{\sigma}{2} \right)^\gamma S\big(B_{\theta_0 \frac{\sigma}{2}}(y)\big) \leq \varepsilon Q + E. \tag{3.26}
$$

For $B_\sigma(y) \subset B$ fixed we consider a finite cover of $B_{\frac{\sigma}{2}}(y)$ of smaller balls of radius $\frac{\theta_0 \sigma}{4}$, i.e., let $\{y_i\}_{i=1}^N \subset B_{\frac{\sigma}{2}}(y)$ be such that

$$
B_{\frac{\sigma}{2}}(y) \subset \bigcup_{i=1}^N B_{\frac{\theta_0 \sigma}{4}}(y_i).
$$

As before, the number of balls needed for this cover depends only on the ratio of the radii, which depends only on $\theta_0$, and the geometry of the space, i.e., $N(\theta_0, s, d)$. Moreover, since $y_i \in B_{\frac{\sigma}{2}}(y)$ then $B_{\frac{\sigma}{2}}(y_i) \subset B_\sigma(y) \subset B$, and by (3.26) applied to $B_{\frac{\sigma}{2}}(y_i)$ instead of $B_\sigma(y)$ we have that

$$\left(\frac{\sigma}{4}\right)^\gamma S\big(B_{\theta_0 \frac{\sigma}{4}}\big) \leqslant \varepsilon Q + E.$$

Then the *almost* subadditivity (3.12) and the monotonicity of $D$ allows us to conclude that

$$\begin{aligned}\sigma^\gamma S\big(B_{\frac{\sigma}{2}}(y)\big) &\lesssim \sum_{i=1}^N \sigma^\gamma S\Big(B_{\frac{\theta_0 \sigma}{4}}(y_i)\Big) + \sigma^\gamma D\big(B_{\frac{\sigma}{2}}(y)\big) \\ &\lesssim \varepsilon Q + E + D(B),\end{aligned}$$

for an implicit constant independent of $B_\sigma(y)$. Taking supremum over $B_\sigma(y) \subset B$ we conclude that

$$Q = \sup_{B_\sigma(y) \subset B} \sigma^\gamma S(B_{\sigma/2}(y)) \lesssim \varepsilon Q + E + D(B). \tag{3.27}$$

Choosing $\varepsilon > 0$ small enough, depending on the implicit proportionality constant in (3.27), we can absorb the $Q$ term to the LHS and conclude the result. $\square$

### 3.3. Proof of the Schauder estimate

**Proof of Theorem 3.1.** If (3.21) does not hold for $B = B_1$ then we would have ($r_0 = 1$)

$$1 \leqslant \left(\frac{\|U\|_{B_1}}{[U]_{\gamma;B_1}}\right)^{\gamma^{-1}} \iff [U]_{\gamma;B_1} \leqslant \|U\|_{B_1},$$

which in particular implies the result, and therefore from now we will on assume that (3.21) holds.

Fix $\eta \in C_c^\infty(\mathbb{R}^d)$ with $\mathrm{supp}(\eta) \subseteq B_2$, $\eta \equiv 1$ in $B_{3/2}$, $0 \leqslant \eta \leqslant 1$ and $\|\nabla \eta\|_{L^\infty(\mathbb{R}^d)} \leqslant 2$. We extend it to $\eta: \mathbb{R}^{d+1} \to \mathbb{R}$ as $\eta(x) = \eta(x_{1:d})$. We have that

$$|\eta(x) - \eta(y)| = |\eta(x_{1:d}) - \eta(y_{1:d})| \leqslant \|\nabla \eta\|_{L^\infty(\mathbb{R}^d)} |x_{1:d} - y_{1:d}| \lesssim d(x,y).$$

We apply Corollary 3.12 to the germ $\tilde{U}$ defined as $\tilde{U}(x,y) := \eta(x)\eta(y) U(x,y)$ to conclude that for every $\delta > 0$ there exists $C_\delta > 0$ such that ((3.21) also holds for $\tilde{U}$ since $\eta \equiv 1$ in $B_{3/2} \supset B_1$)

$$[\tilde{U}]_{\gamma;B_{1/2}} \leqslant \delta [\tilde{U}]_{\gamma;(-1,0)\times \mathbb{R}^d} + C_\delta \big(\|\mathscr{L}\tilde{U}\|_{\gamma-2s;B_1} + [\tilde{U}]_{\gamma\text{-3pt};B_1} + \|\tilde{U}\|_{B_1}\big). \tag{3.28}$$

Observe that for all $x \in \mathbb{R}^{1+d}$ the diagonal derivative of $\tilde{U}$ satisfies:

$$\tilde{\nu}(x) := (\nabla \tilde{U}(x,\cdot))(x) = \eta(x) \nabla(\eta U(x,\cdot))(x) = \eta(x)(\nabla \eta(x) U(x,x) + \eta(x)\nu(x)) = \eta^2(x)\nu(x), \tag{3.29}$$

since we assumed that $U(x,x) = 0$. We write each of the terms in (3.28) in terms of $U$. Since for all $x, y \in B_{3/2}$ we have $x_{1:d}, y_{1:d} \in B_{3/2} \subset \mathbb{R}^d$ and $\tilde{U}(x,y) = \eta(x)\eta(y) U(x,y) = U(x,y)$, then

$$[\tilde{U}]_{\gamma;B_{1/2}} = [U]_{\gamma;B_{1/2}}, \qquad [\tilde{U}]_{\gamma\text{-3pt};B_1} = [U]_{\gamma\text{-3pt};B_1}, \qquad \|\tilde{U}\|_{B_1} = \|U\|_{B_1}. \tag{3.30}$$

For the first term in the right hand side of (3.28), we will first show that

$$[\tilde{U}]_{\gamma;(-1,0]\times \mathbb{R}^d} \leqslant [U]_{\gamma;(-1,0]\times B_2} + C \|\nu\|_{(-1,0]\times B_2},$$

which then by (3.13) from Lemma 3.11, with $\delta = 1$, will imply

$$[\tilde{U}]_{\gamma;(-1,0]\times\mathbb{R}^d} \leqslant 2\,[U]_{\gamma;(-1,0]\times B_2} + C\,\|U\|_{(-1,0]\times B_2}. \tag{3.31}$$

Equation (3.29) implies

$$\begin{aligned}
[\tilde{U}]_{\gamma;(-1,0]\times\mathbb{R}^d} &= \sup_{x,y\in(-1,0)\times\mathbb{R}^d} \frac{|\eta(x)|\,|\eta(y)\,U(x,y) - \eta(x)\,\nu(x)\cdot(y-x)|}{d(x,y)^\gamma} \\
&\leqslant \sup_{\substack{x\in(-1,0)\times B_2 \\ y\in(-1,0)\times\mathbb{R}^d}} \left(\frac{|\eta(y)|\,|U(x,y) - \nu(x)\cdot(y-x)| + |\eta(y)-\eta(x)|\,|\nu(x)|\,|y-x|}{d(x,y)^\gamma}\right) \\
&\leqslant \sup_{x,y\in(-1,0)\times B_2} \frac{|U(x,y) - \nu(x)\cdot(y-x)|}{d(x,y)^\gamma} \\
&\quad + \sup_{\substack{x\in(-1,0)\times B_2 \\ y\in(-1,0)\times B_3}} \frac{|\eta(y)-\eta(x)|\,|\nu(x)|\,d(x,y)}{d(x,y)^\gamma} + \sup_{\substack{x\in(-1,0)\times B_2 \\ y\notin(-1,0)\times B_3}} \frac{|\eta(y)-\eta(x)|\,|\nu(x)|\,d(x,y)}{d(x,y)^\gamma} \\
&\leqslant [U]_{\gamma;(-1,0]\times B_2} + \sup_{\substack{x\in(-1,0)\times B_2 \\ y\in(-1,0)\times B_3}} \frac{d(x,y)\,|\nu(x)|\,d(x,y)}{d(x,y)^\gamma} + \sup_{\substack{x\in(-1,0)\times B_2 \\ y\notin(-1,0)\times B_3}} \frac{|\eta(x)|\,|\nu(x)|\,d(x,y)}{d(x,y)^\gamma} \\
&\leqslant [U]_{\gamma;(-1,0]\times B_2} + \|\nu\|_{(-1,0]\times B_2}\left(\sup_{\substack{x\in(-1,0)\times B_2 \\ y\in(-1,0)\times B_3}} d(x,y)^{2-\gamma} + \sup_{\substack{x\in(-1,0)\times B_2 \\ y\notin(-1,0)\times B_3}} \frac{1}{d(x,y)^{\gamma-1}}\right) \\
&= [U]_{\gamma;(-1,0]\times B_2} + C\|\nu\|_{(-1,0]\times B_2},
\end{aligned}$$

where we used that $x, y$ are uniformly bounded and $2 - \gamma > 0$ in the first supremum and for the second that $1 < d(x,y)$ and $\gamma - 1 > 0$. For the second term in the right hand side of (3.28), we will show that

$$\|\mathscr{L}\tilde{U}\|_{\gamma-2s;B_1} \lesssim \|\mathscr{L}U\|_{\gamma-2s;B_1} + \|U\|_{B_1\times((-1,0]\times\mathbb{R}^d)}. \tag{3.32}$$

Using the linearity of $\mathscr{L}$ we have

$$\|\mathscr{L}\tilde{U}\|_{\gamma-2s;B_1} \leqslant \|\mathscr{L}U\|_{\gamma-2s;B_1} + \|\mathscr{L}(\tilde{U}-U)\|_{\gamma-2s;B_1}.$$

For the second term we fix $x \in B_1$ (for which $\eta(x) = 1$), and define

$$\rho(x,\cdot) := \tilde{U}(x,\cdot) - U(x,\cdot) = \eta(\cdot)\,U(x,\cdot) - U(x,\cdot).$$

We have that $\rho(x,y) = 0$ for all $y \in B_1$, and therefore the locality of the $\partial_t$ implies that $\partial_t \rho(x,\cdot) = 0$ in $B_1$. Fix $\lambda > 0$ such that $B_\lambda(x) \subset B_1$, then

$$\begin{aligned}
\langle \mathscr{L}\rho(x,\cdot), \psi_x^\lambda \rangle &= \langle (-\Delta)^s \rho(x,\cdot), \psi_x^\lambda \rangle \\
&= \langle \rho(x,\cdot), (-\Delta)^s \psi_x^\lambda \rangle \\
&= \int \rho(x,y)\,((-\Delta)^s \psi_x^\lambda)(y)\,\mathrm{d}y \\
&= c_{d,s}\int \rho(x,y) \int_{\mathbb{R}^d} (\psi_x^\lambda(y+w) + \psi_x^\lambda(y-w) - 2\,\psi_x^\lambda(y))\,\frac{\mathrm{d}w}{|w|^{d+2s}}\,\mathrm{d}y \\
&= c_{d,s}\int_{\mathbb{R}^d}\int \rho(x,y)(\psi_x^\lambda(y+w) + \psi_x^\lambda(y-w) - 2\,\psi_x^\lambda(y))\,\mathrm{d}y\,\frac{\mathrm{d}w}{|w|^{d+2s}} \\
&= c_{d,s}\int_{\mathbb{R}^d}\int_{B_\lambda(x)} (\rho(x,y-w) + \rho(x,y+w) - 2\,\rho(x,y))\,\psi_x^\lambda(y)\,\mathrm{d}y\,\frac{\mathrm{d}w}{|w|^{d+2s}} \\
&= c_{d,s}\int_{\mathbb{R}^d\setminus B_{1/2}}\int_{B_\lambda(x)} (\rho(x,y-w) + \rho(x,y+w))\,\psi_x^\lambda(y)\,\mathrm{d}y\,\frac{\mathrm{d}w}{|w|^{d+2s}},
\end{aligned}$$

where we used that for all $y = (y_0, y_{1:d}) \in B_\lambda(x) \subset B_1 \subset \mathbb{R}^{1+d}$ and $w \in B_{1/2} \subset \mathbb{R}^d$ we have that $(y_0, y_{1:d} + w) \in B_{3/2}$, and therefore $\rho(x, (y_0, y_{1:d} \pm w)) = 0$ since $\eta \equiv 1$ in $B_{3/2}$. We conclude that

$$|\langle \mathscr{L}\rho(x,\cdot), \psi_x^\lambda \rangle| \lesssim \int_{\mathbb{R}^d \setminus B_{1/2}} \int_{B_\lambda(x)} |\rho(x, y-w) + \rho(x, y+w)| \psi_x^\lambda(y) \, dy \, \frac{dw}{|w|^{d+2s}}$$

$$\lesssim \|\rho(x,\cdot)\|_{L^\infty(\mathbb{R}^d)} \int_{\mathbb{R}^d \setminus B_{1/2}} \frac{dw}{|w|^{d+2s}} \int_{B_\lambda(x)} \psi_x^\lambda(y) \, dy$$

$$\lesssim \|U\|_{B_1 \times ((-1,0] \times \mathbb{R}^d)},$$

since $|\cdot|^{-d-2s}$ is integrable outside a neighbourhood of the origin. Moreover, $B_\lambda(x) \subset B_1$ implies that $\lambda < 1$ and since $-(\gamma - 2s) > 0$

$$\lambda^{-(\gamma - 2s)} |\langle \mathscr{L}\rho(x,\cdot), \psi_x^\lambda \rangle| \lesssim \|U\|_{B_1 \times ((-1,0] \times \mathbb{R}^d)},$$

which proves (3.32). Combining (3.30), (3.31) and (3.32) we conclude in (3.28) the bound

$$\begin{aligned}[U]_{\gamma; B_{1/2}} &\leqslant \delta [U]_{\gamma; B_2} + C_\delta (\|\mathscr{L}U\|_{\gamma-2s; B_1} + [U]_{\gamma\text{-3pt}; B_1} + \|U\|_{B_1 \times ((-1,0] \times \mathbb{R}^d)}) \\ &\leqslant \delta [U]_{\gamma; B_2} + C_\delta (\|\mathscr{L}U\|_{\gamma-2s; B_2} + [U]_{\gamma\text{-3pt}; B_2} + \|U\|_{B_2 \times ((-2^{2s},0] \times \mathbb{R}^d)}). \end{aligned} \quad (3.33)$$

To conclude we use that by Lemma 3.5 we have that $B \mapsto [U]_{\gamma; B}$ satisfies the condition (3.12) from Lemma 3.6 for

$$D(B) := [U]_{\gamma\text{-3pt}; B} + [\Lambda]_{(\gamma-1)\text{-3pt}; B},$$

and therefore by Lemma 3.6 with $B = B_2$ $S(A) = [U]_{\gamma; A}$ and $\theta_0 = 1/4$ we obtain the existence of $\varepsilon = \varepsilon(\gamma, s, d) > 0$ such that if for some $E \geqslant 0$, the bound

$$\sigma^\gamma [U]_{\gamma; B_{\sigma/4}(y)} \leqslant \varepsilon \sigma^\gamma [U]_{\gamma; B_\sigma(y)} + E \quad (3.34)$$

is satisfied for all $B_\sigma(y) \subseteq B_2$, then we can conclude that

$$\sigma^\gamma [U]_{\gamma; B_{\theta\sigma}(y)} \leqslant C(E + D(B_2)) = C(E + [U]_{\gamma\text{-3pt}; B_2} + [\Lambda]_{(\gamma-1)\text{-3pt}; B_2})$$

for all balls $B_\sigma(y) \subseteq B_2$ and $\theta \in (0,1)$ for some $C = C(d, \theta, \gamma, s) > 0$. To show (3.34) for given $y \in \mathbb{R}^{1+d}$ and $\sigma > 0$, we apply (3.33) to the germ $U(y + \sigma \cdot, y + \sigma \cdot)$, and using the scaling properties of the seminorms (see Lemma 3.15) and the $2s$ homogeneity of $\mathscr{L}$, we conclude that

$$\sigma^\gamma [U]_{\gamma; B_{\sigma/2}(y)} \leqslant \delta \sigma^\gamma [U]_{\gamma; B_{2\sigma}(y)} + C_\delta (\|\mathscr{L}U\|_{\gamma-2s; B_2} + [U]_{\gamma\text{-3pt}; B_2} + \|U\|_{B_2 \times ((-2^{2s},0] \times \mathbb{R}^d)}),$$

were we used that, since $B_{2\sigma}(y) \subset B_2$, and $\gamma > 0$ then $\sigma^\gamma \leqslant 1$. Choosing $\delta = \varepsilon > 0$ we get a concrete $C_\varepsilon > 0$ such that (3.34) holds for all $B_\sigma(y) \subseteq B_2$ and we can conclude that

$$\begin{aligned}&\sigma^\gamma [U]_{\gamma; B_{\theta\sigma}(y)} \\ &\leqslant C_\varepsilon (\|\mathscr{L}U\|_{\gamma-2s; B_2} + [U]_{\gamma\text{-3pt}; B_2} + \|U\|_{B_2 \times ((-2^{2s},0] \times \mathbb{R}^d)}) + ([U]_{\gamma\text{-3pt}; B_2} + [\Lambda]_{(\gamma-1)\text{-3pt}; B_2}) \\ &\lesssim \|\mathscr{L}U\|_{\gamma-2s; B_2} + [U]_{\gamma\text{-3pt}; B_2} + [\Lambda]_{(\gamma-1)\text{-3pt}; B_2} + \|U\|_{B_2 \times ((-2^{2s},0] \times \mathbb{R}^d)} \end{aligned}$$

for all $B_\sigma(y) \subseteq B_2$ and $\theta \in (0,1)$ for some constant $C_\varepsilon = C(s, d, \gamma, \theta)$, and in particular taking $\sigma = 2$, $y = 0$ and $\theta = 1/2$ we obtain the claimed result. □

**Corollary 3.14**. *Let $U$ be a germ which satisfies the conditions from Theorem 3.1 and $z \in \mathbb{R}^{1+d}$, $\lambda > 0$ such that $B_\lambda(z) \subset B_1$, then*

$$\lambda^\gamma [U]_{\gamma; B_{\lambda/2}(z)} \lesssim \lambda^\gamma \|\mathscr{L}(U)\|_{\gamma-2s; B_\lambda(z)} + \lambda^\gamma [U]_{\gamma\text{-3pt}; B_\lambda(z)} + \lambda^\gamma [\Lambda]_{(\gamma-1)\text{-3pt}; B_\lambda(z)}$$
$$+ \|U\|_{B_\lambda(z) \times ((z_0 - \lambda^{2s}, z_0] \times \mathbb{R}^d)}.$$

**Proof.** It follows from applying Theorem 3.1 to the germ $(x, y) \mapsto U(z + \sigma x, z + \sigma y)$ and using the scaling properties of the seminorms as in Lemma 3.15. □

We summarise the scaling properties of some of the seminorms in the following result which proof we omit since it follows directly from the definitions.

**Lemma 3.15**. *Fix $z \in \mathbb{R}^{1+d}$ and $\sigma > 0$. Let $U: \mathbb{R}^{1+d} \times \mathbb{R}^{1+d} \to \mathbb{R}$ be a germ with $[U]_{\gamma; B_\sigma(z)} < +\infty$ for some $\gamma \in (1, 2s)$ and let $\nu$ be defined by (3.2). Assume $\Lambda: \mathbb{R}^{1+d} \times \mathbb{R}^{1+d} \to \mathbb{R}$ is such that the 3-point continuity conditions (3.6) and (3.7) hold in $B_\sigma(z)$. Define $U_{z,\sigma}(x,y) := U(z + \sigma x, z + \sigma y)$, $\nu_{z,\sigma}(x) := \nu(z + \sigma x)$ and $\Lambda_{z,\sigma}(x,y) := \Lambda(z + \sigma x, z + \sigma y)$. Then we have that $\sigma \nu_{z,\sigma}$ is the diagonal derivative of $U_{z,\sigma}$, and $(U_{z,\sigma}, \sigma \Lambda_{z,\sigma})$ satisfy (3.6) and (3.7) in $B_1$. Moreover, we have the following scaling properties of the seminorms*

$$[U_{z,\sigma}]_{\gamma; B_1} = \sigma^\gamma [U]_{\gamma; B_\sigma(z)}, \qquad [U_{z,\sigma}]_{\gamma\text{-3pt}; B_1} = \sigma^\gamma [U]_{\gamma\text{-3pt}; B_\sigma(z)}.$$

*In particular we also have that*

$$[\Lambda_{z,\sigma}]_{\gamma\text{-3pt}; B_1} = \sigma^\gamma [\Lambda]_{(\gamma-1)\text{-3pt}; B_\sigma(z)}, \qquad [\nu_{z,\sigma} - \Lambda_{z,\sigma}]_{\gamma-1; B_1} = [\nu - \Lambda]_{\gamma-1; B_\sigma(z)},$$

*where $(\nu - \Lambda)(x, y) := \nu(y) - \nu(x) - \Lambda(x, y)$. At last, if $F$ is a distributional germ such that $\|F\|_{\gamma-2s; B_1} < +\infty$ and $F_{z,\sigma}$ is defined analogously to $U_{z,\sigma}$ (via duality) then we have that*

$$\|F_{z,\sigma}\|_{\gamma-2s; B_1} = \sigma^{\gamma-2s} \|F\|_{\gamma-2s; B_\sigma(z)}.$$

## 4. A Priori Estimates

For this section we assume we are given $(\Pi; \Gamma)$ a 1-periodic weakly admissible model (see Definition 2.11) and $V \in \mathcal{D}^\gamma$ the 1-periodic modelled distribution defined by (2.13) for some $\gamma \in (3 - 2s, 2s)$, where the coefficients $\{\Upsilon.[\tau]\}_{\tau \in \mathcal{T}}$ are given by Lemma 2.14.

### 4.1. Small scales estimates

The main result of this section is a small scales estimate for some Hölder type seminorms of the remainder $V$. Following [MW20, CMW23] Assumption 4.1 is a convenient way to remove the dependence of some seminorms from the model $(\Pi; \Gamma)$ from the main estimate in Theorem 4.2.

**Assumption 4.1**. *For $c \in (0, 1)$ and $t \in (0, 1)$ we assume that*

$$[\Pi; \tau]_K \leqslant c \, \|v\|_{P_t}^{\alpha(\tau)} \qquad \forall \tau \in (\mathcal{T}_{<2s} \setminus \mathcal{P}) \setminus \mathcal{W}, \tag{4.1}$$

$$[\Gamma \mathcal{I}(\tau); \mathbf{1}]_K, [\Gamma \mathcal{I}(\tau); \boldsymbol{X}]_K, \sup_{x \in P} \|\gamma_{.,.}(\mathcal{I}(\tau))\|_{[0,1] \times \mathbb{R}^d} \leqslant c \, \|v\|_{P_t}^{\alpha(\tau)} \qquad \forall \tau \in \mathcal{V}_{0,\gamma}, \tag{4.2}$$

*with $K = [0,1] \times \overline{B_4} \subset \mathbb{R} \times \mathbb{R}^d$, $P_t := (t^{2s}, 1] \times (0, 1]^d$, $P := P_0$, and for $\tau \in \mathcal{T}$*

$$\alpha(\tau) := \mathfrak{l}(\tau) \, s^{-1} (3s + |\Xi|) = \mathfrak{l}(\tau) \, s^{-1} \left(2s - \frac{3}{2} - \kappa\right) > 0, \tag{4.3}$$

with $\mathfrak{l}(\tau)$ the number of leaves in the tree $\tau$.

**Theorem 4.2.** *Let $(\Pi; \Gamma)$ be a 1-periodic weakly admissible model (see Definition 2.11) and $V \in \mathcal{D}^\gamma$ the 1-periodic modelled distribution defined by (2.13) for some $\gamma \in (3-2s, 2s)$, where the coefficients $\{\Upsilon.[\tau]\}_{\tau \in \mathcal{T}}$ are given by Lemma 2.14. In particular, let $v = \langle \mathbf{1}, V \rangle$. Fix $c \in (0,1)$ and $t \in (0,1)$ such that Assumption 4.1 holds and set*

$$\lambda_t := \|v\|_{P_t}^{-s^{-1}}. \tag{4.4}$$

*If $c$ satisfies a smallness condition which depends only on $s, d$ and $\gamma$, and $t + 3\lambda_t < 1$ then*

$$\sup_{\lambda \leqslant \lambda_t} \lambda^\beta \sup_{z \in P_{t+3\lambda}} [V_\beta; \mathbf{1}]_{B_{\lambda/2}(z)} \lesssim \|v\|_{P_t} \qquad \forall \beta \in (0, \gamma],$$

*where $V_\beta$ is the truncation of $V$ as defined in (2.20) and with the implicit proportionality constant depending only on $s, d$ and $\gamma$.*

**Remark 4.3.** We only use 1-periodicity of $V$ on the identities $\|v\|_{(t^{2s},1] \times \overline{B_4}} = \|v\|_{(t^{2s},1] \times [0,1]^d}$, but the result holds true without periodicity (in both $V$ and $(\Pi; \Gamma)$) by replacing in Assumption 4.1 and Theorem 4.2 the $L^\infty$-norms in the set $P_t$ by an enlargement in space of this set.

We start by defining for each $\beta \in (0, \gamma]$ the germ

$$U_\beta(x, y) := v(y) - v(x) - \sum_{\tau \in \mathcal{V}_{0,\beta}} \frac{\Upsilon_x[\tau]}{\tau!} (\Pi_x \mathcal{I}(\tau))(y), \tag{4.5}$$

and for simplicity we denote $U := U_\gamma$. By Lemma 2.24 we have the identity

$$U(x,y) - v_{\mathbf{X}}(x) \cdot (y-x) = v(y) - (\Pi_x V(x))(y) = \langle \mathbf{1}, V(y) - \Gamma_{yx} V(x) \rangle. \tag{4.6}$$

Moreover, since $(\Pi_y \mathcal{I}(\tau))(y) = 0$ for all $\tau \in \mathcal{V}$ ($|\mathcal{I}(\tau)| > 0$) then $U(y,y) = 0$. In particular we have $[U]_{\gamma;B} = [V; \mathbf{1}]_{\gamma;B} < +\infty$ for any bounded set $B \subset \mathbb{R}^{1+d}$ since $V \in \mathcal{D}^\gamma$ by assumption (see (3.1) and (A.3) for the definitions). Since $\gamma > 1$, this implies that

$$v_{\mathbf{X}}(x) = \nabla U(x, \cdot)|_x. \tag{4.7}$$

The idea is to use the Schauder estimate developed in Section 3, however if used directly with the germ $U$ we obtain a non-local in space $L^\infty$-term which we cannot control with our techniques. To avoid this for $z \in \mathbb{R}^{1+d}$ and $\lambda > 0$ we consider the following localisation of the germ $U$ defined as:

$$U_{z,\lambda}^{\mathrm{loc}}(x,y) := v(y) - v(x) - \sum_{\tau \in \mathcal{V}_{0,\beta}} \frac{\Upsilon_x[\tau]}{\tau!} (\Pi_x \mathcal{I}(\tau))(y) \, \eta(\lambda^{-1}(y-z)_{1:d}), \tag{4.8}$$

where $\eta \in C_c^\infty(\mathbb{R}^d)$ is a non-negative function bounded by 1 which satisfies $\eta \equiv 1$ in $B_2 \subset \mathbb{R}^d$ and $\eta \equiv 0$ in $\mathbb{R}^d \setminus B_3$. Applying the Schauder estimate in the form of Corollary 3.14 to this germ we conclude the next lemma. Observe that now the global $L^\infty$ term turns into a local one, at the price of adding the remainder term $\mathscr{L}(U_{z,\lambda}^{\mathrm{loc}} - U)$.

**Lemma 4.4.** *For all $\lambda > 0$ and $z \in \mathbb{R}^{1+d}$ we have*

$$\begin{aligned}
\lambda^\gamma [U]_{\gamma;B_{\lambda/2}(z)} &\lesssim \lambda^\gamma \|\mathscr{L} U\|_{\gamma-2s;B_\lambda(z)} + \lambda^\gamma [U]_{\gamma\text{-3pt};B_\lambda(z)} + \lambda^\gamma [\Lambda]_{(\gamma-1)\text{-3pt};B_\lambda(z)} \\
&\quad + \lambda^\gamma \|\mathscr{L}(U_{z,\lambda} - U)\|_{\gamma-2s;B_\lambda(z)} + \|U\|_{B_{3\lambda}(z)} + \|v\|_{B_{3\lambda}(z)}.
\end{aligned} \tag{4.9}$$

The goal is now to control all the terms in the right hand side of (4.9). We postpone the proofs of all the stated lemmas in this section to Section 4.3.

**Bounding the $\mathscr{L}U$ term**

Applying the model $\mathscr{L}\Pi_x$ to equation (2.14) and combining it with the PDE (2.15) we write the action of the operator $\mathscr{L}$ on the germ $U$.

**Lemma 4.5**. *Consider the germ $U := U_\gamma$ defined by (4.5). Let $\varepsilon > 0$ satisfy (2.19) and (2.18). For $\tau, \tau_1, \tau_2 \in \mathcal{W}$ define $\beta_\tau := \varepsilon - |\mathcal{I}(\tau)|$ and $\beta_{\tau_1, \tau_2} = \varepsilon - |\mathcal{I}(\tau_1)\mathcal{I}(\tau_2)|$. Then for all $x \in \mathbb{R}^{1+d}$*

$$\begin{aligned}\mathscr{L}(U(x, \cdot)) &= -v^3 - 3 \sum_{\tau \in \mathcal{W}} \frac{\Upsilon[\tau]}{\tau!} (\mathcal{R}(V^2(x)\mathcal{I}(\tau)) - \Pi_x(V^2_{\beta_\tau}(x)\mathcal{I}(\tau)))\,\rho \\ &\quad - 3 \sum_{\tau_1, \tau_2 \in \mathcal{W}} \frac{\Upsilon[\tau_1]\Upsilon[\tau_2]}{\tau_1!\,\tau_2!} (\mathcal{R}(V\mathcal{I}(\tau_1)\mathcal{I}(\tau_2)) - \Pi_x(V_{\beta_{\tau_1}, \beta_{\tau_2}}(x)\mathcal{I}(\tau_1)\mathcal{I}(\tau_2)))\,\rho \\ &\quad + \sum_{\tau \in \delta\mathcal{W} \setminus \mathcal{V}_{0,\gamma}} \frac{\Upsilon[\tau]}{\tau!} \mathbf{\Pi}\tau. \end{aligned} \qquad (4.10)$$

We use the Reconstruction theorem to bound some of the terms in (4.10).

**Lemma 4.6**. *Let $\varepsilon > 0$ satisfy (2.19) and (2.18). Given $\tau \in \mathcal{W}$ let $\beta_\tau := \varepsilon - |\mathcal{I}(\tau)|$, then for all $x \in \mathbb{R}^{1+d}$ and $r > 0$ we have*

$$\begin{aligned}&|\langle \mathcal{R}(\mathcal{I}(\tau)V^2) - \Pi_x(\mathcal{I}(\tau)V^2_{\beta_\tau}(x))\,\rho, \psi^r_x\rangle| \\ &\lesssim r^{\beta_\tau + |\mathcal{I}(\tau)|} \sum_{\sigma_1, \sigma_2 \in \mathcal{T}} [\Pi; \sigma_1\sigma_2\mathcal{I}(\tau)]_{B_{2r}(x)} [V^2_{\beta_\tau}; \sigma_1\sigma_2]_{B_{2r}(x)}. \end{aligned} \qquad (4.11)$$

*Given $\tau_1, \tau_2 \in \mathcal{W}$ let $\beta_{\tau_1, \tau_2} = \varepsilon - |\mathcal{I}(\tau_1)\mathcal{I}(\tau_2)|$, then for all $x \in \mathbb{R}^{1+d}$ and $r > 0$ we have*

$$\begin{aligned}&\langle \mathcal{R}(V\mathcal{I}(\tau_1)\mathcal{I}(\tau_2)) - \Pi_x(V_{\beta_{\tau_1,\tau_2}}(x)\mathcal{I}(\tau_1)\mathcal{I}(\tau_2))\,\rho, \psi^r_x\rangle \\ &\lesssim r^{\beta_{\tau_1,\tau_2} + |\mathcal{I}(\tau_1)\mathcal{I}(\tau_2)|} \sum_{\sigma \in \mathcal{T}} [\Pi; \sigma\mathcal{I}(\tau_1)\mathcal{I}(\tau_2)]_{B_{2r}(x)} [V_{\beta_{\tau_1,\tau_2}}; \sigma]_{B_{2r}(x)}. \end{aligned} \qquad (4.12)$$

*Moreover, all the implicit proportionality constants depend on $s, d$ and $\gamma$.*

In the previous lemma the only terms that actually contribute to the corresponding right hand side are the ones for which $[V_\beta; \sigma], [V^2_\beta; \sigma_1\sigma_2]$ are not zero. By Lemma 2.19 for $V$ this corresponds to $\mathbf{1}, \mathbf{X}$ (depending on $\beta$) and on planted trees. On the other hand, by Lemma 2.21 for $V^2$ this corresponds to coefficients at $\mathbf{1}$, planted trees, or product of two planted trees. Moreover, again by Lemma 2.21, this coefficients of $V^2$ can be written as coefficients of some truncation of $V$ at a different coefficient. Using the explicit representations of $[V_\beta; \sigma]$ given by Lemma 2.19 in terms of expansions of $v, v^2$ and $v_\mathbf{X}$, and combined with Assumption 4.1 on the model we have the following estimates.

**Lemma 4.7**. *If $c, t \in (0,1)$ satisfy Assumption 4.1 and $t + 3\lambda_t < 1$ with $\lambda_t$ is defined by (4.4), then for all $\lambda \in (0, \lambda_t]$, $z \in P_{t+3\lambda}$ and $\beta \in (0, \gamma]$ we have*

$$\lambda^\beta [V_\beta; \mathbf{1}]_{B_{2\lambda}(z)} \lesssim \|v\|_{P_t} + (\mathbb{1}_{\beta \leqslant 1} + c)\,\lambda\,\|v_\mathbf{X}\|_{B_{2\lambda}(z)} + \lambda^\gamma\,[U]_{\gamma; B_{2\lambda}(z)}, \qquad (4.13)$$

$$\lambda^\beta [V_\beta; \mathbf{X}]_{B_{2\lambda}(z)} \lesssim \lambda^{-1} (\|v\|_{P_t} + c\,\lambda\,\|v_\mathbf{X}\|_{B_{2\lambda}(z)} + \lambda^\gamma\,[V; \mathbf{X}]_{B_{2\lambda}(z)}), \qquad (4.14)$$

and for all $\tau \in \mathcal{V}_{0,\beta}$

$$\lambda^{\beta-|\mathcal{I}(\tau)|}\,[V_\beta;\mathcal{I}(\tau)]_{B_{2\lambda}(z)}$$
$$\lesssim \begin{cases} \|v\|_{P_t} + (\mathbb{1}_{\beta-|\mathcal{I}(\tau)|\leqslant 1} + c)\,\lambda\,\|v_{\boldsymbol{X}}\|_{B_{2\lambda}(z)} + \lambda^\gamma\,[U]_{\gamma;B_{2\lambda}(z)} & \text{if } \Upsilon[\tau] \sim v \\ \|v\|_{P_t}(\|v\|_{P_t} + \lambda\,\|v_{\boldsymbol{X}}\|_{B_{2\lambda}(z)} + \lambda^\gamma\,[U]_{\gamma;B_{2\lambda}(z)}) & \text{if } \Upsilon[\tau] \sim v^2 \\ \lambda^{-1}(\|v\|_{P_t} + c\,\lambda\,\|v_{\boldsymbol{X}}\|_{B_{2\lambda}(z)} + \lambda^\gamma\,[V;\boldsymbol{X}]_{B_{2\lambda}(z)}) & \text{if } \Upsilon[\tau] \sim v_{\boldsymbol{X}} \end{cases}\;. \quad (4.15)$$

*Moreover, all the implicit proportionality constants depend only on $s, d$ and $\gamma$.*

Combining Lemma 4.7 with Assumption 4.1 we control the terms on the right hand side of the reconstructions in Lemma 4.6.

**LEMMA 4.8**. *If $c, t \in (0,1)$ satisfy Assumption 4.1 and $t + 3\lambda_t < 1$ with $\lambda_t$ is defined by (4.4). Then for all $\lambda \in (0, \lambda_t]$, $z \in P_{t+3\lambda}$, $\beta \in (0, \gamma)$, $\tau \in \mathcal{W}$ and $\sigma_1, \sigma_2 \in \mathcal{T}$ we have the bound*

$$\lambda^{\beta_\tau + |\mathcal{I}(\tau)| + 2s}[\Pi; \sigma_1\sigma_2\mathcal{I}(\tau)]_{B_{2\lambda}(z)}[V^2_{\beta_\tau};\sigma_1\sigma_2]_{B_{2\lambda}(z)} \lesssim \|v\|_{P_t} + c\,\lambda\,\|v_{\boldsymbol{X}}\|_{B_{2\lambda}(z)} + c\,\lambda^\gamma\,[U]_{\gamma;B_{2\lambda}(z)}.$$

*Similarly, for all $\tau_1, \tau_2 \in \mathcal{W}$ and $\sigma \in \mathcal{T}$ such that $|\mathfrak{n}(\sigma)| = 0$ we have*

$$\lambda^{\beta_{\tau_1,\tau_2} + |\mathcal{I}(\tau_1)\mathcal{I}(\tau_2)| + 2s}[\Pi;\sigma\mathcal{I}(\tau_1)\mathcal{I}(\tau_2)]_{B_{2\lambda}(z)}[V_{\beta_{\tau_1,\tau_2}};\sigma]_{B_{2\lambda}(z)} \lesssim \|v\|_{P_t} + c\,\lambda\,\|v_{\boldsymbol{X}}\|_{B_{2\lambda}(z)} + c\,\lambda^\gamma\,[U]_{\gamma;B_{2\lambda}(z)}.$$

*and for all $\tau_1, \tau_2 \in \mathcal{W}$ and $\sigma \in \mathcal{T}$ such that $|\mathfrak{n}(\sigma)| = 1$ we have*

$$\lambda^{\beta_{\tau_1,\tau_2} + |\mathcal{I}(\tau_1)\mathcal{I}(\tau_2)| + 2s}[\Pi;\sigma\mathcal{I}(\tau_1)\mathcal{I}(\tau_2)]_{B_{2\lambda}(z)}[V_{\beta_{\tau_1,\tau_2}};\sigma]_{B_{2\lambda}(z)} \lesssim \|v\|_{P_t} + c\,\lambda\,\|v_{\boldsymbol{X}}\|_{B_{2\lambda}(z)} + c\,\lambda^\gamma\,[V;\boldsymbol{X}]_{B_{2\lambda}(z)},$$

*where $[V_\beta;\boldsymbol{X}] := \max_{j \in \{1,\dots,d\}}[V_\beta;\boldsymbol{X}_j]$. Moreover, all of the previous implicit proportionality constants depend only on $s, d$ and $\gamma$.*

After bounding the missing cubic term and the noises $\boldsymbol{\Pi}\tau$ in (4.10) we can conclude the following.

**LEMMA 4.9**. *If $c, t \in (0,1)$ satisfy Assumption 4.1 and $t + 3\lambda_t < 1$ with $\lambda_t$ is defined by (4.4), then for all $\lambda \in (0, \lambda_t]$, $z \in P_{t+3\lambda}$ we have*

$$\lambda^\gamma\,\|\mathscr{L}U\|_{\gamma-2s;B_\lambda(z)} \lesssim \|v\|_{P_t} + c\,\lambda\,\|v_{\boldsymbol{X}}\|_{B_{2\lambda}(z)} + c\,\lambda^\gamma\,[U]_{\gamma;B_{2\lambda}(z)} + c\,\lambda^\gamma\,[V;\boldsymbol{X}]_{B_{2\lambda}(z)},$$

*with the implicit proportionality constants depending only on $s, d$ and $\gamma$.*

**Bounding the 3 point continuity term**

Using the underlying model distribution in the definition of the germ $U$ the following result writes the 3-point continuity condition in terms of the model $(\Pi; \Gamma)$ and the modelled distribution $V$.

**LEMMA 4.10**. *Consider the germ $U = U_\gamma$ defined by (4.5) and the vector-valued germ $\Lambda = \{\Lambda_j\}_{j=1}^d$ defined as*

$$\Lambda_j(x, y) := \sum_{\sigma \in \mathcal{V}_{1,\gamma}} \frac{\Upsilon_x[\sigma]}{\sigma!}\,\gamma_{yx}(\mathcal{I}_j(\sigma)). \quad (4.16)$$

*Then $U$ and $\Lambda$ satisfy the 3-point continuity conditions (3.6)-(3.8) with $A := \{|\mathcal{I}(\tau)|: \tau \in \mathcal{V}_{0,\gamma}\}$. Moreover, for any set $B \subset \mathbb{R}^{1+d}$ we have:*

$$[U]_{\gamma\text{-3pt};B} \lesssim \sum_{\tau \in \mathcal{V}_{0,\gamma}} [V;\mathcal{I}(\tau)]_B\,[\Gamma\mathcal{I}(\tau);\mathbf{1}]_B, \quad (4.17)$$

and

$$[\Lambda]_{(\gamma-1)\text{-3pt};B} \lesssim \sum_{\tau \in \mathcal{V}_{1,\gamma}} [V; \mathcal{I}(\tau)]_B \, [\Gamma \mathcal{I}(\tau); \boldsymbol{X}]_B, \qquad (4.18)$$

with the implicit proportionality constants depending only on $s, d$ and $\gamma$.

Combining the previous lemma with Assumption 4.1 we conclude the following bound.

**LEMMA 4.11**. *If $c, t \in (0,1)$ satisfy Assumption 4.1 and $t + 3\lambda_t < 1$ with $\lambda_t$ defined by (4.4), then for all $\lambda \in (0, \lambda_t]$, $z \in P_{t+3\lambda}$, $r \leqslant 3\lambda$ we have*

$$\lambda^\gamma [U]_{\gamma\text{-3pt};B_r(z)} \lesssim \|v\|_{P_t} + c\,\lambda\,\|v_{\boldsymbol{X}}\|_{B_r(z)} + c\,\lambda^\gamma [U]_{\gamma;B_r(z)} + c\,\lambda^\gamma [V; \boldsymbol{X}]_{B_r(z)},$$

and

$$\lambda^\gamma [\Lambda]_{(\gamma-1)\text{-3pt};B_r(z)} \lesssim \|v\|_{P_t} + c\,\lambda\,\|v_{\boldsymbol{X}}\|_{B_r(z)} + c\,\lambda^\gamma [U]_{\gamma;B_r(z)} + c\,\lambda^\gamma [V; \boldsymbol{X}]_{B_r(z)},$$

with the implicit proportionality constants depending only on $s, d$ and $\gamma$.

**Bounding the rest of the terms**

The next lemma bounds the remaining terms in (4.9).

**LEMMA 4.12**. *If $c, t \in (0,1)$ satisfy Assumption 4.1 and $t + 3\lambda_t < 1$ with $\lambda_t$ defined by (4.4), then for all $\lambda \in (0, \lambda_t]$, $z \in P_{t+3\lambda}$ and $r \leqslant 3\lambda$ we have*

$$\|U_\beta\|_{B_r(z)} \lesssim \|v\|_{P_t} + \mathbb{1}_{\beta>1} c\,\lambda\,\|v_{\boldsymbol{X}}\|_{B_r(z)} \qquad \forall \beta \in (0, \gamma] \qquad (4.19)$$

and

$$\lambda^\gamma \|\mathscr{L}(U_{z,\lambda} - U)\|_{\gamma-2s;B_\lambda(z)} \lesssim \|v\|_{P_t} + c\,\lambda\,\|v_{\boldsymbol{X}}\|_{B_\lambda(z)}, \qquad (4.20)$$

with the implicit proportionality constants depending only on $s, d$ and $\gamma$.

Using Lemma 4.9, Lemma 4.11 and Lemma 4.12 we can bound all the terms in (4.9) in terms of $\|v\|$, $\|v_{\boldsymbol{X}}\|$, $[U]_\gamma$ and $[V; \boldsymbol{X}]$. Under an smallness assumption on $c \in (0,1)$, depending only on $s, d$ and $\gamma$, we can bound the gradient terms $\|v_{\boldsymbol{X}}\|, [V; \boldsymbol{X}]$ in terms of the other two.

**LEMMA 4.13**. *If $c, t \in (0,1)$ satisfy Assumption 4.1, $t + 3\lambda_t < 1$ with $\lambda_t$ defined by (4.4) and $c$ satisfies a smallness condition depending only on $d, s$ and $\gamma$, then for all $\lambda \in (0, \lambda_t]$, $z \in P_{t+3\lambda}$, $r \leqslant 3\lambda$ we have*

$$\lambda \|v_{\boldsymbol{X}}\|_{B_r(z)} \lesssim \|v\|_{P_t} + \lambda^\gamma [U]_{\gamma;B_r(z)}, \qquad (4.21)$$

and

$$\lambda^\gamma [V; \boldsymbol{X}]_{B_r(z)} \lesssim \|v\|_{P_t} + \lambda^\gamma [U]_{\gamma;B_r(z)}, \qquad (4.22)$$

with the implicit proportionality constants depending only on $s, d$ and $\gamma$.

With the previous lemmas we can now prove Theorem 4.2.

**Proof of Theorem 4.2.** From Lemma 4.9, Lemma 4.11, Lemma 4.12 and Lemma 4.13 we can conclude on Lemma 4.4:

$$\lambda^\gamma [U]_{\gamma;B_{\lambda/2}(z)} \lesssim \|v\|_{P_t} + c\,\lambda\,\|v_{\boldsymbol{X}}\|_{B_{3\lambda}(z)} + c\,\lambda^\gamma [U]_{\gamma;B_{2\lambda}(z)} + c\,\lambda^\gamma [V;\boldsymbol{X}]_{B_{2\lambda}(z)}.$$

Using the smallness assumption on $c$ we conclude by Lemma 4.13 that

$$\lambda^\gamma [U]_{\gamma;B_{\lambda/2}(z)} \lesssim \|v\|_{P_t} + c\,\lambda^\gamma [U]_{\gamma;B_{3\lambda}(z)}. \tag{4.23}$$

We want to absorb the terms $[U]_{\gamma;B_{3\lambda}(z)}$ on the right hand side into the left hand side using the abstract absorption Lemma 3.6. However, by using this lemma as stated one runs into the issue that, by Lemma 4.11, the three point continuity seminorms appearing in the final estimate would contain again the term $c\,\lambda^\gamma [U]_{\gamma;B_{3\lambda}(z)}$ on the right hand side which turns into a circular argument. However, since this term comes again multiplied with a smallness constant we can prove a slight modification of Lemma 3.6 that avoids this circular argument which we state as Lemma 4.25.

Let $S(B):=[U]_{\gamma;B}$, $D(B):=\lambda^\gamma [U]_{\gamma\text{-3pt};B} + \lambda^\gamma [v_{\boldsymbol{X}_j}-\Lambda]_{(\gamma-1)\text{-3pt};B}$, $B:=P_t$ and $\theta_0=1/6$. We recall that by Lemma 3.5 $S$ satisfies the almost subadditivity condition (4.65) for this choice of $D$, and therefore by Lemma 4.25 there exists $\varepsilon>0$ such that if for all $\lambda\in(0,\lambda_t]$ and $z\in P_{t+3\lambda}$ we have the bound

$$\max\left\{\lambda^\gamma [U]_{\gamma;B_{\lambda/2}(z)}, \lambda^\gamma [U]_{\gamma\text{-3pt};B_{3\lambda}(z)} + \lambda^\gamma [v_{\boldsymbol{X}_j}-\Lambda]_{(\gamma-1)\text{-3pt};B_{3\lambda}(z)}\right\} \leqslant \varepsilon\,\lambda^\gamma [U]_{\gamma;B_{3\lambda}(z)} + E, \tag{4.24}$$

then one can conclude for all $\lambda\in(0,\lambda_t]$ and $z\in P_{t+3\lambda}$ the bound

$$\lambda^\gamma [U]_{\gamma;B_{\lambda/2}(z)} \leqslant C\,E, \tag{4.25}$$

for some constant $0<C=C(d,\gamma,s)$. From (4.23), which holds uniformly on $t\in(0,1)$ satisfying the conditions of this theorem, $\lambda\in(0,\lambda_t]$ and $z\in P_{t+3\lambda}$ as long as $c\in(0,1)$ is small enough so that Lemma 4.13 hold, we can conclude by imposing a smallness condition on $c\in(0,1)$, which depends on the implicit proportionality constant and $\varepsilon$, we obtain the first part of (4.24) for $E=\|v\|_{P_t}$. For the second part of (4.17) we use Lemma 4.13 to improve the conclusions of Lemma 4.11 and obtain

$$\lambda^\gamma [U]_{\gamma\text{-3pt};B_{3\lambda}(z)} + \lambda^\gamma [v_{\boldsymbol{X}_j}-\Lambda]_{(\gamma-1)\text{-3pt};B_{3\lambda}(z)} \lesssim c\,\lambda^\gamma [U]_{\gamma;B_{3\lambda}(z)} + \|v\|_{P_t},$$

and by imposing another smallness condition on $c\in(0,1)$, which depends on the implicit proportionality constant and $\varepsilon$, and therefore on $s,d$ and $\gamma$, one can conclude that (4.24) holds for the second term, which allows us to conclude (4.25). $\square$

## 4.2. Large scales estimates

The first step in the proof of Theorem 1.1 is to show that for fixed $c\in(0,1)$ we can consider that Assumption 4.1 holds, since we can obtain (see Lemma 4.15) an uniform estimate for the times where it does not hold.

**Definition 4.14.** *For $c\in(0,1)$ we define $T_c:=\min\{T_{1,c},T_{2,c},T_{3,c}\}\in[0,1]$ where*

$$T_{1,c} := \inf\left\{t\in[0,1]\,\big|\,\exists\tau\in(\mathcal{T}_{<2s}\setminus\mathcal{P})\setminus\mathcal{W}\ \text{such that}\ c\,\|v\|_{P_t}^{\alpha(\tau)} < [\Pi;\tau]_K\right\},$$

$$T_{2,c} := \inf\left\{t\in[0,1]\,\big|\,\exists\tau\in\mathcal{V}_{0,2s},k\in\{0,e_1,\ldots,e_d\}\ \text{such that}\ c\,\|v\|_{P_t}^{\alpha(\tau)} < [\Gamma\,\mathcal{I}(\tau);\boldsymbol{X}^k]_K\right\},$$

$$T_{3,c} := \inf\left\{t\in[0,1]\,\bigg|\,\exists\tau\in\mathcal{V}_{0,2s}\ \text{such that}\ c\,\|v\|_{P_t}^{\alpha(\tau)} < \sup_{x\in P}\|\Pi_x\mathcal{I}(\tau)\|_{(0,1]\times\mathbb{R}^d}\right\}.$$

**Lemma 4.15**. *Fix $c \in (0,1)$ and let $T_c$ as in Definition 4.14, then for all $t \in [T_c, 1]$ we have*

$$\|v\|_{P_t} \lesssim \max\left\{ \max_{\tau \in (\mathcal{T}_{<2s} \setminus \mathcal{P}) \setminus \mathcal{W}} [\Pi;\tau]_K^{\frac{1}{\alpha(\tau)}}, \max_{\substack{\tau \in \mathcal{V}_{0,2s} \\ |\mathfrak{e}(\tau)|=0}} \max_{\substack{k \in \mathbb{N}^{1+d} \\ |k|<2s}} \left( [\Gamma \mathcal{I}(\tau), \boldsymbol{X}^k]_K^{\frac{1}{\alpha(\tau)}} \vee \|\gamma_{\cdot,\cdot}(\mathcal{I}(\tau))\|_{([0,1] \times \mathbb{R}^d) \times P}^{\frac{1}{\alpha(\tau)}} \right) \right\},$$

*for an implicit proportionality constant which depends only on $c, d, s$ and $\gamma$.*

**Remark 4.16**. Lemma 4.15 tells us that Theorem 1.1 holds true for $t \in [T_c, 1]$, and therefore it is only left to prove the case $t \in (0, T_c]$. Since $\|v\|_{P_t} \leqslant \|v\|_{P_{t,T_c}} + \|v\|_{P_{T_c}}$ with $P_{t,T_c} := (t^{2s}, T^{2s}] \times (0, 1]^d = P_t \setminus P_{T_c}$, and the bound in Lemma 4.15 depends only on $s, d$ and $\gamma$, it will be enough to bound $\|v\|_{P_{t,T_c}}$. Since we will choose $c \in (0,1)$ depending only on $s, d$ and $\gamma$, we will assume that $T_c = 1$ to simplify notation, in which case we have $\|v\|_{P_{t,T_c}} = \|v\|_{P_t}$.

Consider a scale $\lambda$ (to be chosen later and which will depend on $v$) and a regularisation of $v$ given by $(v)_\lambda(x) := (v * \tilde{\psi}^\lambda)(x) = \langle v, \psi_x^\lambda \rangle$ where $\tilde{\psi} = \psi(-\cdot)$ is some fixed non-negative $\psi \in \mathcal{B}_r$ with integral 1 and symmetric in space. Applying this convolution to equation (2.15), and since $\mathscr{L}$ commutes with the convolution, we conclude that $(v)_\lambda$ solves the following PDE:

$$\begin{aligned}
\mathscr{L}(v)_\lambda &= -((v)_\lambda)^3 - ((v^3)_\lambda - ((v)_\lambda)^3) - 3 \sum_{\tau \in \mathcal{W}} \frac{\Upsilon(\tau)}{\tau!} (\mathcal{R}(\mathcal{I}(\tau) V^2))_\lambda - \\
&\quad - 3 \sum_{\tau_1, \tau_2 \in \mathcal{W}} \frac{\Upsilon(\tau_1) \Upsilon(\tau_2)}{\tau_1! \tau_2!} (\mathcal{R}(\mathcal{I}(\tau_1) \mathcal{I}(\tau_2) V))_\lambda + \sum_{\tau \in \partial \mathcal{W}} \frac{\Upsilon(\tau)}{\tau!} (\boldsymbol{\Pi}\tau)_\lambda, \quad (4.26)
\end{aligned}$$

where we have introduced a commutator between the cubic non-linearity and the convolution to preserve the damping cubic term which allows us to use the next result, consequence of the maximum principle of $(-\Delta)^s$ for $s \in (0, 1]$.

**Lemma 4.17**. *Let $s \in (0,1)$ and $\mathscr{L} = (\partial_t + (-\Delta)^s)$. If $u: [0,1] \times \mathbb{R}^d \to \mathbb{R}$ is a smooth function, 1-periodic in space which satisfies pointwise the following PDE:*

$$(\mathscr{L}u)(x) = -u^3(x) + g(x),$$

*where $g: \mathbb{R}^{1+d} \to \mathbb{R}$ is a bounded function, then*

$$|u(x)| \lesssim \max\left\{ \frac{1}{\sqrt{x_0}}, \|g\|^{\frac{1}{3}} \right\} \qquad \forall\, x \in (0, 1] \times \mathbb{R}^d,$$

*for an implicit proportionality constant depending only on $s$ and $d$.*

Applying Lemma 4.17 to $(v)_\lambda$ which solves (4.26) we conclude the following lemma.

**Lemma 4.18**. *For all $t, R, R' \in (0,1)$ such that $t + R' < 1$ and $\lambda \in (0, t)$ we have*

$$\begin{aligned}
\|(v)_\lambda\|_{P_{t+R}} &\lesssim \max_{\substack{\tau_1, \tau_2 \in \mathcal{W} \\ \tau \in \delta \mathcal{W}}} \left\{ (R-R')^{-s}, \|(v^3)_\lambda - (v)_\lambda^3\|_{P_{t+R'}}^{\frac{1}{3}}, \|\mathcal{R}(\mathcal{I}(\tau_1) V^2)_\lambda\|_{P_{t+R'}}^{\frac{1}{3}}, \right. \\
&\qquad\qquad \left. \|\mathcal{R}(\mathcal{I}(\tau_1) \mathcal{I}(\tau_2) V)_\lambda\|_{P_{t+R'}}^{\frac{1}{3}}, \|(\boldsymbol{\Pi}\tau)_\lambda\|_{P_{t+R'}}^{\frac{1}{3}} \right\}, \quad (4.27)
\end{aligned}$$

*for a proportionality constant that depends only on $s$ and $d$.*

The following Lemmas 4.19 and 4.20 allow us to remove the regularisation of the right hand side of (4.27) and control the terms on the right hand side. Both rely strongly on Theorem 4.2 and the improvement this has on the rest of the lemmas from Section 4.1.

**LEMMA 4.19**. *Fix $c \in (0,1)$ and $t \in (0,1)$ such that Assumption 4.1 holds and $t + 3\lambda_t < 1$ for $\lambda_t$ defined as in (4.4). Let $\beta \in (0,1)$ be such that $\beta < |\mathcal{I}(\tau)|$ for all $\tau \in \mathcal{V}_{0,\gamma}$. For $7 \leqslant k \in \mathbb{Z}^+$ we define $\tilde{\lambda}_t := k^{-1} \lambda_t$, then for all $R \in (\lambda_t, 1)$*

$$\|v - (v)_{\tilde{\lambda}_t}\|_{P_{t+R}} \lesssim (k-1)^{-\beta} \|v\|_{P_t}, \qquad (4.28)$$

$$\|(v^3)_{\tilde{\lambda}_t} - (v)^3_{\tilde{\lambda}_t}\|_{P_{t+\lambda_t}} \lesssim k^{-\beta} \|v\|^3_{P_t}, \qquad (4.29)$$

*with the proportionality constant depending only on $s, d$ and $\gamma$.*

**LEMMA 4.20**. *Fix $c \in (0,1)$ and $t \in (0,1)$ such that Assumption 4.1 holds and $t + 3\lambda_t < 1$ for $\lambda_t$ defined as in (4.4). For $k \in \mathbb{Z}^+$ we define $\tilde{\lambda}_t := k^{-1} \lambda_t$, then*

$$\|(\mathbf{\Pi}\tau)_{\tilde{\lambda}_t}\|_{P_{t+\lambda_t}} \lesssim c \|v\|^3_{P_t} \qquad \forall \tau \in \delta\mathcal{W}, \qquad (4.30)$$

$$\|(\mathcal{R}(V^2 \mathcal{I}(\tau)))_{\tilde{\lambda}_t}\|_{P_{t+\lambda_t}} \lesssim c \|v\|^3_{P_t} \qquad \forall \tau \in \mathcal{W}, \qquad (4.31)$$

$$\|(\mathcal{R}(V \mathcal{I}(\tau_1) \mathcal{I}(\tau_2)))_{\tilde{\lambda}_t}\|_{P_{t+\lambda_t}} \lesssim c \|v\|^3_{P_t} \qquad \forall \tau_1, \tau_2 \in \mathcal{W}, \qquad (4.32)$$

*with the proportionality constant depending on $k, s, d$ and $\gamma$.*

In Lemma 4.19 the parameter $k \in \mathbb{Z}^+$ gives us an smallness constant on (4.28) and (4.29). On the other hand, in Lemma 4.20 we have the smallness constant $c \in (0,1)$ in the terms (4.30),(4.31) and (4.32) that allows us to ignore the constant $k$. The proofs of the previous lemmas can be found in Section 4.3. We prove now our main result Theorem 1.1.

**Proof of Theorem 1.1.** Fix $c = c(d,s,\gamma) \in (0,1)$ small enough such that Theorem 4.2 holds, later we will impose another smallness condition only depending on $d, s$ and $\gamma$. Taking Remark 4.16 and Lemma 4.15 into account we assume that $T_c = 1$, which in particular tells us that Assumption 4.1 holds for all $t \in (0,1]$. On the other hand, observe for all $t \in (0, 1/2]$ such that $t + 3\lambda_t \geqslant 1$ we have that

$$\|v\|_{P_t} \lesssim (1-t)^{-s^{-1}} \lesssim 1 \qquad \forall t \in (0, 1/2],$$

which implies (1.2). Therefore, without loss of generality we assume that $t + 3\lambda_t < 1$ for all $t \in (0,1]$. We can now apply Lemma 4.19 and Lemma 4.20 to (4.27) from Lemma 4.18 with $R' = \lambda_t$ to obtain that for all $R \in (\lambda_t, 1)$

$$\|v\|_{P_{t+R}} \leqslant C \max\left\{(R-\lambda_t)^{-s}, \left((k-1)^{-\beta} + k^{-\beta/3} + c^{1/3} C'(k)\right) \|v\|_{P_t}\right\} \qquad (4.33)$$

for some $C = C(s,d,\gamma) > 0$ and $C'(k) > 0$. Choose $k = k(C) \in \mathbb{Z}^+$ big enough and $c = c(C, C'(k)) = c(s,d,\gamma) \in (0,1)$ small enough such that (4.33) turns into

$$\|v\|_{P_{t+R}} \leqslant \max\left\{C (R - \lambda_t)^{-s}, \frac{1}{2} \|v\|_{P_t}\right\} \qquad \forall R \in (\lambda_t, 1).$$

From here the argument is analogous to [MW20, Section 4.5], but we include it for completeness. If we consider $R \geqslant 2\lambda_t$ then we have $0 < R - 2\lambda_t$ and $R \leqslant 2(R - \lambda_t)$ from which we conclude

$$\|v\|_{P_{t+R}} \leqslant \max\left\{C\, 2^s R^{-s}, \frac{1}{2} \|v\|_{P_t}\right\}.$$

To guarantee that both terms inside the maximum are equal we need to chose $R$ such that

$$C\, 2^s\, R^{-s} = \frac{1}{2}\|v\|_{P_t} \iff R := 2^{1+s^{-1}} C^{s^{-1}} \|v\|_{P_t}^{-s^{-1}} = 2^{1+s^{-1}} C^{s^{-1}} \lambda_t, \tag{4.34}$$

and for this particular $R$ we have the bound

$$\|v\|_{P_{t+R}} \leqslant \frac{1}{2}\|v\|_{P_t}. \tag{4.35}$$

Observe that definition (4.34) of $R$ is consistent with the condition $R \geqslant 2\lambda_t$ since $C > 1$ implies that $2^{s^{-1}} C^{s^{-1}} > 1$. We define a finite sequence $0 = R_0 < \cdots < R_N = \frac{1}{2}$ for some $N \in \mathbb{N}$ by setting

$$R_{n+1} - R_n := 2^{1+s^{-1}} C^{s^{-1}} \|v\|_{P_{R_n}}^{-s^{-1}}, \tag{4.36}$$

as long as $R_{n+1}$ defined this way is strictly less than $\frac{1}{2}$. Observe that the map

$$n \mapsto 2^{1+s^{-1}} C^{s^{-1}} \|v\|_{P_{R_n}}^{-s^{-1}}$$

is increasing in $n$ which guarantees that the sequence $0 = R_0 < \cdots < R_N = \frac{1}{2}$ is indeed finite. Considering (4.35) for $t = R_{n-1}$ and $R = R_n - R_{n-1}$ for $n \in \{0, \ldots, N-1\}$ we obtain

$$\|v\|_{P_{R_n}} \leqslant \frac{1}{2}\|v\|_{P_{R_{n-1}}}. \tag{4.37}$$

Now we show that the bound $\|v\|_{P_t} \lesssim t^{-s}$ holds for $t \in \{R_0, \ldots, R_N\}$. First we have that (4.37) implies $\|v\|_{P_{R_n}} \leqslant 2^{-(n-k)} \|v\|_{P_{R_k}}$ for $0 \leqslant k \leqslant n \leqslant N-1$ and therefore

$$R_n = \sum_{k=0}^{n-1} R_{k-1} - R_k = \sum_{k=0}^{n-1} 2^{1+s^{-1}} C^{s^{-1}} \|v\|_{P_{R_n}}^{-s^{-1}} \leqslant 2^{1+s^{-1}} C^{s^{-1}} \|v\|_{P_{R_n}}^{-s^{-1}} \sum_{k=0}^{n-1} 2^{-(n-k)} \lesssim \|v\|_{P_{R_n}}^{-s^{-1}}, \tag{4.38}$$

which implies the desired bound $\|v\|_{P_{R_n}} \lesssim R_n^{-s}$. For the endpoint $\frac{1}{2} = R_N = R_N - R_{N-1} + R_{N-1}$ we have that either $R_{N-1} \geqslant \frac{1}{4}$ or $R_N - R_{N-1} \geqslant \frac{1}{4}$. In the first situation we use (4.38) for $n = N - 1$ and that $t \mapsto \|v\|_{P_t}^{-s^{-1}}$ is increasing to obtain $\frac{1}{4} \leqslant R_{N-1} \lesssim \|v\|_{P_{R_{N-1}}}^{-s^{-1}} \leqslant \|v\|_{P_{R_N}}^{-s^{-1}}$ which implies $\|v\|_{P_{R_N}} \lesssim 1$, similarly for the second case we obtain now by definition of $N$ that $\frac{1}{4} \leqslant R_N - R_{N-1} \leqslant 2^{1+s^{-1}} C^{s^{-1}} \|v\|_{P_{R_{N-1}}}^{-s^{-1}}$ which implies again $\|v\|_{P_{R_N}} \lesssim 1$. Moreover, since $t \mapsto \|v\|_{P_t}$ is decreasing this implies

$$\|v\|_{P_t} \leqslant \|v\|_{P_{R_N}} \lesssim 1 \quad \forall t \in [R_N, 1].$$

At last if $t \in (R_n, R_{n+1})$ for some $n \in \{0, N-1\}$ we use definition (4.36) of $R_{n+1} - R_n$, (4.38) and that $t \mapsto \|v\|_{P_t}^{-s^{-1}}$ is increasing to conclude

$$t \leqslant R_{n+1} = R_{n+1} - R_n + R_n \lesssim \|v\|_{P_{R_n}}^{-s^{-1}} + \|v\|_{P_{R_n}}^{-s^{-1}} \lesssim \|v\|_{P_t}^{-s^{-1}},$$

which implies the desired bound $\|v\|_{P_t} \lesssim t^{-s}$ and concludes the proof of Theorem 1.1. $\square$

### 4.3. Proof of Lemmas

**Lemmas from Section 4.1**

The following two remarks will be used through this section.

**Remark 4.21.** Since $2s > 1$ then $(t-\lambda)^{2s} + \lambda^{2s} \leqslant t^{2s}$ for all $t > 0$ and $\lambda \in (0, t)$, and therefore

$$(t^{2s}, 1] \oplus (-\lambda^{2s}, 0) = (t^{2s} - \lambda^{2s}, 1] \subseteq ((t-\lambda)^{2s}, 1] \tag{4.39}$$

where $\oplus$ denotes the Minkowski sum of sets. In particular for $t \in (0, 1)$ such that $t + 4\lambda_t < 1$ we have that $\lambda_t < 1$ and therefore for all $\lambda \in (0, \lambda_t]$ and $z \in P_{t+3\lambda}$ we have by (4.39) that for all $r \leqslant 3\lambda$

$$B_r(z) \subset P_{t+3\lambda} \oplus B_r \subset ((t + 3\lambda - r)^{2s}, 1] \times ((0, 1]^d \oplus B_1) \subset (t^{2s}, 1] \times B_4 \subset K,$$

and $[\Pi; \tau]_{B_r(z)} \leqslant [\Pi; \tau]_K$ and $[\Gamma \mathcal{I}(\tau); \sigma]_{B_r(z)} \leqslant [\Gamma \mathcal{I}(\tau); \sigma]_K$.

**Remark 4.22.** We will use that if $\tau \in \mathcal{V}$ then for all $x \in \mathbb{R}^{1+d}$ and $y \in B_\lambda(x)$ we have that

$$|(\Pi_x \mathcal{I}(\tau))(y)| \leqslant [\Gamma \mathcal{I}(\tau); \mathbf{1}]_{B_\lambda(x)} d(x, y)^{|\mathcal{I}(\tau)|}.$$

This follows from Lemma 2.9, the identity $\gamma(\mathcal{I}(\tau)) = \langle \mathbf{1}, \Gamma \mathcal{I}(\tau) \rangle$ and definition (A.2), and in particular $(\Pi_y \mathcal{I}(\tau))(y) = 0$.

**Proof of Lemma 4.4.** For $z \in \mathbb{R}^{1+d}$ and $\lambda > 0$ fixed, we apply Corollary 3.14 to the germ $U_{z,\lambda}^{\mathrm{loc}}$ as defined in (4.8) and using that $U_{z,\lambda} \equiv 0$ in $B_\lambda(z) \times ((z_0 - \lambda^{2s}, z_0] \times \mathbb{R}^d \setminus B_{3\lambda}(z_{1:d}))$ (since $\eta(\lambda^{-1}(\cdot - z_{1:d})) \equiv 0$ in $\mathbb{R}^d \setminus B_{3\lambda}(z_{1:d})$) we obtain the estimate

$$\begin{aligned}\lambda^\gamma [U_{z,\lambda}^{\mathrm{loc}}]_{\gamma; B_{\lambda/2}(z)} &\lesssim \lambda^\gamma \|\mathscr{L}(U_{z,\lambda}^{\mathrm{loc}})\|_{\gamma-2s; B_\lambda(z)} + \lambda^\gamma [U_{z,\lambda}^{\mathrm{loc}}]_{\gamma\text{-3pt}; B_\lambda(z)} + \lambda^\gamma [\Lambda_{z,\lambda}^{\mathrm{loc}}]_{\gamma\text{-3pt}; B_\lambda(z)} \\ &\quad + \|U_{z,\lambda}^{\mathrm{loc}}\|_{B_\lambda(z) \times B_{3\lambda}(z)},\end{aligned}$$

where $\Lambda_{z,\lambda}^{\mathrm{loc}}$ is a germ such that $U_{z,\lambda}^{\mathrm{loc}}$ satisfies the 3-point continuity (3.6). On the other hand, since $U \equiv U_{z,\lambda}^{\mathrm{loc}}$ in $B_\lambda(z) \times B_{2\lambda}(z)$ we have that $\Lambda_{z,\lambda}^{\mathrm{loc}}$ can be chosen as $\Lambda$, and both the Hölder and 3-point continuity seminorms are the same for $U$ and $U_{z,\lambda}$. Therefore we obtain the bound

$$\begin{aligned}\lambda^\gamma [U]_{\gamma; B_{\lambda/2}(z)} &\lesssim \lambda^\gamma \|\mathscr{L}(U_{z,\lambda}^{\mathrm{loc}})\|_{\gamma-2s; B_\lambda(z)} + \lambda^\gamma [U]_{\gamma\text{-3pt}; B_\lambda(z)} + \lambda^\gamma [\Lambda]_{\gamma\text{-3pt}; B_\lambda(z)} \\ &\quad + \|U_{z,\lambda}^{\mathrm{loc}}\|_{B_\lambda(z) \times B_{3\lambda}(z)}.\end{aligned}$$

Moreover, for all $(x, y) \in B_\lambda(z) \times B_{3\lambda}(z)$

$$\begin{aligned}|U_{z,\lambda}^{\mathrm{loc}}(x, y)| &= \left| v(y) - v(x) - \sum_\tau \frac{\Upsilon_x[\tau]}{\tau!} (\Pi_x \mathcal{I}(\tau))(y)\, \eta(\lambda^{-1}(y-z)_{1:d}) \right| \\ &= |U(x, y)\, \eta(\lambda^{-1}(y-z)_{1:d}) + (v(y) - v(x))(1 - \eta(\lambda^{-1}(y-z)_{1:d}))| \\ &\leqslant |U(x, y)|\, \eta(\lambda^{-1}(y-z)_{1:d}) + |v(y) - v(x)| \\ &\leqslant \|U\|_{B_\lambda(z) \times B_{3\lambda}(z)} + \|v\|_{B_\lambda(z)} + \|v\|_{B_{3\lambda}(z)} \\ &\leqslant \|U\|_{B_{3\lambda}(z)} + \|v\|_{B_{3\lambda}(z)},\end{aligned}$$

and therefore

$$\|U_{z,\lambda}^{\mathrm{loc}}\|_{B_\lambda(z) \times B_{3\lambda}(z)} \leqslant \|U\|_{B_{3\lambda}(z)} + \|v\|_{B_{3\lambda}(z)}.$$

We conclude the result using that $\mathscr{L}(U_{z,\lambda}^{\mathrm{loc}}) = \mathscr{L}(U_{z,\lambda}^{\mathrm{loc}} - U) + \mathscr{L}U$. $\square$

**Proof of Lemma 4.5.** Recall that $V$ satisfies the algebraic equation (2.14) in $\mathcal{D}^\gamma$, and therefore applying the model $\Pi_x$ and the operator $\mathscr{L}$ to this equation we obtain

$$\begin{aligned}
\mathscr{L}(\Pi_x(V(x))(\cdot)) &= -3 \sum_{\tau \in \mathcal{W}} \frac{\Upsilon[\tau]}{\tau!} \mathscr{L}(\Pi_x \mathcal{I}(\mathcal{I}(\tau) V_{\beta_\tau}^2(x))) \\
&\quad -3 \sum_{\tau_1, \tau_2 \in \mathcal{W}} \frac{\Upsilon[\tau_1] \Upsilon[\tau_2]}{\tau_1! \tau_2!} \mathscr{L}(\Pi_x \mathcal{I}(\mathcal{I}(\tau_1) \mathcal{I}(\tau_2) V_{\beta_{\tau_1,\tau_2}}(x))) \\
&\quad + \sum_{\tau \in \delta\mathcal{W} \cap \mathcal{V}_{0,\gamma}} \frac{\Upsilon[\tau]}{\tau!} \mathscr{L}(\Pi_x \mathcal{I}(\tau)) + \mathscr{L}(\Pi_x(v(x)\mathbf{1} + v_{\mathbf{X}}(x) \cdot \mathbf{X})).
\end{aligned}$$

Since $\Pi_x(v(x)\mathbf{1} + v_{\mathbf{X}}(x) \cdot \mathbf{X})(\cdot)$ is a linear in space polynomial it gets annihilated by $\mathscr{L}$. We use the weak admissibility of the model (Definition 2.11) to simplify the other terms. Using the expansion (2.25) of $V_{\beta_\tau}^2$ one has that

$$\mathcal{I}(\tau) V_{\beta_\tau}^2(x) := v^2(x)\mathbf{1} + \sum_{\sigma \in \mathcal{V}_{0,\beta_\tau}} \frac{2 v(x) \Upsilon_x[\sigma]}{\sigma!} \mathcal{I}(\sigma) \mathcal{I}(\tau)$$
$$+ \sum_{\substack{\mathcal{I}(\sigma_1)\mathcal{I}(\sigma_2) \in \mathcal{T}^+ \\ |\mathcal{I}(\sigma_1)\mathcal{I}(\sigma_2)| < \beta_\tau}} \frac{2 \Upsilon_x[\sigma_1] \Upsilon_x[\sigma_2]}{(\mathcal{I}(\sigma_1)\mathcal{I}(\sigma_2))!} \mathcal{I}(\sigma_1) \mathcal{I}(\sigma_2) \mathcal{I}(\tau),$$

and by Lemma 2.7 it is clear that none of the trees in this expression belong to $\mathcal{W}$. Therefore, for each tree $\sigma$ such that $\langle \sigma, \mathcal{I}(\tau) V_{\beta_\tau}^2 \rangle \neq 0$ the weak admissibility condition from Definition 2.11 takes the form $\mathscr{L}(\Pi_x \mathcal{I}(\sigma)) = \rho \, \Pi_x \sigma$ and therefore

$$\mathscr{L}(\Pi_x(\mathcal{I}(\mathcal{I}(\tau) V_{\beta_\tau}^2(x)))(\cdot)) = \rho \, \Pi_x(\mathcal{I}(\tau) V_{\beta_\tau}^2(x)).$$

Analogously we have that

$$\mathscr{L}(\Pi_x(\mathcal{I}(\mathcal{I}(\tau_1) \mathcal{I}(\tau_2) V_{\beta_{\tau_1,\tau_2}}))(\cdot)) = \rho \, \Pi_x(\mathcal{I}(\tau_1) \mathcal{I}(\tau_2) V_{\beta_{\tau_1,\tau_2}}).$$

On the other hand, we have that for $\tau \in \delta\mathcal{W}$ the weak admissibility condition does not include a localisation $\rho$ and therefore, since $\Gamma$ acts trivially on $\tau$, we have $\mathscr{L}(\Pi_x \mathcal{I}(\tau)) = \Pi_x \tau = \mathbf{\Pi} \tau$. We conclude that

$$\begin{aligned}
\mathscr{L}(\Pi_x(V(x))(\cdot)) &= -3 \sum_{\tau \in \mathcal{W}} \frac{\Upsilon[\tau]}{\tau!} \Pi_x(\mathcal{I}(\tau) V_{\beta_\tau}^2(x)) \rho \\
&\quad -3 \sum_{\tau_1, \tau_2 \in \mathcal{W}} \frac{\Upsilon[\tau_1] \Upsilon[\tau_2]}{\tau_1! \tau_2!} \Pi_x(\mathcal{I}(\tau_1) \mathcal{I}(\tau_2) V_{\beta_{\tau_1,\tau_2}}) \rho + \sum_{\tau \in \delta\mathcal{W} \cap \mathcal{V}_{0,\gamma}} \frac{\Upsilon[\tau]}{\tau!} \mathbf{\Pi} \tau.
\end{aligned}$$

Recalling that $v = \mathcal{R}V$ solves equation (2.15), identity (4.6) and using that $\mathscr{L}$ annihilates linear in space polynomials we obtain that

$$\begin{aligned}
&\mathscr{L}(U(x,\cdot)) \\
&= \mathscr{L}(v(\cdot) - \Pi_x(V(x))(\cdot)) \\
&= -v^3 - 3 \sum_{\tau \in \mathcal{W}} \frac{\Upsilon[\tau]}{\tau!} \{\mathcal{R}(V^2(x) \mathcal{I}(\tau)) - \Pi_x(V_{\beta_\tau}^2(x) \mathcal{I}(\tau))\} \rho \\
&\quad -3 \sum_{\tau_1, \tau_2 \in \mathcal{W}} \frac{\Upsilon[\tau_1] \Upsilon[\tau_2]}{\tau_1! \tau_2!} \{\mathcal{R}(V \mathcal{I}(\tau_1) \mathcal{I}(\tau_2)) - \Pi_x(V_{\beta_{\tau_1},\beta_{\tau_2}}(x) \mathcal{I}(\tau_1) \mathcal{I}(\tau_2))\} \rho + \sum_{\tau \in \partial \mathcal{W} \setminus \mathcal{V}_{0,\gamma}} \frac{\Upsilon[\tau]}{\tau!} \mathbf{\Pi} \tau \\
&\quad +3 \sum_{\tau \in \mathcal{W}} \frac{\Upsilon[\tau]}{\tau!} \mathcal{R}(V^2(x) \mathcal{I}(\tau))(1-\rho) + 3 \sum_{\tau_1, \tau_2 \in \mathcal{W}} \frac{\Upsilon[\tau_1] \Upsilon[\tau_2]}{\tau_1! \tau_2!} \mathcal{R}(V \mathcal{I}(\tau_1) \mathcal{I}(\tau_2))(1-\rho).
\end{aligned}$$

Since $\rho \equiv 1$ in $B_1$ and we are interested in $\mathscr{L}(U(x, \cdot))$ as a distribution in $B_1$, we can ignore the last two terms and conclude the result. $\square$

**Proof of Lemma 4.6.** For (4.11) we define $F_\varepsilon(x) := \mathcal{Q}_\varepsilon(V^2(x) \mathcal{I}(\tau)) = V_\beta^2(x) \mathcal{I}(\tau)$, where $\mathcal{Q}_\varepsilon$ is the projection into symbols of homogeneity less than $\varepsilon$. By (2.19) we have that $\beta \in (0,1)$, and by (2.25) we can write $F_\varepsilon$ as

$$F_\varepsilon(x) = v^2(x) \mathcal{I}(\tau) + 2v(x) \sum_{\sigma \in \mathcal{V}_\beta} \frac{\Upsilon_x[\sigma]}{\sigma!} \mathcal{I}(\sigma) \mathcal{I}(\tau) + \sum_{\substack{\sigma = \mathcal{I}(\sigma_1) \mathcal{I}(\sigma_2) \\ \sigma_1, \sigma_2 \in \mathcal{V}_\beta \\ |\mathcal{I}(\sigma_1) \mathcal{I}(\sigma_2)| < \beta}} \frac{2 \Upsilon_x[\sigma_1] \Upsilon_x[\sigma_2]}{(\mathcal{I}(\sigma_1) \mathcal{I}(\sigma_2))!} \mathcal{I}(\sigma_1) \mathcal{I}(\sigma_2) \mathcal{I}(\tau).$$

Since the structure group acts trivially on $\mathcal{I}(\mathcal{W})$ (Lemma 2.6) we have that

$$F_\varepsilon(y) - \Gamma_{yx} F_\varepsilon(x) = (V_\beta^2(y) - \Gamma_{yx} V_\beta^2(x)) \mathcal{I}(\tau_1),$$

and therefore for $\sigma_1, \sigma_2 \in \mathcal{T}$

$$[F_\varepsilon; \sigma_1 \sigma_2 \mathcal{I}(\tau_1)]_\varepsilon = \sup_{x,y} \frac{|\langle \sigma_1 \sigma_2 \mathcal{I}(\tau_1), F_\varepsilon(y) - \Gamma_{yx} F_\varepsilon(x) \rangle|}{d(x,y)^{\varepsilon - |\sigma_1 \sigma_2 \mathcal{I}(\tau_1)|}} = \sup_{x,y} \frac{|\langle \sigma_1 \sigma_2, V_\beta(y) - \Gamma_{yx} V_\beta(x) \rangle|}{d(x,y)^{\beta - |\sigma_1 \sigma_2|}} = [V_\beta^2; \sigma_1 \sigma_2],$$

which implies that $F_\varepsilon \in \mathcal{D}^\varepsilon(\Gamma)$. By the Reconstruction theorem (Theorem A.5) we have

$$|\langle \mathcal{R} F_\varepsilon - \Pi_x F_\varepsilon(x), \psi_x^r \rangle| \lesssim r^\varepsilon \sum_{\sigma_1, \sigma_2 \in \mathcal{T}} [\Pi; \sigma_1 \sigma_2 \mathcal{I}(\tau_1)]_{B_{2r}(x)} [F; \sigma_1 \sigma_2 \mathcal{I}(\tau_1)]_{B_{2r}(x)}$$

$$= r^{\beta + |\mathcal{I}(\tau_1)|} \sum_{\sigma_1, \sigma_2 \in \mathcal{T}} [\Pi; \sigma_1 \sigma_2 \mathcal{I}(\tau_1)]_{B_{2r}(x)} [V_\beta^2; \sigma_1 \sigma_2]_{B_{2r}(x)},$$

and since $\varepsilon > 0$ we have by the uniqueness of reconstruction that $\mathcal{R} F_\varepsilon = \mathcal{R}(\mathcal{I}(\tau_1) V^2)$ since $F_\varepsilon$ is a truncation of $\mathcal{I}(\tau_1) V^2$. In order to incorporate the localisation $\rho$ observe that

$$\langle \mathcal{R}(\mathcal{I}(\tau) V^2) - \Pi_x(\mathcal{I}(\tau) V_\beta^2(x)) \rho, \psi_x^r \rangle = \langle \mathcal{R} F_\varepsilon - \Pi_x F_\varepsilon(x), \rho \psi_x^r \rangle,$$

and it is enough to argue that $\rho \psi_x^r$ is a test function centred at $x$ and localised at scale $r$. To do this we write

$$\varphi_x^\lambda(y) = \rho(y) \psi_x^\lambda(y) = \rho(y) \lambda^{-(d+2s)} \psi(\lambda^{-1}(y-x)) \iff \varphi(y) = \rho(x + \lambda y) \psi(y).$$

It is clear that $\varphi$ defined like that (which depends on $x$ and on $\lambda$) is smooth and it is compactly supported in $B_1(0)$. To see that $\varphi \in \mathcal{B}_r$ it is enough to see that some norms of $\varphi$ remain bounded uniformly on $x \in \mathbb{R}^{1+d}$ and $\lambda \in (0,1)$. For this we check that

$$\partial^n \varphi(y) = \sum_{n_1 + n_2 = n} \partial^{n_1}(\rho(x + \lambda \cdot))(y) (\partial^{n_2} \psi)(y) = \sum_{n_1 + n_2 = n} \lambda^{|n_1|} \partial^{n_1}(\rho)(x + \lambda y) (\partial^{n_2} \psi)(y),$$

and therefore for all $|n| \leqslant r = 2$ we have that

$$\|\partial^n \varphi\|_\infty \leqslant \sum_{n_1 + n_2 = n} \lambda^{|n_1|} \|(\partial^{n_1} \rho)\|_\infty \|(\partial^{n_2} \psi)\|_\infty \leqslant \sum_{n_1 + n_2 = n} \lambda^{|n_1|} \lesssim 1.$$

To show (4.12) we define $F_\varepsilon(x) := \mathcal{Q}_\varepsilon(V(x)\mathcal{I}(\tau_1)\mathcal{I}(\tau_2)) = V_\beta(x)\mathcal{I}(\tau_1)\mathcal{I}(\tau_2)$. By (2.20) we have

$$F_\varepsilon(x) = v(x)\mathcal{I}(\tau_1)\mathcal{I}(\tau_2) + \mathbb{1}_{\beta>1}v_{\boldsymbol{X}}(x)\cdot \boldsymbol{X}\mathcal{I}(\tau_1)\mathcal{I}(\tau_2) + \sum_{\sigma\in\mathcal{V}_{0,\beta}}\frac{\Upsilon_x[\sigma]}{\sigma!}\mathcal{I}(\sigma)\mathcal{I}(\tau_1)\mathcal{I}(\tau_2).$$

Analogously, by Lemma 2.6 we have that $F_\varepsilon(y) - \Gamma_{yx}F_\varepsilon(x) = (V_\beta(y) - \Gamma_{yx}V_\beta(x))\mathcal{I}(\tau_1)\mathcal{I}(\tau_2)$, and therefore $[F_\varepsilon; \sigma\mathcal{I}(\tau_1)\mathcal{I}(\tau_2)] = [V_\beta; \sigma]$ for all $\sigma\in\mathcal{T}$, which implies that $F_\varepsilon\in\mathcal{D}^\varepsilon(\Gamma)$. By the Reconstruction theorem (Theorem A.5) we have

$$|\langle \mathcal{R}F_\varepsilon - \Pi_x F_\varepsilon(x), \psi_x^r\rangle| \lesssim r^\varepsilon \sum_{\sigma\in\mathcal{T}}[\Pi; \sigma\mathcal{I}(\tau_1)\mathcal{I}(\tau_2)]_{B_{2r}(x)}[F_\varepsilon; \sigma\mathcal{I}(\tau_1)\mathcal{I}(\tau_2)]_{B_{2r}(x)}$$
$$= r^{\beta+|\mathcal{I}(\tau_1)\mathcal{I}(\tau_2)|} \sum_{\sigma\in\mathcal{T}}[\Pi; \sigma\mathcal{I}(\tau_1)\mathcal{I}(\tau_2)]_{B_{2r}(x)}[V_\beta; \sigma]_{B_{2r}(x)},$$

and since $\varepsilon > 0$ we have by the uniqueness of reconstruction that $\mathcal{R}F_\varepsilon = \mathcal{R}(\mathcal{I}(\tau_1)\mathcal{I}(\tau_2)V)$ since $F_\varepsilon$ is a truncation of $\mathcal{I}(\tau_1)\mathcal{I}(\tau_2)V^2$. The same argument as before allows us to include the localiser $\rho$ and conclude. $\square$

We will need the following representation of the exponents $\alpha$ defined in (4.3).

**LEMMA 4.23**. *For all $\tau\in(\mathcal{T}\setminus\mathcal{P})\setminus\{\Xi\}$ we have that*

$$\alpha(\tau) = s^{-1}(|\tau| - |\mathfrak{n}(\tau)| + |\mathfrak{e}(\tau)|) - \mathfrak{m}(\tau) + 3. \tag{4.40}$$

**Proof.** We define $\alpha': \mathcal{T}\to\mathbb{R}$ as the right hand side of (4.40). Given $\tau\in\mathcal{T}$ we can easily check directly by the definition of $\alpha'$ that

$$\alpha'(\mathcal{I}(\tau)) = \alpha'(\tau), \qquad \alpha'(\tau\boldsymbol{X}^k) = \alpha'(\tau), \qquad \alpha'(\mathcal{I}_j(\tau)) = \alpha(\mathcal{I}(\tau)), \tag{4.41}$$

i.e., $\alpha'$ is invariant under planting, addition of both polynomial or edge decorations. Moreover, $\alpha'$ satisfies the following identities on sub-ternary trees:

$$\alpha'(\mathcal{I}(\tau_1)\mathcal{I}(\tau_2)) = \alpha'(\tau_1) + \alpha'(\tau_2), \qquad \alpha'(\mathcal{I}(\tau_1)\mathcal{I}(\tau_2)\mathcal{I}(\tau_3)) = \alpha'(\tau_1) + \alpha'(\tau_2) + \alpha'(\tau_3). \tag{4.42}$$

Since $\alpha$ as defined in (4.3) only depends on the number of leaves $\mathfrak{l}(\tau)$ then it is clear that it also satisfies properties (4.41) and (4.42). This implies that if $\tau_1,\tau_2,\tau_3\in\mathcal{T}$ are such that $\alpha(\tau_i) = \alpha'(\tau_i)$ then $\alpha(\sigma) = \alpha'(\sigma)$ for $\sigma\in\{\mathcal{I}(\tau_1)\boldsymbol{X}^k, \mathcal{I}_j(\tau_1)\boldsymbol{X}^k, \mathcal{I}(\tau_1)\mathcal{I}(\tau_2)\boldsymbol{X}^k, \mathcal{I}(\tau_1)\mathcal{I}(\tau_2)\mathcal{I}(\tau_3)\}$, and therefore by the recursive definition (2.1) of $\mathcal{T}$ it is only left to show that the result holds for some trees to be used as the basis for the induction. An explicit computation shows that $\alpha(\Xi)\neq\alpha'(\Xi)$, which is why it was left out in the statement. We show that the result holds true for $\mathcal{I}(\Xi), \mathcal{I}(\Xi)\mathcal{I}(\Xi)$, $\mathcal{I}(\Xi)\mathcal{I}(\Xi)\mathcal{I}(\Xi)\in\mathcal{T}$. Since $(\mathfrak{l}(\mathcal{I}(\Xi)), |\mathfrak{n}(\mathcal{I}(\Xi))|, |\mathfrak{e}(\mathcal{I}(\Xi))|, \mathfrak{m}(\mathcal{I}(\Xi))) = (1,0,0,2)$, we have that

$$\alpha'(\mathcal{I}(\Xi)) = s^{-1}|\mathcal{I}(\Xi)| + 1 = s^{-1}\left(s - \frac{3}{2} - \kappa\right) + 1 = s^{-1}\left(2s - \frac{3}{2} - \kappa\right) = \alpha(\mathcal{I}(\Xi)).$$

Analogous computations show that $\alpha'(\mathcal{I}(\Xi)\mathcal{I}(\Xi)) = \alpha(\mathcal{I}(\Xi)\mathcal{I}(\Xi))$ and $\alpha'(\mathcal{I}(\Xi)\mathcal{I}(\Xi)\mathcal{I}(\Xi)) = \alpha(\mathcal{I}(\Xi)\mathcal{I}(\Xi)\mathcal{I}(\Xi))$. By invariance under polynomial decorations the result is also true for $\mathcal{I}(\Xi)\boldsymbol{X}^k$, $\mathcal{I}(\Xi)\mathcal{I}(\Xi)\boldsymbol{X}^k$, which is enough to use as the basis of induction and conclude result. $\square$

Before proving Lemma 4.7 we have need the following lemma which uses Assumption 4.1.

**LEMMA 4.24**. *If $c, t \in (0,1)$ satisfy Assumption 4.1 and $t + 3\lambda_t < 1$ with $\lambda_t$ is defined by (4.4), then for all $\lambda \in (0, \lambda_t]$, $z \in P_{t+4\lambda}$, $r \leqslant 3\lambda$, $\tau \in \mathcal{V}_{0,\gamma}$ and $k \in \mathbb{N}^{1+d}$*

$$\lambda^{|\mathcal{I}(\tau)|} [\Gamma \mathcal{I}(\tau); \boldsymbol{X}^k]_{B_r(z)} \|\Upsilon.[\tau]\|_{B_r(z)} \leqslant \begin{cases} \|v\|_{P_t} & \text{if } \Upsilon[\tau] \nsim v_{\boldsymbol{X}} \\ c\lambda \|v_{\boldsymbol{X}}\|_{B_r(z)} & \text{if } \Upsilon[\tau] \sim v_{\boldsymbol{X}} \end{cases},$$

*with the implicit proportionality constants depending only on $s, d$ and $\gamma$.*

**Proof.** In the first case we have $\Upsilon.[\tau] \sim v^{\mathfrak{m}(\tau)}$, and using that $|\mathcal{I}(\tau)| > 0$ we have that

$$\begin{aligned}
\lambda^{|\mathcal{I}(\tau)|} [\Gamma \mathcal{I}(\tau); \boldsymbol{X}^k]_{B_r(z)} \|\Upsilon.[\tau]\|_{B_r(z)} &\lesssim [\Gamma \mathcal{I}(\tau); \boldsymbol{X}^k]_K \lambda^{|\mathcal{I}(\tau)|} \|v\|_{B_r(z)}^{\mathfrak{m}(\tau)} \\
&\leqslant [\Gamma \mathcal{I}(\tau); \boldsymbol{X}^k]_K \lambda_t^{|\mathcal{I}(\tau)|} \|v\|_{B_{3\lambda}(z)}^{\mathfrak{m}(\tau)} \\
&\leqslant c \|v\|_{P_t}^{\alpha(\tau)} \|v\|_{P_t}^{-s^{-1}|\mathcal{I}(\tau)|} \|v\|_{P_t}^{\mathfrak{m}(\tau)} \\
&\leqslant \|v\|_{P_t}^{\alpha(\tau) - s^{-1}|\mathcal{I}(\tau)| + \mathfrak{m}(\tau)},
\end{aligned} \quad (4.43)$$

where we used Remark 4.21 and Assumption (4.2) to bound uniformly the terms $[\Gamma \mathcal{I}(\tau); \boldsymbol{X}^k]_{B_\lambda(z)}$, and Remark 4.21 with the space periodicity of $v$ to bound $\|v\|_{B_{3\lambda}(z)} \leqslant \|v\|_{(t^{2s},1] \times \overline{B_4}} = \|v\|_{P_t}$. Recall that by Lemma 2.14 all the trees $\tau \in \mathcal{V}_{0,\gamma}$ such that $\Upsilon[\tau] \neq 0$ have no edge decorations, i.e., $|\mathfrak{e}(\tau)| = 0$, and in this case ($\Upsilon.[\tau] \sim v^{\mathfrak{m}(\tau)}$) Lemma 2.14 also implies that $\tau$ has no polynomial decorations, i.e., $|\mathfrak{n}(\tau)| = 0$. By Lemma 4.23 we have

$$\alpha(\tau) = s^{-1}|\tau| - \mathfrak{m}(\tau) + 3 = s^{-1}|\mathcal{I}(\tau)| - \mathfrak{m}(\tau) + 1 \quad (4.44)$$

which implies in (4.43) the claimed bound. For the case $\Upsilon.[\tau] \sim v_{\boldsymbol{X}}$ we have that $(\mathfrak{m}(\tau), |\mathfrak{n}(\tau)|) = (1, 1)$ by Lemma 2.14 and Lemma 2.7, and therefore $|\mathcal{I}(\tau)| > 1$ by Lemma 2.7 and Assumption (4.1) we have

$$\begin{aligned}
\lambda^{|\mathcal{I}(\tau)|} [\Gamma \mathcal{I}(\tau); \boldsymbol{X}^k]_{B_r(z)} \|\Upsilon.[\tau]\|_{B_r(z)} &\lesssim [\Gamma \mathcal{I}(\tau); \boldsymbol{X}^k]_K \lambda^{|\mathcal{I}(\tau)|} \|v_{\boldsymbol{X}}\|_{B_r(z)} \\
&\leqslant [\Gamma \mathcal{I}(\tau); \boldsymbol{X}^k]_K \lambda_t^{|\mathcal{I}(\tau)|-1} \lambda \|v_{\boldsymbol{X}}\|_{B_r(z)} \\
&\leqslant c \|v\|_{P_t}^{\alpha(\tau)} \|v\|_{P_t}^{-s^{-1}(|\mathcal{I}(\tau)|-1)} \lambda \|v_{\boldsymbol{X}}\|_{B_r(z)},
\end{aligned}$$

and the claimed bound follows since by Lemma 4.23

$$\alpha(\tau) = s^{-1}(|\tau| - 1) - 1 + 3 = s^{-1}(|\mathcal{I}(\tau)| - 1). \quad (4.45) \quad \square$$

**Proof of Lemma 4.7.** Since

$$\begin{aligned}
U_\beta(x,y) &= v(y) - v(x) - \sum_{\tau \in \mathcal{V}_{0,\beta}} \frac{\Upsilon_x[\tau]}{\tau!} (\Pi_x \mathcal{I}(\tau))(y) \\
&= v(y) - v(x) - \sum_{\tau \in \mathcal{V}_{0,\gamma}} \frac{\Upsilon_x[\tau]}{\tau!} (\Pi_x \mathcal{I}(\tau))(y) + \sum_{\tau \in \mathcal{V}_{\beta,\gamma}} \frac{\Upsilon_x[\tau]}{\tau!} (\Pi_x \mathcal{I}(\tau))(y) \\
&= U(x,y) + \sum_{\tau \in \mathcal{V}_{\beta,\gamma}} \frac{\Upsilon_x[\tau]}{\tau!} (\Pi_x \mathcal{I}(\tau))(y),
\end{aligned}$$

then we have for all $x, y \in B_{2\lambda}(z)$ with $y_0 \leqslant x_0$

$$
\begin{aligned}
&|U_\beta(x,y) - \mathbb{1}_{\beta>1} v_{\boldsymbol{X}}(x) \cdot (y-x)| \\
&\leqslant |U(x,y) - v_{\boldsymbol{X}}(x) \cdot (y-x)| + \mathbb{1}_{\beta \leqslant 1} |v_{\boldsymbol{X}}(x) \cdot (y-x)| + \sum_{\tau \in \mathcal{V}_{\beta,\gamma}} \frac{|\Upsilon_x[\tau]|}{\tau!} |(\Pi_x \mathcal{I}(\tau))(y)| \\
&\lesssim [U]_{\gamma;B_{2\lambda}(z)} d(x,y)^\gamma + \mathbb{1}_{\beta \leqslant 1} |v_{\boldsymbol{X}}(x)| d(x,y) + \sum_{\tau \in \mathcal{V}_{\beta,\gamma}} \|\Upsilon_\cdot[\tau]\|_{B_{2\lambda}(z)} [\Gamma \mathcal{I}(\tau); \mathbf{1}]_{B_{2\lambda}(z)} d(x,y)^{|\mathcal{I}(\tau)|} \\
&\leqslant d(x,y)^\beta \bigg( [U]_{\gamma;B_{2\lambda}(z)} \lambda^{\gamma-\beta} + \mathbb{1}_{\beta \leqslant 1} \|v_{\boldsymbol{X}}\|_{B_{2\lambda}(z)} \lambda^{1-\beta} + \sum_{\tau \in \mathcal{V}_{\beta,\gamma}} \|\Upsilon_\cdot[\tau]\|_{B_{2\lambda}(z)} [\Gamma \mathcal{I}(\tau); \mathbf{1}]_{B_{2\lambda}(z)} \lambda^{|\mathcal{I}(\tau)|-\beta} \bigg),
\end{aligned}
$$

where we used that $d(x,y) \lesssim \lambda$ and by definition $|\mathcal{I}(\tau)| > \beta$ for $\tau \in \mathcal{V}_{\beta,\gamma}$. We conclude that

$$
\begin{aligned}
{[U_\beta]_{\beta;B_{2\lambda}(z)}} \lesssim{} & [U]_{\gamma;B_{2\lambda}(z)} \lambda^{\gamma-\beta} + \mathbb{1}_{\beta \leqslant 1} \|v_{\boldsymbol{X}}\|_{B_{2\lambda}(z)} \lambda^{1-\beta} \\
& + \sum_{\tau \in \mathcal{V}_{\beta,\gamma}} \|\Upsilon_\cdot[\tau]\|_{B_{2\lambda}(z)} [\Gamma \mathcal{I}(\tau); \mathbf{1}]_{B_{2\lambda}(z)} \lambda^{|\mathcal{I}(\tau)|-\beta}.
\end{aligned} \tag{4.46}
$$

Multiplying (4.46) by $\lambda^\beta$ and using the first part of Lemma 4.24 to bound the terms in the sum we conclude (4.13) since $[U_\beta]_\beta = [V_\beta; \mathbf{1}]$. On the other hand, by (2.22) we have that

$$
\begin{aligned}
\langle \boldsymbol{X}_j, V_\beta(y) - \Gamma_{yx} V_\beta(x) \rangle &= v_{\boldsymbol{X}_j}(y) - v_{\boldsymbol{X}_j}(x) - \sum_{\mu \in \mathcal{V}_{1,\beta}} \frac{\Upsilon_x[\mu]}{\mu!} \gamma_{yx}(\mathcal{I}_j(\mu)) \\
&= v_{\boldsymbol{X}_j}(y) - v_{\boldsymbol{X}_j}(x) - \sum_{\mu \in \mathcal{V}_{1,\gamma}} \frac{\Upsilon_x[\mu]}{\mu!} \gamma_{yx}(\mathcal{I}_j(\mu)) + \sum_{\mu \in \mathcal{V}_{\beta,\gamma}} \frac{\Upsilon_x[\mu]}{\mu!} \gamma_{yx}(\mathcal{I}_j(\mu)) \\
&= \langle \boldsymbol{X}_j, V(y) - \Gamma_{yx} V(x) \rangle + \sum_{\mu \in \mathcal{V}_{\beta,\gamma}} \frac{\Upsilon_x[\mu]}{\mu!} \langle \boldsymbol{X}_j, \Gamma \mathcal{I}(\tau) \rangle,
\end{aligned}
$$

where the identity $\langle \boldsymbol{X}_j, \Gamma \mathcal{I}(\mu) \rangle = \gamma_{yx}(\mathcal{I}_j(\mu))$ follows from applying $(\mathrm{Id} \otimes \gamma)$ to (A.14). Therefore, given $x, y \in B$ we have

$$
\begin{aligned}
&|\langle \boldsymbol{X}_j, V_\beta(y) - \Gamma_{yx} V_\beta(x) \rangle| \\
&\lesssim |\langle \boldsymbol{X}_j, V(y) - \Gamma_{yx} V(x) \rangle| + \sum_{\mu \in \mathcal{V}_{\beta,\gamma}} |\Upsilon_x[\mu]| |\langle \boldsymbol{X}_j, \Gamma \mathcal{I}(\tau) \rangle| \\
&\leqslant [V; \boldsymbol{X}_j]_B d(x,y)^{\gamma-1} + \sum_{\mu \in \mathcal{V}_{\beta,\gamma}} \|\Upsilon_\cdot[\mu]\|_B [\Gamma \mathcal{I}(\mu); \boldsymbol{X}_j]_B d(x,y)^{|\mathcal{I}_j(\mu)|} \\
&\leqslant d(x,y)^{\beta-1} \bigg( [V; \boldsymbol{X}_j]_B d(x,y)^{\gamma-\beta} + \sum_{\mu \in \mathcal{V}_{\beta,\gamma}} \|\Upsilon_\cdot[\mu]\|_B [\Gamma \mathcal{I}(\mu); \boldsymbol{X}_j]_B d(x,y)^{|\mathcal{I}(\mu)|-\beta} \bigg),
\end{aligned}
$$

which implies

$$
[V_\beta; \boldsymbol{X}_j]_{B_{2\lambda}(z)} \lesssim \lambda^{\gamma-\beta} [V; \boldsymbol{X}_j]_{B_{2\lambda}(z)} + \sum_{\mu \in \mathcal{V}_{\beta,\gamma}} \lambda^{|\mathcal{I}(\mu)|-\beta} \|\Upsilon_\cdot[\mu]\|_{B_{2\lambda}(z)} [\Gamma \mathcal{I}(\mu); \boldsymbol{X}_j]_{B_{2\lambda}(z)}. \tag{4.47}
$$

The terms in the sum in (4.47) can be bounded by Lemma 4.24, and multiplying with $\lambda^\beta$ we conclude (4.14). To show (4.15) let $\tau \in \mathcal{V}_{1,\gamma}$. If $\Upsilon[\tau] \sim v$ then by Lemma 2.24 $[V_\beta; \mathcal{I}(\tau)] \sim [V_{\beta-|\mathcal{I}(\tau)|}; \mathbf{1}]$ and the corresponding case follows from (4.13). If $\Upsilon[\tau] \sim v_{\boldsymbol{X}}$ then by Lemma 2.24 $[V_\beta; \sigma] = [V_{\beta-(|\sigma|-1)}; \boldsymbol{X}_j]$ and the corresponding case follows from (4.14).

For $\Upsilon[\tau] \sim v^2$ first we consider some $\beta \in (0, 1]$. By (2.25) we can write

$$
\begin{aligned}
& v^2(y) - (\Pi_x V_\beta^2(x))(y) \\
= {}& v(y)\,(v(y) - (\Pi_x V_\beta(x))(y)) + v(y)\,(\Pi_x V_\beta(x))(y) - (\Pi_x V_\beta^2(x))(y) \\
= {}& v(y)\,(v(y) - (\Pi_x V_\beta(x))(y)) + v(y)\left( v(x) + \sum_{\sigma \in \mathcal{V}_{0,\beta}} \frac{\Upsilon_x[\sigma]}{\sigma!} (\Pi_x \mathcal{I}(\sigma))(y) \right) \\
& - v^2(x) - 2\,v(x) \sum_{\sigma \in \mathcal{V}_{0,\beta}} \frac{\Upsilon_x[\sigma]}{\sigma!} (\Pi_x \mathcal{I}(\sigma))(y) - \sum_{\substack{\sigma = \mathcal{I}(\sigma_1)\mathcal{I}(\sigma_2) \\ \sigma_1, \sigma_2 \in \mathcal{V}_{0,\beta} \\ |\mathcal{I}(\sigma_1)\mathcal{I}(\sigma_2)| < \beta}} \frac{\Upsilon_x[\sigma_1]\,\Upsilon_x[\sigma_2]}{\sigma_1!\,\sigma_2!} \Pi_x(\mathcal{I}(\sigma_1)\,\mathcal{I}(\sigma_2))(y) \\
= {}& v(y)\,(v(y) - (\Pi_x V_\beta(x))(y)) + v(x)\left( v(y) - v(x) - \sum_{\tau \in \mathcal{V}_{0,\beta}} \frac{\Upsilon_x[\sigma]}{\sigma!} (\Pi_x \mathcal{I}(\sigma))(y) \right) \\
& + \sum_{\sigma \in \mathcal{V}_{0,\beta}} \frac{\Upsilon_x[\sigma]}{\sigma!} (\Pi_x \mathcal{I}(\sigma))(y) \left( v(y) - v(x) - \sum_{\sigma_2 \in \mathcal{V}_{0,\beta-|\mathcal{I}(\tau)|}} \frac{\Upsilon_x[\sigma_2]}{\sigma_2!} (\Pi_x \mathcal{I}(\sigma_2))(y) \right) \\
= {}& (v(y) + v(x))\,(v(y) - (\Pi_x V_\beta(x))(y)) \\
& + \sum_{\sigma \in \mathcal{V}_{0,\beta}} \frac{\Upsilon_x[\sigma]}{\sigma!} (\Pi_x \mathcal{I}(\sigma))(y)\,(v(y) - (\Pi_x V_{\beta-|\mathcal{I}(\sigma)|}(x))(y)).
\end{aligned}
\tag{4.48}
$$

Given $\beta \in (0, \gamma]$ and $\tau \in \mathcal{V}_\beta$ such that $\Upsilon.[\tau] \sim v^2$ we claim that $\beta - |\mathcal{I}(\tau)| \leq 1$. Arguing by contradiction we have that

$$
1 < \beta - |\mathcal{I}(\tau)| \Leftrightarrow |\mathcal{I}(\tau \boldsymbol{X}_j)| = 2\,s + |\tau| + 1 = |\mathcal{I}(\tau)| + 1 < \beta \leq \gamma,
$$

which implies that $\tau \boldsymbol{X}_j \in \mathcal{V}_\gamma$ and $(\mathfrak{m}(\tau \boldsymbol{X}_j), \mathfrak{n}(\tau \boldsymbol{X}_j)) = (\mathfrak{m}(\tau), \mathfrak{n}(\tau) + e_j) = (2, e_j)$, a contradiction to Lemma 2.7. On the other hand, by Lemma 2.24 we have:

$$
\langle \mathcal{I}(\tau), V_\beta(y) - \Gamma_{yx} V_\beta(x) \rangle \sim v^2(y) - (\Pi_x V_{\beta-|\mathcal{I}(\tau)|}^2(x))(y),
$$

and since $\beta - |\mathcal{I}(\tau)| \leq 1$ we can use (4.48) to conclude for all $x, y \in B$ that

$$
\begin{aligned}
& |\langle \mathcal{I}(\tau), V_\beta(y) - \Gamma_{zx} V_\beta(x) \rangle| \\
\lesssim {}& |v^2(y) - (\Pi_x V_{\beta-|\mathcal{I}(\tau)|}^2(x))(y)| \\
\lesssim {}& 2(|v(y)| + |v(x)|)\,|v(y) - (\Pi_x V_{\beta-|\mathcal{I}(\tau)|}(x))(y)| \\
& + \sum_{\sigma \in \mathcal{V}_{0,\beta-|\mathcal{I}(\tau)|}} |\Upsilon_x[\sigma]|\,|(\Pi_x \mathcal{I}(\sigma))(y)|\,|v(y) - (\Pi_x V_{\beta-|\mathcal{I}(\tau)|-|\mathcal{I}(\sigma)|}(x))(y)| \\
\lesssim {}& \|v\|_B\,|\langle \boldsymbol{1}, V_{\beta-|\mathcal{I}(\tau)|}(y) - \Gamma_{yx} V_{\beta-|\mathcal{I}(\tau)|}(x) \rangle| \\
& + \sum_{\sigma \in \mathcal{V}_{0,\beta-|\mathcal{I}(\tau)|}} |\Upsilon_x[\sigma]|\,|(\Pi_x \mathcal{I}(\sigma))(y)|\,|\langle \boldsymbol{1}, V_{\beta-|\mathcal{I}(\tau)\mathcal{I}(\sigma)|}(y) - \Gamma_{yx} V_{\beta-|\mathcal{I}(\tau)\mathcal{I}(\sigma)|}(x) \rangle| \\
\lesssim {}& \|v\|_B\,[V_{\beta-|\mathcal{I}(\tau)|}; \boldsymbol{1}]_B\, d(x,y)^{\beta-|\mathcal{I}(\tau)|} \\
& + \sum_{\sigma \in \mathcal{V}_{0,\beta-|\mathcal{I}(\tau)|}} \|\Upsilon.[\sigma]\|_B [\Gamma \mathcal{I}(\tau); \boldsymbol{1}]_{B_{8\theta}(z)} d(x,y)^{|\mathcal{I}(\sigma)|} [V_{\beta-|\mathcal{I}(\tau)\mathcal{I}(\sigma)|}; \boldsymbol{1}]_B\, d(x,y)^{\beta-|\mathcal{I}(\tau)\mathcal{I}(\sigma)|} \\
= {}& \|v\|_B\,[V_{\beta-|\mathcal{I}(\tau)|}; \boldsymbol{1}]_B\, d(x,y)^{\beta-|\mathcal{I}(\tau)|} \\
& + \sum_{\sigma \in \mathcal{V}_{0,\beta-|\mathcal{I}(\tau)|}} \|\Upsilon.[\sigma]\|_B\,[\Gamma \mathcal{I}(\tau); \boldsymbol{1}]_B\,[V_{\beta-|\mathcal{I}(\tau)\mathcal{I}(\sigma)|}; \boldsymbol{1}]_B\, d(x,y)^{\beta-|\mathcal{I}(\tau)|},
\end{aligned}
$$

and therefore

$$
[V_\beta; \mathcal{I}(\tau)]_B \lesssim \|v\|_B\,[V_{\beta-|\mathcal{I}(\tau)|}; \boldsymbol{1}]_B + \sum_{\sigma \in \mathcal{V}_{0,\beta-|\mathcal{I}(\tau)|}} \|\Upsilon.[\sigma]\|_B\,[\Gamma \mathcal{I}(\tau); \boldsymbol{1}]_B\,[V_{\beta-|\mathcal{I}(\tau)\mathcal{I}(\sigma)|}; \boldsymbol{1}]_B.
$$

Since $\beta - |\mathcal{I}(\tau)| \leqslant 1$ then $\mathcal{V}_{0,\beta-|\mathcal{I}(\tau)|} \subset \mathcal{V}_{0,1}$ and $\Upsilon.[\sigma] \approx v_{\mathbf{X}}$ for all $\sigma \in \mathcal{V}_{0,\beta-|\mathcal{I}(\tau)|}$ by Lemma 2.7, which allows us to conclude the result with Lemma 4.24 and the previously shown (4.13) as

$$\begin{aligned}
&\lambda^{\beta-|\mathcal{I}(\tau)|}[V_\beta; \mathcal{I}(\tau)]_{B_{2\lambda}(z)} \\
&\lesssim \|v\|_{B_{2\lambda}(z)} \lambda^{\beta-|\mathcal{I}(\tau)|}[V_{\beta-|\mathcal{I}(\tau)|}; \mathbf{1}]_{B_{2\lambda}(z)} \\
&\quad + \sum_{\sigma \in \mathcal{V}_{0,\beta-|\mathcal{I}(\tau)|}} \lambda^{|\mathcal{I}(\sigma)|}[\Gamma \mathcal{I}(\tau); \mathbf{1}]_{B_{2\lambda}(z)} \|\Upsilon.[\sigma]\|_{B_{2\lambda}(z)} \lambda^{\beta-|\mathcal{I}(\tau)\mathcal{I}(\sigma)|}[V_{\beta-|\mathcal{I}(\tau)\mathcal{I}(\sigma)|}; \mathbf{1}]_{B_{2\lambda}(z)} \\
&\lesssim \|v\|_{P_t} \lambda^{\beta-|\mathcal{I}(\tau)|}[V_{\beta-|\mathcal{I}(\tau)|}; \mathbf{1}]_{B_{2\lambda}(z)} + \sum_{\sigma \in \mathcal{V}_{0,\beta-|\mathcal{I}(\tau)|}} \|v\|_{P_t} \lambda^{\beta-|\mathcal{I}(\tau)\mathcal{I}(\sigma)|}[V_{\beta-|\mathcal{I}(\tau)\mathcal{I}(\sigma)|}; \mathbf{1}]_{B_{2\lambda}(z)} \\
&\lesssim \|v\|_{P_t}(\|v\|_{P_t} + (\mathbb{1}_{\beta-|\mathcal{I}(\tau)|\leqslant 1} + c)\lambda\|v_{\mathbf{X}}\|_{B_{2\lambda}(z)} + \lambda^\gamma[U]_{\gamma; B_{2\lambda}(z)}) \\
&\quad + \sum_{\sigma \in \mathcal{V}_{0,\beta-|\mathcal{I}(\tau)|}} \|v\|_{P_t}(\|v\|_{P_t} + (\mathbb{1}_{\beta-|\mathcal{I}(\tau)\mathcal{I}(\sigma)|\leqslant 1} + c)\lambda\|v_{\mathbf{X}}\|_{B_{2\lambda}(z)} + \lambda^\gamma[U]_{\gamma; B_{2\lambda}(z)}) \\
&\lesssim \|v\|_{P_t}(\|v\|_{P_t} + \lambda\|v_{\mathbf{X}}\|_{B_{2\lambda}(z)} + \lambda^\gamma[U]_{\gamma; B_{2\lambda}(z)}),
\end{aligned}$$

which concludes the proof of (4.15). $\square$

**Proof of Lemma 4.8.** Given $\tau \in \mathcal{W}$ the condition $[V^2_{\beta_\tau}; \sigma_1\sigma_2] \neq 0$ implies that $0 \leqslant |\sigma_1|, |\sigma_2| < \beta_\tau \leqslant 1$, and therefore $\sigma_1, \sigma_2$ have no polynomial decorations by Lemma 2.7. By Lemma 2.21 we have that $[V^2_{\beta_\tau}; \sigma_1\sigma_2] = [V_{\beta_\tau+|\mathcal{I}(\mathcal{I}(\Xi))|}; \mathcal{I}(\sigma_1\sigma_2\mathcal{I}(\Xi))]$, and since $\sigma_1\sigma_2\mathcal{I}(\Xi)$ has no polynomial decorations then $\Upsilon.[\sigma_1\sigma_2\mathcal{I}(\Xi)] \approx v_{\mathbf{X}}$. Lemma 2.21 also implies that $\sigma_1\sigma_2$ has no edge decorations. Using (4.15) from Lemma 4.7 (not the case $\sim v_{\mathbf{X}}$ since $\beta_\tau \leqslant 1$) we have

$$\begin{aligned}
\lambda^{\beta_\tau-|\sigma_1\sigma_2|}[V^2_{\beta_\tau}; \sigma_1\sigma_2]_{B_{2\lambda}(z)} &= \lambda^{\beta_\tau+|\mathcal{I}(\mathcal{I}(\Xi))|-|\mathcal{I}(\sigma_1\sigma_2\mathcal{I}(\Xi))|}[V_{\beta_\tau+|\mathcal{I}(\mathcal{I}(\Xi))|}; \mathcal{I}(\sigma_1\sigma_2\mathcal{I}(\Xi))]_{B_{2\lambda}(z)} \\
&\lesssim \|v\|_{P_t}^{\mathfrak{m}(\sigma_1\sigma_2\mathcal{I}(\Xi))-1}(\|v\|_{P_t} + \lambda\|v_{\mathbf{X}}\|_{B_{2\lambda}(z)} + \lambda^\gamma[U]_{\gamma; B_{2\lambda}(z)}). \quad (4.49)
\end{aligned}$$

On the other hand, by Remark 4.21 and Assumption (4.1)

$$\begin{aligned}
\lambda^{|\mathcal{I}(\tau)|+2s+|\sigma_1\sigma_2|}[\Pi; \sigma_1\sigma_2\mathcal{I}(\tau)]_{B_{2\lambda}(z)} &= \lambda^{|\sigma_1\sigma_2\mathcal{I}(\tau)|+2s}[\Pi; \sigma_1\sigma_2\mathcal{I}(\tau)]_{B_{2\lambda}(z)} \\
&\leqslant c\,\lambda_t^{2s+|\sigma_1\sigma_2\mathcal{I}(\tau)|}\|v\|_{P_t}^{\alpha(\sigma_1\sigma_2\mathcal{I}(\tau))} \\
&= c\,\|v\|_{P_t}^{-s^{-1}(2s+|\sigma_1\sigma_2\mathcal{I}(\tau)|)+\alpha(\sigma_1\sigma_2\mathcal{I}(\tau))} \\
&= c\,\|v\|_{P_t}^{-2-s^{-1}|\sigma_1\sigma_2\mathcal{I}(\tau)|+\alpha(\sigma_1\sigma_2\mathcal{I}(\tau))},
\end{aligned}$$

where we used that $|\sigma_1\sigma_2\mathcal{I}(\tau)| + 2s = |\sigma_1\sigma_2| + |\mathcal{I}(\mathcal{I}(\tau))| \geqslant |\mathcal{I}(\mathcal{I}(\tau))| > 0$ since $\mathcal{I}(\tau) \in \mathcal{V}$ (even though $\tau_1 \in \mathcal{W}$). Moreover, we have that since $\tau$ is full it has no decorations and therefore $\sigma_1\sigma_2\mathcal{I}(\tau)$ also has no decorations, and by Lemma 4.23 we have that

$$\begin{aligned}
\alpha(\sigma_1\sigma_2\mathcal{I}(\tau)) &= s^{-1}(|\sigma_1\sigma_2\mathcal{I}(\tau)| - |\mathfrak{n}(\sigma_1\sigma_2\mathcal{I}(\tau))| + |\mathfrak{e}(\sigma_1\sigma_2\mathcal{I}(\tau))|) - \mathfrak{m}(\sigma_1\sigma_2\mathcal{I}(\tau)) + 3 \\
&= s^{-1}|\sigma_1\sigma_2\mathcal{I}(\tau)| - \mathfrak{m}(\sigma_1\sigma_2\mathcal{I}(\tau)) + 3,
\end{aligned}$$

and therefore

$$\lambda^{|\mathcal{I}(\tau)|+2s+|\sigma_1\sigma_2|}[\Pi; \sigma_1\sigma_2\mathcal{I}(\tau)]_{B_{2\lambda}(z)} \leqslant c\,\|v\|_{P_t}^{1-\mathfrak{m}(\sigma_1\sigma_2\mathcal{I}(\tau))}. \quad (4.50)$$

Since $\tau$ and $\Xi$ are full, then $\mathfrak{m}(\sigma_1\sigma_2\mathcal{I}(\tau)) = \mathfrak{m}(\sigma_1\sigma_2\mathcal{I}(\Xi))$ and we can combine bounds (4.49) and (4.50) to conclude that

$$\begin{aligned}
&\lambda^{\beta_\tau+|\mathcal{I}(\tau)|+2s}[\Pi; \sigma_1\sigma_2\mathcal{I}(\tau)]_{B_{2\lambda}(z)}[V^2_{\beta_\tau}; \sigma_1\sigma_2]_{B_{2\lambda}(z)} \\
&= \lambda^{|\mathcal{I}(\tau)|+2s+|\sigma_1\sigma_2|}[\Pi; \sigma_1\sigma_2\mathcal{I}(\tau)]_{B_{2\lambda}(z)} \lambda^{\beta_\tau-|\sigma_1\sigma_2|}[V^2_{\beta_\tau}; \sigma_1\sigma_2]_{B_{2\lambda}(z)} \\
&\leqslant c\,\|v\|_{P_t}^{1-\mathfrak{m}(\sigma_1\sigma_2\mathcal{I}(\tau))}\|v\|_{P_t}^{\mathfrak{m}(\sigma_1\sigma_2\mathcal{I}(\Xi))-1}(\|v\|_{P_t} + \lambda\|v_{\mathbf{X}}\|_{B_{2\lambda}(z)} + \lambda^\gamma[U]_{\gamma; B_{2\lambda}(z)}) \\
&\leqslant \|v\|_{P_t} + c\lambda\|v_{\mathbf{X}}\|_{B_{2\lambda}(z)} + \lambda^\gamma[U]_{\gamma; B_{2\lambda}(z)},
\end{aligned}$$

which is precisely the first stated bound.

Fix $\tau_1, \tau_2 \in \mathcal{W}$ and let $\beta := \beta_{\tau_1,\tau_2}$, it is enough to consider $\sigma \in \{\mathbf{1}, \mathbf{X}_j\}_{j=1}^d \cup \mathcal{I}(\mathcal{V}_{0,\beta})$ for the second bound. If $\sigma = \mathbf{1}$ then by Lemma 4.23 we have

$$\begin{aligned}\alpha(\mathcal{I}(\tau_1)\mathcal{I}(\tau_2)) &= s^{-1}\left(|\mathcal{I}(\tau_1)\mathcal{I}(\tau_2)| - |\mathfrak{n}(\mathcal{I}(\tau_1)\mathcal{I}(\tau_2))|\right) - \mathfrak{m}(\mathcal{I}(\tau_1)\mathcal{I}(\tau_2)) + 3 \\ &= s^{-1}|\mathcal{I}(\tau_1)\mathcal{I}(\tau_2)| + 2,\end{aligned}$$

where we used that since $\tau_1, \tau_2 \in \mathcal{W}$ then $\mathcal{I}(\tau_1)\mathcal{I}(\tau_2)$ has no decoration, and by definition of $\mathfrak{m}$ (see (2.2)) we have $\mathfrak{m}(\mathcal{I}(\tau_1)\mathcal{I}(\tau_2)) = 1 + \mathfrak{m}(\tau_1) + \mathfrak{m}(\tau_2) = 1$ since $\tau_1, \tau_2 \in \mathcal{W}$ are full. By Remark 4.21 and Assumption (4.1) we have

$$\begin{aligned}\lambda^{|\mathcal{I}(\tau_1)\mathcal{I}(\tau_2)|+2s}\left[\Pi; \mathcal{I}(\tau_1)\mathcal{I}(\tau_2)\right]_{B_{2\lambda}(z)} &\lesssim c\lambda_t^{2s+|\mathcal{I}(\tau_1)\mathcal{I}(\tau_2)|}\|v\|_{P_t}^{\alpha(\mathcal{I}(\tau_1)\mathcal{I}(\tau_2))} \\ &\leqslant c\|v\|_{P_t}^{-2-s^{-1}|\mathcal{I}(\tau_1)\mathcal{I}(\tau_2)|+\alpha(\mathcal{I}(\tau_1)\mathcal{I}(\tau_2))} = c,\end{aligned} \quad (4.51)$$

and combining this with (4.13) from Lemma 4.7 we conclude

$$\begin{aligned}&\lambda^{\beta+|\mathcal{I}(\tau_1)\mathcal{I}(\tau_2)|+2s}\left[\Pi; \mathcal{I}(\tau_1)\mathcal{I}(\tau_2)\right]_{B_{2\lambda}(z)}\left[V_{\beta_\tau}; \mathbf{1}\right]_{B_{2\lambda}(z)} \\ &= \lambda^{|\mathcal{I}(\tau_1)\mathcal{I}(\tau_2)|+2s}\left[\Pi; \mathcal{I}(\tau_1)\mathcal{I}(\tau_2)\right]_{B_{2\lambda}(z)}\lambda^{\beta_\tau}\left[V_{\beta_\tau}; \mathbf{1}\right]_{B_{2\lambda}(z)} \\ &\leqslant \|v\|_{P_t} + c\lambda\|v_{\mathbf{X}}\|_{B_{2\lambda}(z)} + \lambda^\gamma[U]_{\gamma;B_{2\lambda}(z)},\end{aligned}$$

which shows the case $\sigma = \mathbf{1}$. If $\sigma = \mathbf{X}_j$ for some $j \in \{1, \ldots d\}$, then we have by Lemma 4.23 that

$$\begin{aligned}\alpha(\mathbf{X}_j\mathcal{I}(\tau_1)\mathcal{I}(\tau_2)) &= s^{-1}\left(|\mathbf{X}_j\mathcal{I}(\tau_1)\mathcal{I}(\tau_2)| - |\mathfrak{n}(\mathbf{X}_j\mathcal{I}(\tau_1)\mathcal{I}(\tau_2))|\right) - \mathfrak{m}(\mathbf{X}_j\mathcal{I}(\tau_1)\mathcal{I}(\tau_2)) + 3 \\ &= s^{-1}\left(|\mathbf{X}_j\mathcal{I}(\tau_1)\mathcal{I}(\tau_2)| - 1\right) - \mathfrak{m}(\tau_1) - \mathfrak{m}(\tau_2) + 2 \\ &= s^{-1}|\mathcal{I}(\tau_1)\mathcal{I}(\tau_2)| + 2,\end{aligned} \quad (4.52)$$

and therefore by (4.14) from Lemma 4.7, Remark 4.21 and Assumption (4.1) we obtain

$$\begin{aligned}&\lambda^{\beta+|\mathcal{I}(\tau_1)\mathcal{I}(\tau_2)|+2s}\left[\Pi; \mathbf{X}_j\mathcal{I}(\tau_1)\mathcal{I}(\tau_2)\right]_{B_{2\lambda}(z)}\left[V_\beta; \mathbf{X}_j\right]_{B_{2\lambda}(z)} \\ &= \lambda^{2s+|\mathcal{I}(\tau_1)\mathcal{I}(\tau_2)|}\left[\Pi; \mathbf{X}_j\mathcal{I}(\tau_1)\mathcal{I}(\tau_2)\right]_{B_{2\lambda}(z)}\lambda^\beta\left[V_\beta; \mathbf{X}_j\right]_{B_{2\lambda}(z)} \\ &\leqslant c\lambda_t^{2s+|\mathcal{I}(\tau_1)\mathcal{I}(\tau_2)|}\|v\|_{P_t}^{\alpha(\mathbf{X}_j\mathcal{I}(\tau_1)\mathcal{I}(\tau_2))}\left(\|v\|_{P_t} + c\lambda\|v_{\mathbf{X}}\|_{B_{2\lambda}(z)} + \lambda^\gamma[V;\mathbf{X}]_{B_{2\lambda}(z)}\right) \\ &\leqslant c\|v\|_{P_t}^{-2-s^{-1}|\mathcal{I}(\tau_1)\mathcal{I}(\tau_2)|+\alpha(\mathbf{X}_j\mathcal{I}(\tau_1)\mathcal{I}(\tau_2))}\left(\|v\|_{P_t} + c\lambda\|v_{\mathbf{X}}\|_{B_{2\lambda}(z)} + \lambda^\gamma[V;\mathbf{X}]_{B_{2\lambda}(z)}\right) \\ &\leqslant \|v\|_{P_t} + c\lambda\|v_{\mathbf{X}}\|_{B_{2\lambda}(z)} + c\lambda^\gamma[V;\mathbf{X}]_{B_{2\lambda}(z)},\end{aligned}$$

which shows the case $\sigma = \mathbf{X}_j$. It is only left to consider $\mu \in \mathcal{V}_{0,\beta}$ and $\sigma = \mathcal{I}(\mu)$ which we split in two cases. If $|\mathfrak{n}(\sigma)| = 0$, then $|\mathfrak{n}(\mu)| = 0$, and therefore $\Upsilon[\mu] \sim v^{\mathfrak{m}(\sigma)}$. By Lemma 4.7 we have

$$\lambda^{\beta-|\mathcal{I}(\mu)|}[V_\beta;\mathcal{I}(\mu)]_{B_{2\lambda}(z)} \lesssim \|v\|_{P_t}^{\mathfrak{m}(\mu)-1}\left(\|v\|_{P_t} + (\mathbb{1}_{\beta-|\mathcal{I}(\mu)|\leqslant 1}+c)\lambda\|v_{\mathbf{X}}\|_{B_{2\lambda}(z)} + \lambda^\gamma[U]_{\gamma;B_{2\lambda}(z)}\right). \quad (4.53)$$

On the other hand, by Lemma 4.23

$$\begin{aligned}\alpha(\mathcal{I}(\mu)\mathcal{I}(\tau_1)\mathcal{I}(\tau_2)) &= s^{-1}\left(|\mathcal{I}(\mu)\mathcal{I}(\tau_1)\mathcal{I}(\tau_2)| - |\mathfrak{n}(\mathcal{I}(\mu)\mathcal{I}(\tau_1)\mathcal{I}(\tau_2))|\right) - \mathfrak{m}(\mathcal{I}(\mu)\mathcal{I}(\tau_1)\mathcal{I}(\tau_2)) + 3 \\ &= s^{-1}|\mathcal{I}(\mu)\mathcal{I}(\tau_1)\mathcal{I}(\tau_2)| - \mathfrak{m}(\mu) + 3,\end{aligned}$$

and by Remark 4.21 and Assumption (4.1) we conclude

$$\begin{aligned}\lambda^{|\mathcal{I}(\mu)|+|\mathcal{I}(\tau_1)\mathcal{I}(\tau_2)|+2s}\left[\Pi; \mathcal{I}(\mu)\mathcal{I}(\tau_1)\mathcal{I}(\tau_2)\right]_{B_{2\lambda}(z)} &\lesssim c\lambda_t^{2s+|\mathcal{I}(\mu)\mathcal{I}(\tau_1)\mathcal{I}(\tau_2)|}\|v\|_{P_t}^{\alpha(\mathcal{I}(\mu)\mathcal{I}(\tau_1)\mathcal{I}(\tau_2))} \\ &\leqslant c\|v\|_{P_t}^{-2-s^{-1}|\mathcal{I}(\mu)\mathcal{I}(\tau_1)\mathcal{I}(\tau_2)|+\alpha(\mathcal{I}(\mu)\mathcal{I}(\tau_1)\mathcal{I}(\tau_2))} \\ &= c\|v\|_{P_t}^{1-\mathfrak{m}(\mu)}.\end{aligned} \quad (4.54)$$

Combining (4.53) and (4.54) we conclude this case. On the other hand, if $|\mathfrak{n}(\sigma)|=1$, then $|\mathfrak{n}(\mu)|=1$, and therefore $\Upsilon[\mu] \sim v_{\boldsymbol{X}}$. By Lemma 2.24 we have $[V_\beta;\sigma]=[V_{\beta-(|\sigma|-1)};\boldsymbol{X}_j]$ for some $j \in \{1,\ldots,d\}$ and by Lemma 4.23 we have that

$$\begin{aligned}
\alpha(\sigma\mathcal{I}(\tau_1)\mathcal{I}(\tau_2)) &= s^{-1}(|\sigma\mathcal{I}(\tau_1)\mathcal{I}(\tau_2)| - |\mathfrak{n}(\sigma\mathcal{I}(\tau_1)\mathcal{I}(\tau_2))|) - \mathfrak{m}(\sigma\mathcal{I}(\tau_1)\mathcal{I}(\tau_2)) + 3 \\
&= s^{-1}(|\sigma\mathcal{I}(\tau_1)\mathcal{I}(\tau_2)| - |\mathfrak{n}(\sigma)|) - \mathfrak{m}(\sigma) + 3 \\
&= s^{-1}(|\sigma\mathcal{I}(\tau_1)\mathcal{I}(\tau_2)| - 1) + 2.
\end{aligned} \tag{4.55}$$

With this we can conclude that

$$\begin{aligned}
& \lambda^{\beta+|\mathcal{I}(\tau_1)\mathcal{I}(\tau_2)|+2s} [\Pi;\sigma\mathcal{I}(\tau_1)\mathcal{I}(\tau_2)]_{B_{2\lambda}(z)} [V_\beta;\sigma]_{B_{2\lambda}(z)} \\
&= \lambda^{2s+|\mathcal{I}(\tau_1)\mathcal{I}(\tau_2)|+|\sigma|-1} [\Pi;\sigma\mathcal{I}(\tau_1)\mathcal{I}(\tau_2)]_{B_{2\lambda}(z)} \lambda^{\beta-(|\sigma|-1)} [V_{\beta-(|\sigma|-1)};\boldsymbol{X}_j]_{B_{2\lambda}(z)} \\
&\leqslant c \lambda_t^{2s+|\mathcal{I}(\tau_1)\mathcal{I}(\tau_2)|-1} \|v\|_{P_t}^{\alpha(\sigma\mathcal{I}(\tau_1)\mathcal{I}(\tau_2))} (\|v\|_{P_t} + c\lambda\|v_{\boldsymbol{X}}\|_{B_{2\lambda}(z)} + \lambda^\gamma[V;\boldsymbol{X}]_{B_{2\lambda}(z)}) \\
&\leqslant c \|v\|_{P_t}^{-2-s^{-1}(|\sigma\mathcal{I}(\tau_1)\mathcal{I}(\tau_2)|-1)+\alpha(\sigma\mathcal{I}(\tau_1)\mathcal{I}(\tau_2))} (\|v\|_{P_t} + c\lambda\|v_{\boldsymbol{X}}\|_{B_{2\lambda}(z)} + \lambda^\gamma[V;\boldsymbol{X}]_{B_{2\lambda}(z)}) \\
&= \|v\|_{P_t} + c\lambda\|v_{\boldsymbol{X}}\|_{B_{2\lambda}(z)} + c\lambda^\gamma[V;\boldsymbol{X}]_{B_{2\lambda}(z)},
\end{aligned}$$

which concludes the proof. $\square$

**Proof of Lemma 4.9.** For $\tau,\tau_1,\tau_2 \in \mathcal{W}$ define the distributional germs

$$\begin{aligned}
F_\tau(x) &:= (\mathcal{R}(V^2(x)\mathcal{I}(\tau)) - \Pi_x(V^2_{\beta_\tau}(x)\mathcal{I}(\tau)))\rho \\
F_{\tau_1,\tau_2}(x) &:= (\mathcal{R}(V\mathcal{I}(\tau_1)\mathcal{I}(\tau_2)) - \Pi_x(V_{\beta_{\tau_1},\beta_{\tau_2}}(x)\mathcal{I}(\tau_1)\mathcal{I}(\tau_2)))\rho.
\end{aligned}$$

By Lemma 4.5 we have that

$$\begin{aligned}
& \|\mathscr{L}U\|_{\gamma-2s;B_\lambda(z)} \\
&\leqslant \|v^3\|_{\gamma-2s;B_\lambda(z)} + \sum_{\tau\in\mathcal{W}} \frac{|\Upsilon[\tau]|}{\tau!} \|F_\tau\|_{\gamma-2s;B_\lambda(z)} + \sum_{\tau_1,\tau_2} \frac{|\Upsilon[\tau_1]\Upsilon[\tau_2]|}{\tau_1!\tau_2!} \|F_{\tau_1,\tau_2}\|_{\gamma-2s;B_\lambda(z)} \\
&\quad + \sum_{\tau\in\delta\mathcal{W}\setminus\mathcal{V}_{0,\gamma}} \frac{\Upsilon[\tau]}{\tau!} \|\Pi\tau\|_{\gamma-2s;B_\lambda(z)} \\
&\lesssim \|v^3\|_{\gamma-2s;B_\lambda(z)} + \sum_{\tau\in\mathcal{W}} \|F_\tau\|_{\gamma-2s;B_\lambda(z)} + \sum_{\tau_1,\tau_2} \|F_{\tau_1,\tau_2}\|_{\gamma-2s;B_\lambda(z)} + \sum_{\tau\in\delta\mathcal{W}\setminus\mathcal{V}_{0,\gamma}} \|\Pi\tau\|_{\gamma-2s;B_\lambda(z)},
\end{aligned}$$

where the term with $v^3$ is interpreted as a germ constant on the base point. The bounds for the terms $F_\tau$ and $F_{\tau_1,\tau_2}$ follow immediately by Lemma 4.6 and Lemma 4.8, so we only need to bound the first and last terms. We proceed to bound each term using definition (3.4) of the seminorm:

i. Given $x \in B_\lambda(z)$ and $r > 0$ such that $B_r(x) \subset B_\lambda(z)$ we have

$$|\langle v^3, \psi_x^r\rangle| \leqslant \|v\|_{B_r(x)}^3 \|\psi_x^r\|_{L^1(B_r(x))} = \|v\|_{B_r(x)}^3 \leqslant \|v\|_{B_\lambda(z)}^3 \leqslant \|v\|_{P_t}^3,$$

where in the last bound we used Remark 4.21 and the 1-periodicity in space of $v$. Since $-(\gamma-2s) > 0$ and $0 < r \leqslant \lambda \leqslant \lambda_t$ we conclude that

$$\lambda^\gamma \|v^3\|_{\gamma-2s;B_\lambda(z)} \leqslant \lambda^\gamma \|v\|_{P_t}^3 \lambda^{-(\gamma-2s)} \leqslant \lambda_t^{2s}\|v\|_{P_t}^3 = (\|v\|_{P_t}^{-s^{-1}})^{2s} \|v\|_{P_t}^3 = \|v\|_{P_t}.$$

ii. For $\tau \in \mathcal{W}$ by Lemma 4.6 and Lemma 4.8 we have

$$\begin{aligned}
\lambda^\gamma \|F_\tau\|_{\gamma-2s;B_\lambda(z)} &= \lambda^\gamma \sup_{x\in B_\lambda(z)} \sup_{\psi\in\mathcal{B}_r} \sup_{\substack{r>0\\B_r(x)\subset B_\lambda(z)}} |\langle F_\tau, \psi_x^r\rangle| r^{-(\gamma-2s)} \\
&\lesssim \lambda^{\beta_\tau+|\mathcal{I}(\tau)|+2s} \sum_{\sigma_1,\sigma_2\in\mathcal{T}} [\Pi;\sigma_1\sigma_2\mathcal{I}(\tau)]_{B_{2\lambda}(z)} [V_\beta^2;\sigma_1\sigma_2]_{B_{2\lambda}(z)} \\
&\lesssim \|v\|_{P_t} + c\lambda\|v_{\boldsymbol{X}}\|_{B_{2\lambda}(z)} + c\lambda^\gamma[U]_{\gamma;B_{2\lambda}(z)}.
\end{aligned}$$

iii. Analogously, for $\tau_1, \tau_2 \in \mathcal{W}$ by Lemma 4.6 and Lemma 4.8 we have

$$\lambda^\gamma \|F_{\tau_1,\tau_2}\|_{\gamma-2s; B_\lambda(z)} \lesssim \|v\|_{P_t} + c\,\lambda\,\|v_{\boldsymbol{X}}\|_{B_{2\lambda}(z)} + c\,\lambda^\gamma [U]_{\gamma; B_{2\lambda}(z)}.$$

iv. For $\tau \in \delta\mathcal{W} \setminus \mathcal{V}_{0,\gamma}$, since $\Gamma\tau = \tau$ we have that $\boldsymbol{\Pi}\tau = \Pi_x \tau$ for any $x \in \mathbb{R}^{1+d}$, and therefore

$$\lambda^\gamma \|\boldsymbol{\Pi}\tau\|_{\gamma-2s; B_\lambda(z)} = \lambda^\gamma \sup_{x \in B_\lambda(z)} \sup_{\psi \in \mathcal{B}_r} \sup_{\substack{r>0 \\ B_r(x) \subset B_\lambda(z)}} |\langle \Pi_x\tau, \psi_x^r\rangle|\, r^{-(\gamma-2s)}$$
$$\leqslant \lambda^\gamma [\Pi;\tau]_{B_\lambda(z)}\, r^{|\tau|-\gamma+2s} = \lambda^{|\mathcal{I}(\tau)|} [\Pi;\tau]_{B_\lambda(z)},$$

where we used that since $\tau \in \delta\mathcal{W} \setminus \mathcal{V}_{0,\gamma}$ then $2s + |\tau| = |\mathcal{I}(\tau)| > \gamma$ and $r^{|\mathcal{I}(\tau)|-\gamma} \leqslant \lambda^{|\mathcal{I}(\tau)|-\gamma}$. We conclude that

$$\lambda^\gamma \|\boldsymbol{\Pi}\tau\|_{\gamma-2s; B_\lambda(z)} \lesssim [\Pi;\tau]_K\, \lambda^{|\mathcal{I}(\tau)|} \leqslant c\, \|v\|_{P_t}^{\alpha(\tau)}\, \lambda_t^{|\tau|+2s} \leqslant \|v\|_{P_t}^{\alpha(\tau)-s^{-1}|\tau|-2} = \|v\|_{P_t},$$

where the last line follows since $\tau \in \delta\mathcal{W}$ is full which implies

$$\alpha(\tau) = s^{-1}(|\tau| - |\mathfrak{n}(\tau)| + |\mathfrak{e}(\tau)|) - \mathfrak{m}(\tau) + 3 = s^{-1}|\tau| + 3. \tag{4.56}$$

Combining (i)-(iv) we conclude the proof. $\square$

**Proof of Lemma 4.10.** Since $V$ is a function-like modelled distribution then $(\Pi_y V(y))(y) = v(y)$ and by (4.6) then $U(y, y) = 0$. Therefore, for all $x, y, z \in B \subset \mathbb{R}^{1+d}$

$$\begin{aligned}
&U(x,z) - U(x,y) - U(y,z) + U(y,y) \\
&= v(z) - (\Pi_x V(x))(z) + v_{\boldsymbol{X}}(x)\cdot(z-x) - (v(y) - (\Pi_x V(x))(y) + v_{\boldsymbol{X}}(x)\cdot(y-x)) \\
&\quad - (v(z) - (\Pi_y V(y))(z) + v_{\boldsymbol{X}}(y)\cdot(z-y)) \\
&= (\Pi_y V(y))(z) - (\Pi_x V(x))(z) - (v(y) - (\Pi_x V(x))(y)) \\
&\quad + v_{\boldsymbol{X}}(x)\cdot(z-x) - v_{\boldsymbol{X}}(x)\cdot(y-x) - v_{\boldsymbol{X}}(y)\cdot(z-y) \\
&= (\Pi_y V(y))(z) - (\Pi_x V(x))(z) - (v(y) - (\Pi_x V(x))(y)) + v_{\boldsymbol{X}}(x)\cdot(z-y) - v_{\boldsymbol{X}}(y)\cdot(z-y) \\
&= (\Pi_y V(y))(z) - (\Pi_x V(x))(z) - (v(y) - (\Pi_x V(x))(y)) - (v_{\boldsymbol{X}}(y) - v_{\boldsymbol{X}}(x))\cdot(z-y).
\end{aligned}$$

By Lemma 2.20 we have that $v(y) - (\Pi_x V(x))(y) = \langle \mathbf{1}, V(y) - \Gamma_{yx} V(x)\rangle$, and using the relationship of the model $\Pi_y \Gamma_{yx} = \Pi_x$ we obtain

$$\begin{aligned}
&(\Pi_y V(y))(z) - (\Pi_x V(x))(z) - (v(y) - (\Pi_x V(x))(y)) \\
&= (\Pi_y V(y))(z) - (\Pi_x V(x))(z) - \langle \mathbf{1}, V(y) - \Gamma_{yx} V(x)\rangle \\
&= (\Pi_y(V(y) - \Gamma_{yx} V(x)))(z) - \langle \mathbf{1}, V(y) - \Gamma_{yx} V(x)\rangle \\
&= \sum_{\tau \in \mathcal{T}} \frac{\langle \tau, V(y) - \Gamma_{yx} V(x)\rangle}{\tau!} (\Pi_y \tau)(z) - \langle \mathbf{1}, V(y) - \Gamma_{yx} V(x)\rangle (\Pi_y \mathbf{1})(y) \\
&= \sum_{\tau \in \mathcal{T}\setminus\{\mathbf{1}\}} \frac{\langle \tau, V(y) - \Gamma_{yx} V(x)\rangle}{\tau!} (\Pi_y \tau)(z) \\
&= \sum_{j=1}^{d} \langle \boldsymbol{X}_j, V(y) - \Gamma_{yx} V(x)\rangle\, (\Pi_y \boldsymbol{X}_j)(z) + \sum_{\tau \in \mathcal{V}_{0,\gamma}} \frac{\langle \mathcal{I}(\tau), V(y) - \Gamma_{yx} V(x)\rangle}{\tau!} (\Pi_y \mathcal{I}(\tau))(z) \\
&= \langle \boldsymbol{X}, V(y) - \Gamma_{yx} V(x)\rangle \cdot (z-y) + \sum_{\tau \in \mathcal{V}_{0,\gamma}} \frac{\langle \mathcal{I}(\tau), V(y) - \Gamma_{yx} V(x)\rangle}{\tau!} (\Pi_y \mathcal{I}(\tau))(z),
\end{aligned}$$

where we used the admissibility of the model for the polynomial symbol, the representation of $V$ in (2.13), and that the span of polynomials and planted trees is a sector. By (2.22) we have that

$$\langle \boldsymbol{X}_j, V(y) - \Gamma_{yx} V(x)\rangle = v_{\boldsymbol{X}_j}(y) - v_{\boldsymbol{X}_j}(x) - \sum_{\sigma \in \mathcal{V}_{1,\gamma}} \frac{\Upsilon_x[\sigma]}{\sigma!} \gamma_{yx}(\mathcal{I}_j(\sigma)) = v_{\boldsymbol{X}_j}(y) - v_{\boldsymbol{X}_j}(x) - \Lambda_j(x,y),$$

and conclude that

$$U(x,z) - U(x,y) - U(y,z) + \Lambda(x,y) \cdot (z-y) = \sum_{\tau \in \mathcal{V}_{0,\gamma}} \frac{\langle \mathcal{I}(\tau), V(y) - \Gamma_{yx} V(x) \rangle}{\tau!} (\Pi_y \mathcal{I}(\tau))(z).$$

This identity implies then

$$\begin{aligned}
&|U(x,z) - U(x,y) - U(y,z) + U(y,y) + \Lambda(x,y) \cdot (z-y)| \\
&\lesssim \sum_{\tau \in \mathcal{V}_{0,\gamma}} |\langle \mathcal{I}(\tau), V(y) - \Gamma_{yx} V(x) \rangle| |(\Pi_y \mathcal{I}(\tau))(z)| \\
&\leqslant \sum_{\tau \in \mathcal{V}_{0,\gamma}} [V; \mathcal{I}(\tau)]_B \, d(x,y)^{\gamma - |\mathcal{I}(\tau)|} [\Gamma \mathcal{I}(\tau); \mathbf{1}]_B \, d(y,z)^{|\mathcal{I}(\tau)|} \\
&\leqslant \left( \sum_{\tau \in \mathcal{V}_{0,\gamma}} [V; \mathcal{I}(\tau)]_B [\Gamma \mathcal{I}(\tau); \mathbf{1}]_B \right) \left( \sum_{\tau \in \mathcal{V}_{0,\gamma}} d(x,y)^{\gamma - |\mathcal{I}(\tau)|} d(y,z)^{|\mathcal{I}(\tau)|} \right),
\end{aligned} \qquad (4.57)$$

which shows (4.17). For the second part, since $\gamma_{yx}(\mathcal{I}_j(\tau)) = \langle \mathbf{X}_j, \Gamma_{yx} \mathcal{I}(\tau) \rangle$ (which follows from (A.14)) we have for $\mathcal{I}_j(\tau) \in \mathcal{T}^+ \setminus \{0\}$ that $|\mathcal{I}_j(t)| > 0$ and therefore by (A.2) $\gamma_{yy}(\mathcal{I}(\tau)) = 0$, which implies that $\Lambda(y,y) = 0$. By (2.37) we can write for all $x, y \in B \subset \mathbb{R}^{1+d}$

$$V^{(j)}(x) = v_{\mathbf{X}_j}(x) \mathbf{1} + \sum_{\tau \in \mathcal{V}_{1,\gamma}} \frac{\Upsilon_x[\tau]}{\tau!} \mathcal{I}_j(\tau),$$

and therefore

$$\Lambda_j(x,y) = \sum_{\tau \in \mathcal{V}_{1,\gamma}} \frac{\Upsilon_x[\tau]}{\tau!} \gamma_{yx}(\mathcal{I}_j(\tau)) = \gamma_{yx}\left( \sum_{\tau \in \mathcal{V}_{1,\gamma}} \frac{\Upsilon_x[\tau]}{\tau!} \mathcal{I}_j(\tau) \right) = \gamma_{yx}(V^{(j)}(x)) - v_{\mathbf{X}_j}(x).$$

By Lemma 2.23 we have for all $x, y, z \in B$

$$\begin{aligned}
&\Lambda_j(x,z) - \Lambda_j(x,y) - \Lambda_j(y,z) + \Lambda_j(y,y) \\
&= \gamma_{zx}(V^{(j)}(x)) - v_{\mathbf{X}_j}(x) - \Lambda_j(x,y) - \gamma_{zy}(V^{(j)}(y)) + v_{\mathbf{X}_j}(y) \\
&= (v_{\mathbf{X}_j}(y) - v_{\mathbf{X}_j}(x) - \Lambda_j(x,y)) - \gamma_{zy}(V^{(j)}(y) - \Gamma_{yx}(V \wedge (j)(x))) \\
&= (v_{\mathbf{X}} - \Lambda)_j(x,y) - \sum_{\mu \in \mathcal{T}^+} \frac{\langle \mu, V^{(j)}(y) - \Gamma_{yx} V^{(j)}(x) \rangle}{\mu!} \gamma_{zy}(\mu) \\
&= (v_{\mathbf{X}} - \Lambda)_j(x,y) - \langle \mathbf{X}_j, V(y) - \Gamma_{yx} V(x) \rangle \gamma_{zy}(\mathbf{1}) - \sum_{\tau \in \mathcal{V}_{1,\gamma}} \frac{\langle \mathcal{I}(\tau), V(y) - \Gamma_{yx} V(x) \rangle}{\mathcal{I}_j(\tau)!} \gamma_{zy}(\mathcal{I}_j(\tau)) \\
&= -\sum_{\tau \in \mathcal{V}_{1,\gamma}} \frac{\langle \mathcal{I}(\tau), V(y) - \Gamma_{yx} V(x) \rangle}{\tau!} \gamma_{zy}(\mathcal{I}_j(\tau)),
\end{aligned}$$

the last identity following from (2.22) and $\gamma(\mathbf{1}) = 1$. By (A.24) $\gamma_{zy}(\mathcal{I}_j(\tau)) = \langle \mathbf{X}_j, \Gamma_{zy} \mathcal{I}(\tau) \rangle$ and

$$\begin{aligned}
&|\Lambda_j(x,z) - \Lambda_j(x,y) - \Lambda_j(y,z) + \Lambda_j(y,y)| \\
&\lesssim \sum_{\tau \in \mathcal{V}_{1,\gamma}} |\langle \mathcal{I}(\tau), V(y) - \Gamma_{yx} V(x) \rangle| |\gamma_{zy}(\mathcal{I}_j(\tau))| \\
&\leqslant \sum_{\tau \in \mathcal{V}_{1,\gamma}} [\mathcal{I}(\tau); V]_B \, d(x,y)^{\gamma - |\mathcal{I}(\tau)|} [\Gamma \mathcal{I}(\tau); \mathbf{X}]_B \, d(y,z)^{|\mathcal{I}_j(\tau)|} \\
&\leqslant \left( \sum_{\tau \in \mathcal{V}_{1,\gamma}} [\mathcal{I}(\tau); V]_B [\Gamma \mathcal{I}(\tau); \mathbf{X}]_B \right) \left( \sum_{\tau \in \mathcal{V}_{1,\gamma}} d(x,y)^{(\gamma-1) - |\mathcal{I}_j(\tau)|} d(y,z)^{|\mathcal{I}_j(\tau)|} \right),
\end{aligned}$$

where we used that $|\mathcal{I}_j(\tau)| = |\mathcal{I}(\tau)| - 1 > 0$ since $\tau \in \mathcal{V}_{1,\gamma}$, and which shows (4.18). □

**Proof of Lemma 4.11.** By (4.17) of Lemma 4.10 we have to bound $\lambda^\gamma [V; \mathcal{I}(\tau)]_{B_r(z)} [\Gamma \mathcal{I}(\tau); \mathbf{1}]_{B_r(z)}$ for $\tau \in \mathcal{V}_{0,\gamma}$. If is $\tau \in \mathcal{V}_{0,\gamma}$ such that $\Upsilon[\tau] \not\sim v_{\boldsymbol{X}}$, then by Lemma 4.7 we have

$$\lambda^{\gamma - |\mathcal{I}(\tau)|} [V; \mathcal{I}(\tau)]_{B_r(z)} \lesssim \|v\|_{P_t}^{\mathfrak{m}(\tau)-1} (\|v\|_{P_t} + \lambda \|v_{\boldsymbol{X}}\|_{B_r(z)} + \lambda^\gamma [U]_{\gamma; B_r(z)}). \tag{4.58}$$

On the other hand, by (4.44), Remark 4.21, and Assumption (4.2) we have that

$$\lambda^{|\mathcal{I}(\tau)|} [\Gamma \mathcal{I}(\tau); \mathbf{1}]_{B_r(z)} \leqslant c \, \lambda_t^{|\mathcal{I}(\tau)|} \|v\|_{P_t}^{\alpha(\tau)} = c \, \|v\|_{P_t}^{-s^{-1}|\mathcal{I}(\tau)| + \alpha(\tau)} = c \, \|v\|_{P_t}^{1 - \mathfrak{m}(\tau)},$$

which combined with (4.58) allows us to conclude that

$$\begin{aligned}
\lambda^\gamma [V; \mathcal{I}(\tau)]_{B_r(z)} [\Gamma \mathcal{I}(\tau); \mathbf{1}]_{B_r(z)} &\lesssim c \, \|v\|_{P_t}^{1 - \mathfrak{m}(\tau)} \|v\|_{P_t}^{\mathfrak{m}(\tau) - 1} (\|v\|_{P_t} + \lambda \|v_{\boldsymbol{X}}\|_{B_r(z)} + \lambda^\gamma [U]_{\gamma; B_r(z)}) \\
&\leqslant \|v\|_{P_t} + c \, \lambda \|v_{\boldsymbol{X}}\|_{B_r(z)} + c \, \lambda^\gamma [U]_{\gamma; B_r(z)}.
\end{aligned}$$

If $\tau \in \mathcal{V}_{0,\gamma}$ such that $\Upsilon[\tau] \sim v_{\boldsymbol{X}}$ we proceed analogously. By Lemma 4.7, Remark 4.21 and Assumption (4.2) we have

$$\begin{aligned}
\lambda^\gamma [V; \mathcal{I}(\tau)]_{B_r(z)} [\Gamma \mathcal{I}(\tau); \mathbf{1}]_{B_r(z)} &= \lambda^{\gamma - |\mathcal{I}(\tau)| + 1} [V; \mathcal{I}(\tau)]_{B_r(z)} \lambda^{(|\mathcal{I}(\tau)| - 1)} [\Gamma \mathcal{I}(\tau); \mathbf{1}]_{B_r(z)} \\
&\leqslant (\|v\|_{P_t} + c \, \lambda \|v_{\boldsymbol{X}}\|_{B_r(z)} + \lambda^\gamma [V; \boldsymbol{X}]_{B_r(z)}) \, c \, \|v\|_{P_t}^{-s^{-1}(|\mathcal{I}(\tau)| - 1) + \alpha(\tau)} \\
&\leqslant \|v\|_{P_t} + c \, \lambda \|v_{\boldsymbol{X}}\|_{B_r(z)} + c \, \lambda^\gamma [V; \boldsymbol{X}]_{B_r(z)},
\end{aligned}$$

since in this case $\alpha(\tau) = s^{-1}(|\mathcal{I}(\tau)| - 1)$ by (4.45). Putting both cases together and using (4.17) we conclude the first part. For the second part, by (4.18) of Lemma 4.10 we have to bound $\lambda^\gamma [V; \mathcal{I}(\tau)]_{B_r(z)} [\Gamma \mathcal{I}(\tau); \boldsymbol{X}]_{B_r(z)}$ for $\tau \in \mathcal{V}_{1,\gamma}$. If $\tau \in \mathcal{V}_{1,\gamma}$ is such that $\Upsilon[\tau] \not\sim v_{\boldsymbol{X}}$ then by Remark 4.21, Assumption (4.2) and (4.44) we have

$$\lambda^{|\mathcal{I}(\tau)|} [\Gamma \mathcal{I}(\tau); \boldsymbol{X}]_{B_r(z)} \leqslant c \, \lambda_t^{|\mathcal{I}(\tau)|} \|v\|_{P_t}^{\alpha(\tau)} = c \, \|v\|_{P_t}^{-s^{-1}|\mathcal{I}(\tau)| + \alpha(\tau)} = c \, \|v\|_{P_t}^{1 - \mathfrak{m}(\tau)},$$

which combined with Lemma 4.7 allows us to conclude

$$\begin{aligned}
\lambda^\gamma [\mathcal{I}(\tau); V]_{B_r(z)} [\Gamma \mathcal{I}(\tau); \boldsymbol{X}]_{B_r(z)} &= \lambda^{\gamma - |\mathcal{I}(\tau)|} [V; \mathcal{I}(\tau)]_{B_r(z)} \lambda^{|\mathcal{I}(\tau)|} [\Gamma \mathcal{I}(\tau); \boldsymbol{X}]_{B_r(z)} \\
&\lesssim \|v\|_{P_t}^{\mathfrak{m}(\tau) - 1} (\|v\|_{P_t} + \lambda \|v_{\boldsymbol{X}}\|_{B_r(z)} + \lambda^\gamma [U]_{\gamma; B_r(z)}) \, c \, \|v\|_{P_t}^{1 - \mathfrak{m}(\tau)} \\
&\leqslant \|v\|_{P_t} + c \, \lambda \|v_{\boldsymbol{X}}\|_{B_r(z)} + c \, \lambda^\gamma [U]_{\gamma; B_r(z)}.
\end{aligned}$$

If $\tau \in \mathcal{V}_{1,\gamma}$ is such that $\Upsilon[\tau] \sim v_{\boldsymbol{X}}$, then we have by Lemma 4.7, Remark 4.21, Assumption (4.2), and (4.45) that

$$\begin{aligned}
&\lambda^\gamma [V; \mathcal{I}(\tau)]_{B_r(z)} [\Gamma \mathcal{I}(\tau); \boldsymbol{X}]_{B_r(z)} \\
&\leqslant \lambda^{\gamma - |\mathcal{I}(\tau)| + 1} [V; \mathcal{I}(\tau)]_{B_r(z)} \lambda^{(|\mathcal{I}(\tau)| - 1)} [\Gamma \mathcal{I}(\tau); \boldsymbol{X}]_{B_r(z)} \\
&\leqslant (\|v\|_{P_t} + c \, \lambda \|v_{\boldsymbol{X}}\|_{B_r(z)} + \lambda^\gamma [V; \boldsymbol{X}]_{B_r(z)}) \, c \, \|v\|_{P_t}^{-s^{-1}(|\mathcal{I}(\tau)| - 1) + \alpha(\tau)} \\
&\leqslant \|v\|_{P_t} + c \, \lambda \|v_{\boldsymbol{X}}\|_{B_r(z)} + c \, \lambda^\gamma [V; \boldsymbol{X}]_{B_r(z)}.
\end{aligned}$$

Putting both cases together and using (4.18) we conclude the result. $\square$

**Proof of Lemma 4.12.** To show (4.19) let $x, y \in B_r(z)$, then

$$\begin{aligned}
|U_\beta(x, y)| &\leqslant |v(y)| + |v(x)| + \sum_{\tau \in \mathcal{V}_{0,\beta}} \frac{|\Upsilon_x[\tau]|}{\tau!} |\Pi_x(\mathcal{I}(\tau))(y)| \\
&\lesssim \|v\|_{P_t} + \sum_{\tau \in \mathcal{V}_{0,\beta}} \lambda^{|\mathcal{I}(\tau)|} [\Gamma \mathcal{I}(\tau); \mathbf{1}]_{B_r(z)} \|\Upsilon_\cdot [\tau]\|_{B_r(z)},
\end{aligned}$$

and therefore by Lemma 4.24 we have

$$\begin{aligned}
\|U_\beta\|_{B_r(z)} &\lesssim \|v\|_{P_t} + \sum_{\tau \in \mathcal{V}_{0,\beta}} (\|v\|_{P_t} + \mathbb{1}_{\beta > 1} c \, \lambda \|v_{\boldsymbol{X}}\|_{B_r(z)}) \\
&\leqslant \|v\|_{P_t} + \mathbb{1}_{\beta > 1} c \, \lambda \|v_{\boldsymbol{X}}\|_{B_r(z)},
\end{aligned}$$

where we used that if $\beta \leqslant 1$ then $\Upsilon[\tau] \not\propto v_{\boldsymbol{X}}$ for all $\tau \in \mathcal{V}_{0,\beta}$ by Lemma 2.7 and Lemma 2.14.

To show (4.20) we have by definition (3.4) of the seminorm we have

$$
\begin{aligned}
&\lambda^{\gamma} \|\mathcal{L}(U_{z,\lambda}^{\mathrm{loc}} - U)\|_{\gamma - 2s; B_\lambda(z)} \\
&= \lambda^{\gamma} \sup_{\substack{x \in B_\lambda(z) \\ B_r(x) \subset B_\lambda(z)}} \sup_{r > 0} |\langle \mathcal{L}(U_{z,\lambda}^{\mathrm{loc}} - U)(x,\cdot), \psi_x^r \rangle| r^{-(\gamma - 2s)} \\
&\lesssim \lambda^{2s} \sup_{\substack{x \in B_\lambda(z) \\ B_r(x) \subset B_\lambda(z)}} \sup_{r > 0} \|(-\Delta)^s (U_{z,\lambda}^{\mathrm{loc}}(x,\cdot) - U(x,\cdot))\|_{B_r(x)},
\end{aligned}
\tag{4.59}
$$

where we used that $\|\psi_x^r\|_{L^1} = 1$ and that $\partial_t (U_{z,\lambda}^{\mathrm{loc}}(x,\cdot) - U(x,\cdot)) = 0$ in $B_\lambda(z)$ by locality of $\partial_t$ since $U_{z,\lambda}^{\mathrm{loc}}(x,\cdot) \equiv U(x,\cdot)$ in $B_{2\lambda}(z)$. On the other hand, we can write the difference $U_{z,\lambda}^{\mathrm{loc}} - U$ as

$$
U_{z,\lambda}(x,\cdot) - U(x,\cdot) = \sum_{\tau \in \mathcal{V}_{0,\gamma}} \frac{\Upsilon_x[\tau]}{\tau!} \Pi_x \mathcal{I}(\tau)(1 - \eta(\lambda^{-1}(\cdot - z)_{1:d})),
$$

and therefore

$$
(-\Delta)^s (U_{z,\lambda}(x,\cdot) - U(x,\cdot))(y) = \sum_{\tau \in \mathcal{V}_{0,\gamma}} \frac{\Upsilon_x[\tau]}{\tau!} (-\Delta)^s (\Pi_x \mathcal{I}(\tau)(1 - \eta(\lambda^{-1}(\cdot - z)_{1:d})))(y). \tag{4.60}
$$

We proceed to bound the $L^\infty$-norm of (4.60) in $B_r(x) \subset B_\lambda(z)$. First we make the observation that if $f: \mathbb{R}^d \to \mathbb{R}$ is such that $f = 0$ in a neighbourhood around $y_{1:d} \in \mathbb{R}^d$, then (assuming $f$ has growth as in (1.6)) the singular integral (1.4) which defines the action of $(-\Delta)^s$ can be written as

$$
((-\Delta)^s f)(y_{1:d}) = 2 c_{d,s} \int_{\mathbb{R}^d} f(w) \frac{\mathrm{d}w}{|w - y_{1:d}|^{d+2s}}.
$$

By choice of $\eta$ we have that $1 - \eta(\lambda^{-1}(\cdot - z_{1:d})) = 0$ in $B_{2\lambda}(z_{1:d})$ and therefore for all $y \in B_\lambda(z)$ we have that

$$
\begin{aligned}
&|(-\Delta)^s (\Pi_x \mathcal{I}(\tau)(1 - \eta(\lambda^{-1}(\cdot - z)_{1:d})))(y)| \\
&= 2 c_{d,s} \left| \int_{\mathbb{R}^d \setminus B_{2\lambda}(z_{1:d})} (\Pi_x \mathcal{I}(\tau)(1 - \eta(\lambda^{-1}(\cdot - z_{1:d}))))(y_0, w) \frac{\mathrm{d}w}{|w - y_{1:d}|^{d+2s}} \right| \\
&\lesssim \int_{B_{1+\lambda}(z_{1:d}) \setminus B_{2\lambda}(z_{1:d})} |(\Pi_x \mathcal{I}(\tau))(y_0, w)| \frac{\mathrm{d}w}{|w - y_{1:d}|^{d+2s}} \\
&\quad + \int_{\mathbb{R}^d \setminus B_{1+\lambda}(z_{1:d})} |(\Pi_x \mathcal{I}(\tau))(y_0, w)| \frac{\mathrm{d}w}{|w - y_{1:d}|^{d+2s}}
\end{aligned}
\tag{4.61}
$$

We consider the first term in (4.61). Since $y \in B_r(x) \subset B_\lambda(z)$, then $d(x, (y_0, w)) \leqslant \lambda + |w - y_{1:d}|$, and since $0 < |\mathcal{I}(\tau)| < \gamma$ for $\tau \in \mathcal{V}_{0,\gamma}$, we have that $d(x, (y_0, w))^{|\mathcal{I}(\tau)|} \lesssim \lambda^{|\mathcal{I}(\tau)|} + |w - y_{1:d}|^{|\mathcal{I}(\tau)|}$, for a constant depending only on $\gamma$. On the other hand, we have that $\mathbb{R}^d \setminus B_{2\lambda}(z_{1:d}) \subset \mathbb{R}^d \setminus B_\lambda(y_{1:d})$ since $|w - y_{1:d}| \geqslant |w - z_{1:d}| - d(z, y) \geqslant |w - z_{1:d}| - \lambda$, and therefore

$$
\begin{aligned}
&\int_{B_{1+\lambda}(z_{1:d}) \setminus B_{2\lambda}(z_{1:d})} |(\Pi_x \mathcal{I}(\tau))(y_0, w)| \frac{\mathrm{d}w}{|w - y_{1:d}|^{d+2s}} \\
&\leqslant [\Gamma \mathcal{I}(\tau); \mathbf{1}]_K \int_{B_{1+\lambda}(z_{1:d}) \setminus B_{2\lambda}(z_{1:d})} \frac{d(x, (y_0, w))^{|\mathcal{I}(\tau)|}}{|w - y_{1:d}|^{d+2s}} \mathrm{d}w \\
&\leqslant [\Gamma \mathcal{I}(\tau); \mathbf{1}]_K \int_{\mathbb{R}^d \setminus B_\lambda(y_{1:d})} \frac{(\lambda^{|\mathcal{I}(\tau)|} + |w - y_{1:d}|^{|\mathcal{I}(\tau)|})}{|w - y_{1:d}|^{d+2s}} \mathrm{d}w \\
&\lesssim [\Gamma \mathcal{I}(\tau); \mathbf{1}]_K \int_\lambda^\infty \frac{(\lambda^{|\mathcal{I}(\tau)|} + r^{|\mathcal{I}(\tau)|})}{r^{1+2s}} \mathrm{d}r = \lambda^{|\tau|} [\Gamma \mathcal{I}(\tau); \mathbf{1}]_K,
\end{aligned}
$$

where we used that $x, (y_0, w) \in K$ for all $w \in B_{1+\lambda}(z_{1:d}) \subset \mathbb{R}^d$. For the second term in (4.61), since the integral $\int_1^\infty r^{-1-2s}\,\mathrm{d}r$ is finite, we have that

$$
\int_{\mathbb{R}^d \setminus B_{1+\lambda}(z_{1:d})} |(\Pi_x \mathcal{I}(\tau))(y_0, w)| \frac{\mathrm{d}w}{|w-y|^{d+2s}} \leqslant \int_{\mathbb{R}^d \setminus B_1(y)} \|\gamma_{\cdot,x}(\mathcal{I}(\tau))\|_{(0,1] \times \mathbb{R}^d} \frac{\mathrm{d}w}{|w-y|^{d+2s}}
$$
$$
\lesssim \|\gamma_{\cdot,x}(\mathcal{I}(\tau))\|_{(0,1] \times \mathbb{R}^d} \int_1^\infty \frac{r^{d-1}\,\mathrm{d}r}{r^{d+2s}}
$$
$$
\lesssim \|\gamma_{\cdot,x}(\mathcal{I}(\tau))\|_{(0,1] \times \mathbb{R}^d}.
$$

We conclude on (4.61) that

$$
\|(-\Delta)^s(\Pi_x \mathcal{I}(\tau)(1-\eta_{z,\lambda}))(y)\|_{B_r(x)} \lesssim \lambda^{|\tau|} [\Gamma \mathcal{I}(\tau); \mathbf{1}]_K + \|\gamma_{\cdot,x}(\mathcal{I}(\tau))\|_{(0,1] \times \mathbb{R}^d},
$$

and therefore by (4.60) we conclude on (4.59)

$$
\lambda^\gamma \|\mathscr{L}(U^{\mathrm{loc}}_{z,\lambda} - U)\|_{\gamma - 2s; B_\lambda(z)} \lesssim \lambda^{2s} \sum_{\tau \in \mathcal{V}_{0,\gamma}} \|\Upsilon.[\tau]\|_{B_\lambda(z)} \left( \lambda^{|\tau|} [\Gamma \mathcal{I}(\tau); \mathbf{1}]_K + \sup_{x \in P} \|\gamma_{\cdot,x}(\mathcal{I}(\tau))\|_{P \times \mathbb{R}^d} \right)
$$
$$
\leqslant \sum_{\tau \in \mathcal{V}_{0,\gamma}} \|\Upsilon.[\tau]\|_{B_\lambda(z)} \left( \lambda^{|\mathcal{I}(\tau)|} [\Gamma \mathcal{I}(\tau); \mathbf{1}]_K + \lambda^{2s} \sup_{x \in P} \|\gamma_{\cdot,x}(\mathcal{I}(\tau))\|_{P \times \mathbb{R}^d} \right)
$$
$$
\leqslant \sum_{\tau \in \mathcal{V}_{0,\gamma}} \lambda^{|\mathcal{I}(\tau)|} \|\Upsilon.[\tau]\|_{B_\lambda(z)} \left( [\Gamma \mathcal{I}(\tau); \mathbf{1}]_K + \sup_{x \in P} \|\gamma_{\cdot,x}(\mathcal{I}(\tau))\|_{P \times \mathbb{R}^d} \right)
$$
$$
\leqslant \sum_{\tau \in \mathcal{V}_{0,\gamma}} c\, \lambda^{|\mathcal{I}(\tau)|} \|\Upsilon.[\tau]\|_{B_\lambda(z)} \|v\|_{P_t}^{\alpha(\tau)},
$$

where we used Assumption 4.1 in the last bound. The same argument as Lemma 4.24 allows us to conclude the result. □

**Proof of Lemma 4.13.** For (4.21) consider $\delta > 0$ such that there is no tree in $\mathcal{V}_\gamma$ with homogeneity in the set $(1, 1+\delta)$ (can be done since by subcriticality $\mathcal{V}_\gamma$ is finite). On the other hand, we have that the family $\{B_r(z)\}_{z \in P_{t+4\lambda}}$ satisfies a $(\beta, r)$-spatial cone condition (see Definition 3.7) uniformly for $\beta = \frac{1}{\sqrt{2}}$, and then by Lemma 3.11 (using (4.7)) we conclude that

$$
\lambda \|v_\mathbf{X}\|_{B_r(z)} \lesssim \lambda^{1+\delta} [U_{1+\delta}]_{1+\delta; B_r(z)} + \|U_{1+\delta}\|_{B_r(z)}, \tag{4.62}
$$

with the proportionality constant depending only on $\beta$ and $\delta$. Moreover, since for all $x, y \in B_r(z)$

$$
U_{1+\delta}(x, y) = v(y) - \sum_{\tau \in \mathcal{V}_{0,1+\delta}} \frac{\Upsilon_x[\tau]}{\tau!} \Pi_x(\mathcal{I}(\tau))(y) = v(y) - \Pi_x(V_{1+\delta}(x))(y) = \langle \mathbf{1}, V_{1+\delta}(y) - \Gamma_{yx} V_{1+\delta}(x) \rangle,
$$

then $[U_{1+\delta}]_{1+\delta; B_r(z)} = [V_{1+\delta}; \mathbf{1}]_{B_r(z)}$, and by Lemma 4.7 we conclude that

$$
\lambda^{1+\delta} [U_{1+\delta}]_{1+\delta; B_r(z)} \lesssim \|v\|_{P_t} + c\, \lambda \|v_\mathbf{X}\|_{B_r(z)} + \lambda^\gamma [U]_{\gamma; B_r(z)}.
$$

On the other hand, $U_{1+\delta} = U_1$ by choice of $\delta$, and therefore by Lemma 4.12 we have that

$$
\|U_{1+\delta}\|_{B_r(z)} = \|U_1\|_{B_r(z)} \lesssim \|v\|_{P_t} + \mathbb{1}_{1>1} \lambda \|v_\mathbf{X}\|_{B_r(z)} = \|v\|_{P_t},
$$

which allows us to conclude on (4.62) that

$$
\lambda \|v_\mathbf{X}\|_{B_r(z)} \lesssim \|v\|_{P_t} + c\, \lambda \|v_\mathbf{X}\|_{B_r(z)} + \lambda^\gamma [U]_{\gamma; B_r(z)}. \tag{4.63}
$$

By imposing another smallness condition on $c \in (0,1)$ depending on the implicit proportionality constant in (4.63), which itself depends only on $d, s$ and $\gamma$, we can absorb the gradient term into the left-hand side and conclude (4.21).

For (4.20), recall that by (2.22) and (4.16) we have for all $x, y \in B \subset \mathbb{R}^{1+d}$

$$\langle \boldsymbol{X}_j, V(y) - \Gamma_{yx} V(y) \rangle = v_{\boldsymbol{X}_j}(y) - v_{\boldsymbol{X}_j}(x) - \sum_{\mu \in \mathcal{V}_{\beta,\gamma}} \frac{\Upsilon_x[\mu]}{\mu!} \gamma_{yx}(\mathcal{I}_j(\mu)) = (v_{\boldsymbol{X}_j} - \Lambda_j)(x, y),$$

and therefore $[V; \boldsymbol{X}]_B = [v_{\boldsymbol{X}} - \Lambda]_{\gamma-1;B}$. By (4.7) $v_{\boldsymbol{X}_j}$ plays the role of the generalised gradient $\nu$ for the germ $U$ and by Lemma 4.10 the germ $\Lambda$ is such that (3.6) holds, and by Lemma 3.9 we can conclude that

$$[V; \boldsymbol{X}]_{B_r(z)} \lesssim [U]_{\gamma\text{-3pt};B_r(z)} + [U]_{\gamma;B_r(z)}.$$

Combining this with Lemma 4.11 we have that

$$\lambda^\gamma [V; \boldsymbol{X}]_{B_r(z)} \lesssim \|v\|_{P_t} + c\,\lambda\,\|v_{\boldsymbol{X}}\|_{B_r(z)} + \lambda^\gamma [U]_{\gamma;B_r(z)} + c\,\lambda^\gamma [V; \boldsymbol{X}]_{B_r(z)}. \tag{4.64}$$

By imposing a smallness condition on $c \in (0, 1)$, which depends on the implicit proportionality constant which depends only on $s, d$ and $\gamma$, we can absorb the last term into the right hand side, and then we conclude combining and conclude that

$$\lambda^\gamma [V; \boldsymbol{X}]_{B_r(z)} \lesssim \|v\|_{P_t} + c\,\lambda\,\|v_{\boldsymbol{X}}\|_{B_r(z)} + \lambda^\gamma [U]_{\gamma;B_r(z)},$$

and with this (4.22) follows from (4.21). $\square$

The following modification of Lemma 3.6 was needed in the proof of Theorem 4.2.

**Lemma 4.25**. *Let $S$ and $D$ be non-negative functions defined on the convex subsets of a convex set $B \subset \mathbb{R}^{1+d}$ with $\mathrm{diam}(B) \leqslant 1$. Assume $D$ is monotone and that $S$ satisfies that for all half-parabolic balls $B_0, B_1, \ldots B_n \subset B$ such that $B_0 \subset \bigcup_{i=1}^n B_j$ one has the almost subadditivity property:*

$$S(B_0) \leqslant \sum_{i=1}^n S(B_j) + n\,D(B_0). \tag{4.65}$$

*Then for any given constants $\theta_0 \in (0, 1/2], \gamma > 0$, there exists $\varepsilon = \varepsilon(\theta_0, \gamma, s, d,) \in (0, 1)$ such that if for some $E \geqslant 0$ the following bounds*

$$\max\{\sigma^\gamma D(B_{\sigma_i}(y)), \sigma^\gamma S(B_{\theta_0 \sigma}(y))\} \leqslant \varepsilon\,\sigma^\gamma S(B_{\sigma_i}(y_i)) + E$$

*are satisfied for all balls in the family $\{B_\sigma(y) \subset B : \sigma \leqslant \bar{\sigma}, y \in B'\}$ for some $\bar{\sigma} \leqslant \mathrm{diam}(B)$ and $B' \subset B$. Then for each $\theta \in (0, 1)$ there exists a constant $C = C(d, \theta_0, \theta, \gamma, s) > 0$ (in particular independent of $\bar{\sigma}$ and $B'$) such that for all balls in the same family it holds*

$$\sigma^\gamma S(B_{\theta\sigma}(y_i)) \leqslant CE.$$

**Proof.** The proof is analogue to Lemma 3.6 by instead defining

$$Q := \sup_{\substack{\sigma \leqslant \bar{\sigma}, y \in B' \\ B_\sigma(y) \subset B}} \sigma^\gamma S(B_{\sigma/2}(y)) \leqslant \mathrm{diam}(B)^\gamma S(B) < +\infty,$$

and with the only difference being in the following bound

$$\sigma^\gamma S\big(B_{\frac{\sigma}{2}}(y)\big) \leqslant \sigma^\gamma \left( \sum_{i=1}^n S\big(B_{\frac{\theta_0 \sigma}{4}}(y_i)\big) + n\, D\big(B_{\frac{\sigma}{2}}(y)\big) \right)$$
$$\leqslant \sum_{i=1}^n \sigma^\gamma S\big(B_{\frac{\theta_0 \sigma}{4}}(y_i)\big) + n\, \sigma^\gamma D\big(B_{\frac{\sigma}{2}}(y)\big)$$
$$\lesssim \varepsilon Q + E,$$

which one uses to conclude analogously. $\square$

**Lemmas from Section 4.2**

**Proof of Lemma 4.15.** We show the result for $T = T_1$, the other cases being completely analogous. For any $t \in (T_1^{2s}, 1]$ we have that there exists $s \in [T_1, t^{1/2s}]$ and $\tau \in \mathcal{T}_{<2s} \setminus \mathcal{P}$ such that $c \, \|v\|_{P_s}^{\alpha(\tau)} \leqslant [\Pi; \tau]$. Since $r^{2s} < t$ we have that $(t, x_{1:d}) \in P_r$ for all $x_{1:d} \in (0, 1]^d$ and therefore

$$|v(t, x_{1:d})| \leqslant \|v\|_{P_r} \leqslant c^{-\alpha(\tau)^{-1}} [\Pi; \tau]^{\alpha(\tau)^{-1}} \lesssim c^{-\alpha(\tau)^{-1}} \max_{\tau \in (\mathcal{T}_{<2s} \setminus \mathcal{P}) \setminus \mathcal{W}} \{[\Pi; \tau]^{\alpha(\tau)^{-1}}\},$$

where we used that $\alpha$ is strictly positive precisely on the subcritical regime (see (4.3)). By taking supremum over $t \in (T_1^{2s}, 1]$ and $x_{1:d} \in (0, 1]^d$ we conclude

$$\|v\|_{P_{T_1}} \lesssim c^{-\alpha(\tau)^{-1}} \max_{\tau \in (\mathcal{T}_{<2s} \setminus \mathcal{P}) \setminus \mathcal{W}} \{[\Pi; \tau]^{\alpha(\tau)^{-1}}\}, \tag{4.66}$$

Analogously we can conclude that

$$\|v\|_{P_{T_2}} \lesssim c^{-\alpha(\tau)^{-1}} \max_{k \in \{0, e_1, \ldots, e_j\}} \max_{\tau \in \mathcal{V}_{0,2s}} \{[\Gamma \mathcal{I}(\tau); \boldsymbol{X}^k]^{\alpha(\tau)^{-1}}\}, \tag{4.67}$$

and

$$\|v\|_{P_{T_3}} \lesssim c^{-\alpha(\tau)^{-1}} \max_{\tau \in \mathcal{V}_{0,2s}} \|\gamma_{\cdot,\cdot}(\mathcal{I}(\tau))\|_{([0,1] \times \mathbb{R}^d) \times P}. \qquad \square$$

The proof of the next lemma follows the same ideas as [MW20, Lemma 2.7]. The differences are the use of the maximum principle for the fractional Laplacian $(-\Delta)^s$ for $s \in (0, 1)$ and the corresponding simplifications for periodic boundary conditions in space.

**Proof of Lemma 4.17.** By space periodicity, it is enough to check the result on $P = (0, 1] \times (0, 1]^d$. First we will deal with the case where $u$ is non-negative. Let $\eta \in C^2([0, \infty); \mathbb{R})$ be strictly positive on $(0, \infty)$ and zero at the boundary, i.e, $\eta(0) = 0$. We will later impose extra condition on $\eta$ that will allow us to bound the function $u$. We extend the function to $\eta \colon [0, \infty) \times \mathbb{R}^d \to \mathbb{R}$ as $\eta(x_0, x_{1:d}) := \eta(x_0)$. Let $z = \arg\max_{x \in [0,1] \times [0,1]^d} (u\,\eta)(x)$ be where the maximum of the continuous function $u\,\eta$ is attained on the compact $[0, 1] \times [0, 1]^d$. If $u(z) = 0$ or $z_0 = 0$, then we have that $\max_{[0,1] \times \mathbb{R}^d} u\,\eta = u(z)\,\eta(z) = 0$, since in the case $z_0 = 0$ we have $\eta(z) = \eta(z_0) = \eta(0) = 0$ by assumption on $\eta$. Since $\eta > 0$ on $(0, 1] \times [0, 1]^d$ we conclude that $u \leqslant 0$ on $(0, 1] \times [0, 1]^d$ and by continuity it extends to $[0, 1] \times [0, 1]^d$. Assume that $u(z) \neq 0$ and $z_0 \neq 0$, then we have that $z \in (0, 1] \times [0, 1]^d$. We claim that

$$(-\Delta)^s (u\,\eta)(z) \geqslant 0. \tag{4.68}$$

Since $u$ is smooth and periodic then we can use Bochner's representation (1.5) of the fractional Laplacian, which combined with the maximum principle for the heat semigroup $\{e^{t\Delta}\}_{t \geqslant 0}$ allows us to conclude that $(e^{t\Delta}(u\,\eta))(z) \leqslant (u\,\eta)(z)$ and therefore

$$-(-\Delta)^s(u\,\eta)(z) = \frac{1}{|\Gamma(-s)|} \int_0^\infty \{e^{t\Delta}(u\,\eta)(z) - (u\,\eta)(z)\} \frac{\mathrm{d}t}{t^{1+s}} \leqslant 0,$$

which shows (4.68). If $z \in (0,1) \times [0,1]^d$ then it is a maximum at an interior (in time) point and therefore $\partial_t(u\,\eta)(z_0) = 0$. If $z_0 \in \{1\} \times [0,1]^d$ then the maximum is obtained at the (time) boundary and therefore $\partial_t(u\,\eta)(z_0) \geq 0$ since the function must be increasing to its maximum from the left. In either case $\partial_t(u\,\eta)(z_0) \geq 0$ which combined with (4.68) allows us to conclude that

$$0 \leq \mathscr{L}(u\,\eta)(z) = (\eta\,\mathscr{L}(u))(z) + (u\,\partial_t\,\eta)(z) = -\eta\,(u^3 - g(u(\cdot),\cdot))(z) + (u\,\partial_t\,\eta)(z),$$

where we used that $\eta$ is constant in space and therefore $(-\Delta)^s(u\,\eta) = \eta\,(-\Delta)^s u$. Since $\eta > 0$ in $(0,\infty)$ and $z_0 \neq 0$, then

$$(u^3 - g(u(\cdot),\cdot))(z) \leq \left(u\,\frac{\partial_t\,\eta}{\eta}\right)(z). \tag{4.69}$$

If we assume the following condition on $\eta$:

$$\frac{\partial_t \eta}{\eta} \leq \frac{1}{2\,\eta^2}, \tag{4.70}$$

then we can conclude on (4.69) that

$$(u^3 - g)(z) \leq \left(\frac{u}{2\,\eta^2}\right)(z),$$

and since $u(z) \neq 0$, then

$$u^2(z) \leq \frac{1}{2\,\eta^2(z)} + \frac{g(z)}{u(z)} \leq 2\max\left\{\frac{1}{2\,\eta^2(z)}, \frac{\|g\|_{\mathbb{R}^{1+d}}}{u(z)}\right\}. \tag{4.71}$$

First suppose that the maximum in (4.71) is obtained by the first term, then by (4.71) we have that $u^2(z) \leq \frac{1}{\eta^2(z)}$, and therefore $(u\,\eta)(z) \leq 1$. In the other case, when case the maximum is obtained by the second term in (4.71) we have

$$\frac{1}{2\,\eta^2(z)} \leq \frac{\|g\|_{\mathbb{R}^{1+d}}}{u(z)} \Longrightarrow u\,\eta(z) \leq 2\,\eta^3(z)\,\|g\|_{\mathbb{R}^{1+d}}.$$

If we impose on $\eta$ the following extra condition:

$$\eta \leq \|g\|_{\mathbb{R}^{1+d}}^{-\frac{1}{3}}, \tag{4.72}$$

this we conclude that $u\,\eta(z) \leq 2$. Since $u\,\eta$ attains its maximum at $z$ on the compact $[0,1] \times [0,1]^d$ then we can conclude that

$$u(x) \leq \frac{2}{\eta(x)} \qquad \forall\,x \in [0,1] \times [0,1]^d. \tag{4.73}$$

$u \leq \frac{2}{\eta}$ on the set $(0,1] \times [0,1]^d$. We propose $\eta$ to be of the form

$$\eta_\lambda \colon (0,\infty) \to \mathbb{R} \qquad t \mapsto \frac{\lambda}{\lambda\,\|g\|_{L^\infty(\mathbb{R}\times\mathbb{R}^{1+d})}^{\frac{1}{3}} + t^{-\frac{1}{2}}},$$

where $\lambda > 0$ is a parameter to be chosen later. It is clear that in this case $\eta_\lambda$ is $C^2$, strictly positive and $\lim_{t \to 0^+} \eta_\lambda(t) = 0$ which allows it to be extended to $[0,\infty)$ with $\eta_\lambda(0) = 0$. Now we show that conditions (4.70) and (4.72) hold. Since

$$\eta_\lambda(z) = \frac{\lambda}{\lambda\,\|g\|_{\mathbb{R}^{1+d}}^{\frac{1}{3}} + z_0^{-\frac{1}{2}}} \leq \frac{\lambda}{\lambda\,\|g\|_{\mathbb{R}^{1+d}}^{\frac{1}{3}}} = \|g\|_{\mathbb{R}^{1+d}}^{-\frac{1}{3}},$$

which is precisely (4.72). To verify (4.70) we have that

$$\partial_t \eta_\lambda = \frac{\lambda^2}{2\lambda t^{\frac{3}{2}}(\lambda \|g\|_{\mathbb{R}^{1+d}}^{\frac{1}{3}} + t^{-\frac{1}{2}})^2} = \frac{1}{2\lambda t^{\frac{3}{2}}} \eta_\lambda^2,$$

and therefore

$$\eta_\lambda \partial_t \eta_\lambda = \frac{1}{2\lambda t^{\frac{3}{2}}} \eta_\lambda^3 = \frac{1}{2}\lambda^2 \left(\frac{1}{\lambda t^{-\frac{1}{2}}}\right)^3 \eta_\lambda^3 = \frac{1}{2}\lambda^2 \left(\frac{1}{\eta_\lambda} - \|g\|_{\mathbb{R}^{1+d}}^{\frac{1}{3}}\right)^3 \eta_\lambda^3 \leqslant \frac{1}{2}\lambda^2 \left(\frac{1}{\eta_\lambda}\right)^3 \eta_\lambda^3 \leqslant \frac{1}{2}\lambda^2,$$

and for any $\lambda \in (0,1]$ we conclude that $\eta_\lambda \partial_t \eta_\lambda \leqslant \frac{1}{2}$ which shows that (4.70) holds for any $\eta_\lambda$ and $\lambda \in (0,1]$. We can conclude now by (4.73) that for all $x \in [0,1] \times [0,1]^d$ we have the bound

$$u(x) \leqslant \frac{2}{\eta_\lambda(x)} = 2\left(\|g\|_{\mathbb{R}^{1+d}}^{-1} + \lambda^{-1} t^{-\frac{1}{2}}\right) \lesssim \max\left\{t^{-\frac{1}{2}}, \|g\|_{\mathbb{R}^{1+d}}^{\frac{1}{3}}\right\},$$

which shows the result for $u$ non-negative. The case where $u$ is non-positive follows by symmetry. For the general case we consider $z = \arg\max_{x \in [0,1] \times [0,1]^d}(|u|\eta_\lambda)(x)$. Without loss of generality suppose that $u(z) \geqslant 0$ since the other case follows by symmetry. By the same argument as before we conclude that $(|u|\eta_\lambda)(z) \leqslant 2$, and by definition of $z$ that $(|u|\eta_\lambda)(z) \leqslant 2$ for all $z \in (0,1] \times [0,1]^d$, and the desired bound follows from here. $\square$

**Proof of Lemma 4.18.** Applying Lemma 4.17 to the smooth and 1-periodic in space function $(v)_\lambda$ we conclude for every $R' \in (0,1)$, $R \in (R', 1)$ such that $t + R < 1$ and $x \in P_{t+R} \subset P_{t+R'}$:

$$|(v)_\lambda(x)| \lesssim \max_{\substack{\tau_1, \tau_2 \in \mathcal{W} \\ \tilde{\tau} \in \delta\mathcal{W}}} \left\{x_0^{-1/2}, \|(v^3)_\lambda - (v)_\lambda^3\|_{P_{t+R'}}^{\frac{1}{3}}, \|\mathcal{R}(\mathcal{I}(\tau_1) V^2)_\lambda\|_{P_{t+R'}}^{\frac{1}{3}}, \right.$$
$$\left. \|\mathcal{R}(\mathcal{I}(\tau_1)\mathcal{I}(\tau_2) V)_\lambda\|_{P_{t+R'}}^{\frac{1}{3}}, \|(\mathbf{\Pi}\tilde{\tau})_\lambda\|_{P_{t+R'}}^{\frac{1}{3}}\right\}.$$

Moreover, given $x \in P_{t+R} \subset P_{t+R'}$ we have that by definition of $P_{t+R}$

$$0 < (t+R')^{2s} < (t+R)^{2s} < x_0 \leqslant 1 \implies 0 < t+R' < t+R < (x_0)^{\frac{1}{2s}} \leqslant 1$$
$$\implies 0 < R - R' = t+R - (t+R') < (x_0)^{\frac{1}{2s}} \leqslant 1$$
$$\implies (R-R')^s \leqslant x_0^{1/2}$$
$$\implies x_0^{-1/2} \leqslant (R-R')^{-s},$$

which is enough to conclude the result. $\square$

**Proof of Lemma 4.19.** For all $t \in (0,1)$ and $\lambda \in (0,t)$ we have

$$\|v - (v)_\lambda\|_{P_t} \leqslant \sup_{z \in P_t} \int_{B_\lambda(x)} |(v(x) - v(y))|\, \psi_x^\lambda(y)\, \mathrm{d}y \leqslant \sup_{z \in P_t} [v]_{\beta; B_\lambda(z)} \lambda^\beta. \qquad (4.74)$$

By definition:

$$\lambda_t - \tilde{\lambda}_t = (k-1)\tilde{\lambda}_t = (1-k^{-1})\lambda_t \leqslant \lambda_t, \qquad (4.75)$$

and in particular $(\lambda_t - \tilde{\lambda}_t)/3 \leqslant \lambda_t$. By Theorem 4.2 we obtain:

$$\left(\frac{\lambda_t - \tilde{\lambda}_t}{3}\right)^\beta \sup_{z \in P_{t+\lambda_t - \tilde{\lambda}_t}} [V_\beta; \mathbf{1}]_{B_{\frac{\lambda_t - \tilde{\lambda}_t}{6}}(z)} \lesssim \|v\|_{P_t}.$$

Since $\tilde{\lambda}_t > 0$ and $R \in (\lambda_t, 1)$, we have that $\lambda_t - \tilde{\lambda}_t < R$ and $P_{t+R} \subset P_{t+\lambda_t - \tilde{\lambda}_t}$. On the other hand, by (4.75) we have that $\tilde{\lambda}_t = \frac{\lambda_t - \tilde{\lambda}_t}{k-1} \leqslant \frac{\lambda_t - \tilde{\lambda}_t}{6}$ holds as long as $k \geqslant 7$, and therefore $B_{\tilde{\lambda}_t}(z) \subset B_{\frac{\lambda_t - \tilde{\lambda}_t}{6}}(z)$. Putting this together, we can conclude the bound

$$\tilde{\lambda}^\beta \sup_{z \in P_{t+R}} [V_\beta; \mathbf{1}]_{B_{\tilde{\lambda}_t}(z)} \leqslant \tilde{\lambda}^\beta \sup_{z \in P_{t+\lambda_t - \tilde{\lambda}_t}} [V_\beta; \mathbf{1}]_{B_{\frac{\lambda_t - \tilde{\lambda}_t}{6}}(z)} \lesssim \tilde{\lambda}_t^\beta (\lambda_t - \tilde{\lambda}_t)^{-\beta} \|v\|_{P_t}. \tag{4.76}$$

Moreover, the choice of $\beta$ is such that $V_\beta(x) = v(x) \mathbf{1}$, and therefore $[v]_\beta = [V_\beta; \mathbf{1}]$. By (4.75) $\tilde{\lambda}_t^\beta (\lambda_t - \tilde{\lambda}_t)^{-\beta} = (k-1)^{-\beta}$, and combining (4.74) and (4.76) we obtain

$$\|v - (v)_{\tilde{\lambda}_t}\|_{P_{t+R}} \lesssim (k-1)^{-\beta} \|v\|_{P_t},$$

which proves (4.28).

To show (4.29) first we prove that for all $t \in (0,1)$ and $\lambda \in (0,t)$ we have

$$\|(v^3)_\lambda - (v)_\lambda^3\|_{P_t} \leqslant \lambda^\beta \|v\|_{P_{t-\lambda}}^2 \sup_{z \in P_t} [v]_{\beta; B_\lambda(z)}. \tag{4.77}$$

For $x \in \mathbb{R}^{1+d}$ we have that

$$(v^3)_\lambda(x) - (v)_\lambda^3(x) = (v^3)_\lambda(x) - v^3(x) + v^3(x) - (v)_\lambda^3(x), \tag{4.78}$$

by (4.74) we can bound the first difference as

$$|(v^3)_\lambda(x) - v^3(x)| \leqslant \lambda^\beta \sup_{z \in P_t} [v^3]_{\beta; B_\lambda(z)} \leqslant \lambda^\beta \sup_{z \in P_t} (\|v\|_{B_\lambda(z)}^2 [v]_{\beta; B_\lambda(z)}) \leqslant \lambda^\beta \|v\|_{P_{t-\lambda}}^2 \sup_{z \in P_t} [v]_{\beta; B_\lambda(z)}, \tag{4.79}$$

where we used that $B_\lambda(z) \subset ((t-\lambda)^{2s}, 1] \times \mathbb{R}^d$ by (4.21) and the 1-periodicity in space of $v$ to bound $\sup_{z \in P_t} \|v\|_{B_\lambda(z)}^2 \leqslant \|v\|_{P_{t-\lambda}}$. On the other hand, the second difference in (4.78) can be bounded as

$$\begin{aligned}|v^3(x) - (v)_\lambda^3(x)| &= \left|\int_0^1 3\left(\lambda (v)_\lambda(x) + (1-\lambda) v(x)\right)^2 ((v)_\lambda(x) - v(x)) \, d\lambda\right| \\ &\lesssim |(v)_\lambda(x) - v(x)| \int_0^1 (\lambda (v)_\lambda^2(x) + (1-\lambda) v^2(x)) \, d\lambda \\ &\leqslant \lambda^\beta \sup_{z \in P_t} [v]_{\beta; B_\lambda(z)} (\|(v)_\lambda^2\|_{P_t} + \|v^2\|_{P_t}) \\ &\leqslant \lambda^\beta \|v\|_{P_{t-\lambda}}^2 \sup_{z \in P_t} [v]_{\beta; B_\lambda(z)}, \end{aligned} \tag{4.80}$$

where we used that

$$|(v)_\lambda(x)| = \left|\int_{B_\lambda(x)} v(y) \, \psi_x^\lambda(y) \, dy\right| \leqslant \|v\|_{B_\lambda(x)} \left|\int_{B_\lambda(x)} \psi_x^\lambda(y) \, dy\right| \leqslant \|v\|_{P_{t-\lambda}}.$$

Combining (4.79) and (4.80) we obtain (4.77). On the other hand, observe that $B_{\tilde{\lambda}_t}(z) \subset B_{\lambda_t/6}(z)$ since $\tilde{\lambda}_t = k^{-1} \lambda_t \leqslant \lambda_t/6$ for $k \geqslant 7$. Proceeding as with (4.76) we have that

$$\tilde{\lambda}_t^\beta \|v\|_{P_{t+\lambda_t}}^2 \sup_{z \in P_{t+\lambda_t}} [v]_{\beta; B_{\tilde{\lambda}_t}(z)} \leqslant \tilde{\lambda}_t^\beta \|v\|_{P_t}^2 \sup_{z \in P_{t+\lambda_t}} [v]_{\beta; B_{\lambda_t/6}(z)} \lesssim \tilde{\lambda}_t^\beta \|v\|_{P_t}^2 \lambda_t^{-\beta} \|v\|_{P_t} = k^{-\beta} \|v\|_{P_t}^3,$$

which combined with (4.77) concludes the proof of (4.29). $\square$

**Proof of Lemma 4.20.** Using the analytic bounds of the model and the definition of $\tilde{\lambda}_t$, we have

$$\|(\mathbf{\Pi}\tau)_{\tilde{\lambda}_t}\|_{P_{t+\lambda_t}} = \sup_{x \in P_{t+\lambda_t}} |\langle \Pi_x \tau, \psi_x^{\tilde{\lambda}_t}\rangle| \leqslant [\Pi;\tau]_K \tilde{\lambda}_t^{|\tau|} \leqslant c \, \|v\|_{P_t}^{\alpha(\tau)} (k^{-1}\lambda_t)^{|\tau|} = c \, k^{-|\tau|} \|v\|_{P_t}^{\alpha(\tau)-s^{-1}|\tau|},$$

and (4.30) follows from the identity (4.56) for $\alpha(\tau)$.

For the other terms, we recall that by Lemma 2.14 we only look at tress with $\mathfrak{e}(\tau) = 0$. To show (4.31) consider $\tau \in \mathcal{W}$, then

$$\|(\mathcal{R}(\mathcal{I}(\tau)\,V^2))_{\tilde{\lambda}_t}\|_{P_{t+\lambda_t}} \leqslant \sup_{x \in P_{t+\lambda_t}} \left|\langle \mathcal{R}(\mathcal{I}(\tau)\,V_{\beta_\tau}^2) - \Pi_x(V_{\beta_\tau}^2(x)\,\mathcal{I}(\tau)), \psi_x^{\tilde{\lambda}_t}\rangle\right| \\ + \sup_{x \in P_{t+\lambda_t}} \left|\langle \Pi_x(V_{\beta_\tau}^2(x)\,\mathcal{I}(\tau)), \psi_x^{\tilde{\lambda}_t}\rangle\right|. \qquad (4.81)$$

By the Reconstruction Theorem in the form Lemma 4.6 we have

$$\left|\langle \mathcal{R}(\mathcal{I}(\tau)\,V^2) - \Pi_x(V_{\beta_\tau}^2(x)\,\mathcal{I}(\tau)), \psi_x^{\tilde{\lambda}_t}\rangle\right| \lesssim \tilde{\lambda}_t^{\beta_\tau + |\mathcal{I}(\tau)|} \sum_{\sigma_1,\sigma_2 \in \mathcal{T}} [\Pi; \sigma_1\sigma_2\mathcal{I}(\tau_1)]_{B_{2\tilde{\lambda}_t}(x)} [V_{\beta_\tau}^2; \sigma_1\sigma_2]_{B_{2\tilde{\lambda}_t}(x)}.$$

As in the proof of Lemma 4.19 we can use Theorem 4.2 to strengthen the results of Lemma 4.13 and Lemma 4.8 and conclude, since $\tilde{\lambda}_t \leqslant \lambda_t$, that

$$\tilde{\lambda}_t^{\beta_\tau + |\mathcal{I}(\tau)| + 2s} [\Pi; \sigma_1\sigma_2\mathcal{I}(\tau_1)]_{B_{2\tilde{\lambda}_t}(x)} [V_\beta^2; \sigma_1\sigma_2]_{B_{2\tilde{\lambda}_t}(x)} \lesssim c \, \|v\|_{P_t}.$$

By definition of $\tilde{\lambda}_t$ we have that

$$\tilde{\lambda}_t^{-2s} = (k^{-1}\lambda_t)^{-2s} = k^{2s}\left(\|v\|_{P_t}^{-s^{-1}}\right)^{2s} = k^{2s}\,\|v\|_{P_t}^2, \qquad (4.82)$$

which allows us to conclude that

$$\sup_{x \in P_{t+\lambda_t}} \left|\langle \mathcal{R}(\mathcal{I}(\tau)\,V_{\beta_\tau}^2) - \Pi_x(V_{\beta_\tau}^2(x)\,\mathcal{I}(\tau)), \psi_x^{\tilde{\lambda}_t}\rangle\right| \lesssim c\, k^{2s}\,\|v\|_{P_t}^3. \qquad (4.83)$$

To bound the term $\Pi_x(V_{\beta_\tau}^2(x)\,\mathcal{I}(\tau))$ in (4.81) we obtain by (2.25),(4.50) and (4.82):

$$\left|\langle \Pi_x(V_{\beta_\tau}^2(x)\,\mathcal{I}(\tau)), \psi_x^{\tilde{\lambda}_t}\rangle\right|$$
$$\lesssim \|v\|_{P_t}^2\, \tilde{\lambda}_t^{|\mathcal{I}(\tau)|}\, [\Pi;\mathcal{I}(\tau)]_{B_{2\tilde{\lambda}_t}(x)} + \sum_\sigma \|v\|\,\|\Upsilon.[\sigma]\|\, \tilde{\lambda}_t^{|\mathcal{I}(\sigma)\mathcal{I}(\tau)|}\, [\Pi;\mathcal{I}(\sigma)\,\mathcal{I}(\tau)]_{B_{2\tilde{\lambda}_t}(x)}$$
$$+ \sum_{\sigma_1,\sigma_2} \|\Upsilon.[\sigma_1]\|\,\|\Upsilon.[\sigma_2]\|\, \tilde{\lambda}_t^{|\mathcal{I}(\sigma_1)\mathcal{I}(\sigma_2)\mathcal{I}(\tau)|}\, [\Pi;\mathcal{I}(\sigma_1)\,\mathcal{I}(\sigma_2)\,\mathcal{I}(\tau)]_{B_{2\tilde{\lambda}_t}(x)}$$
$$\lesssim c\left(k^{-|\mathcal{I}(\tau)|} + \sum_\sigma k^{-|\mathcal{I}(\sigma)\mathcal{I}(\tau)|} + \sum_{\sigma_1,\sigma_2} k^{-|\mathcal{I}(\sigma_1)\mathcal{I}(\sigma_2)\mathcal{I}(\tau)|}\right)\|v\|_{P_t}^3, \qquad (4.84)$$

where we also used that since $\beta_\tau \in (0,1)$ then $\Upsilon[\sigma] \sim v^{\mathfrak{m}(\sigma)}$ for all the trees in (2.25). Combining (4.84) with (4.83) on (4.81) we conclude (4.31).

To show (4.32) consider $\tau_1, \tau_2 \in \mathcal{W}$, then

$$\|(\mathcal{R}(\mathcal{I}(\tau_1)\,\mathcal{I}(\tau_2)\,V))_{\tilde{\lambda}_t}\|_{P_{t+\lambda_t}} \leqslant \sup_{x \in P_{t+\lambda_t}} \left|\langle \mathcal{R}(\mathcal{I}(\tau_1)\,\mathcal{I}(\tau_2)\,V_{\beta_{\tau_1,\tau_2}}) - \Pi_x(\mathcal{I}(\tau_1)\,\mathcal{I}(\tau_2)\,V_{\beta_{\tau_1,\tau_2}}), \psi_x^{\tilde{\lambda}_t}\rangle\right| \\ + \sup_{x \in P_{t+\lambda_t}} \left|\langle \Pi_x(\mathcal{I}(\tau_1)\,\mathcal{I}(\tau_2)\,V_{\beta_{\tau_1,\tau_2}}), \psi_x^{\tilde{\lambda}_t}\rangle\right|. \qquad (4.85)$$

For simplicity denote $\beta = \beta_{\tau_1, \tau_2}$. By the Reconstruction Theorem in the form of Lemma 4.6 we have

$$\left|\langle \mathcal{R}(\mathcal{I}(\tau_1)\mathcal{I}(\tau_2)V_\beta) - \Pi_x(\mathcal{I}(\tau_1)\mathcal{I}(\tau_2)V_\beta), \psi_x^{\tilde{\lambda}_t}\rangle\right| \lesssim \tilde{\lambda}_t^{\beta + |\mathcal{I}(\tau_1)\mathcal{I}(\tau_2)|} \sum_\sigma [\Pi; \sigma\mathcal{I}(\tau_1)\mathcal{I}(\tau_2)]_{B_{2\tilde{\lambda}_t}(x)} [V_\beta; \sigma]_{B_{2\tilde{\lambda}_t}(x)}.$$

As in the previous case we can use Theorem 4.2 to strengthen the results of Lemma 4.13 and Lemma 4.8 and conclude, since $\tilde{\lambda}_t \leqslant \lambda_t$, that

$$\tilde{\lambda}_t^{\beta + |\mathcal{I}(\tau_1)\mathcal{I}(\tau_2)| + 2s} [\Pi; \sigma\mathcal{I}(\tau_1)\mathcal{I}(\tau_2)]_{B_{2\tilde{\lambda}_t}(x)} [V_\beta; \sigma]_{B_{2\tilde{\lambda}_t}(x)} \lesssim c \, \|v\|_{P_t},$$

which combined with (4.82) allows us to conclude that

$$\sup_{x \in P_{t+\lambda_t}} \left|\langle \mathcal{R}(\mathcal{I}(\tau_1)\mathcal{I}(\tau_2)V_\beta) - \Pi_x(\mathcal{I}(\tau_1)\mathcal{I}(\tau_2)V_\beta), \psi_x^{\tilde{\lambda}_t}\rangle\right| \lesssim c \, k^{2s} \|v\|_{P_t}^3. \tag{4.86}$$

To bound the term $\Pi_x(\mathcal{I}(\tau_1)\mathcal{I}(\tau_2)V_\beta)$ in (4.85) we obtain by (2.20)

$$\begin{aligned}
&\left|\langle \Pi_x(\mathcal{I}(\tau_1)\mathcal{I}(\tau_2)V_\beta), \psi_x^{\tilde{\lambda}_t}\rangle\right| \\
&\lesssim \|v\|_{P_t} \tilde{\lambda}_t^{|\mathcal{I}(\tau_1)\mathcal{I}(\tau_2)|} [\Pi; \mathcal{I}(\tau_1)\mathcal{I}(\tau_2)]_{B_{2\tilde{\lambda}_t}(x)} + \|v_{\boldsymbol{X}}\|_{B_{2\tilde{\lambda}_t}(x)} \tilde{\lambda}_t^{|\mathcal{I}(\tau_1)\mathcal{I}(\tau_2)\boldsymbol{X}|} [\Pi; \mathcal{I}(\tau_1)\mathcal{I}(\tau_2)\boldsymbol{X}]_{B_{2\tilde{\lambda}_t}(x)} \\
&\quad + \sum_\sigma \|\Upsilon.[\sigma]\| \tilde{\lambda}_t^{|\mathcal{I}(\sigma)\mathcal{I}(\tau_1)\mathcal{I}(\tau_2)|} [\Pi; \mathcal{I}(\tau_1)\mathcal{I}(\tau_2)\mathcal{I}(\sigma)]_{B_{2\tilde{\lambda}_t}(x)}.
\end{aligned} \tag{4.87}$$

For the first term and the terms $\sigma$ in the sum of (4.87) such that $\mathfrak{n}(\sigma) = 0$ we use (4.51) and (4.54) respectively combined with (4.82) to conclude

$$\|v\|_{P_t} \tilde{\lambda}_t^{|\mathcal{I}(\tau_1)\mathcal{I}(\tau_2)|} [\Pi; \mathcal{I}(\tau_1)\mathcal{I}(\tau_2)]_{B_{2\tilde{\lambda}_t}(x)} \lesssim c \, k^{-|\mathcal{I}(\tau_1)\mathcal{I}(\tau_2)|} \|v\|_{P_t}^3 \tag{4.88}$$

and

$$\|v\|_{P_t} \tilde{\lambda}_t^{|\mathcal{I}(\sigma)\mathcal{I}(\tau_1)\mathcal{I}(\tau_2)|} [\Pi; \mathcal{I}(\sigma)\mathcal{I}(\tau_1)\mathcal{I}(\tau_2)]_{B_{2\tilde{\lambda}_t}(x)} \lesssim c \, k^{-|\mathcal{I}(\sigma)\mathcal{I}(\tau_1)\mathcal{I}(\tau_2)|} \|v\|_{P_t}^3. \tag{4.89}$$

For the term with $v_{\boldsymbol{X}}$ in (4.87) we use that by Lemma 4.13 improved by Theorem 4.2

$$\|v_{\boldsymbol{X}}\|_{B_{2\tilde{\lambda}_t}(x)} \lesssim \lambda_t^{-1} \|v\|_{P_t} = \|v\|_{P_t}^{s^{-1}+1}, \tag{4.90}$$

and therefore

$$\begin{aligned}
&\|v_{\boldsymbol{X}}\|_{B_{2\tilde{\lambda}_t}(x)} \tilde{\lambda}_t^{|\mathcal{I}(\tau_1)\mathcal{I}(\tau_2)\boldsymbol{X}|} [\Pi; \mathcal{I}(\tau_1)\mathcal{I}(\tau_2)\boldsymbol{X}]_{B_{2\tilde{\lambda}_t}(x)} \\
&\lesssim c \, k^{-|\boldsymbol{X}\mathcal{I}(\tau_1)\mathcal{I}(\tau_2)|} \|v\|_{P_t}^{1+s^{-1}-s^{-1}|\boldsymbol{X}\mathcal{I}(\tau_1)\mathcal{I}(\tau_2)| + \alpha(\boldsymbol{X}\mathcal{I}(\tau_1)\mathcal{I}(\tau_2))} \\
&\leqslant c \, k^{-|\boldsymbol{X}\mathcal{I}(\tau_1)\mathcal{I}(\tau_2)|} \|v\|_{P_t}^3
\end{aligned} \tag{4.91}$$

by (4.52). For the terms in the sum of (4.87) such that $\mathfrak{n}(\sigma) = 1$, since $\Upsilon[\sigma] \sim v_{\boldsymbol{X}}$, we can conclude analogously using (4.55) that

$$\|\Upsilon.[\sigma]\|_{B_{2\tilde{\lambda}_t}(x)} \tilde{\lambda}_t^{|\mathcal{I}(\tau_1)\mathcal{I}(\tau_2)\boldsymbol{X}|} [\Pi; \mathcal{I}(\tau_1)\mathcal{I}(\tau_2)\boldsymbol{X}]_{B_{2\tilde{\lambda}_t}(x)} \lesssim c \, k^{-|\mathcal{I}(\sigma)\mathcal{I}(\tau_1)\mathcal{I}(\tau_2)|} \|v\|_{P_t}^3. \tag{4.92}$$

Putting (4.88),(4.89),(4.91) and (4.92) together we conclude on (4.87) the bound

$$\begin{aligned}
&\left|\langle \Pi_x(\mathcal{I}(\tau_1)\mathcal{I}(\tau_2)V_\beta), \psi_x^{\tilde{\lambda}_t}\rangle\right| \\
&\lesssim c \left( k^{-|\mathcal{I}(\tau_1)\mathcal{I}(\tau_2)|} + k^{-|\boldsymbol{X}\mathcal{I}(\tau_1)\mathcal{I}(\tau_2)|} + \sum_\sigma k^{-|\mathcal{I}(\sigma)\mathcal{I}(\tau_1)\mathcal{I}(\tau_2)|} + \right) \|v\|_{P_t}^3.
\end{aligned} \tag{4.93}$$

Combining (4.93) with (4.86) on (4.85) we conclude (4.32). □

## Appendix A.  Regularity Structures

We recall some notions of the theory of regularity structures. We follow mostly [FH20, Chapter 13], but our definitions have some minor modifications to the ones in there.

**Definition A.1**. ([FH20, Definition 13.1]) *A regularity structure $(T, G)$ consists of*

- *A structure space given as a graded vector space $T = \bigoplus_{\alpha \in A} T_\alpha$ where each $T_\alpha$ is a Banach space, with index set $A \subset \mathbb{R}$ bounded from below and locally finite. Elements of $T_\alpha$ are said to have homogeneity $\alpha$, and we write $|\tau| = \alpha$ for $\tau \in T_\alpha$.*

- *A structure group $G$ of continuous linear operators acting on $T$ such that for every $\Gamma \in G$, every $\alpha \in A$, and every $\tau \in T_\alpha$, one has*

$$\Gamma \tau - \tau \in T_{<\alpha} := \bigoplus_{\beta < \alpha} T_\beta.$$

*A sector $S$ of $\mathcal{T}$ is a linear subspace $S = \bigoplus_{\alpha \in A} S_\alpha \subset T$, with closed linear subspaces $S_\alpha \subset T_\alpha$, which is invariant under $G$, such that $(S, G|_S)$ is a regularity structure on its own.*

For our purposes we can assume that $T$ is a finite dimensional Hilbert space with inner product $\langle \cdot, \cdot \rangle_T$, and that the subspaces $\{T_\alpha\}_{\alpha \in A}$ are mutually orthogonal subspaces. Moreover, we let $\mathcal{T}$ denote an orthogonal basis such that $\{\tau \in \mathcal{T} : |\tau| = \alpha\}$ is an orthogonal basis of $T_\alpha$ for each $\alpha \in A$.

**Definition A.2**. ([FH20, Definition 13.5]) *A model on the regularity structure $(T, G)$ consists of a pair $(\Pi, \Gamma)$ of maps*

$$\Pi : \mathbb{R}^{1+d} \to \mathcal{L}(T, \mathscr{D}'(\mathbb{R}^{1+d})) \qquad \Gamma : \mathbb{R}^{1+d} \times \mathbb{R}^{1+d} \to G$$
$$x \mapsto \Pi_x, \qquad\qquad\qquad (x, y) \mapsto \Gamma_{xy},$$

*which satisfy the non-linear relationships for all $x, y, z \in \mathbb{R}^{1+d}$:*

$$\Gamma_{xy} \Gamma_{yz} = \Gamma_{xz} \qquad \text{and} \qquad \Pi_x \Gamma_{xy} = \Pi_y.$$

*Moreover, $\Pi$ has to satisfy for every compact set $D \subset \mathbb{R}^{1+d}$ and $\tau \in \mathcal{T}$, the analytic bound*

$$[\Pi; \tau]_D := \sup_{x \in D} \sup_{\psi \in \mathcal{B}_r} \sup_{\substack{\lambda > 0 \\ B_\lambda(x) \subset D}} |\langle \Pi_x \tau, \psi_x^\lambda \rangle| \lambda^{-|\tau|} < +\infty, \tag{A.1}$$

*where $r \in \mathbb{Z}^+$ is the smallest integer such that $r > |\min A| \geqslant 0$, $\mathcal{B}_r$ is defined (as in (3.5))*

$$\mathcal{B}_r := \{\psi \in \mathscr{D}(B_1) : \operatorname{supp}(\psi) \subset B_1, \|\psi\|_{C^r} \leqslant 1\},$$

*$\langle \cdot, \cdot \rangle$ is the dual pairing between distributions in $\mathscr{D}'(\mathbb{R}^{1+d})$ and test functions $\mathscr{D}(\mathbb{R}^{1+d})$. Similarly, $\Gamma$ has to satisfy for all $\sigma \in \mathcal{T}$ with $|\sigma| < |\tau|$ the analytic bound*

$$[\Gamma \tau; \sigma]_D := \sup_{x \in D} \sup_{\substack{y \in D \\ y_0 \leqslant x_0}} \frac{|\langle \sigma, \Gamma_{xy} \tau \rangle_T|}{d(x, y)^{|\tau| - |\sigma|}} < +\infty. \tag{A.2}$$

**Remark A.3.** In our definition (A.2) we have the restriction that the second variable $y$ has to be in the past of $x$, i.e., $y_0 \leqslant x_0$. A similar restriction is implicit in (A.1) since by definition for $\psi \in \mathcal{B}_r$ we have that the support of $\psi$ is contained in the half-parabolic ball $B_1$ as defined in Section 1.1, and therefore $\psi_x$ has support in the past of $x$.

**Definition A.4.** ([FH20, Definition 13.10]) *Let $(\Pi, \Gamma)$ be a model on the regularity structure $(T, G)$. For $\gamma \in \mathbb{R}$ we define the set of modelled distributions $\mathcal{D}^\gamma(\Gamma)$ as the set of functions $F$: $\mathbb{R}^{1+d} \to T_{<\gamma}$ such that for every $\tau \in \mathcal{T}$ and bounded set $D \subset \mathbb{R}^{1+d}$:*

$$[F; \tau]_{\gamma; D} := \sup_{x \in D} \sup_{\substack{y \in D \\ y_0 \leqslant x_0}} \frac{|\langle \tau, F(x) - \Gamma_{xy} F(y) \rangle_T|}{d(x, y)^{\gamma - |\tau|}} < +\infty. \tag{A.3}$$

*Given a sector $S$ of $(T, G)$ we denote by $\mathcal{D}^\gamma(S, \Gamma)$ the set of modelled distributions taking values in $S$. When the context is clear we will only write $\mathcal{D}^\gamma$.*

One of the main results of the theory of regularity structure is the reconstruction theorem which we state below.

**Theorem A.5.** ([FH20, Theorem 13.26]) *Let $(\Pi, \Gamma)$ be a model for a regularity structure $(T, G)$ on $\mathbb{R}^{1+d}$. Then for $\gamma > 0$ there exists a unique linear map $\mathcal{R} : \mathcal{D}^\gamma(\Gamma) \to \mathscr{D}'(\mathbb{R}^{1+d})$ such that*

$$|\langle \mathcal{R} F - \Pi_x F(x), \psi_x^r \rangle| \lesssim r^\gamma \sum_{\tau \in \mathcal{T}} [\Pi; \tau]_{B_{2r}(x)} [F; \tau]_{\gamma; B_{2r}(x)}, \tag{A.4}$$

*uniformly over $\psi \in \mathcal{B}_r$, $r \in (0, 1]$. The implicit proportionality constant depends only on $d$ and $\gamma$.*

**Remark A.6.** The more explicit constants appearing in the right hand side of (A.4) can be obtained by following the proof in [FH20, Theorem 13.26]. The fact that a bigger ball appears on the right hand side follows from the auxiliary result [FH20, Theorem 13.24].

## A.1. Periodicity

Following [Hai14, Section 3.6] we consider $\mathbb{Z}^d$ acting on $\mathbb{R}^{1+d}$ by spatial translations, i.e., for each $x \in \mathbb{R}^{1+d}$ and $k \in \mathbb{Z}^d$ we have the action $\tau_k x = (x_0, x_{1:d} + k)$.

**Definition A.7.** ([Hai14, Definition 3.33]) *We say that a model $(\Pi, \Gamma)$ is 1-periodic if it is adapted to the action of $\mathbb{Z}^d$, i.e.:*

  i. *For every $\psi \in \mathscr{D}'(\mathbb{R}^{1+d})$, $x \in \mathbb{R}^{1+d}$, $\tau \in \mathcal{T}$ and $k \in \mathbb{Z}^d$ one has $\langle \Pi_{x+k} \tau, \psi_k \rangle = \langle \Pi_x \tau, \psi \rangle$.*

  ii. *For every $x, y \in \mathbb{R}^{1+d}$ and $k \in \mathbb{Z}^d$ one has $\Gamma_{x+k, y+k} = \Gamma_{x,y}$.*

*A modelled distribution $F \in \mathcal{D}^\gamma$ is called 1-periodic if $F(x) = F(x + k)$ for all $x \in \mathbb{R}^{1+d}$ and $k \in \mathbb{Z}^d$.*

**Remark A.8.** A model being 1-periodic does not imply that $\Pi_x \tau$ are 1-periodic distributions. This is easily seen by considering the polynomial model for which we have that $\Pi_x \boldsymbol{X}^k = (\cdot - x)^k$ which is not periodic for $k \neq 0$.

With this definitions the reconstruction operator preserves periodicity.

**Lemma A.9.** ([Hai14, Proposition 3.38] ) *If $(\Pi, \Gamma)$ is a 1-periodic model and $F \in \mathcal{D}^\gamma$ is a 1-periodic modelled distribution for some $\gamma > 0$, then the reconstruction of $F$ is a 1-periodic distribution, i.e., for all $k \in \mathbb{Z}^d$ we have $\langle \mathcal{R} F, \varphi_k \rangle = \langle \mathcal{R} F, \varphi \rangle$.*

## A.2. Decorated trees

The basis $\mathcal{T}$ of our regularity structure will consist of decorated trees which we introduce next. Our definitions follow [BB21, Section 2], the main difference being that in there the authors allow noises on all internal nodes while in our case noises are only allowed on leave nodes. The difference comes from the noise in our PDE (1.1) acting additively, while in their general framework they consider SPDEs with the noise acting affinely. Moreover, they consider noises as decorations on edges, while for us noises are decorations on leaves.

A decorated tree consist of a tuple $(\tau, \rho_\tau, \mathfrak{l}, \mathfrak{n}, \mathfrak{e})$ where $\tau$ is a non-planar rooted tree with node set $N_\tau$, root $\rho_\tau \in N_\tau$, edge set $E_\tau$, node decorations $\mathfrak{l} \colon N_\tau \to \{0,1\}$ and $\mathfrak{n} \colon N_\tau \to \mathbb{N}^{d+1}$ and an an edge decoration $\mathfrak{e} \colon E_\tau \to \mathbb{N}^{d+1}$. The decoration $\mathfrak{l}$ encodes noises, whilst $\mathfrak{n}$ will encode multiplication by polynomials and $\mathfrak{e}$ derivatives of kernels. We denote by $\Xi$ the tree consisting of a single node, its root $\rho_\Xi$, with decorations $\mathfrak{l}(\rho_\Xi) = 1$ and $\mathfrak{n}(\rho_\Xi) = 0$. The tree $\Xi$ is an abstract representation of the noise $\xi$ appearing in equation (1.1). Let $L_\tau \subset N_\tau \setminus \{\rho_\tau\}$ be the set of leaves of $\tau$, i.e., the non-root nodes which belong exactly to one edge. Since the noise we consider in (1.1) is additive, we only need to consider trees $\tau$ such that $\mathfrak{l}(v) = 1$ if and only if $v \in L_\tau$. On the other hand, we denote by $\boldsymbol{X}^k$ the tree consisting of a single node, with decorations $\mathfrak{l}(\rho_{\boldsymbol{X}^k}) = 0$, $\mathfrak{n}(\rho_{\boldsymbol{X}^k}) = k \in \mathbb{N}^{1+d}$. Given a decorated tree $\tau$ and $m \in \mathbb{N}^{1+d}$ we denote by $\mathcal{I}_m(\tau)$ the decorated tree obtained by rooting the decorated tree $\tau$ into a new root with an edge decorated by $m \in \mathbb{N}^{d+1}$, and we set the $\mathfrak{n}$ and $\mathfrak{l}$ decoration of this new root to 0. For simplicity we denote $\mathcal{I} := \mathcal{I}_0$.

Given two decorated trees $\tau_1$ and $\tau_2$ different from $\Xi$, we define their product $\tau_1 \tau_2$ as the decorated tree obtained by identifying their roots into a single root, adding the decorations $\mathfrak{n}$ at the root, i.e., $\mathfrak{n}(\rho_{\tau_1 \tau_2}) = \mathfrak{n}(\rho_{\tau_1}) + \mathfrak{n}(\rho_{\tau_2})$ and preserving all the other decorations. This product is commutative, and the tree $\mathbf{1} := \boldsymbol{X}^0$ is a unit for this product.

Consider $\mathcal{T}$ the set of trees such that $\Xi \in \mathcal{T}$ and which is closed under recursive applications of planting a tree, finite product of planted trees and product with abstract monomials $\boldsymbol{X}^k$. Then we have that given $\tau \in \mathcal{T} \setminus \{\Xi\}$ there exists a unique finite (possibly empty) collection $\{(\tau_i, m_i, \beta_i)\}_{i \in I} \subset \mathcal{T} \times \mathbb{N}^{d+1} \times \mathbb{N}$ and a unique $k \in \mathbb{N}^{d+1}$ such that

$$\tau = \boldsymbol{X}^k \prod_{i \in I} \mathcal{I}_{m_i}(\tau_i)^{\beta_i}. \tag{A.5}$$

Trees of the form $\mathcal{I}_m(\tau)$ are called planted trees. By definition noise decorations, encoded by $\mathfrak{l}$, are only allowed at leaves for trees in $\mathcal{T} \setminus \{\Xi\}$, since otherwise we would have needed to allow a factor of the form $\Xi^\ell$ in (A.5). Moreover, from (A.5) we see that $\mathcal{T} \setminus \{\Xi\}$ is closed under products.

We define the symmetry factor $\tau!$ recursively by setting $\Xi! := 1$ and extending it to a tree of the form (A.5) as

$$\tau! := k! \prod_{i \in I} (\tau_i!)^{\beta_i} \beta_i!, \tag{A.6}$$

where for the multi-index $k \in \mathbb{N}^{1+d}$ its factorial is given by $k! = k_0! \cdots k_d!$. If we are given a tree $\tau$ not in its unique factorisation (A.5), but as a product

$$\tau = \boldsymbol{X}^k \prod_{i=1}^n \mathcal{I}_{m_i}(\tau_i), \tag{A.7}$$

where repetition between the trees $\mathcal{I}_{m_i}(\tau_i)$ is allowed, then we can write the symmetry factor of $\tau$ as

$$\tau! = k! \frac{n!}{\delta(\mathcal{I}_{m_1}(\tau_1),\ldots,\mathcal{I}_{m_n}(\tau_n))} \prod_{i=1}^{n} \tau_i!, \tag{A.8}$$

where $\delta(\mathcal{I}_{m_1}(\tau_1),\ldots,\mathcal{I}_{m_n}(\tau_n))$ counts the number of different ordered $n$-tuples $(\mathcal{I}_{m_1}(\tau_1),\ldots,\mathcal{I}_{m_n}(\tau_n))$ which correspond to the same unordered collection $\{\mathcal{I}_{m_1}(\tau_1),\ldots,\mathcal{I}_{m_n}(\tau_n)\}$, and therefore the factor $n!/\delta(\mathcal{I}_{m_1}(\tau_1),\ldots,\mathcal{I}_{m_n}(\tau_n))$ counts the order of the subgroup of permutations of $n$ elements that preserve the ordered collection $(\mathcal{I}_{m_1}(\tau_1),\ldots,\mathcal{I}_{m_n}(\tau_n))$.

To define the homogeneity of a decorated tree we consider the scaling associated to the operator $\mathscr{L} = (\partial_t + (-\Delta)^{2s})$, which is given by the scaling vector $\mathfrak{s} := (2s, 1, \ldots, 1) \in \mathbb{R}^{d+1}$. We define the homogeneity of the noise as

$$|\Xi|_{\mathfrak{s}} := -\frac{3+2s}{2} - \kappa, \tag{A.9}$$

for $0 < \kappa \ll 1$ fixed. This value of $|\Xi|_{\mathfrak{s}}$ corresponds to the (negative) Hölder regularity of the space-time white noise when measured with the fractional parabolic scaling $\mathfrak{s}$ for $d=3$ (see [Hai14, Lemma 10.2]). The homogeneity of a monomial is defined as:

$$|\boldsymbol{X}^k|_{\mathfrak{s}} := |k|_{\mathfrak{s}} := 2s\, k_0 + \sum_{i=1}^{d} k_i. \tag{A.10}$$

Definition (A.10) takes into account that under the scaling $\mathfrak{s}$ a time coordinate counts $2s$ times a spatial one. By Schauder theory (e.g. Theorem 3.1) $\mathscr{L}^{-1}$ improves regularity by $2s$ whilst taking a derivative $\partial^k$, $k \in \mathbb{N}^{d+1}$, decreases it by $|k|_{\mathfrak{s}}$. This is reflected in the definition:

$$|\mathcal{I}_m(\tau)|_{\mathfrak{s}} := |\tau|_{\mathfrak{s}} + 2s - |m|_{\mathfrak{s}}. \tag{A.11}$$

At last we extend the homogeneity to an arbitrary decorated tree of the form (A.5) recursively as:

$$|\tau|_{\mathfrak{s}} := |\boldsymbol{X}^k|_{\mathfrak{s}} + \sum_{i \in I} \beta_i\, |\mathcal{I}_{m_i}(\tau_i)|_{\mathfrak{s}}. \tag{A.12}$$

We define the total polynomial and edge decoration as the homogeneity of the correspondent sum

$$|\mathfrak{n}(\tau)|_{\mathfrak{s}} := \sum_{v \in N_\tau} |\mathfrak{n}(v)|_{\mathfrak{s}}, \qquad |\mathfrak{e}(\tau)|_{\mathfrak{s}} := \sum_{v \in E_\tau} |\mathfrak{e}(v)|_{\mathfrak{s}}. \tag{A.13}$$

### A.3. Structure group

Assume we are given a subset of $\mathcal{T} \subset \mathscr{T}$ of decorated trees which contains the monomials $\{\boldsymbol{X}^k\}_{k \in \mathbb{N}^{1+d}}$ and such that $A := \{|\tau|_{\mathfrak{s}} : \tau \in \mathcal{T}\}$ is bounded from below, locally finite and $\mathbf{1} := \boldsymbol{X}^0$ is the only tree such that $|\mathbf{1}|_{\mathfrak{s}} = 0$. Let $T = \langle \mathcal{T} \rangle$ be the $\mathbb{R}$-vector space spanned by $\mathcal{T}$ with grading induced by $A$. Following [FH20, Section 15.3] we describe the construction of the structure group associated to $T$.

Consider the set $\mathcal{T}^+$ of elements of the form (A.5) such that $I$ is a finite set and $\{\tau_i\}_{i \in I} \subset \mathcal{T}$, $\{\rho_i\}_{i \in I} \subset \mathbb{N}^{1+d}$ are such that $|\mathcal{I}_{m_i}(\tau_i)| > 0$. We define $T^+$ as the vector space spanned by $\mathcal{T}^+$ which has a natural commutative algebra structure with unit $\mathbf{1}$. This algebra is freely generated by the trees $\{\boldsymbol{X}^{e_j}\}_{j=0}^{d} \cup \{\mathcal{I}_m(\tau)\}_{\tau \in \mathcal{T}, m \in \mathbb{N}^{1+d}, |\mathcal{I}_m(\tau)| > 0}$. Define a linear map $\Delta: T \to T \otimes T^+$, called *coaction*, by setting $\Delta \Xi := \Xi \otimes \mathbf{1}, \Delta\, \boldsymbol{X}_j := \mathbf{1} \otimes \boldsymbol{X}_j + \boldsymbol{X}_j \otimes \mathbf{1}$, and extend it inductively to planted trees as

$$\Delta\, \mathcal{I}_m(\tau) = (\mathcal{I}_m \otimes \mathrm{Id}_{T^+})\, \Delta \tau + \sum_{k \in \mathbb{N}^{1+d}} \frac{\boldsymbol{X}^k}{k!} \otimes \mathcal{I}_{m+k}(\tau), \tag{A.14}$$

with the convention that $\mathcal{I}_{m+k}(\tau)=0\in\mathcal{T}^+$ if $|\mathcal{I}_{m+k}(\tau)|<0$. Finally, we extend $\Delta$ it to the full basis $\mathcal{T}$ by setting it to act multiplicatively on symbols of the form (A.5). By multiplicativity this definition implies that

$$\Delta \boldsymbol{X}^k = (\Delta \boldsymbol{X})^k = \sum_{m\in\mathbb{N}^{1+d}} \binom{k}{m} \boldsymbol{X}^m \otimes \boldsymbol{X}^{k-m}, \qquad k\in\mathbb{N}^{d+1}. \tag{A.15}$$

We define an algebra homomorphism $\Delta\colon T^+\to T^+\otimes T^+$ which acts on the polynomial symbols as (A.15) and on planted trees $\mathcal{I}_m(\tau)\in\mathcal{T}^+$ as (A.14), where $\Delta\tau$ refers to the *coaction* previously defined for $\tau\in\mathcal{T}$, and with the same convention of setting $\mathcal{I}_{m+k}(\tau)=0\in\mathcal{T}^+$ if $|\mathcal{I}_{m+k}(\tau)|<0$. Since $T^+$ is freely generated by these symbols the action of the algebra homomorphism is determined by these definitions, which extend multiplicatively to symbols of the form (A.5). Observe that for polynomials and planted trees which can be interpreted as both elements of $T$ and of $T^+$ the definitions of the coaction $\Delta$ and the algebra homomorphism $\Delta$ coincide with the suitable identifications.

It is well-known that $T^+$ has a Hopf algebra structure with coproduct given by $\Delta$, and $T$ is a right comodule over $T^+$ with coaction $\Delta$ (see [FH20, Section 15.3 p.300]). We consider $\mathrm{Char}(T^+)$, the set of characters on $T^+$ whose elements are the algebra homomorphisms $\gamma\colon T^+\to\mathbb{R}$. We consider the convolution product in $\mathrm{Char}(T^+)$ defined by $(\gamma_1*\gamma_2)(\tau)=(\gamma_1\otimes\gamma_2)\Delta\tau$, and which turns $(\mathrm{Char}(T^+),*)$ into a group. The inverse elements in this group are given by $\gamma^{-1}=\gamma\circ\mathcal{A}$ where $\mathcal{A}\in\mathrm{End}(T^+)$ is the *antipode* map, which is part of the Hopf algebra structure of $T^+$. This antipode will not play a role in our work, so we omit details on this, but we refer to [FH20, Section 15.3] for details. Given a character $\gamma\in\mathrm{Char}(T^+)$ we define $\Gamma_\gamma\in\mathrm{End}(T)$ as the linear endomorphism:

$$T\ni\tau\mapsto\Gamma_\gamma\tau:=(\mathrm{Id}_T\otimes\gamma)\Delta\tau. \tag{A.16}$$

The structure group is then defined as $G:=\{\Gamma_\gamma\in\mathrm{End}(T):\gamma\in\mathrm{Char}(T^+)\}$. One can check (see [FH20, Section 15.3 p.300]) that $\Delta$ satisfies the *coassociativity* $(\Delta\otimes\mathrm{Id})\Delta=(\mathrm{Id}\otimes\Delta)\Delta$ when interpreted as linear maps both in $\mathcal{L}(T^+,T^+\otimes T^+\otimes T^+)$ and $\mathcal{L}(T,T\otimes T^+\otimes T^+)$. For the first interpretation this identity implies the associativity of the convolution product of the characters, and for the second one it translates to the property $\Gamma_{\gamma_1}\Gamma_{\gamma_2}=\Gamma_{\gamma_1*\gamma_2}$ of the structure group.

From (A.15) we see that $\Gamma$ acts on monomials by translation, and therefore $\mathcal{P}:=\{\boldsymbol{X}^k\}_{k\in\mathbb{N}^{1+d}}$ spans a sector $P:=\mathrm{span}\{\mathcal{P}\}$ isomorphic to the polynomial structure as in [FH20, Section 13.2.1].

For $\Gamma\in G$ we will denote by $\gamma\in\mathrm{Char}(T^+)$ to the character such that $\Gamma=(\mathrm{Id}_T\otimes\gamma)\Delta$.

## A.4. Duality

In [BB21, Section 2] a $\star$-product between decorated trees which satisfies a duality property to the *coaction* $\Delta$ is constructed. We describe their construction specialising it to the decorated trees defined in Section A.2.

Denote by $\langle\mathscr{T}\rangle$ the $\mathbb{R}$-linear vector space spanned by $\mathscr{T}$, and similarly by $\langle\mathscr{T}_\alpha\rangle=\mathrm{span}\{\tau\in\mathscr{T}: |\tau|_\mathfrak{s}=\alpha\}$ for $\alpha\in\mathbb{R}$. We extend the definition of homogeneity $|\cdot|_\mathfrak{s}$ to $\langle\mathscr{T}\rangle_\alpha$ by setting it to $\alpha$. Given some decorated tree $\tau$, some $k\in\mathbb{Z}^{d+1}$ and $v\in N_\tau$ the decorated tree $\uparrow_v^k\tau$ is obtained by adding $k$ to the decoration $\mathfrak{n}(v)$ in the vertex $v\in N_\tau$, and in particular we have $|\uparrow_v^k\tau|_\mathfrak{s}=|\tau|_\mathfrak{s}+|k|_\mathfrak{s}$. We extend this operator linearly to $\langle\mathscr{T}\rangle$. Given two decorated trees $\sigma,\tau\in\mathscr{T}$, a node $v\in N_\tau$ and some $m\in\mathbb{N}^{d+1}$ the decorated tree $\sigma\curvearrowright_m^v\tau$ is given by the grafting of $\sigma$ into the vertex $v\in N_\tau$ of the tree $\tau$ by adding an edge decorated by $m$ between $\rho_\sigma$, the root of $\sigma$, and $v$. With this definition we have that

$$|\sigma\curvearrowright_m^v\tau|_\mathfrak{s}=|\mathcal{I}_m(\sigma)|_\mathfrak{s}+|\tau|_\mathfrak{s}. \tag{A.17}$$

For $\sigma,\tau\in\mathscr{T}$ with $\tau\neq\Xi$ we define

$$\sigma\curvearrowright_m\tau:=\sum_{v\in N_\tau\setminus L_\tau}\sum_{\substack{j\in\mathbb{N}^{d+1}\\ j\leqslant\mathfrak{n}(v)}}\binom{\mathfrak{n}(v)}{j}\sigma\curvearrowright_{m-j}^v(\uparrow_v^{-j}\tau)\in\langle\mathscr{T}\rangle, \tag{A.18}$$

which encodes all the ways to graft $\sigma$ into $\tau$, with *deformations* on the decorations (represented by $j \in \mathbb{N}^{1+d}$ which decreases decorations in the node $v$ of $\tau$). In (A.18) the term $\binom{\mathfrak{n}(v)}{j}$ is the binomial coefficient for multi-indexes. The first sum is restricted to the set $N_\tau \setminus L_\tau$ since by definition of $\mathscr{T} \setminus \{\Xi\}$ noise decorations are only allowed on leaves and grafting over a noise will break this condition. Moreover, the restriction $j \leqslant \mathfrak{n}(v)$ guarantees that the resulting trees have $\mathbb{N}^{1+d}$-valued $\mathfrak{n}$ decorations. Every term in (A.18) satisfies

$$|\sigma \curvearrowright_{m-j}^v (\uparrow_v^{-j} \tau)|_\mathfrak{s} = |\mathcal{I}_{m-j}(\sigma)|_\mathfrak{s} + |\uparrow_v^{-j} \tau|_\mathfrak{s} = |\mathcal{I}(\sigma)|_\mathfrak{s} - (|m|_\mathfrak{s} - |j|_\mathfrak{s}) + |\tau|_\mathfrak{s} - |j|_\mathfrak{s} = |\mathcal{I}_m(\sigma)|_\mathfrak{s} + |\tau|_\mathfrak{s},$$

and therefore $\sigma \curvearrowright_m \tau$ is homogeneous with homogeneity $|\sigma \curvearrowright_m \tau|_\mathfrak{s} = |\mathcal{I}_m(\sigma)|_\mathfrak{s} + |\tau|_\mathfrak{s}$. This grafting allows us to define a product $\curvearrowright$ between planted trees and trees by $\mathcal{I}_m(\sigma) \curvearrowright \tau := \sigma \curvearrowright_m \tau \in \langle \mathscr{T} \rangle$, which can be extended to the finite product of planted trees $(\prod_{i \in I} \mathcal{I}_{m_i}(\sigma_i)) \curvearrowright \tau$ by performing each grafting $\sigma_i \curvearrowright_{m_i} \tau$ independently of each other but grafting exclusively on nodes of $\tau$. To include polynomials decoration we extend the definition of $\uparrow_v^k$ to sets of nodes $B \subset N_\tau \setminus L_\tau$ as:

$$\uparrow_B^k \tau := \sum_{\substack{(k_v) \in (\mathbb{N}^{1+d})^B \\ \sum_{v \in B} k_v = k}} \left( \prod_{v \in B} \uparrow_v^{k_v} \right) \tau \in \langle \mathscr{T} \rangle.$$

This operation results in all the ways to add decorations on the nodes $B$ with total contribution $k \in \mathbb{N}^{1+d}$, and the result is a homogeneous element with $|\uparrow_B^k \tau|_\mathfrak{s} = |\tau|_\mathfrak{s} + |k|_\mathfrak{s}$.

Given decorated trees $\tau, \sigma \in \mathscr{T} \setminus \{\Xi\}$ with $\sigma = \boldsymbol{X}^k \prod_{i \in I} \mathcal{I}_{m_i}(\sigma_i)$ we define their $\star$-product as

$$\sigma \star \tau := \uparrow_{N_\tau \setminus L_\tau}^k \left( \left( \prod_{i \in I} \mathcal{I}_{m_i}(\sigma_i) \right) \curvearrowright \tau \right). \tag{A.19}$$

This product is homogeneous with homogeneity

$$|\sigma \star \tau|_\mathfrak{s} = |k|_\mathfrak{s} + \left| \prod_{i \in I} \mathcal{I}_{m_i}(\sigma_i) \curvearrowright \tau \right|_\mathfrak{s} = |k|_\mathfrak{s} + \sum_{i \in I} |\mathcal{I}_{m_i}(\tau)|_\mathfrak{s} + |\tau|_\mathfrak{s} = |\sigma|_\mathfrak{s} + |\tau|_\mathfrak{s}. \tag{A.20}$$

The following explicit computation of the $\star$-product illustrates the action of the *deformation* in its definition. Since $\boldsymbol{X}^k$ consists only of the root with a decoration $k$, we have

$$\mathcal{I}_m(\tau) \star \boldsymbol{X}^k = \sum_{j \in \mathbb{N}^{1+d}} \binom{\mathfrak{n}(\text{root})}{j} \tau \curvearrowright_{m-j}^{\rho_{\boldsymbol{X}^k}} (\uparrow_{\rho_{\boldsymbol{X}^k}}^{-j} \boldsymbol{X}^k) = \sum_{j \in \mathbb{N}^{1+d}} \binom{k}{j} \mathcal{I}_{m-j}(\tau) \, \boldsymbol{X}^{k-j}. \tag{A.21}$$

In general this is what happens when the grafting $\tau \curvearrowright_m^v \sigma$ occurs on a decorated node with $\mathfrak{n}(v) = k$.

**Remark A.10.** Observe that in the previous expression in order to produce a term without edge decoration at the root we need $k$ to be such that $\binom{k}{j} \neq 0$, or equivalently $k - j \in \mathbb{N}^{1+d}$.

We consider an inner product on $T$ given by

$$\langle \cdot, \cdot \rangle : T \times T \to \mathbb{R}, \qquad \langle \tau, \sigma \rangle := \tau! \, \delta_{\tau, \sigma}, \tag{A.22}$$

where $\delta_{\tau, \sigma}$ is Kronecker's delta function and $\tau!$ is the symmetry factor of the decorated tree as defined in (A.6). We can now state the duality property between $\star$ and $\Delta$ as stated in [BB21, Equation 2.2].

**Lemma A.11.** *For all $\mu \in \mathcal{T}^+$ and $\tau, \sigma \in \mathcal{T}$ we have the following duality property:*

$$\langle \mu \star \tau, \sigma \rangle = \langle \tau \otimes \mu, \Delta \sigma \rangle.$$

The proof of this duality can be found in [BM23, Theorem 4.2], and it is a generalisation of the duality shown in [BCMW22, Theorem 3.5] in the context of rough paths.

**Lemma A.12.** *For every $\sigma, \tau \in \mathcal{T}$ and $\Gamma \in G$ we have*

$$\langle \sigma, \Gamma \tau \rangle = \sum_{\mu \in \mathcal{T}^+} \langle \mu \star \sigma, \tau \rangle \frac{\gamma(\mu)}{\mu!}. \tag{A.23}$$

*Moreover, for every $\mathcal{I}_m(\tau) \in \mathcal{T}$ and $\Gamma \in G$ the nontrivial components of $\Gamma \mathcal{I}_m(\tau)$ are given by*

$$\langle \mathcal{I}_m(\sigma), \Gamma \mathcal{I}_m(\tau) \rangle = \sum_{\mu \in \mathcal{T}^+} \langle \mu \star \sigma, \tau \rangle \frac{\gamma(\mu)}{\mu!}, \qquad \langle \mathbf{X}^k, \Gamma \mathcal{I}_m(\tau) \rangle = \gamma(\mathcal{I}_{m+k}(\tau)). \tag{A.24}$$

**Proof.** Using the duality of Lemma A.11 we can write

$$\Delta \tau = \sum_{\sigma \in \mathcal{T}, \mu \in \mathcal{T}^+} \frac{\langle \sigma \otimes \mu, \Delta \tau \rangle}{\sigma! \, \mu!} \sigma \otimes \mu = \sum_{\sigma \in \mathcal{T}, \mu \in \mathcal{T}^+} \frac{\langle \mu \star \sigma, \tau \rangle}{\sigma! \, \mu!} \sigma \otimes \mu, \tag{A.25}$$

which induces a representation of the action of the structure group given for $\Gamma \in G$ as

$$\Gamma \tau = (\mathrm{Id} \otimes \gamma) \Delta \tau = \sum_{\sigma \in \mathcal{T}, \mu \in \mathcal{T}^+} \frac{\langle \mu \star \sigma, \tau \rangle}{\sigma! \, \mu!} \gamma(\mu) \, \sigma = \sum_{\sigma \in \mathcal{T}} \left\{ \sum_{\mu \in \mathcal{T}^+} \langle \mu \star \sigma, \tau \rangle \frac{\gamma(\mu)}{\mu!} \right\} \frac{\sigma}{\sigma!},$$

where $\gamma$ is the character associated to $\Gamma$ (see (A.16)), and (A.23) follows. Given $\mathcal{I}_m(\tau) \in \mathcal{T}$ and using definition (A.14) and formula (A.25) for $\Delta \tau$ we obtain

$$\begin{aligned}
\Delta \mathcal{I}_m(\tau) &= (\mathcal{I}_m \otimes \mathrm{Id}_{\mathcal{T}^+}) \Delta \tau + \sum_{k \in \mathbb{N}^{1+d}} \frac{\mathbf{X}^k}{k!} \otimes \mathcal{I}_{m+k}(\tau) \\
&= (\mathcal{I}_m \otimes \mathrm{Id}_{\mathcal{T}^+}) \sum_{\sigma \in \mathcal{T}, \mu \in \mathcal{T}^+} \frac{\langle \mu \star \sigma, \tau \rangle}{\sigma! \, \mu!} \sigma \otimes \mu + \sum_{k \in \mathbb{N}^{1+d}} \frac{\mathbf{X}^k}{k!} \otimes \mathcal{I}_{m+k}(\tau) \\
&= \sum_{\sigma \in \mathcal{T}, \mu \in \mathcal{T}^+} \frac{\langle \mu \star \sigma, \tau \rangle}{\sigma! \, \mu!} \mathcal{I}_m(\sigma) \otimes \mu + \sum_{k \in \mathbb{N}^{1+d}} \frac{\mathbf{X}^k}{k!} \otimes \mathcal{I}_{m+k}(\tau), \\
&= \sum_{\sigma \in \mathcal{T}, \mu \in \mathcal{T}^+} \frac{\langle \mathcal{I}_m(\mu \star \sigma), \mathcal{I}_m(\tau) \rangle}{\sigma! \, \mu!} \mathcal{I}_m(\sigma) \otimes \mu + \sum_{k \in \mathbb{N}^{1+d}} \frac{\mathbf{X}^k}{k!} \otimes \mathcal{I}_{m+k}(\tau) \\
&= \sum_{\sigma \in \mathcal{T}, \mu \in \mathcal{T}^+} \frac{\langle \mu \star \mathcal{I}_m(\sigma), \mathcal{I}_m(\tau) \rangle}{\sigma! \, \mu!} \mathcal{I}_m(\sigma) \otimes \mu + \sum_{k \in \mathbb{N}^{1+d}} \frac{\mathbf{X}^k}{k!} \otimes \mathcal{I}_{m+k}(\tau),
\end{aligned}$$

where in the last line we used that the part of $\mu \star \mathcal{I}_m(\sigma)$ that grafts or decorates at the root cannot produce a planted tree. Applying $\Gamma$ we conclude that

$$\Gamma \mathcal{I}_m(\tau) = \sum_{\sigma \in \mathcal{T}} \left\{ \sum_{\mu \in \mathcal{T}^+} \langle \mu \star \sigma, \tau \rangle \frac{\gamma(\mu)}{\mu!} \right\} \frac{\mathcal{I}_m(\sigma)}{\sigma!} + \sum_{k \in \mathbb{N}^{1+d}} \gamma(\mathcal{I}_{m+k}(\tau)) \frac{\mathbf{X}^k}{k!},$$

from where (A.24) follows. □

## A.5. Coherence

In the framework of regularity structures solutions to SPDEs, encoded by modelled distributions, have an analytic and an algebraic component. The analytic part guarantees that the reconstruction of the modelled distribution satisfies the correct (renormalised) SPDE, while the algebraic component relates the coefficients of the rooted trees in terms of the polynomial coefficients. This algebraic relationship is encoded in the *coherence map*.

To describe this coherence map we consider an SPDE with additive noise of the form

$$\mathcal{L}\varphi = F(\varphi, D\varphi, \dots) + \xi.$$

Given some modelled distribution $\Phi \colon \mathbb{R}^{1+d} \to \mathbb{R}$ we denote by $F(\Phi, D\Phi, \dots)$ the *lift* of the non-linearity $F$ which acts on modelled distributions. If $F$ is a polynomial function this lift can be defined for arbitrary modelled distributions as long as the regularity structure $T$ contains the required symbols to described the products. However, if $F$ is an arbitrary smooth function then it can only be lifted to act on modelled distributions which take values on a function-like sector of the regularity structure. We refer the reader to [FH20, Chapter 14] for details on this lift.

It is convenient to consider the non-linearity $F \colon \mathbb{R}^{\mathbb{N}^{1+d}} \to \mathbb{R}$ as a function of all possible derivatives of $\varphi$, and denote by $D_k F \colon \mathbb{R}^{\mathbb{N}^{1+d}} \to \mathbb{R}$ the partial derivative of $F$ with respect to it's $k$-th coordinate for $k \in \mathbb{N}^{1+d}$. Given $j \in \{0, \dots, d\}$ and $e_j \in \mathbb{N}^{1+d}$ the canonical multi-index, we define $\partial^{e_j} F \colon \mathbb{R}^{\mathbb{N}^{1+d}} \to \mathbb{R}$ as the non-linearity given by

$$\partial^{e_j} F := \sum_{k \in \mathbb{N}^{1+d}} \mathcal{X}_{k+e_j} D_{e_j} F, \tag{A.26}$$

where for $k \in \mathbb{N}^{1+d}$ we denote by $\mathcal{X}_k \colon \mathbb{R}^{\mathbb{N}^{1+d}} \to \mathbb{R}$ to the non-linearity $\mathcal{X}_k(\{\partial^n \varphi\}_{n \in \mathbb{N}^{1+d}}) := \partial^k \varphi$. This definition can be extended to $\partial^k F$ for all $k \in \mathbb{N}^{1+d}$ by composition. We can now define the coherence map.

**Definition A.13.** *Given a non-linearity $F \colon \mathbb{R}^{\mathbb{N}^{1+d}} \to \mathbb{R}$ we define $\Upsilon^F[\Xi] := 1$, and given a tree $\tau \in \mathcal{T}$ of the form (A.7) we define recursively a non-linearity $\Upsilon^F[\tau] \colon \mathbb{R}^{\mathbb{N}^{1+d}} \to \mathbb{R}$ as*

$$\Upsilon^F[\tau] = \left(\prod_{i=1}^m \Upsilon^F[\tau_i]\right) \cdot \left(\partial^k \prod_{i=1}^n D_{m_i} F\right). \tag{A.27}$$

*The map $\tau \mapsto \Upsilon^F[\tau]$ is called the coherence map.*

One of the main results in [BCCH20] is that a modelled distribution $\Phi$ solves the equation

$$\Phi = \mathcal{I}(F(\Phi, D\Phi, \dots) + \Xi) + \sum_{k \in \mathbb{N}^{1+d}} \frac{\langle \boldsymbol{X}^k, \Phi \rangle}{k!} \boldsymbol{X}^k,$$

if and only if for all $\tau \in \mathcal{T} \setminus \{\boldsymbol{X}^k\}_{k \in \mathbb{N}^{1+d}}$ the coefficient at $\mathcal{I}(\tau)$ of $\Phi$ is determined by the polynomial part of $\Phi$ and the coherence map via the identity

$$\langle \mathcal{I}(\tau), \Phi \rangle = \Upsilon^F[\tau](\{\langle \boldsymbol{X}^k, \Phi \rangle\}_{k \in \mathbb{N}^{1+d}}).$$

We have the following morphism property between the $\star$-product and the coherence map.

**LEMMA A.14**. ([BB21, PROPOSITION 2].) *For every $\tau \in \mathcal{T}$ and $\mu = \boldsymbol{X}^k \prod_{i=1}^n \mathcal{I}_{m_i}(\sigma_i) \in \mathcal{T}^+$ we have*

$$\Upsilon[\mu \star \tau] = \left( \prod_{i=1}^n \Upsilon[\sigma_i] \right) (\partial^k D_{m_1} \cdots D_{m_n}) \Upsilon[\tau].$$

**Acknowledgments.** HW thanks Scott Smith for explaining the use of Simon's method to prove basepoint dependent Schauder estimates for the local heat operator as in [SS24]. SE and HW are funded by the Deutsche Forschungsgemeinschaft (DFG, German Research Foundation) under Germany's Excellence Strategy EXC 2044-390685587, Mathematics Münster: Dynamics–Geometry–Structure. HW is funded by the European Union (ERC, GE4SPDE, 101045082). SE acknowledges support from CoNaCyT scholarship 2020-000000-01EXTF-00129 CVU 779198 and the EPSRC Centre for Doctoral Training in Statistical Applied Mathematics at Bath (SAMBa), under the project EP/S022945/1. This paper has been written with GNU T<sub>E</sub>X<sub>MACS</sub> (www.texmacs.org).